\documentclass[11pt,letterpaper,twoside,reqno,nosumlimits]{amsart}

\makeatletter

\renewcommand{\l@section}{%
  \@tocline{1}{0.8em plus 0.2em}{1pc}{}{\bfseries}}

\renewcommand{\l@subsection}{%
  \@tocline{2}{0.15em}{1pc}{5pc}{}}

\makeatother

\usepackage{newclude}

\usepackage[usenames,dvipsnames]{xcolor}
\usepackage{fancyhdr}
\usepackage{amsmath,amsfonts,amsbsy,amsgen,amscd,amssymb,amsthm}
\usepackage{url}

\usepackage{mathtools}

\mathtoolsset{showonlyrefs}

\usepackage[font=small,margin=0.25in,labelfont={bf},labelsep={colon}]{caption}

\usepackage{subcaption}

\usepackage{tikz}
\usepackage{microtype}
\usepackage{enumitem}

\definecolor{dark-gray}{gray}{0.3}
\definecolor{dkgray}{rgb}{.4,.4,.4}
\definecolor{dkblue}{rgb}{0,0,.5}
\definecolor{medblue}{rgb}{0,0,.75}
\definecolor{rust}{rgb}{0.5,0.1,0.1}

\usepackage[colorlinks=true]{hyperref}

\hypersetup{urlcolor=rust}
\hypersetup{citecolor=dkblue}
\hypersetup{linkcolor=dkblue}

\usepackage{setspace}

\usepackage{graphicx}
\usepackage{booktabs,longtable,tabu} 
\usepackage{multirow} 

\usepackage{float}

\usepackage[full]{textcomp}

\usepackage[scaled=.98,sups,osf]{XCharter}
\usepackage[scaled=1.04,varqu,varl]{inconsolata}
\usepackage[type1]{cabin}
\usepackage[charter,vvarbb,scaled=1.07]{newtxmath}
\usepackage[cal=boondoxo]{mathalfa}
\usepackage[T1]{fontenc}

\usepackage{bm} 

\graphicspath{{figures/}}

\newtheorem{theorem}{Theorem}[section]
\newtheorem{lemma}[theorem]{Lemma}

\newtheorem{proposition}[theorem]{Proposition}

\newtheorem{corollary}[theorem]{Corollary}

\theoremstyle{definition}

\newtheorem{remark}[theorem]{Remark}

\numberwithin{equation}{section} 
\numberwithin{figure}{section}
\numberwithin{table}{section}

\floatstyle{plaintop}
\newfloat{recipe}{thp}{lor}
\floatname{recipe}{Recipe}
\numberwithin{recipe}{section}

\renewcommand{\phi}{\varphi}

\newcommand{\iunit}{\mathrm{i}}

\providecommand{\mathbbm}{\mathbb} 

\newcommand{\R}{\mathbbm{R}}

\newcommand{\N}{\mathbbm{N}}
\newcommand{\Z}{\mathbbm{Z}}

\newcommand{\diff}{\mathrm{d}}
\newcommand{\idiff}{\, \diff}

\newcommand{\norm}[1]{\left\Vert {#1} \right\Vert}

\newcommand{\triplenorm}[1]{{\left\vert\kern-0.25ex\left\vert\kern-0.25ex\left\vert #1
    \right\vert\kern-0.25ex\right\vert\kern-0.25ex\right\vert}}

\renewcommand{\div}{\text{div}}

\def\XXint#1#2#3{{%
\setbox0=\hbox{$#1{#2#3}{\int}$}
\vcenter{\hbox{$#2#3$}}\kern-.5\wd0}}

\renewcommand{\leq}{\leqslant}
\renewcommand{\geq}{\geqslant}
\renewcommand{\subset}{\subseteq}

\usepackage{placeins,needspace}

\allowdisplaybreaks

\usepackage[numbers]{natbib}
\mathtoolsset{showonlyrefs}

\makeatletter
\renewcommand{\paragraph}{%
  \@startsection{paragraph}{4}{\z@}%
    {3pt}{-\fontdimen2\font}%
    {$\circ$\ \normalfont\bfseries}%
}
\makeatother

\hypersetup{
  pdftitle={A viscous cascade construction of nonunique Leray--Hopf solutions to the unforced 3D Navier--Stokes equations},
  pdfsubject={Nonuniqueness of suitable Leray--Hopf solutions by a viscous cascade},
  pdfkeywords={Navier--Stokes equations, Leray--Hopf solutions, nonuniqueness, viscous cascade}
}

\newcommand{\PP}{\mathbb P}
\newcommand{\tp}{{t_{\mathrm p}}}
\newcommand{\tf}{{t_{\mathrm f}}}
\newcommand{\ts}{{t_\ast}}

\title{A Viscous Cascade Construction
of Nonunique Leray--Hopf Solutions\\
to the Unforced $3$D Navier--Stokes Equations}
\author{Thomas Y. Hou}
\address{Applied and Computational Mathematics,
California Institute of Technology, Pasadena, CA, USA}
\email{hou@cms.caltech.edu}

\author{Xiang Qin}
\address{Applied and Computational Mathematics,
California Institute of Technology, Pasadena, CA, USA}
\email{xqin2@caltech.edu}

\author{Peicong Song}
\address{Applied and Computational Mathematics,
California Institute of Technology, Pasadena, CA, USA}
\email{psong2@caltech.edu}

\author{Wenhan Zhang}
\address{Applied and Computational Mathematics,
California Institute of Technology, Pasadena, CA, USA}
\email{wzhang9@caltech.edu}

\date{}

\begin{document}
\begin{abstract}
We construct a nonzero divergence-free initial velocity on $\R^3$ and a
continuum of pairwise distinct suitable Leray--Hopf solutions of the
unforced three-dimensional Navier--Stokes equations on a common time
interval. The initial datum is axisymmetric, smooth away from the origin, and
belongs to every $L^{\mathfrak q}$ with $1\leq \mathfrak q<3$, but not to $L^3$. All the solutions
attain this datum strongly in $L^2$, are smooth at positive times, and
satisfy both the global and local energy equalities through time zero.
The same family belongs simultaneously to every
$L^{\mathfrak p}_tL^{\mathfrak q}_x$ with $1\leq \mathfrak p\leq\infty$, $2\leq \mathfrak q\leq\infty$, and
$2/\mathfrak p+3/\mathfrak q>1$.
The same solutions also belong to the critical
Lorentz space $L^\infty_tL^{3,\infty}_x$.
This endpoint is sharp for uniqueness in the Leray--Hopf class
within the Lorentz scale.

The proof is entirely analytic and uses a cascade through geometrically
separated Gaussian swirls. We establish growth for the time-dependent
viscous linearization and determine the sign of its complete quadratic
stress, including the interaction of the reference velocity with the
second variation. The integrated stress generates a controlled perturbation
at the next larger spatial scale. At each finite cutoff, a prescribed
terminal coordinate selects an initial perturbation whose physical
$L^2$ norm tends to zero as the cutoff is removed. Uniform energy, moment,
and transfer estimates yield the common strong initial trace while
preserving a nonzero observable at a fixed positive time.
\end{abstract}
\maketitle

\begingroup
\fontsize{10}{10.5}\selectfont
\tableofcontents
\endgroup

\section{Introduction}\label{sec:intro}

\subsection{The Cauchy problem and the energy class}

We consider the unforced incompressible Navier--Stokes equations
\begin{equation}\label{intro:NS}
 \partial_tu+(u\cdot\nabla)u-\Delta u+\nabla p=0,
 \qquad \nabla\cdot u=0,
 \qquad u|_{t=0}=u_0,
 \qquad (t,x)\in(0,T)\times\R^3.
\end{equation}
The viscosity is normalized to one. Leray's construction on the whole space,
and Hopf's construction in domains, provide global weak solutions for
divergence-free initial velocities with finite kinetic energy
\cite{Ler34,Hop51}. The energy inequality controls both the velocity and its
integrated dissipation. It does not, by itself, control the amplification of
a difference between two solutions near a singular initial time. The purpose
of this paper is to construct nonunique solutions in this energy class by
a sequence of viscous interactions between geometrically separated spatial
scales.

We specify the solution class before stating the result. Let $L^2_\sigma$
be the closed subspace of divergence-free vector fields in $L^2(\R^3)$.
A Leray--Hopf solution on $[0,T]$ is a distributional solution of
\eqref{intro:NS} such that
\begin{equation}\label{intro:energy-space}
 u\in C_{\mathrm w}([0,T];L^2_\sigma)
              \cap L^2(0,T;H^1(\R^3)),
 \qquad \lim_{t\downarrow0}\|u(t)-u_0\|_2=0,
\end{equation}
and
\begin{equation}\label{intro:energy-inequality}
 \frac12\|u(t)\|_2^2+\int_s^t\|\nabla u(\tau)\|_2^2\idiff\tau
 \leq\frac12\|u(s)\|_2^2
\end{equation}
for $s=0$, for almost every $s\in(0,T)$, and for every $t\in[s,T]$.
Here $C_{\mathrm w}$ denotes weak continuity. The distributional equation
includes the initial term; equivalently, for every divergence-free
$\psi\in C_c^\infty([0,T)\times\R^3;\R^3)$,
\begin{equation}\label{intro:weak-formulation}
 \int_0^T\int_{\R^3}
       \bigl(u\cdot(\partial_t\psi+\Delta\psi)
                    +(u\otimes u):\nabla\psi\bigr)\idiff x\idiff t
       +\int_{\R^3}u_0\cdot\psi(0)\idiff x=0.
\end{equation}

A solution is suitable if it has a pressure
$p\in L^{3/2}_{\mathrm{loc}}((0,T)\times\R^3)$ and satisfies the local
energy inequality. When the initial time is included, this means that every
nonnegative $\phi\in C_c^\infty([0,T)\times\R^3)$ satisfies
\begin{equation}\label{intro:local-energy}
 \begin{aligned}
 \int_0^T\int_{\R^3}|\nabla u|^2\phi\idiff x\idiff t
 \leq{}&\int_0^T\int_{\R^3}
       \Bigl(\frac{|u|^2}{2}(\partial_t\phi+\Delta\phi)
             +\Bigl(\frac{|u|^2}{2}+p\Bigr)u\cdot\nabla\phi\Bigr)
                   \idiff x\idiff t\\
 &+\frac12\int_{\R^3}|u_0|^2\phi(0)\idiff x.
 \end{aligned}
\end{equation}
The local inequality is the additional condition underlying the partial
regularity theory of Caffarelli, Kohn, and Nirenberg \cite{CKN82}; see also
Lin's proof \cite{Lin98}. The solutions constructed here satisfy equality
in both \eqref{intro:energy-inequality} and \eqref{intro:local-energy}.

Uniqueness holds under stronger integrability. In particular, the
Prodi--Serrin condition
\begin{equation}\label{intro:serrin}
 u\in L^{\mathfrak p}(0,T;L^{\mathfrak q}(\R^3)),\qquad
 \frac2{\mathfrak p}+\frac3{\mathfrak q}\leq1,\qquad \mathfrak q>3,
\end{equation}
implies uniqueness among Leray--Hopf solutions with the same initial datum
\cite{Pro59,Ser63}. The strong spatial $L^3$ endpoint is covered by the
theorem of Escauriaza, Seregin, and \v{S}ver\'ak \cite{ESS03}.
The scaling of \eqref{intro:NS},
\begin{equation}\label{intro:scaling}
 u_\lambda(t,x)=\lambda u(\lambda^2t,\lambda x),
 \qquad p_\lambda(t,x)=\lambda^2p(\lambda^2t,\lambda x),
\end{equation}
explains the boundary $2/\mathfrak p+3/\mathfrak q=1$, where $\mathfrak p$ denotes the temporal exponent
and $\mathfrak q$ the spatial exponent. Our solutions belong simultaneously to
every strong mixed space on its strict supercritical side.

The distinction between finite energy and critical regularity also appears
in the mild-solution theory. Fujita and Kato developed the critical
Sobolev framework \cite{FK64}, and Kato established the critical $L^3$
theory \cite{Kat84}; small data in the larger critical space
$\mathrm{BMO}^{-1}$ are covered by Koch and Tataru \cite{KT01}.
Phuc established regularity in the nonendpoint Lorentz classes
$L^\infty_tL^{3,\mathfrak r}_x$, $3<\mathfrak r<\infty$, and weak--strong
uniqueness for Leray--Hopf solutions under the additional assumption
$u_0\in L^3$ \cite{Phu15}.
We show that this initial-data assumption can be removed. The same
Lorentz condition $L^\infty_tL^{3,\mathfrak r}_x$ implies weak--strong uniqueness among Leray--Hopf
solutions. Proposition~\ref{app:lorentz-uniqueness} gives the strengthened
statement and proof, using the local Lorentz-space theory and relative
energy method \cite{Dub03}.
The solutions constructed here belong to the endpoint class
$L^\infty_tL^{3,\infty}_x$, with a common initial datum in
$L^{3,\infty}\setminus \left(\bigcup_{1\leq \mathfrak r<\infty}L^{3,\mathfrak r}\right)$.

\subsection{The main result}

The initial datum will be an axisymmetric sum of Gaussian swirls. It is smooth
away from the origin, has finite kinetic energy, and has Gaussian tails at
spatial infinity. The radii decrease geometrically toward the origin. A
perturbation concentrated at increasingly small scales can therefore vanish in
the initial energy norm while its effect reaches a fixed larger scale at
positive time.

\begin{theorem}[Nonuniqueness by a viscous cascade]\label{intro:main}
There exist $T>0$, a nonzero real-valued divergence-free initial velocity
$u_0\in L^2(\R^3)$, a solution $\bar u$, and solutions
$u^\vartheta$, $0<\vartheta<1$, of \eqref{intro:NS}, with the following
properties.
\begin{enumerate}[label=\textup{(\roman*)},leftmargin=*]
\item The initial datum is axisymmetric and satisfies
\begin{equation}\label{intro:datum-properties}
 u_0\in C^\infty(\R^3\setminus\{0\})\cap L^{3,\infty}(\R^3)
              \cap\bigcap_{1\leq \mathfrak q<3}L^{\mathfrak q}(\R^3),
 \qquad u_0\notin L^3(\R^3).
\end{equation}
\item Every listed solution belongs to
\begin{equation}\label{intro:solution-regularity}
 C([0,T];L^2_\sigma)\cap L^2(0,T;H^1(\R^3))
                 \cap C^\infty((0,T)\times\R^3)
\end{equation}
and has initial value $u_0$ strongly in $L^2$. It has a pressure in
$L^{5/3}((0,T)\times\R^3)$, satisfies the distributional momentum equation
with the initial term, and satisfies the global energy equality
\begin{equation}\label{intro:energy-equality}
 \frac12\|u(t)\|_2^2+\int_s^t\|\nabla u(\tau)\|_2^2\idiff\tau
 =\frac12\|u(s)\|_2^2,\qquad0\leq s\leq t\leq T,
\end{equation}
as well as local energy equality through time zero. In particular, every
listed solution is suitable and Leray--Hopf.
\item The same solutions belong to every space
\begin{equation}\label{intro:strict-mixed-spaces}
 L^{\mathfrak p}(0,T;L^{\mathfrak q}(\R^3)),\qquad
 1\leq \mathfrak p\leq\infty,\quad2\leq \mathfrak q\leq\infty,
 \quad\frac2{\mathfrak p}+\frac3{\mathfrak q}>1,
\end{equation}
with the convention $1/\infty=0$.
In addition, the same family satisfies
\begin{equation}\label{intro:lorentz-endpoint}
 \bar u,u^\vartheta\in L^\infty(0,T;L^{3,\infty}(\R^3)),
 \qquad 0<\vartheta<1.
\end{equation}
\item The listed solutions are pairwise distinct. For any two of them,
$u\ne v$, and every $0<\varepsilon\leq T$, one has
\begin{equation}\label{intro:immediate-distinction}
 \int_0^\varepsilon\|u(t)-v(t)\|_2^2\idiff t>0.
\end{equation}
Moreover, there are one time $t_{\mathrm w}\in(0,T)$ and one real test
$\varphi\in C_c^\infty(\R^3;\R^3)$ such that
\begin{equation}\label{intro:fixed-witness}
 \langle u^\vartheta(t_{\mathrm w})-\bar u(t_{\mathrm w}),\varphi\rangle>0
 \qquad(0<\vartheta<1).
\end{equation}
\end{enumerate}
\end{theorem}

The proof is completed in Section~\ref{limit:main}. The construction uses
analytic estimates for the viscous evolution, its full quadratic stress, and
the infinite-scale limit. The datum is specified in \eqref{bg:datum} after a
finite ordered choice of parameters. The Gaussian tails allow a finite-energy realization of the entire
array.  The parameter $\vartheta$ is assigned by prescribing a terminal
$L^2$ coordinate of the limiting solution.

The mixed-norm statement is made for one family and one datum, with all
exponents quantified simultaneously. Its strict inequality is consistent
with \eqref{intro:serrin} and the strong $L^3$ endpoint. Positive-time
smoothness does not give an integrable Lipschitz bound down to the initial
time; this is the interval on which the increasingly small scales are active.

\begin{remark}
The endpoint \eqref{intro:lorentz-endpoint} is sharp within the
Lorentz scale at spatial exponent three, for bounds uniform down to
the initial time.  By Proposition~\ref{app:lorentz-uniqueness}, no
listed solution belongs to $L^\infty(0,T;L^{3,\mathfrak r}(\R^3))$ for any
$1\leq \mathfrak r<\infty$: such a bound would force uniqueness among
Leray--Hopf solutions with the same datum.  Independently,
\eqref{limit:datum-not-finite-lorentz} shows that the datum itself
belongs to no finite-index $L^{3,\mathfrak r}$.  This does not exclude
membership at individual positive times.
\end{remark}

\subsection{The mechanism in brief}\label{intro:mechanism}
The essential idea is to make a difference between two solutions travel
from an arbitrarily small spatial scale to a fixed observable scale.
At each step a perturbation grows near an unstable viscous swirl.
Its \emph{signed quadratic stress}, integrated over the growth interval,
produces the perturbation that reaches the next larger scale.  Diffusion
controls the intervening evolution, with viscosity fixed at one.

The amplification of a \emph{difference} is compatible with the
energy inequality for each solution.  Indeed, if $u$ and $\bar u$
are two smooth unforced solutions and $h=u-\bar u$, then
\[
 \frac12\frac{\mathrm d}{\mathrm dt}\|h\|_2^2
    +\|\nabla h\|_2^2
  =-\int_{\R^3}(h\cdot\nabla)\bar u\cdot h\,\idiff x.
\]
The right-hand side has no fixed sign.  When
$\int_0^T\|\nabla\bar u(t)\|_\infty\,\idiff t<\infty$,
Gronwall's inequality gives uniqueness; the increasingly small
swirls make this sufficient uniqueness estimate unavailable uniformly
down to the initial time.

\paragraph{A finite-energy array of viscous swirls.}
Write $x=(x_h,z)\in\R^2\times\R$, $r=|x_h|$, and
$x_h^\perp=(-x_2,x_1,0)$.  Our background begins with the
axisymmetric, divergence-free Gaussian swirl
\begin{equation}\label{intro:brief-swirl}
 Q_\ell(0,x)=\frac{A}{\ell^2}
    \exp\!\left(-\frac{r^2}{4\ell^2}-\frac{z^2}{2b^2\ell^2}\right)x_h^\perp.
\end{equation}
Here $\ell$ is its radial length scale, $A$ is a large but fixed
amplitude parameter, and $b$ is a large axial-width parameter.
Thus each swirl has velocity of order $A/\ell$ in its core and extends
much farther along the symmetry axis than across it.  We choose a
fixed geometric ratio $R>1$ and place such swirls at radii
$1,R^{-1},R^{-2},\ldots$.  The sum has finite kinetic energy because
a direct Gaussian integration gives
\begin{equation}\label{intro:single-swirl-energy}
 \|Q_\ell(0)\|_2^2=4\pi^{3/2}A^2b\ell.
\end{equation}
Thus $\sum_{n\geq0}\|Q_{R^{-n}}(0)\|_2<\infty$ for fixed
$A$, $b$, and $R>1$.

The heat evolutions of these swirls are explicit, but their sum is
not an exact Navier--Stokes solution: the pressure must balance
interactions between \emph{all pairs} of swirls, and axial localization
leaves a residual.  We construct a small correction to the heat sum
so that the resulting reference velocity $\bar u$ is an \emph{exact
unforced, finite-energy} solution.  The estimates for this correction
are uniform in the number of retained scales.  Proposition~\ref{bg:reference} provides the reference, its quantitative
bounds, and its scaling relation between consecutive stages.

\paragraph{An unstable mode of one viscous swirl.}
To see the local amplification mechanism, focus on the swirl of radial
scale one.  In its central axial region, the background is close to
the axially invariant Gaussian swirl
\[
 U_{\mathrm d}(t,r)=A\,r(1+t)^{-2}
      e^{-r^2/(4(1+t))}e_\theta,
 \qquad e_\theta=\frac{x_h^\perp}{r}.
\]
This comparison field solves the viscous Navier--Stokes equations
with an appropriate radial pressure, but has infinite energy in
$\R^3$ because it is independent of $z$.  We use it only to analyze
one local interaction; the actual solution always uses the localized
finite-energy reference $\bar u$.

Linearize around $U_{\mathrm d}$, take the Fourier transform in the
axial coordinate $z$, and let $k>0$ denote the axial Fourier
frequency.  Freeze the coupling at time $t=0$, temporarily omit the heat
operator, and factor out the amplitude $A$.  The remaining projected
\emph{inviscid coupling operator} on the radial and azimuthal
components is denoted by $\mathcal C_k$.  It has a positive eigenvalue
$\sigma_k$.  For a fixed $k$, the full
linearized equation contains both the amplifying term
$A\mathcal C_k$ and the dissipative heat operator.

\begin{figure}[!tbp]
\centering
\begin{minipage}[t]{0.485\textwidth}
    \centering
    \includegraphics[width=\linewidth,height=2.25in,
      keepaspectratio]{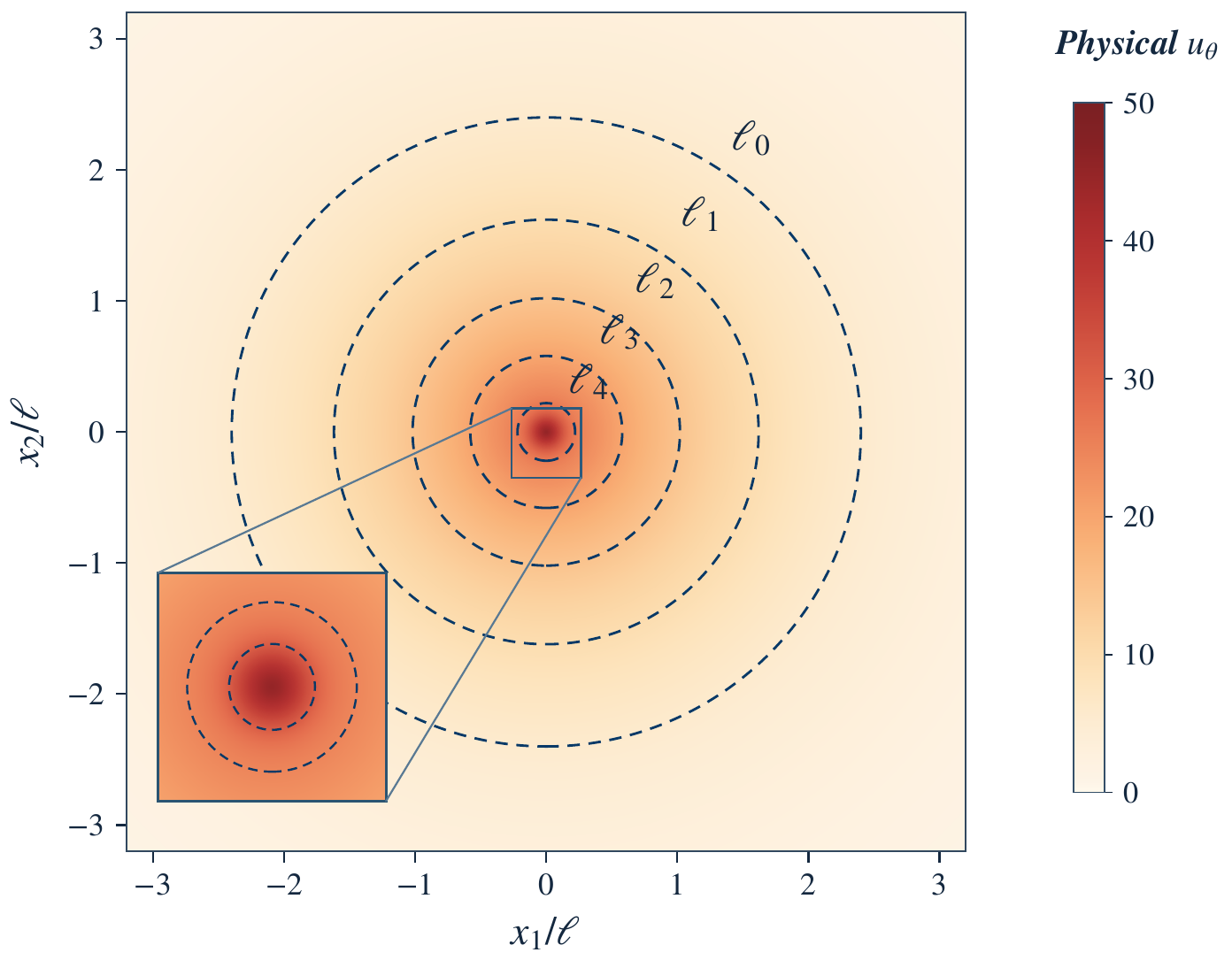}
    \par\smallskip
    {\small (a)}
\end{minipage}\hfill
\begin{minipage}[t]{0.485\textwidth}
    \centering
    \includegraphics[width=\linewidth,height=2.25in,
      keepaspectratio]{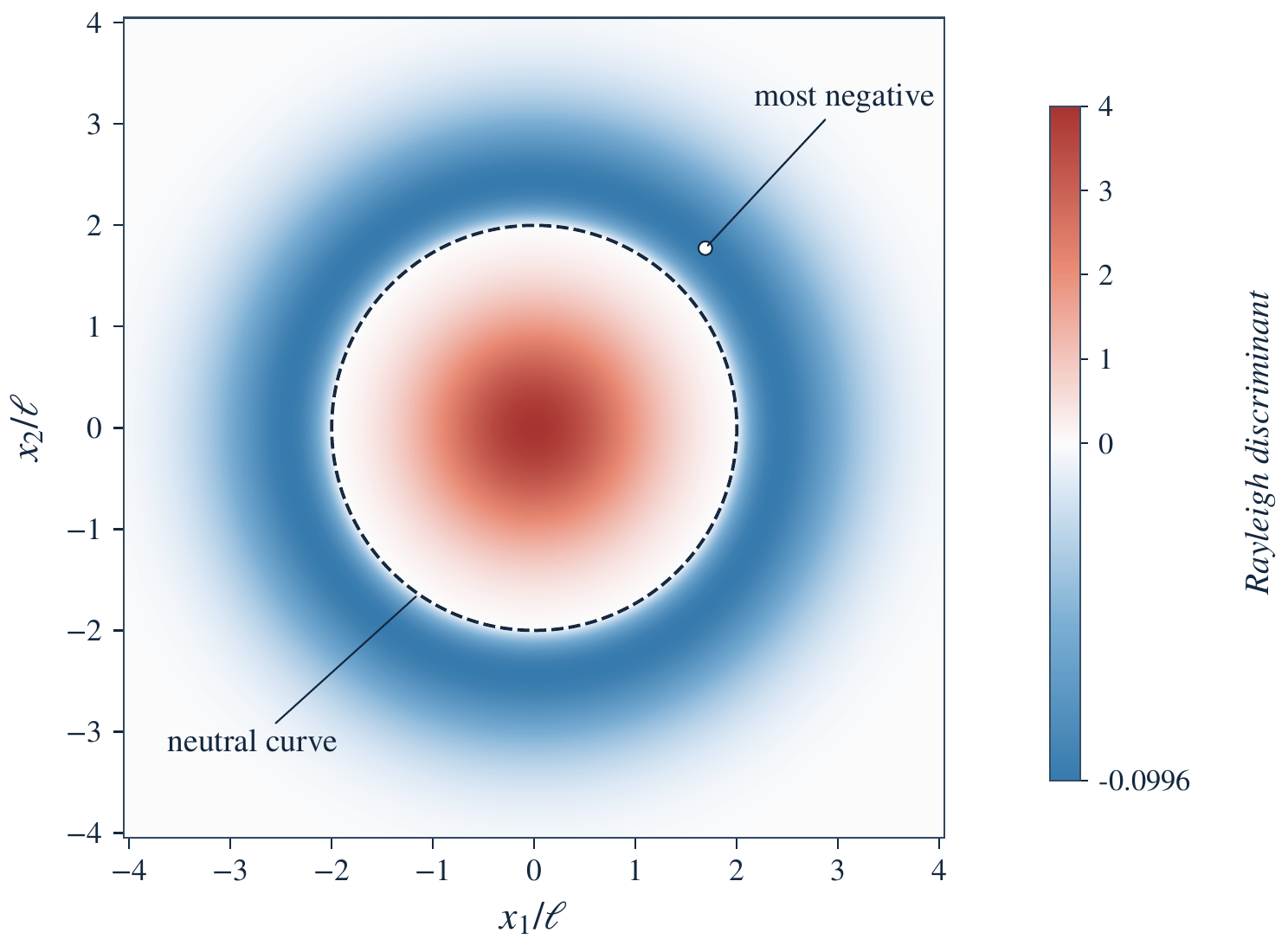}
    \par\smallskip
    {\small (b)}
\end{minipage}

\caption{\textbf{Nested Gaussian swirls and centrifugal instability.}
Both panels show top-down views in the horizontal
$(x_1,x_2)$-plane.
\textup{(a)} The Gaussian-swirl array in
\eqref{intro:brief-swirl}, with characteristic radii
$\ell_n=R^{-n}$. The inset magnifies the innermost
displayed scales.
\textup{(b)} The Rayleigh discriminant of the normalized
frozen swirl $U_\theta(r)=r e^{-r^2/4}$ is
$(4-r^2)e^{-r^2/2}$. It changes sign at $r=2$ and
attains its minimum $-2e^{-3}$ at $r=\sqrt{6}$.
The color scale is asymmetric about zero, with the
positive and negative ranges normalized separately
to make the comparatively small negative values visible.
Thus, color intensities of opposite signs are not
directly comparable in magnitude.
The dashed circles in \textup{(a)} mark characteristic
radii, while the dashed circle in \textup{(b)} marks
the neutral curve $r=2$.}
\label{fig:intro-swirl-rayleigh}
\end{figure}

The positive eigenvalue has an explicit variational explanation.
After eliminating one component of the frozen eigenvector, the
relevant radial potential is
$(r^2-4)e^{-r^2/2}$.  It is positive outside $r=2$ and
has maximum $2e^{-3}$ at $r=\sqrt6$.  The Rayleigh quotient
for $\sigma_k^2$ contains this potential in its numerator and
nonnegative radial-derivative penalties proportional to $k^{-2}$
in its denominator.  Thus
\[
  0<\sigma_k<c_\star,\qquad
  \sigma_k\longrightarrow c_\star=\sqrt2\,e^{-3/2}
  \quad\text{as }k\to\infty.
\]
The high-frequency limit here concerns the \emph{frozen inviscid}
operator; viscosity would suppress arbitrarily large axial
frequencies.  We instead fix a finite frequency band where
$\sigma_k$ is sufficiently close to $c_\star$, and \emph{then} take $A$
large.  The resulting time-dependent viscous response has a growing
component, with amplification on the exponential scale
$e^{c_\star A}$.  Showing that this growth survives diffusion and the
movement of the growing eigenspace is one of the principal estimates
in Section~\ref{sec:seed}.  No numerical eigenvalue computation is
needed: the unstable eigenvalue is obtained from the radial
variational problem. 
The nested Gaussian-swirl geometry and the sign structure of
the frozen Rayleigh discriminant are illustrated in
Figure~\ref{fig:intro-swirl-rayleigh}.

\paragraph{From linear growth to a persistent signed stress.}
The existence of an unstable mode is not yet a cascade: the next
scale must receive the \emph{particular} perturbation created by the
preceding one, rather than a freely chosen eigenvector.  We prove
that this prescribed input has a nonzero component along the
growing mode.  After a short heat preparation, the nonlinear calculation uses the complete linear response at all Fourier frequencies.

Let $U$ be the reference in a single transition, and write a small
perturbed solution formally as
$u^\varepsilon=U+\varepsilon v+\varepsilon^2w+\cdots$.
The first variation $v$ solves the equation linearized about $U$;
the second variation $w$ is driven by the quadratic product
$v\otimes v$.  When the linearized terms involving $U$ are moved
to the right side of the heat equation for $w$, the source becomes
the \emph{complete quadratic stress}
\begin{equation}\label{intro:brief-stress}
 \mathsf S=v\otimes v+U\otimes w+w\otimes U,
 \qquad (\partial_t-\Delta)w=-\PP\operatorname{div}\mathsf S,
\end{equation}
where $\PP$ is the Leray projection onto divergence-free fields.
The mixed terms retain an interaction of the same quadratic order:
the product of first variations creates a nonzero axial mean of the
azimuthal component of $w$, and this mean continues to interact
with the swirl after the rapidly varying first-order perturbation
has been damped by viscosity.  The elementary identity
$\cos^2(kz)=\frac12(1+\cos(2kz))$ illustrates how an oscillatory
quantity can produce a slowly varying quadratic component.

By axisymmetry, the spatial integral of the symmetric stress is
diagonal with equal horizontal entries.  Write it as
$\int_{\R^3}\mathsf S(t,x)\,\idiff x
=\operatorname{diag}(s_h(t),s_h(t),s_z(t))$.
The quantity that determines transfer is the \emph{time-integrated
anisotropy}
\begin{equation}\label{intro:brief-mass}
 M=-\int_{t_{\mathrm p}}^{t_{\mathrm f}}
                    (s_h(t)-s_z(t))\,\idiff t,
\end{equation}
where $t_{\mathrm p}$ is the short preparation time and
$t_{\mathrm f}$ is the later time at which the first variation has
decayed.  A central part of our analysis proves that $M>0$ and
that it is sufficiently large to overcome the loss incurred when
passing to the next spatial scale.  The proof combines a sharp
viscous energy estimate with an exact calculation for the azimuthal
mean.  It controls the full range of Fourier frequencies and then
shows that the sign survives passage from the isolated comparison
swirl to the exact finite-energy background.

\paragraph{A viscous transfer to the next larger scale.}
The signed stress has already been accumulated when the first
variation becomes small; its effect is recorded in the second
variation.  Seen from a much larger radius, the lowest nontrivial
spatial-frequency term of this accumulated stress acts like a
\emph{derivative of a point mass}.  More precisely, the relevant
shape is the projected distribution
$\mathcal I=\PP\partial_z(e_z\delta_0)$, where
$e_z=(0,0,1)$ and $\delta_0$ is the \emph{three-dimensional}
Dirac mass at $x=0$.  The isotropic part of the stress disappears
under $\PP\operatorname{div}$ because it is a gradient; the signed
difference between its horizontal and vertical components produces
this prescribed distribution.  Heat evolution at the larger scale
turns it into a smooth perturbation, with a controlled approximation
error.

This transfer requires more than a formal Taylor expansion.
The outgoing perturbation must pass through the previously
active, smaller swirls. We compute their leading interaction
explicitly and retain its contribution to the outgoing amplitude.
To estimate the remaining interactions, we propagate a smooth
divergence-free terminal test backward through the adjoint
linearized evolution. Axisymmetry forces its horizontal components
to vanish linearly at the symmetry axis, providing an additional
scale gain near the concentrated old swirls.

After extracting the leading linear terms, the spatial Taylor
remainder of the adjoint test begins at cubic order.
Incompressibility and the odd axial parity of the incoming
vertical component provide further small factors.
After diffusion and rescaling into
the larger swirl's coordinates, the outgoing perturbation
recovers essentially the same prepared shape, allowing
the cascade to continue.

To describe the leading amplitude, write the incoming difference
as $-af$, where $a>0$ and $f$ is a
normalized smooth perturbation shape close to the prepared one.
The outgoing difference has the form $-a_{\mathrm{out}}f_{\mathrm{out}}$,
with $f_{\mathrm{out}}$ close to the same prepared shape.  At the
level of scalar amplitudes, one step behaves schematically as
\begin{equation}\label{intro:brief-amplitude}
  a_{\mathrm{out}}\approx K_A(f)\,a^2,
  \qquad K_A(f)>0.
\end{equation}
The gain $K_A(f)$ depends on the incoming shape through its signed
stress, as well as on the scale ratio and the heat-prepared response;
it is bounded in the manner needed for every admissible input shape.
The analysis supplies strict error bounds and shows that the gain
can compensate for the scale loss.  Proposition~\ref{nl:transition} gives the exact transition map and
the remainder bounds corresponding to \eqref{intro:brief-amplitude}.
It also controls the output in $H^4$ and in a weighted $L^2$ norm,
so that the next stage receives a shape in the same admissible class.

\paragraph{A finite cascade and its infinite-depth limit.}
For a fixed cascade depth $N$, with $N+1$ retained radii and $N$
transitions, we construct an exact perturbed
Navier--Stokes evolution by choosing \emph{one} initial amplitude
at the smallest radius.  A continuous shooting argument selects
that amplitude so the difference has a prescribed, nonzero
observation at a fixed later time.
The evolution determines all intermediate perturbations.
All analytic parameters and error
margins are fixed before $N$ is increased.

\begin{figure}[!tbp]
\centering
\begin{minipage}[t]{0.485\textwidth}
    \centering
    \includegraphics[width=\linewidth,height=2.12in,
      keepaspectratio]{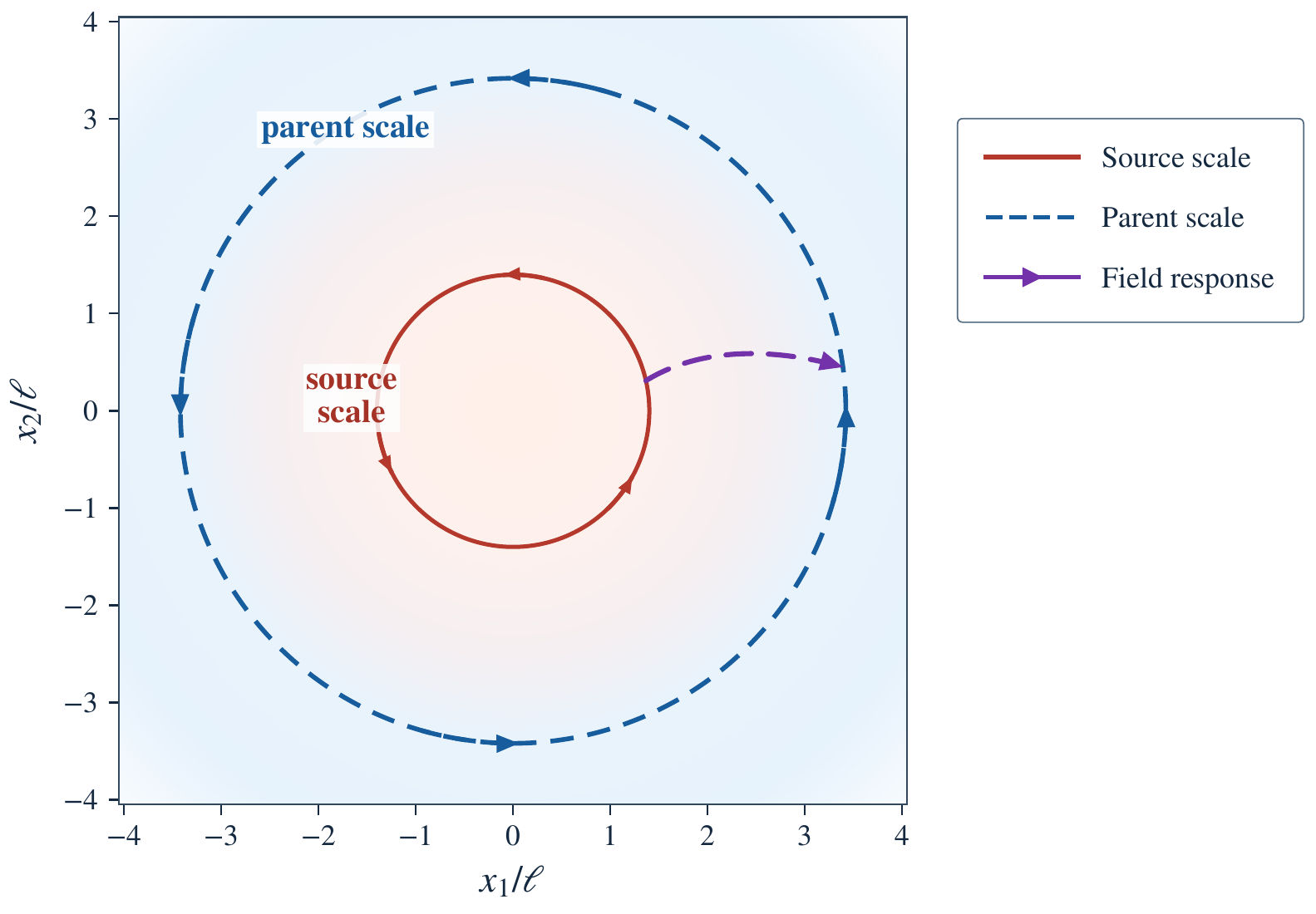}
    \par\smallskip
    {\small (a)}
\end{minipage}\hfill
\begin{minipage}[t]{0.485\textwidth}
    \centering
    \includegraphics[width=\linewidth,height=2.12in,
      keepaspectratio]{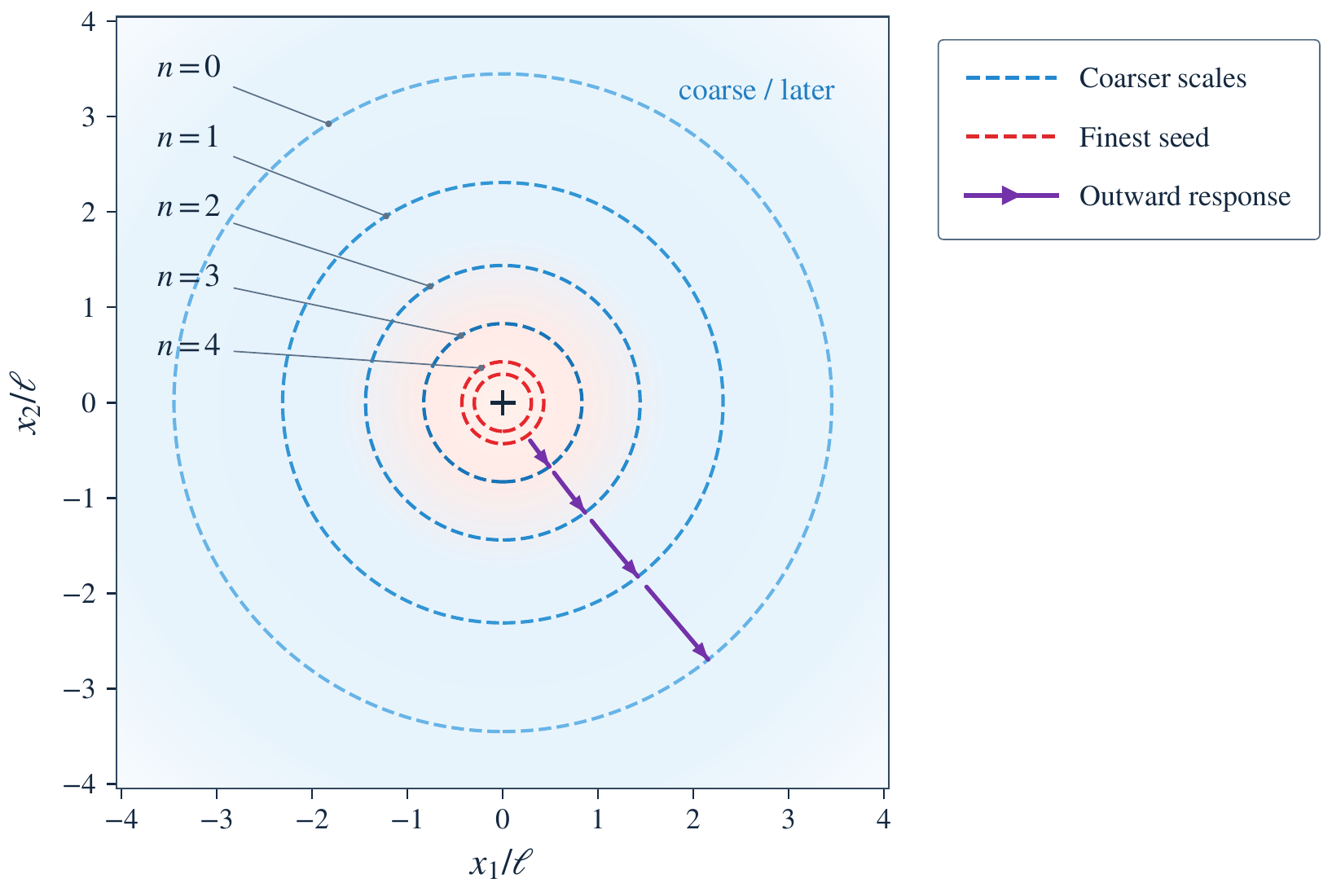}
    \par\smallskip
    {\small (b)}
\end{minipage}

\caption{\textbf{Viscous transfer and the outward cascade.}
Both panels are schematic top-down views in the horizontal
$(x_1,x_2)$-plane.
\textup{(a)} A perturbation at the source scale $\ell_n$
excites an unstable viscous response. The time-integrated
anisotropy of its complete quadratic stress produces,
after Leray projection and diffusion, a prepared perturbation
at the parent scale $\ell_{n-1}=R\ell_n$.
\textup{(b)} Iterating this transition propagates a
discrepancy from the finest level $n=N$ toward the
outermost level $n=0$, with successive transition times
$T_N<T_{N-1}<\cdots<T_0$.
At finite depth, the initial perturbation is selected
to retain a nonzero observation at a fixed outer scale.
The arrows represent the transfer of a discrepancy,
rather than fluid-particle trajectories or an outward
flux of kinetic energy. The geometry is not drawn
to scale.}
\label{fig:intro-transfer-cascade}
\end{figure}

The initial perturbation is a multiple of the regularized impulse
$\mathcal I_{\eta_{\mathrm{reg}}}=e^{\eta_{\mathrm{reg}}\Delta}\mathcal I$, with $\eta_{\mathrm{reg}}>0$ fixed
before the cascade depth. If its radius is $\ell_N=R^{-N}$,
Navier--Stokes scaling gives
\[
 \bigl\|\ell_N^{-1}\mathcal I_{\eta_{\mathrm{reg}}}(\cdot/\ell_N)\bigr\|_{L^2(\R^3)}
       =\ell_N^{1/2}\|\mathcal I_{\eta_{\mathrm{reg}}}\|_{L^2(\R^3)}.
\]
The selected amplitudes are bounded independently of $N$.
Consequently the initial differences converge to zero strongly in
$L^2$ even though their prescribed positive-time observations stay
apart.  Compactness of the finite solutions then produces distinct
weak solutions from the \emph{same} limiting initial velocity.
Uniform energy estimates, compactness, and convergence of the
nonlinear and pressure fluxes give the local energy inequality;
additional positive-time estimates yield smoothness for $t>0$, and
the common strong initial trace lets us extend the energy
\emph{equalities} through $t=0$.  Propositions~\ref{limit:compactness} and \ref{limit:mixed} give
the strong trace, positive-time regularity, and both energy identities
for these limits.
The single-scale transfer and its iteration across the
Gaussian-swirl array are illustrated in
Figure~\ref{fig:intro-transfer-cascade}.

Section~\ref{sec:strategy} states the estimates that connect the
reference construction, linear growth, quadratic stress, and nonlinear
iteration.  The proofs occupy Sections~\ref{sec:background}--\ref{sec:limit}.

\subsection{Relation to earlier work}

The present paper builds on our recent constructions of nonunique unforced two-dimensional Euler evolutions \cite{HQPZ26b,HQPZ26a}. The common mechanism is an outward cascade of a discrepancy from smaller inner scales to larger outer scales. At each finite cutoff, two initial configurations differ only at the finest, innermost scale. This discrepancy vanishes in the limiting initial data, but its effect propagates outward to a fixed larger scale, where the two limiting evolutions remain distinct at positive times. We place this mechanism alongside the instability,
convex-integration, and cascade approaches to nonuniqueness.

\paragraph{Centrifugal instability.}
The frozen instability is of the classical axisymmetric centrifugal
type; see Howard and Gupta \cite{HG62} for the inviscid theory of
swirling flows.  For the unit-amplitude frozen profile
$U_\theta(r)=r e^{-r^2/4}$, the Rayleigh discriminant is
\[
 \frac1{r^3}\frac{\diff}{\diff r}(rU_\theta(r))^2
       =(4-r^2)e^{-r^2/2}.
\]
Its negative is the potential in the variational problem of
Lemma~\ref{seed:rayleigh}.  The cascade requires quantitative control
of the corresponding time-dependent viscous response, its excitation
by the prescribed point input, and the sign of its complete quadratic
stress.  These are the estimates developed in Section~\ref{sec:seed}.

\paragraph{Instability, self-similarity, and forcing.}
Jia and \v{S}ver\'ak constructed forward self-similar solutions for
$(-1)$-homogeneous data \cite{JS14} and formulated spectral conditions
for nonuniqueness, including a localization argument for finite-energy
data \cite{JS15}.  In similarity variables, an unstable trajectory can
approach a reference profile as logarithmic time tends to $-\infty$,
producing a distinct evolution with the same initial trace.  Guillod
and \v{S}ver\'ak obtained numerical evidence of a symmetry-breaking
bifurcation \cite{GS23}.

Vishik's construction for the forced two-dimensional Euler equations
uses a linearly unstable radial vortex with integrable, unbounded
vorticity \cite{Vis18a,Vis18b}.  Albritton, Bru\'e, and Colombo adapted
this instability to an axisymmetric vortex ring and constructed two
distinct suitable Leray--Hopf solutions with zero initial velocity and
the same external force \cite{ABC22}.  Their force belongs to
$L^1(0,T;L^2(\R^3))$, is smooth at positive times, and is singular
at the initial space--time origin.  It makes the prescribed unstable
background an exact solution.  Their subsequent gluing argument gives
forced nonuniqueness in bounded domains and on the torus \cite{ABC23}.
Here the exact reference is obtained by correcting a Gaussian heat
array, and the common initial datum supplies the energy for the cascade.

\paragraph{Unforced Navier--Stokes constructions.}
Hou, Wang, and Yang give a computer-assisted construction of an exact
self-similar profile and an unstable eigenpair, yielding nonunique
suitable Leray--Hopf solutions of the unforced equations \cite{HWY25}.
Their solutions are smooth at positive times, have finite-energy
initial data smooth away from the origin, and belong to the strict
supercritical mixed spaces in \eqref{intro:strict-mixed-spaces}.
The present proof uses analytic estimates for explicit, evolving
Gaussian swirls, followed by signed stress transfer across an infinite
array of spatial scales.  Ionescu, Jia, and Palasek construct numerical
approximate unstable self-similar profiles in the axisymmetric
swirl-free class \cite{IJP26}.  Their result supplies numerical evidence
and a framework for studying nonuniqueness in that class.  The active
swirls considered here have nonzero azimuthal velocity, and the local
instability occurs within the axisymmetric class.

\paragraph{Convex integration and admissibility.}
Buckmaster and Vicol proved nonuniqueness of finite-kinetic-energy weak
Navier--Stokes solutions by convex integration \cite{BV19}.
Buckmaster, Colombo, and Vicol constructed weak solutions that are
smooth outside a set of times of Hausdorff dimension less than one
and can join prescribed smooth flows \cite{BCV22}.  These constructions
do not impose the Leray--Hopf energy inequality and integrated
dissipation bound.  Cheskidov and Luo obtained sharp nonuniqueness
in $L^{\mathfrak p}_tL^\infty_x$ for $\mathfrak p<2$ in periodic dimensions at least
two \cite{CL22}, and nonuniqueness in $C_tL^{\mathfrak q}_x$ for $\mathfrak q<2$ in two
dimensions \cite{CL23}.  The energy inequality, dissipation bound,
and suitability are additional requirements in the present problem.

\paragraph{Cascade models and critical data.}
Tao constructed finite-time blowup for an averaged Navier--Stokes
equation whose modified nonlinearity preserves the energy cancellation
\cite{Tao16}.  For a dyadic Navier--Stokes model, Palasek proved
Leray--Hopf nonuniqueness by coupling finite-dimensional systems across
scales \cite{Pal25}.  For the full equations, Coiculescu and Palasek
constructed nonunique global solutions from $\mathrm{BMO}^{-1}$ data,
smooth at every positive time \cite{CP26}.  Their result complements
the small-data theory of Koch and Tataru \cite{KT01}.
Cheskidov, Dai, and Palasek constructed an inverse frequency cascade
with instantaneous Type~I blowup and nonuniqueness from smooth data
\cite{CDP26}.  Those solutions belong to a class that permits energy
to enter from arbitrarily high frequencies; weak--strong uniqueness
excludes this behavior for Leray--Hopf solutions during the classical
existence interval of the smooth initial datum.

\paragraph{Relation to $2$D Euler nonuniqueness.}
Yudovich's uniqueness theory for bounded, integrable vorticity
\cite{Yud63}, and its extension to certain unbounded vorticities
\cite{Yud95}, provide the classical reference points for the
two-dimensional Euler problem.  Bru\'e, Colombo, and Kumar established
nonuniqueness with vorticity in $L^\infty_tL^{\mathfrak p}_x$ for some $\mathfrak p>1$
by convex integration \cite{BCK26}.

Our two recent papers develop a novel mechanism of propagating
instability across spatial scales.
In \cite{HQPZ26b}, for each
prescribed $1<\mathfrak p<\infty$, we construct two unforced Lagrangian
evolutions from the same initial datum, with vorticities strongly
continuous in $L^1\cap L^{\mathfrak p}$.  Unstable four-vortex configurations
formed by concentrated signed patches separate into dipoles; the
induced strain then excites the next larger configuration.  Relative
patch estimates control errors in proportion to the transmitted
discrepancy.  A finite-cutoff selection and a fixed positive-time
observation yield distinct limits with a common datum.  Independently,
Choi \cite{Choi26} develops a related four-vortex cascade for one
compactly supported, zero-mean datum belonging to every finite $L^{\mathfrak p}$.

In \cite{HQPZ26a}, the reference consists of stationary arrays of
disjoint neutral radial vortices.  Nonuniqueness is driven by the
unstable second angular mode of a two-layer vortex.  Exact shielding
also permits a matching uniqueness theorem, giving a sharp criterion
within the neutral-array class in terms of the growth of the
vorticity $L^{\mathfrak p}$ norms.  
Although the local instabilities and reference
evolutions in these two Euler constructions are therefore different,
both retain a discrepancy at a fixed outer scale as the
initial perturbation recedes to increasingly fine scales. Finally, we refer the readers to two other recent works on $2$D Euler nonuniqueness with different methods \cite{Hexagon2026,Coiculescu2026}. 

The Navier--Stokes construction uses the time-integrated anisotropy
of the complete quadratic stress as the transmitted quantity.
The first variation grows and is then damped by diffusion; its
quadratic response retains the accumulated stress.  After rescaling,
the latter produces the point derivative $\mathcal I$ and evolves
into the next prepared input.  Controlling this passage requires
estimates for diffusion, the nonlocal pressure, and the strong
initial energy trace.  The finite-cutoff selection and the fixed
positive-time observation then provide the common-data limit.

Recent work on Euler blowup and Navier--Stokes blowup with smooth
forcing also uses multiscale amplification
\cite{OpenAIEuler26,OpenAINS26}.  The questions addressed there concern
singularity formation rather than uniqueness in the energy class.

\subsection{Organization and conventions}
The remainder of the paper is organized as follows.
Section~\ref{sec:strategy} develops the proof strategy and defines the
specialized quantities for later use.  Section~\ref{sec:background}
constructs an exact finite-energy reference from Gaussian swirls.
Section~\ref{sec:seed} establishes the isolated instability and the
signed complete stress; Section~\ref{sec:history} transfers those results
to the actual array.  Section~\ref{sec:decay} controls the first
variation after diffusion without losing its spatial moment.
Section~\ref{sec:transfer} proves the parent-scale transfer,
Section~\ref{sec:nonlinear} closes the exact nonlinear iteration, and
Section~\ref{sec:limit} takes the limit and separates the solutions.
Appendix~\ref{app:lorentz} proves the complementary uniqueness
criterion for finite Lorentz second exponents.

Function-space indices are denoted by $\mathfrak p,\mathfrak q,\mathfrak r$ throughout.
We write $x=(x_h,z)\in\R^2\times\R$, $r=|x_h|$, and
$x_h^\perp=(-x_2,x_1,0)$.  The cylindrical unit vectors are
$e_r,e_\theta,e_z$.  The Leray projection is $\PP$,
$\langle x\rangle=(1+|x|^2)^{1/2}$, and
$M_\mu(f)=\|\langle x\rangle^\mu f\|_2$.
We use $(v\otimes w)_{ij}=v_iw_j$ and
$(\operatorname{div}F)_i=\sum_j\partial_jF_{ij}$.
We write $\|f\|_{H^s}=\|(1-\Delta)^{s/2}f\|_2$;
negative Sobolev norms are defined by the same Fourier multiplier on
tempered distributions.
Unless another domain is specified, $\|\cdot\|_2$ denotes
$L^2(\R^3)$; radial and horizontal-fiber norms are specified when introduced.
Complex inner products are linear in the first factor.
The full spatial Fourier transform is
$\widehat f(\xi)=\int_{\R^3}e^{-ix\cdot\xi}f(x)\idiff x$;
the unitary axial transform is denoted separately by $\widehat f^{\,z}$
in \eqref{bg:axial-transform}.  We also use the unnormalized horizontal
transform and its Plancherel identity,
\begin{equation}\label{intro:horizontal-transform}
 \widehat q^{\,h}(\xi_h)=\int_{\R^2}e^{-ix_h\cdot\xi_h}q(x_h)\,\idiff x_h,
 \qquad \|q\|_{L^2_{x_h}}=(2\pi)^{-1}\|\widehat q^{\,h}\|_{L^2_{\xi_h}}.
\end{equation}
Thus $\widehat{\widehat f^{\,z}}^{\,h}=(2\pi)^{-1/2}\widehat f$.
The projection and Riesz-transform conventions are
\begin{equation}\label{intro:operators}
 \widehat{(\PP f)}_i(\xi)
   =\sum_{j=1}^3\left(\delta_{ij}-\frac{\xi_i\xi_j}{|\xi|^2}\right)
       \widehat f_j(\xi),\qquad
 \widehat{\mathcal R_j f}(\xi)=\frac{i\xi_j}{|\xi|}\widehat f(\xi).
\end{equation}
Thus $\PP$ removes gradients and retains the pressure effect of the
nonlinear term.  The matrix is defined away from the single frequency
$\xi=0$, whose value does not affect the $L^2$ operator.
The detailed Fourier and operator conventions are recalled where first
needed.  Generic constants may change from line to line.
Before the final choice of $A$, $CA^P$ denotes polynomial dependence on $A$:
$C$ and $P$ are independent of $A$, the scale cutoffs, and the cascade depth,
with other permitted dependencies stated locally.  After $A$ and the
parameters chosen from it are fixed, $C(A)$ denotes an unrestricted finite
constant independent of the cutoffs and depth.  In particular, a bound
written $C(A)$ is not used as a polynomial prefactor in an exponential
comparison.

\section{Overview of the proof strategy}\label{sec:strategy}
The proof has two quantitative parts.  We first establish growth and
a signed quadratic stress for a prescribed perturbation of one swirl.
We then transfer the resulting field to the next spatial scale, with
enough control of its amplitude, shape, and tails to repeat the
construction.  This section introduces the objects in those estimates
and explains how the exact nonlinear transition yields solutions with
a common initial trace.

\subsection{The reference solution and the transition scales}

The estimates below implement the mechanism of
Section~\ref{intro:mechanism} uniformly in the number of retained radii.
The reference solves the unforced equations exactly.  At each scale,
the incoming field excites a viscous instability, and its complete
quadratic stress has a definite sign.  We control the amplitude and
shape of the resulting field after it has passed through the smaller
swirls, together with all higher nonlinear terms.  These bounds allow
the finite solutions to converge while preserving a fixed positive-time
observation.

{
We first fix the spatial and temporal conventions for one
transition. Let $A\geq1$ be the swirl amplitude, $R>1$ the
ratio of adjacent spatial scales, and $b>1$ the relative
axial width. At level $n$, the active swirl has radius
$\ell_n=R^{-n}$. We call this the source scale; its
immediately larger neighbor, of radius
$\ell_{n-1}=R\ell_n$, is the parent scale.

A transition is analyzed in source-scale variables, in
which the active swirl has radius one. If $t$ and $x$
denote these normalized variables, the corresponding
physical variables are $(\ell_n^2t,\ell_n x)$.

We use three times in source-scale coordinates:
$t_{\mathrm p}$ is the preparation time,
$t_{\mathrm f}$ marks the end of the main source history,
and $t_\ast$ is the time at which the perturbation has
reached the prepared state of the parent scale.
Their values, together with the scale parameters, are
\begin{equation}\label{strategy:brown-scales}
 \begin{gathered}
 \ell_n=R^{-n},\qquad
 R=e^{\varrho A},\qquad b=e^{\beta A},
 \qquad \varrho=\frac{c_\star}{8},\\
 t_{\mathrm p}=\frac{2+\Theta\log A}{A},
 \qquad t_{\mathrm f}=t_{\mathrm p}+16A^2,
 \qquad t_\ast=R^2t_{\mathrm p},\\
 T_n=\ell_n^2t_{\mathrm p},
 \qquad \ell_n^2t_\ast=T_{n-1}.
 \end{gathered}
\end{equation}
Here $c_\star$ is the limiting growth rate defined in
Section~\ref{intro:mechanism}, while $\Theta>0$ and
$\beta>0$ are fixed preparation and width parameters.
The last identity ensures that consecutive transitions
meet at the same physical time. We also fix a growth deficit $\omega>0$ and write
$G=c_\star-\omega$.  The frequency band and the preparation
constant $\Theta$ are selected first, then the width exponent
$\beta$, and finally $A$.  The quantitative restrictions are
specified in \eqref{nl:parameters}; the cascade depth varies only
after these choices have been fixed.
}

\subsection{A single viscous transition: instability, stress, and diffusion}

Recall the finite-energy Gaussian swirls $Q_\ell(0,x)$ in
\eqref{intro:brief-swirl}.  Let $Q(t,x)=\sum_{\ell\in\mathscr A}Q_\ell(t,x)$
be their explicit heat evolution for a finite set of geometrically
separated radii $\mathscr A$.  This is not an exact solution of
Navier--Stokes.  Every ordered pair of swirls contributes to the
radial pressure $P_Q$; when that pressure is included, the remaining
residual is purely axial:
\begin{equation}\label{strategy:reference-residual}
 Q_t-\Delta Q+(Q\cdot\nabla)Q+\nabla P_Q
    =(\partial_zP_Q)e_z,\qquad
 \|\partial_zP_Q(t)\|_2\leq C A^2b^{-1/2}t^{-3/4}.
\end{equation}
The factor $b^{-1/2}$ comes from the large axial width, while
$C$ is independent of the number of radii.
The correction cancels this residual and gives the exact finite-energy solution
$\bar u=Q+e$ of the unforced equations.

This residual is
\emph{not} an external force in the final construction: solving an
unforced perturbation equation removes it and gives the exact
finite-energy reference $\bar u=Q+e$.  

{
For a smooth divergence-free background $V$, define the
linearized operator on solenoidal fields by
\begin{equation}\label{strategy:linearized-generator}
 \mathcal L_Vh
 :=\Delta h-\PP\operatorname{div}
       (V\otimes h+h\otimes V).
\end{equation}
We denote by $\mathcal U_V(t,s)$ the evolution operator
associated with $\partial_t h=\mathcal L_Vh$.
In particular, $\mathcal U_Q$ is the linear evolution
about the Gaussian heat array $Q$.

Set $\bar u=Q+e$. Substituting this decomposition into
the unforced Navier--Stokes equations gives the equation
for the correction $e$, with zero initial correction data
$e(0)=0$:
\begin{equation}\label{strategy:brown-correction}
 \begin{cases}
 (\partial_t-\mathcal L_Q)e
   =F_Q-\PP\operatorname{div}(e\otimes e),\\
 e(0)=0,
 \end{cases}
 \qquad
 F_Q=-\PP\bigl((\partial_zP_Q)e_z\bigr).
\end{equation}
Equivalently, the correction satisfies
\[
 e(t)=\int_0^t\mathcal U_Q(t,s)
    \Bigl(F_Q(s)
      -\PP\operatorname{div}(e\otimes e)(s)\Bigr)
       \idiff s.
\]
The fixed-point argument uses the norm
\begin{equation}\label{strategy:brown-correction-norm}
 \|e\|_{\mathcal X}
 =\sup_{0<t\leq \ell_{\max}^2}
 \left(
 t^{-1/4}\|e(t)\|_2
 +t^{1/4}\|\nabla e(t)\|_2
 +t^{1/2}\|e(t)\|_\infty
 \right),
\end{equation}
where $\ell_{\max}$ is the largest radius in the retained
array. Proposition~\ref{bg:reference} establishes
\(
 \|e\|_{\mathcal X}
 \leq C_\varrho A^P e^{C_\varrho A}b^{-1/2}.
\)
The weights give
$\|e(t)\|_2\leq t^{1/4}\|e\|_{\mathcal X}$ and
$\|e(t)\|_\infty\|\nabla e(t)\|_2
 \leq t^{-3/4}\|e\|_{\mathcal X}^2$.
Thus the quadratic source has an integrable time singularity, and
the factor $b^{-1/2}$ offsets the amplitude-dependent losses.
The bound is uniform in the number of retained scales.

}

The estimate for $\mathcal U_Q$ must retain a time singularity
integrable against the residual in \eqref{strategy:reference-residual}.
Lemma~\ref{bg:propagator} gives this bound, and
Proposition~\ref{bg:reference} uses it to construct $e$ and establish
the scaling covariance of the reference.

In the notation just introduced, the homogeneous linearized equation is
\begin{equation}\label{intro:linearized-general}
 \partial_t v-\Delta v
  +\PP\operatorname{div}(V\otimes v+v\otimes V)=0,
 \qquad \nabla\cdot v=0.
\end{equation}
For the Gaussian heat array, its evolution is $\mathcal U_Q$.

For $0<s\leq t$, the estimate proved in Section~\ref{sec:background}
for this propagator has the form on $L^2_{\sigma,\mathrm{axi}}$, the closed axisymmetric
subspace of $L^2_\sigma$,
\begin{equation}\label{intro:heat-array-propagator}
 \|\mathcal U_Q(t,s)\|_{L^2_{\sigma,\mathrm{axi}}\to L^2_{\sigma,\mathrm{axi}}}
 \leq K_Q\left(\frac ts\right)^{\epsilon_Q},
 \qquad 0<\epsilon_Q<\frac14.
\end{equation}
The constants $K_Q$ and $\epsilon_Q$ are explicitly defined in
Lemma~\ref{bg:propagator}; $K_Q$ can depend substantially on $A$,
but neither constant depends on the number of retained swirls.
The axial width $b$ can be chosen so that $\epsilon_Q$ is as small
as required in the construction.  The bound is not obtained by simply
estimating $\|\nabla Q(t)\|_\infty$: that estimate is of order $A/t$
and would give a much larger, unacceptable power of $t/s$.

For an axially invariant Gaussian swirl, the linearized
convection terms, after applying the Leray projection,
have a special cancellation at zero axial frequency.
Indeed, incompressibility forces the radial component
of an axisymmetric finite-energy Fourier fiber to vanish
at $k=0$, while the remaining radial forcing is a
horizontal gradient and is therefore annihilated by
the Leray projection. Consequently, the linearized
coupling vanishes at $k=0$ and is of order $|k|$
for small axial frequencies. This cancellation,
combined with axial heat diffusion, is the key to
obtaining the refined propagator estimate.  Combining the
$O(|k|)$ coupling with the axial heat factor $e^{-k^2(t-s)}$
before summing the radii gives the bounded frequency multiplier
in Lemma~\ref{bg:propagator}.
For the axially localized array, the remaining error is made small by
taking $b$ large.  Consequently the pressure residual of the heat
array, whose $L^2$ size is bounded by a constant times $t^{-3/4}$,
has an integrable response at $s=0$.  Indeed, since
$\epsilon_Q<1/4$, its Duhamel integral satisfies
\begin{equation}\label{strategy:integrable-response}
 \int_0^t\left(\frac ts\right)^{\epsilon_Q}s^{-3/4}\,\idiff s
  =\frac{t^{1/4}}{1/4-\epsilon_Q},
\end{equation}
which yields the strong zero initial trace of the
linear response. This is the reason the correction $e(t)$ tends to zero in $L^2$ at the
initial time while the exact reference retains its finite-energy datum.
The calculation uses the $L^2$ propagator bound; the remaining
regularity and contraction estimates are proved in Section~\ref{sec:background}.

\paragraph{Instability of the isolated swirl.}
We now examine \emph{one} interaction in coordinates where the active
swirl has radial scale one.  Temporarily remove the axial cutoff and
consider
\[
 U_{\mathrm d}(t,r)=A W_{1+t}(r)e_\theta,
 \qquad W_T(r)=rT^{-2}e^{-r^2/(4T)}, \qquad T=1+t.
\]
This flow is independent of $z$ and is an exact viscous
Navier--Stokes solution, with a radial pressure balancing its
centrifugal force.  It has infinite three-dimensional energy because
of its axial invariance, so it is \emph{not} one of the finite-energy
solutions in our theorem.  Its role is to make the local linearized
interaction explicit; Sections~\ref{sec:history} and
\ref{sec:transfer} then compare this model with the actual reference.

{
We next specify precisely the operator whose eigenvalue will be
used below.  Linearize about $U_{\mathrm d}$ using
\eqref{intro:linearized-general}.
We take the unitary Fourier transform in the axial
variable $z$, with frequency $k>0$. An axisymmetric
divergence-free velocity fiber can be represented by
\[
 \widehat v^{\,z}(r,k)
 =\bigl(kf_r(r),\,f_\theta(r),\,iD_rf_r(r)\bigr),
 \qquad D_rf_r=f_r'+r^{-1}f_r.
\]
The relation between the first and third components
expresses incompressibility. Thus the pair
$(f_r,f_\theta)$ determines the full velocity fiber.
Let $L_k$ be the positive radial operator
\[
 L_k=-\partial_r^2-r^{-1}\partial_r+r^{-2}+k^2
\]
with its Friedrichs realization in
$L^2((0,\infty),r\,\idiff r)$.
The projected inviscid coupling, frozen at $t=0$
and divided by the swirl amplitude $A$, is
\begin{equation}\label{intro:frozen-coupling}
 \mathcal C_k(0)
 =\begin{pmatrix}
 0&2kL_k^{-1}\Omega_1\\
 -kB_1&0
 \end{pmatrix},
\end{equation}
where
$\Omega_1(r)=e^{-r^2/4}$ and
$B_1(r)=(2-r^2/2)e^{-r^2/4}$.
The matrix
$\mathcal C_k(0)$ acts on the pair $(f_r,f_\theta)$; the full
\emph{viscous} equation has generator
$A\mathcal C_k(t)-\operatorname{diag}(L_k,L_k)$.
The positive eigenvalue \(\sigma_k\) of \(\mathcal C_k(0)\) identifies the growing direction; Lemma~\ref{seed:metric} controls its evolution in the presence of viscosity.
}

An eigenvector $(h,f_\theta)$ with positive eigenvalue $\sigma$
satisfies $\sigma f_\theta=-kB_1h$.  Eliminating $f_\theta$ gives
\begin{equation}\label{intro:radial-eigenvalue}
 \sigma^2L_kh=-k^2\Phi(r)h,
 \qquad \Phi(r)=(4-r^2)e^{-r^2/2}.
\end{equation}
This scalar form is the reason the unstable eigenvalue can be obtained
analytically.  On the regular radial form domain, the \emph{largest
positive} eigenvalue $\sigma_k$ is characterized by the Rayleigh
quotient
\begin{equation}\label{intro:rayleigh-explained}
 \sigma_k^2=\sup_{h\ne0}
 \frac{\displaystyle\int_0^\infty
       (r^2-4)e^{-r^2/2}|h(r)|^2r\,\idiff r}
      {\displaystyle\int_0^\infty\left(|h(r)|^2
          +k^{-2}|h'(r)|^2+k^{-2}r^{-2}|h(r)|^2\right)
                           r\,\idiff r}.
\end{equation}
The supremum is positive because the numerator is positive for
functions supported in $r>2$.  The potential in the numerator is
largest at $r^2=6$, where its value is $2e^{-3}$.  For every finite
$k$, the derivative terms in the denominator impose a positive
penalty, so $\sigma_k^2<2e^{-3}$.  Conversely, choose a smooth test
function concentrated near $r=\sqrt6$.  As $k\to\infty$, the
$k^{-2}$ derivative penalty for this \emph{fixed test function}
tends to zero.  By choosing the support increasingly close to
$r=\sqrt6$ only after that limit, we obtain the matching lower bound.
Therefore
\[
 0<\sigma_k<c_\star,\qquad
 \sigma_k\nearrow c_\star=\sqrt2\,e^{-3/2}
 \quad\text{as }k\to\infty.
\]
Lemma~\ref{seed:rayleigh} proves the existence, simplicity, and
monotonicity of this eigenvalue, together with a quantitative trial
function estimate.  The associated instability is related to the
centrifugal sign of the swirl: with $T=1+t$,
\begin{equation}\label{intro:centrifugal-discriminant}
 \frac1{r^3}\partial_r(rU_{{\mathrm d},\theta})^2
 =A^2T^{-4}\left(4-\frac{r^2}{T}\right)e^{-r^2/(2T)}.
\end{equation}
The sign changes at $r=2\sqrt T$; the negative part of this
quantity gives the same radial potential as in
\eqref{intro:radial-eigenvalue}.  Thus the negative centrifugal discriminant identifies the
annulus responsible for the positive Rayleigh quotient.

The limit $k\to\infty$ concerns only the \emph{frozen
inviscid} operator.  For the full viscous equation, arbitrarily high
$k$ are strongly damped by axial diffusion.  In the proof we first
choose a \emph{fixed} frequency interval on which $\sigma_k$ is close
to $c_\star$, and only then increase $A$.  The resulting viscous
first variation can grow at a rate close to $e^{c_\star A}$.

{

The leading exponent can be read directly from the dilation of the
frozen operator. Formula~\eqref{seed:growing-rate} gives
\begin{equation}\label{strategy:brown-growth-clock}
 \lambda_k(t)=(1+t)^{-2}\sigma_{k\sqrt{1+t}},\qquad
 A\int_\tp^t\lambda_k(s)\idiff s
 \leq c_\star A\left(\frac1{1+\tp}-\frac1{1+t}\right)
 =\Gamma_{\mathrm{src}}(t).
\end{equation}
The last quantity is the growth clock in \eqref{nl:growth}.
The lower estimate on the chosen frequency band and the diffusive
complementary estimate are what justify retaining growth close to
this upper exponent for the actual viscous response.
}

A sharp physical energy estimate bounds its growth from above at the
same leading rate.  Matching these two rates leaves the narrow but
essential exponential margins needed to carry errors through the
cascade.  Following the growing eigenspace in time also requires
controlling the other spectral directions while preserving diffusion;
that is the purpose of the complementary metric in
Lemma~\ref{seed:metric}.

To apply the spectral growth estimate, we must compute the projection
of the input delivered by the preceding scale onto the growing
eigenspace.  This input is the point derivative defined below.

Let $\delta_0=\delta_{(0,0,0)}$ denote the \emph{three-dimensional}
Dirac mass at the spatial origin of $\R^3$.
Let $\mathcal U_{\mathrm d}(t,s)$ denote the evolution of
\eqref{intro:linearized-general} about the isolated swirl
$U_{\mathrm d}$, and set
\begin{equation}\label{intro:point-input}
 \mathcal I=\PP\partial_z(e_z\delta_0),\qquad
 S_A=\mathcal U_{\mathrm d}(\tp,0)\mathcal I,\qquad
 n_A=\|S_A\|_2,\qquad g=\frac{S_A}{n_A}.
\end{equation}
The derivative is distributional: for a smooth compactly supported
vector test $\phi$, one has
$\langle\partial_z(e_z\delta_0),\phi\rangle
=-\partial_z\phi_z(0)$.  The projection $\PP$ makes this
vector distribution divergence free and also gives it a nonlocal
pressure-related tail.  In particular, $\mathcal I$ is \emph{not}
a compactly supported velocity field.
After the unitary Fourier transform in $z$, the three-dimensional
Dirac mass gives $(2\pi)^{-1/2}\delta_{0,h}$ on each axial fiber,
where $\delta_{0,h}$ is the \emph{two-dimensional} Dirac mass at
$x_h=0$.  Thus
\[
 \widehat{\mathcal I}^{\,z}(k)
  =\frac{ik}{\sqrt{2\pi}}\PP_k(e_z\delta_{0,h}).
\]
Here $\PP_k$ is the Fourier-fiber version of the same Leray
projection.  These different Dirac masses and the factor $ik$ are
important in the point-input pairing of Lemma~\ref{seed:point-growth}.

Since $\mathcal I$
is a distribution, we first evolve it to the \emph{positive}
preparation time
\[
 \tp=\frac{2+\Theta\log A}{A}>0,
\]
where $\Theta>0$ is a fixed preparation constant chosen to obtain the
needed regularity and alignment.  Heat smoothing then gives the smooth,
finite-energy profile $S_A$, and we normalize it to $g$.
The proof shows that the pairing of the original distribution with
the growing adjoint eigenfunction is \emph{strictly positive}; in
particular the contribution of the horizontal Dirac mass at the axis
must be retained.  The short preparation interval controls the
remaining spectral components without changing the leading growth.
The subsequent nonlinear calculation uses the \emph{entire} prepared response $S_A$.

\paragraph{The complete quadratic stress and its sign.}
Consider for the moment an exact reference solution $U$ and a small
initial perturbation.  The formal expansion
\[
 u^\varepsilon=U+\varepsilon v+\varepsilon^2w+O(\varepsilon^3)
\]
shows the distinction between a first-order perturbation $v$ and the
second-order response $w$.  The field $v$ solves the homogeneous
linearized equation \eqref{intro:linearized-general}, and $w$,
initialized at zero at the preparation time, solves
\[
 (\partial_t-\mathcal L_U)w
   =-\PP\operatorname{div}(v\otimes v),\qquad
 \mathcal L_U w
 =\Delta w-\PP\operatorname{div}(U\otimes w+w\otimes U).
\]
There is no factorial in the definition of $w$.  If the linear terms
involving $U$ are moved to the source side of the \emph{free heat}
equation, then $w$ is driven by
\begin{equation}\label{intro:complete-stress}
 \mathsf S=v\otimes v+U\otimes w+w\otimes U,
 \qquad (\partial_t-\Delta)w=-\PP\operatorname{div}\mathsf S.
\end{equation}
This identity explains the word \emph{complete}.  Keeping only
$v\otimes v$ would omit the effects of the second-order response
interacting with the reference swirl; those effects are of the same
relevant order and influence the sign of the transfer.

The extra terms have a concrete origin.  Quadratic interactions of the
first variation produce an azimuthal component of $w$ with a
\emph{nonzero axial mean} $\int_\R w_\theta(r,z)\,\idiff z$.
This mean continues to diffuse after the oscillatory part of $v$
has been damped, and it still interacts with the background.
Lemma~\ref{seed:stress-history} computes this mean and its
contribution to the time-integrated stress for the full Fourier
response.

For a sufficiently localized axisymmetric variation, rotational
symmetry forces the spatial integral of the symmetric stress to have
the form
\[
 \int_{\R^3}\mathsf S(t,x)\,\idiff x
  =\operatorname{diag}\bigl(s_h(t),
                              s_h(t),s_z(t)\bigr).
\]
As in \eqref{intro:brief-mass}, the scalar that controls the outgoing perturbation is the signed \emph{time-integrated
anisotropy}
\[
 M=-\int_{\tp}^{\tf}
        \bigl(s_h(t)-s_z(t)\bigr)\,\idiff t,
 \qquad \tf=\tp+16A^2.
\]
{

The normalization is important: the exponential lower bound quoted
below is for the $L^2$-normalized prepared shape $g$ and its admissible
neighbors, as in \eqref{nl:mass}; the unnormalized point response
has its corresponding mass before division by $\|S_A\|_2^2$.
The sign identity can also be made explicit. Since
\begin{equation}\label{strategy:brown-stress-decomposition}
 \int\mathsf S\,\idiff x
   =s_hI_3+(s_z-s_h)e_z\otimes e_z,
 \qquad
 -\int_\tp^\tf\PP\operatorname{div}
       \left(\Big(\int\mathsf S\,\idiff x\Big)\delta_0\right)\idiff t
   =-M\mathcal I,
\end{equation}
positive $M$ sends a negative multiple of the same prescribed input
to the next scale. This fixes the sign in the convention $h=-af$.
}

The first variation initially grows, but the relatively long interval
up to $\tf$ allows axial diffusion to damp it.  

{
For an incoming shape $f$ in the admissible class $\mathcal C_A$
defined in \eqref{strategy:brown-shape-class} below, let $v_f,w_f$
be the first and second variations about the exact reference, with
$v_f(\tp)=f$ and $w_f(\tp)=0$.  Denote their positive complete-stress
mass by $\mathcal M_f$.}
{

For the actual first variation $v_f$, the required endpoint estimate
is stronger than a statement that oscillatory modes decay:
\begin{equation}\label{strategy:brown-linear-endpoint}
 \begin{aligned}
 M_{s_\star}(v_f(\tf))
     &\leq CA^Pe^{2\varrho A}M_{s_\star}(f),\\
 \|v_f(\tf)\|_{H^4}
     &\leq CA^P\|f\|_2
       +CA^Pe^{K_{s_\star}A}b^{-1/2}M_{s_\star}(f),
       \qquad s_\star=19/8.
 \end{aligned}
\end{equation}
These are \eqref{nl:linear-decay}; after the width choice, the
unweighted endpoint norm is only polynomial on the admissible class.
}

The stress created
during its growth is nevertheless retained in the second variation.
Sections~\ref{sec:seed} and \ref{sec:history} prove that the full
history has the \emph{correct sign} $M>0$ and, for the prepared
response, satisfies
\[
 A^{-P}e^{2GA}\lesssim M\lesssim A^Pe^{2c_\star A}
\]
for a fixed $G<c_\star$ close to $c_\star$.  The lower bound uses a
Gaussian-weighted sign estimate, including the azimuthal mean;
low axial frequencies are estimated separately rather than discarded.
The next section of the proof transfers these statements from the
isolated comparison flow to the exact localized reference.

\paragraph{Diffusive transfer to the next scale.}
The significance of $M$ can be seen directly from the small-frequency
expansion of the heat response to the stress.  Write
$\int\mathsf S=\operatorname{diag}(s_h,s_h,s_z)$.
Its isotropic part is invisible to $\PP\operatorname{div}$, since the
Leray projection annihilates gradients.  The remaining anisotropy,
when multiplied by $-\PP\operatorname{div}$, produces a scalar
multiple of $\PP\partial_z(e_z\delta_0)=\mathcal I$.
Indeed the Fourier-symbol identity
\[
 -i\PP(\xi)\operatorname{diag}(1,1,0)\xi
    =i\xi_z\PP(\xi)e_z
\]
shows the precise sign.  Let $m=-M<0$ and let $R\gg1$ be the
ratio between the active radius and its larger \emph{parent} radius.
After rescaling the integrated second response $w$ into parent
coordinates, its leading, distributional approximation is
\[
 \frac{R^4}{m}\,w(\tf,R\,\cdot)\ \simeq\ \mathcal I.
\]

{

For a shape $f$, this statement has the following precise fiberwise
meaning.  Write $y=x/R$ for parent coordinates; parent axial frequency
$k$ corresponds to source frequency $k/R$.  Put $m_f=-\mathcal M_f$, let
$F_{f,R}(y)=R^4w_f(\tf,Ry)/m_f$, and denote its unitary axial
Fourier transform in $y_3$ by $F_{f,k}$.  Proposition~\ref{hist:generic} gives
\begin{equation}\label{strategy:brown-point-error}
 \begin{gathered}
 \mathcal I_k=\frac{ik}{\sqrt{2\pi}}\PP_k(e_z\delta_{0,h}),\\
 \|F_{f,k}-\mathcal I_k\|_{H^{-17/4}_{y_h}}
   \leq CA^P(1+|k|)^{11/4}R^{-7/4}e^{2\omega A},
       \qquad k\in\R.
 \end{gathered}
\end{equation}
This is \eqref{hist:point}. Since $R=e^{\varrho A}$ and
$2\omega<\varrho/32$, the error decays exponentially for each
bounded frequency range. The approximation is made before parent
heat smoothing, in a negative horizontal Sobolev norm.
}

The factor $R^4$ is the scaling
of a first derivative of a three-dimensional point mass.  In contrast,
velocity rescaling between Navier--Stokes solutions uses only a
factor $R$; the difference in these factors is responsible for the
$R^{-3}$ in the nonlinear amplitude law below.

During transfer, the second variation interacts with the smaller
swirls that are still present.  Some of these contributions have a
large norm even though their effect on the smooth parent-scale
observation is small.  Section~\ref{sec:transfer} first computes the
leading old-column interaction and then estimates the remaining
ones by testing against a backward-evolved smooth function.  The
leading correction changes the outgoing coefficient by a real factor
$a_f=1+O(A^{-1})>0$, which must be retained in the amplitude.

{

After propagation through the old columns and parent smoothing,
\eqref{nl:quadratic-transfer} becomes
\begin{equation}\label{strategy:brown-quadratic-transfer}
 \left\|R\bigl(\mathcal U_{\bar u}(\ts,\tf)w_f(\tf)\bigr)(R\,\cdot)
       +a_f\mathcal M_fR^{-3}S_A\right\|_{H^4}
   \leq CA^P\mathcal M_fR^{-3}e^{-\kappa_{\mathrm r}A},
 \qquad 0<\kappa_{\mathrm r}<\varrho.
\end{equation}
The factor $a_f$ therefore belongs to the leading scalar coefficient;
the residual shape error in \eqref{strategy:brown-quadratic-transfer}
is measured after this coefficient has been extracted.  This formula
propagates only $w_f(\tf)$: the full second variation continues to
receive quadratic forcing on $[\tf,\ts]$, as controlled in
Lemma~\ref{trans:late} and Proposition~\ref{trans:quadratic}.
Higher-order terms of the exact nonlinear difference are treated in
Section~\ref{sec:nonlinear}.
}

After time for diffusion and preparation in the parent coordinates,
the output again resembles the same smooth profile $g$ that was
introduced in \eqref{intro:point-input}.  The parent observation time
in the source-scale coordinates is
$\ts=R^2\tp$; consequently $\ts/R^2=\tp$ exactly, so the next
transition begins at the correct prepared heat age.  The leading
amplitude law is best described as follows.  If the incoming
perturbation at time $\tp$ is $-af$, where $a>0$ is the amplitude and
$f$ is a normalized shape close to $g$, then
\[
 R\bigl(u(\ts,R\,\cdot)-\bar u(\ts,R\,\cdot)\bigr)
 \simeq-a_f\,\|S_A\|_2\,\mathcal M_f R^{-3}a^2g.
\]
Here $\mathcal M_f$ is the positive complete-stress mass generated by that shape.
This formula is a guide to the mechanism: the propositions in
Sections~\ref{sec:transfer} and \ref{sec:nonlinear} state the actual
quantitative remainder bounds needed for iteration. Lemmas
\ref{trans:invariant-entrance} and \ref{trans:actual-comparison} first
control the entrance error in $L^2$. The parent-scale estimate in
Lemma~\ref{trans:smoothing} then gives the required $H^4$ bound.

\subsection{From one transition to an exact nonlinear cascade}

If the outgoing field had the correct scalar amplitude but the wrong
spatial shape, the next swirl might not amplify it.  Thus the output
must return to a controlled neighborhood of the prepared profile
$g=S_A/\|S_A\|_2$.  We write the incoming perturbation as $-af$,
where $a>0$ is the scalar amplitude and $f$ is normalized by
$\langle f,g\rangle_{L^2}=1$.  This fixes what we mean by the
amplitude and keeps errors parallel to $g$ separate from errors in
shape.

{

For the real, centrally odd, axisymmetric, divergence-free fields,
we control the shape in $H^4$ and its spatial tail by
$M_{s_\star}$ from \eqref{bg:moment}, where $s_\star=19/8$.
The resulting class is the affine slice
\begin{equation}\label{strategy:brown-shape-class}
 \mathcal C_A=\left\{g+\zeta:
   \langle\zeta,g\rangle=0,\quad
   \|\zeta\|_{H^4}\leq e^{-3\varrho A/8},\quad
   M_{19/8}(\zeta)\leq e^{\varrho A/4}\right\}.
\end{equation}
See \eqref{nl:shape-class}. In particular $\langle f,g\rangle=1$
does not require $\|f\|_2=1$.
}

Two kinds of estimates are needed.  The $H^4$ norm measures the
\emph{strong shape error} and controls the nonlinear products.  The
ordinary spatial moment
\[
 M_s(h)=\bigl\|\langle x\rangle^s h\bigr\|_{L^2(\R^3)},
 \qquad \langle x\rangle=(1+|x|^2)^{1/2},
\]
measures how much of a perturbation is spread far from its core.  We
use $s=s_\star=19/8$.  The moment is needed to justify the stress
integrals, control the nonlocal pressure tails, and approximate the
outgoing field by the point distribution.  Its tolerance is much
larger than the allowed $H^4$ shape error: a small $L^2$ perturbation
can be spread across a very large region and have a large moment.
Moreover the Leray projection can create algebraic velocity tails, so
we cannot assume that every admissible perturbation has Gaussian decay.

The key nonlinear estimates are therefore \emph{tame}: one places the
spatial weight on a well-controlled factor in each product and uses a
small unweighted norm for the other factor.  This prevents a large
incoming moment from being squared at every transition.  

{

The full perturbation is constructed, on the source interval, by
the convergent series
\begin{equation}\label{strategy:brown-nonlinear-series}
 \begin{gathered}
 h(\tp)=-af,\qquad h(t)=\sum_{m\geq1}(-a)^mV_m(t),\\
 V_1=v_f,\quad V_2=w_f,\quad V_1(\tp)=f,\quad V_m(\tp)=0\ (m\geq2),\\
 (\partial_t-\mathcal L_{\bar u})V_m
   =-\PP\operatorname{div}\sum_{j=1}^{m-1}V_j\otimes V_{m-j}\quad(m\geq2).
 \end{gathered}
\end{equation}
This is the recurrence \eqref{nl:series}.  The degree-independent,
unweighted $H^4$ source estimate stated in the proof of
Proposition~\ref{nl:transition} makes the majorant series converge;
the weighted estimate \eqref{nl:early-source} separately prevents the
incoming moment tolerance from being squared.
}

Diffusive smoothing near the end of the transition upgrades the
outgoing $L^2$ approximation to the $H^4$ control needed for the next
transition.  Proposition~\ref{nl:transition} establishes the resulting
invariant \emph{amplitude--shape class} for the exact nonlinear
Navier--Stokes equation.

\subsection{The exact transition, finite shooting, and limiting solutions}
The nonlinear transition is continuous in its input amplitude and shape.
It has strict crossings at the two amplitude endpoints, controls the
solution at intermediate times, and uses a reference that scales into
the next stage.  We now explain how these properties select a finite
cascade and preserve its terminal coordinate in the limit.

{
Write $\bar u^{I,J}$ for the exact reference associated with the radii $R^{-I},\ldots,1,\ldots,R^J$; thus $I$ and $J$ count the retained smaller and larger scales.}
In coordinates in which the active swirl has radius one, the input is
$\bar u^{I,J}(\tp)-af$.  Here $a>0$ is a scalar and $f$ is a normalized
shape close to $g$, with an ordinary spatial moment; the exact class is
\eqref{nl:shape-class}.  Proposition~\ref{nl:transition} constructs the
strong solution through $\ts=R^2\tp$ and writes its rescaled difference as
\[
 H(y):=R\bigl(u(\ts,Ry)-\bar u^{I,J}(\ts,Ry)\bigr)
       =-\mathcal A^{I,J}(a,f)\,\mathcal F^{I,J}(a,f)(y).
\]

{

The scalar and shape are defined from this exact output $H$ by
\begin{equation}\label{strategy:brown-exact-map}
 \mathcal A^{I,J}=-\langle H,g\rangle>0,\qquad
 \mathcal F^{I,J}=-H/\mathcal A^{I,J},\qquad
 \langle\mathcal F^{I,J},g\rangle=1.
\end{equation}
The three conclusions of Proposition~\ref{nl:transition} are
continuity, $\mathcal F^{I,J}(a,f)\in\mathcal C_A$, and
\begin{equation}\label{strategy:brown-crossings}
 \begin{gathered}
 a_-=e^{-(13\varrho+\varrho/32)A},\qquad
 a_+=e^{-12\varrho A},\\
 A^{-P}e^{(2G-3\varrho)A}a^2
    \leq\mathcal A^{I,J}(a,f)
    \leq A^Pe^{(2c_\star-3\varrho)A}a^2,\\
 \mathcal A^{I,J}(a_-,f)<a_-/2,\qquad
 \mathcal A^{I,J}(a_+,f)>2a_+.
 \end{gathered}
\end{equation}
For example, $2c_\star=16\varrho$ makes the upper bound for
$\mathcal A(a_-,f)/a_-$ equal to
$A^Pe^{-\varrho A/32}$, whereas the lower bound for
$\mathcal A(a_+,f)/a_+$ is
$A^{-P}e^{(\varrho-2\omega)A}$. Thus the crossings follow by
taking $A$ large, uniformly in the shape and the cutoffs.
}

The intermediate-time bound \eqref{nl:interstage} is equally important:
control at observation times alone would not suffice for the limiting PDE.
The reference covariance \eqref{bg:covariance} identifies the output
background with $\bar u^{I+1,J-1}(\tp)$.
The leading scalar law is
\[
 \mathcal A^{I,J}(a,f)
   \simeq a_f n_A \mathcal M_fR^{-3}a^2,
 \qquad a_f=1+O(A^{-1})>0.
\]
The symbol $\simeq$ in this guide summarizes the quantitative comparisons
in Sections~\ref{sec:transfer} and \ref{sec:nonlinear}.
The positive stress mass $\mathcal M_f$ comes from the
\emph{complete} stress.  Its large size compensates for the factor $R^{-3}$.

\paragraph{Finite-depth shooting.}
Choose a terminal value $a_*\in(a_-,a_+)$ and a depth $N$.  The smooth
preparation in Proposition~\ref{nl:preparation}, followed by the continuous
selection in Lemma~\ref{nl:shooting}, produces $N$ genuine transitions
ending at $a_*$.  The selection varies one finest-scale initial parameter.
It uses the endpoint crossings rather than monotonicity of the full map,
and requires neither a fixed point nor compactness of the shape class.
For the initial preparation we fix
$\eta_{\mathrm{reg}}=e^{-\varrho A}$ and set
$\mathcal I_{\eta_{\mathrm{reg}}}=e^{\eta_{\mathrm{reg}}\Delta}\mathcal I$.
This smooth, divergence-free field is used only at the finest scale.
Proposition~\ref{nl:preparation} controls its evolution from time
zero to $\tp$, with $\eta_{\mathrm{reg}}$ fixed before the depth $N$ is varied.

{
For the array with radii $1,R,\ldots,R^N$, let $B_N(\lambda)$ and $F_N(\lambda)$ denote the amplitude and normalized shape at $t_{\mathrm p}$ obtained from the initial perturbation $-\lambda\mathcal I_{\eta_{\mathrm{reg}}}/n_A$. The shooting argument selects a parameter $\lambda_N$.}
{

Concretely, the orbit selected in \eqref{nl:orbit} satisfies
\begin{equation}\label{strategy:brown-finite-orbit}
 \begin{gathered}
 a_0=B_N(\lambda_N),\quad f_0=F_N(\lambda_N),\qquad
 a_-/2<\lambda_N<2a_+,\\
 (a_{j+1},f_{j+1})=
    (\mathcal A^{j,N-j}(a_j,f_j),\mathcal F^{j,N-j}(a_j,f_j)),
       \quad 0\leq j<N,\qquad a_N=a_*.
 \end{gathered}
\end{equation}
The clipping argument in Lemma~\ref{nl:shooting} gives a continuous
terminal-amplitude path from the lower endpoint to the upper endpoint.
Once it hits either clipped boundary it stays there, so a parameter
with an interior terminal value has undergone only genuine transitions.
}

Set $\ell_n=R^{-n}$ and $T_n=\ell_n^2\tp$.  At the finest radius,
Proposition~\ref{nl:preparation} is started at physical time zero; subsequent
stages join at the times $T_N,T_{N-1},\ldots,T_0$.  Navier--Stokes scaling
keeps the viscosity equal to one throughout.
{
Let $U_N$ be the physical solution obtained from the selected depth-$N$ cascade,
and let $\bar u_N$ be its corresponding reference solution.}
The initial difference from the finite reference is exactly
\[
 U_N(0,x)-\bar u_N(0,x)
  =-\frac{\lambda_N}{n_A}\ell_N^{-1}\mathcal I_{\eta_{\mathrm{reg}}}(x/\ell_N),
 \qquad
 \|U_N(0)-\bar u_N(0)\|_2
  =\frac{|\lambda_N|}{n_A}\ell_N^{1/2}\|\mathcal I_{\eta_{\mathrm{reg}}}\|_2\longrightarrow0.
\]
The convergence follows from the depth-independent bound on
$\lambda_N$ and the fixed value of $\eta_{\rm reg}$.
The finite reference data
converge strongly in $L^2$ to the explicit Gaussian series
\eqref{bg:datum}.  At the largest observation time, by contrast,
\[
 \langle U_N(T_0)-\bar u_N(T_0),-g\rangle=a_*.
\]
The task of Section~\ref{sec:limit} is to preserve both statements in the
same subsequential limit.

\paragraph{Compactness and the energy trace.}
Each finite evolution satisfies the unforced energy identity.  Its initial
energy is uniformly bounded, so there are uniform energy and dissipation
bounds.  Local compactness alone gives strong convergence only on
bounded spatial domains; kinetic energy could otherwise
escape to spatial infinity as the cascade depth increases.
To exclude this possibility, we use the local energy
identity with an exterior spatial cutoff to show that
the kinetic energy outside a sufficiently large ball is
uniformly small, independently of the cascade depth
and time. This spatial tightness, combined with local
compactness, yields strong convergence in
$L^2((0,T_0)\times\R^3)$. Uniform weak convergence in time preserves the
fixed observation, including its endpoint value.  Strong convergence of
the initial data and the energy inequality give the strong common initial
trace.

{

The estimates and convergence modes needed here are, explicitly,
\begin{equation}\label{strategy:brown-compactness}
 \begin{gathered}
 \sup_N\left(\sup_{0\leq t\leq T_0}\|U_N(t)\|_2^2
       +2\int_0^{T_0}\|\nabla U_N(t)\|_2^2\idiff t\right)<\infty,\\
 \lim_{L\to\infty}\sup_N\sup_{0\leq t\leq T_0}
       \int_{|x|>L}|U_N(t,x)|^2\idiff x=0,\\
 U_N\longrightarrow U\quad\text{strongly in }
       L^2((0,T_0)\times\R^3),\qquad
 U_N\longrightarrow U\quad\text{in }C_{\mathrm w}([0,T_0];L^2).
 \end{gathered}
\end{equation}
The last convergence means uniform convergence of pairings against
each fixed $L^2$ test. In particular, since $\bar u_N\to\bar u$
in $C([0,T_0];L^2)$, the terminal coordinate survives:
\begin{equation}\label{strategy:brown-surviving-coordinate}
 \langle U(T_0)-\bar u(T_0),-g\rangle
  =\lim_{N\to\infty}\langle U_N(T_0)-\bar u_N(T_0),-g\rangle
  =a_*>0.
\end{equation}
These are the separate assertions of
Proposition~\ref{limit:compactness} and \eqref{limit:terminal-coordinate}.
}

The bounds between adjacent scales also yield, for the limiting difference,
$\|U(t)-\bar u(t)\|_2\leq C_At^{1/4}$ and
$\|U(t)-\bar u(t)\|_\infty\leq C_At^{-1/2}$. 

{

Interpolation for $h=U-\bar u$ gives the calculation
\begin{equation}\label{strategy:brown-mixed-norms}
 \begin{aligned}
 \|h(t)\|_{\mathfrak q}
 &\leq\|h(t)\|_2^{2/\mathfrak q}\|h(t)\|_\infty^{1-2/\mathfrak q}
   \leq C_{A,\mathfrak q}t^{\frac{3}{2\mathfrak q}-\frac12},
       &&2\leq \mathfrak q\leq\infty,\\
 \int_0^{T_0}t^{-\frac{\mathfrak p}{2}(1-3/\mathfrak q)}\idiff t<\infty
 &\quad\Longleftrightarrow\quad
       \frac2{\mathfrak p}+\frac3{\mathfrak q}>1,
       &&\mathfrak q>3,\quad \mathfrak p<\infty.
 \end{aligned}
\end{equation}
For the full velocity, the reference bounds are bounded for $\mathfrak q<3$,
logarithmic at $\mathfrak q=3$, and of the same power for $\mathfrak q>3$; see
\eqref{limit:mixed-bounds}. This proves the full strict mixed-norm
range, including $\mathfrak p=\infty$ when $\mathfrak q<3$, for the same limiting family.
Together these estimates also give positive-time regularity.
{
The weak-$L^3$ bound for the Gaussian datum and the $\mathfrak q=3$ bounds
for the correction and the difference also yield
\eqref{intro:lorentz-endpoint}, as proved in
Proposition~\ref{limit:mixed}.}
}

Energy \emph{equalities} are
then recovered from the smooth positive-time identities and the strong
trace; compactness by itself gives only inequalities.  Different terminal
values produce different solutions with the same datum.  This completes
the logical route from the transition estimate to Theorem~\ref{intro:main},
subject to the detailed arguments in Section~\ref{sec:limit}.

\paragraph{Transition phases and norms.}\phantomsection\label{nl:phase-guide}
The source column has radius one.  Smooth preparation on $[0,\tp]$
initializes only the finest retained scale; the source history on
$[\tp,\tf]$ and transfer on $[\tf,\ts]$ are repeated at each transition.
During transfer the source and parent coordinates satisfy
\begin{equation}\label{strategy:parent-clocks}
 y=x/R,\qquad s=(t-\tf)/R^2,\qquad d=\tf/R^2,\qquad
 s_\ast=\tp-d,\qquad \ts=R^2\tp.
\end{equation}
The parent heat age is $1+d+s$, so its terminal value is exactly
$1+\tp$.  Thus consecutive transitions meet at the same physical time,
without resetting the heat clock.

The $H^4$ norm measures shape accuracy and controls nonlinear products.
The ordinary moment $M_{s_\star}$ controls the pressure-generated tails
and permits iteration.  The negative horizontal Sobolev norm in the
transfer estimate measures the incoming approximation to the point
distribution, before smoothing.  These controls are not interchangeable.
In particular, $R^4$ normalizes the point derivative, whereas $R$
rescales a velocity; their difference produces the factor $R^{-3}$
in the leading amplitude law.

\section{A finite-energy reference solution}
\label{sec:background}

The cascade takes place near an exact unforced solution assembled from
Gaussian swirls at geometrically separated scales.  An axially invariant
swirl solves the equation explicitly but has infinite three-dimensional
energy.  We therefore give each swirl a long Gaussian axial cutoff and
solve for the resulting nonlinear correction.  This is done before any
growing perturbation is introduced.

The output is a finite-energy background with estimates independent of
both scale cutoffs, in particular the smallest radius that will tend to
zero.  The weighted operator lemma specifies the admissible algebraic
tails; a frequency-dependent energy estimate retains axial diffusion
before summing the scales; and a contraction removes the full
axial-localization residual.  These steps avoid a growth loss proportional
to the number of retained swirls.

\subsection{Weighted operators and symmetries}

The weighted estimates below control polynomial spatial moments.
Such moments are needed because application of
$\PP\operatorname{div}$ to a compactly supported tensor source
produces an algebraic pressure tail.
The threshold $\mu<5/2$ follows from the first nonvanishing
far-field term of an odd velocity: a kernel of order $|x|^{-4}$ has
weighted square integral proportional to
$\int_1^\infty r^{2\mu-6}\idiff r$.
In contrast, for the pressure of a general tensor we use the
unrestricted double Riesz transform only at the smaller moment
$\lambda<3/2$.  This difference determines where the weight is
placed in the later nonlinear energy inequalities.

We write $x=(x_h,z)\in\R^2\times\R$, $r=|x_h|$, and
$x_h^\perp=(-x_2,x_1,0)$.  For $r>0$, the cylindrical unit vectors are
$e_r=(x_1,x_2,0)/r$, $e_\theta=x_h^\perp/r$, and $e_z=(0,0,1)$.
A vector field is axisymmetric if it is covariant under rotations about the
$z$-axis, and it is centrally odd if $v(-x)=-v(x)$.  In particular, a
centrally odd axisymmetric field has $v_r,v_\theta$ even in $z$ and $v_z$
odd in $z$.  Its smooth Cartesian representation is used at the axis.

For $\mu\geq0$, we use the ordinary spatial moment
\begin{equation}\label{bg:moment}
 M_\mu(v)=\norm{\langle x\rangle^\mu v}_2,
 \qquad \langle x\rangle=(1+|x|^2)^{1/2}.
\end{equation}
For tensors we use the Frobenius norm, so
$M_\mu(F)^2=\int_{\R^3}\langle x\rangle^{2\mu}
\sum_{i,j}|F_{ij}|^2\,\idiff x$.
The axial Fourier transform
is unitary and is given by
\begin{equation}\label{bg:axial-transform}
 \widehat v^{\,z}(x_h,k)
 =(2\pi)^{-1/2}\int_\R e^{-ikz}v(x_h,z)\idiff z.
\end{equation}
We denote the solenoidal axisymmetric physical velocity space at frequency
$k$ by $\mathcal H_k$.  Its norm includes all three velocity components and
the cylindrical measure $2\pi r\idiff r$.  The Leray projection is $\PP$,
and its restriction to this frequency is $\PP_k$.

\begin{lemma}[Weighted projection and heat estimates]\label{bg:weighted}
Let $0\leq\mu<5/2$.  On centrally odd vector fields one has
\begin{equation}\label{bg:odd-operators}
 M_\mu(\PP v)+M_\mu(e^{h\Delta}v)
 \leq C_\mu M_\mu(v),\qquad h\geq0.
\end{equation}
The projection estimate also holds on horizontal axisymmetric fields,
meaning axisymmetric fields with zero vertical component.
For any tensor $F$ and $h>0$, one has
\begin{equation}\label{bg:tensor-heat}
 M_\mu(e^{h\Delta}\PP\div F)
 \leq C_\mu(1+h^{-1/2})M_\mu(F).
\end{equation}
For $0\leq\lambda<3/2$, each double Riesz transform obeys
\begin{equation}\label{bg:riesz}
 M_\lambda(\mathcal R_i\mathcal R_j f)
 \leq C_\lambda M_\lambda(f).
\end{equation}
The constants remain uniform when the polynomial weight is capped at an
arbitrary radius.
\end{lemma}

\begin{proof}
Partition space into $D_0=\{|x|<2\}$ and
$D_j=\{2^j\leq|x|<2^{j+1}\}$ for $j\geq1$.  Adjacent shells are handled
by the unweighted operator bounds.  The nonlocal kernel $K$ of $\PP$
satisfies $|\partial^\alpha K(x)|\leq C_\alpha|x|^{-3-|\alpha|}$.
If $j\geq k+3$, central oddness gives
$\int_{D_k}v=0$.  Subtracting $K(x)$ in the integral kernel gives
\begin{equation}\label{bg:outward-shell}
 \norm{1_{D_j}\PP(1_{D_k}v)}_2
 \leq C2^{-5(j-k)/2}\norm{1_{D_k}v}_2.
\end{equation}

Here is the scaling in the outward estimate.  For $x\in D_j$ and
$y\in D_k$, with $j\geq k+3$, the mean-value theorem gives
\[
 |K(x-y)-K(x)|\leq C|y|\,|x|^{-4}\leq C2^k2^{-4j}.
\]
Cauchy--Schwarz for the input integral and then the output volume give
\[
 \begin{aligned}
 \|1_{D_j}\PP(1_{D_k}v)\|_2
 &\leq C|D_j|^{1/2}2^{k-4j}|D_k|^{1/2}
                                      \|1_{D_k}v\|_2\\
 &\leq C2^{3j/2+k-4j+3k/2}\|1_{D_k}v\|_2
  =C2^{-5(j-k)/2}\|1_{D_k}v\|_2.
 \end{aligned}
\]
The subtraction is legitimate because the integral of an odd field on
one centered shell is zero.  Without it, the kernel has order $-3$ and
the exponent is only $3/2$.  Multiplying by the ratio of shell weights
$2^{\mu(j-k)}$ explains the two thresholds $\mu<5/2$ and
$\lambda<3/2$ that reappear throughout the paper.  The estimates are
summed as a discrete convolution on the shell $L^2$ norms; absolute
summability of its coefficients gives the weighted operator bound.

For horizontal axisymmetric fields the angular integral vanishes at every
$(r,z)$; subtraction of $K(x-(0,0,z))$ gives the same bound.  In the
opposite direction, without a cancellation hypothesis, one has
\begin{equation}\label{bg:inward-shell}
 \norm{1_{D_j}\PP(1_{D_k}v)}_2
 \leq C2^{-3(k-j)/2}\norm{1_{D_k}v}_2,
 \qquad k\geq j+3.
\end{equation}
The weighted shell coefficients are consequently
$C2^{-(5/2-\mu)(j-k)}$ and $C2^{-(3/2+\mu)(k-j)}$.  Their sums converge
in the asserted range.  The identity part of the projection contributes
only on adjacent shells.

Uniformly for $h>0$, the heat kernel satisfies
$|G_h(x)|\leq C|x|^{-3}$ and $|\nabla G_h(x)|\leq C|x|^{-4}$.
The same cancellation argument, together with its $L^2$ contraction
property, proves the heat estimate in \eqref{bg:odd-operators}.
The kernel of $e^{h\Delta}\PP\div$ is bounded by
$C(|x|+\sqrt h)^{-4}$ and has $L^2$ operator norm at most $Ch^{-1/2}$.
Its outward shell coefficient is
$C2^{-5(j-k)/2}2^{-k}$, without cancellation of $F$; the inward
coefficient is bounded by $C2^{-3(k-j)/2}$.  Summation proves
\eqref{bg:tensor-heat}.

The degree $-3$ kernel of a double Riesz transform has outward shell
coefficient $C2^{-3(j-k)/2}$ without cancellation.  This gives
\eqref{bg:riesz} for $\lambda<3/2$.  Capping $2^{\mu j}$ at $L^\mu$
decreases each outward weight ratio, while the inward ratio is at most
one.  The same summations are therefore uniform in $L$.
\end{proof}

\subsection{The heat array and its linear evolution}

A direct sum of large strain bounds would obscure the smallness of the
coupling near zero axial frequency.  We first retain that frequency
dependence for one column and then sum the scales.  A bounded Fourier
multiplier uses axial dissipation to control the invariant part; only a
small localization error remains in the energy estimate.

The retained radius set $\mathscr A\subset R^{\mathbb Z}$ may be finite or
infinite.  Its largest member is finite whenever the set is nonempty and
bounded above, since the allowed radii are discrete and geometric.  A
bound ``independent of the retained set'' must therefore be proved using
geometric sums, rather than a constant multiplied by the number of
columns.  For an infinite inner array, we first construct finite
approximations and use the positive-time convergence proved below.
Uniformity as the smallest radius tends to zero is needed later when
the cascade depth increases; $A$ remains fixed throughout that limit.

Fix $A\geq1$, $R\geq2$, and $b\geq1$.  For a heat age $T>0$, define
\begin{equation}\label{bg:unit-profiles}
 \Omega_T(r)=T^{-2}e^{-r^2/(4T)},\qquad
 W_T(r)=r\Omega_T(r),\qquad
 B_T(r)=2\left(1-\frac{r^2}{4T}\right)\Omega_T(r).
\end{equation}
These unit profiles contain neither the amplitude nor the axial cutoff;
their subscript always denotes heat age.  For a radius $\ell\in R^\Z$, set
\begin{equation}\label{bg:Qa}
 \begin{split}
 T_{\ell}(t)&=\ell^2+t,\qquad D_{\ell}(t)=b^2\ell^2+2t,\\
 \eta_{\ell}(t,z)&=\frac{b\ell}{\sqrt{D_{\ell}(t)}}
              \exp\left(-\frac{z^2}{2D_{\ell}(t)}\right),\\
 Q_{\ell}(t,x)&=A\ell^2\eta_{\ell}(t,z)
              \Omega_{T_{\ell}(t)}(r)x_h^\perp.
 \end{split}
\end{equation}
Each $Q_{\ell}$ is exactly the three-dimensional heat evolution of $Q_{\ell}(0)$.
For a retained set of radii $\mathscr A\subset R^\Z$, write
$Q_{\mathscr A}=\sum_{\ell\in\mathscr A}Q_{\ell}$.
The two-sided sum and its spatial derivatives converge locally uniformly
for $t>0$.  Splitting at $\ell=\sqrt t$ gives
\begin{equation}\label{bg:heat-sup}
 \norm{\nabla^j Q_{\mathscr A}(t)}_\infty
 \leq C_jAt^{-(j+1)/2},\qquad j\geq0.
\end{equation}
All estimates in this subsection are independent of the retained set.

The uniform derivative bounds do not assume a lower bound on $b\ell$.
Indeed $D_{\ell}\geq T_{\ell}$ and $b\ell/\sqrt{D_{\ell}}\leq1$, so Gaussian
differentiation gives
\[
 \|\partial_z^j\eta_{\ell}(t)\|_\infty
       \leq C_jD_{\ell}(t)^{-j/2}\leq C_jT_{\ell}(t)^{-j/2}.
\]
In particular the axial shear is bounded by
$CA\ell^2T_{\ell}^{-3/2}T_{\ell}^{-1/2}=CA\ell^2T_{\ell}^{-2}$.
For estimates that need the additional width gain, one instead uses
$\|\eta_{\ell}'\|_\infty\leq C/(b\ell)$ and
$\|\eta_{\ell}''\|_\infty\leq C/(b\ell\sqrt{D_{\ell}})$.
These two types of bounds serve different purposes and are both
uniform when old axial widths tend to zero.

The cancellation at small axial frequencies will be encoded in
\begin{equation}\label{bg:clock-def}
 c_{\ell,\mathrm{coup}}(t)=\frac{2A\ell^2}{T_{\ell}(t)^2},\qquad
 L_{\ell}(t)=\frac{\sqrt2 A\ell^2}{T_{\ell}(t)^{3/2}},\qquad
 b_0(t,k)=\sum_{\ell\in R^\Z}\min\bigl(c_{\ell,\mathrm{coup}}(t),L_{\ell}(t)|k|\bigr).
\end{equation}

The coefficient $c_{\ell,\mathrm{coup}}(t)$ in \eqref{bg:clock-def} is a bound that ignores
axial-frequency cancellation.  The alternative $L_{\ell}(t)|k|$ uses
incompressibility and the fact that a horizontal pressure gradient is
annihilated by the Leray projection at $k=0$.  We retain their
\emph{minimum} before summing the radii.  The notation $b_0(t,k)$ is a
frequency-dependent coupling majorant; it is unrelated to the axial
width $b$.  Its integral bound will be combined with the dissipative
factor $e^{-k^2(t-s)}$, not used as a replacement for diffusion.

\begin{lemma}[Frequency-dependent coupling bound]\label{bg:clock}
The function $b_0(t,k)$ decreases in $t$ and satisfies
\begin{equation}\label{bg:clock-bounds}
 b_0(t,k)\leq\frac{16A}{3t},\qquad
 \int_s^{s+h}b_0(t,k)\idiff t
 \leq12A+\frac{2A}{\log R}\log_+(|k|\sqrt h).
\end{equation}
The full axisymmetric linear coupling about the axially invariant sum of
the same swirls has operator norm at most $b_0(t,k)$ on $\mathcal H_k$.
\end{lemma}

\begin{proof}
For one invariant swirl the coupling is
\begin{equation}\label{bg:column-coupling}
 \mathcal B_{\ell}(t,k)v=A\ell^2\bigl(\PP_k(2\Omega_Tv_\theta e_r)
                    -B_Tv_re_\theta\bigr),
\end{equation}
where $T=T_{\ell}(t)$ and the profiles are defined in
\eqref{bg:unit-profiles}.  The two outputs are orthogonal, and their multiplier
bounds give $c_{\ell,\mathrm{coup}}(t)$.  For the first output, define
$p_v(r)=-\int_r^\infty2\Omega_T(\rho)v_\theta(\rho)\idiff\rho$.
The gradient cancellation and the Hilbert--Schmidt calculation give
\begin{equation}\label{bg:primitive}
 \PP_k(\nabla_hp_v,0)=-ik\PP_k(p_ve_z),\qquad
 \norm{v_\theta\mapsto p_v}_{\mathrm{HS}}^2
 =2\int_0^\infty\rho\Omega_T(\rho)^2\idiff\rho=2T^{-3}.
\end{equation}
For the second output, incompressibility gives
\begin{equation}\label{bg:radial-from-vertical}
 v_r(r)=-\frac{ik}{r}\int_0^r\rho v_z(\rho)\idiff\rho,
 \qquad |v_r(r)|\leq\frac{|k|}{\sqrt2}
                    \norm{v_z}_{L^2(r\idiff r)}.
\end{equation}
Since $\norm{B_T}_{L^2(r\idiff r)}=\sqrt2T^{-3/2}$, these bounds
give $L_{\ell}(t)|k|$.

The geometric sums on either side of $\sqrt t$ yield the first bound in
\eqref{bg:clock-bounds}.  By monotonicity in time, the integral of one
summand over an interval of length $h$ is bounded both by
$2Ah/(\ell^2+h)$ and by $2\sqrt2A\ell|k|$.  If $|k|\sqrt h\geq1$, split
the sum into $\ell\leq|k|^{-1}$, $|k|^{-1}<\ell\leq\sqrt h$, and
$\ell>\sqrt h$.  The three bounds are, respectively,
\begin{equation}\label{bg:three-scale-sums}
 \frac{2\sqrt2 A}{1-R^{-1}},\qquad
 2A\left(1+\frac{\log(|k|\sqrt h)}{\log R}\right),\qquad
 \frac{2A}{1-R^{-2}}.
\end{equation}
Their constant parts have sum less than $12A$.  If $|k|\sqrt h<1$,
splitting at $\sqrt h$ and using the same two geometric tails proves
the remaining case of \eqref{bg:clock-bounds}.
\end{proof}

For the heat array $Q=Q_{\mathscr A}$, let $\mathcal U_Q(t,s)$
denote the evolution of the homogeneous linearized equation on
axisymmetric solenoidal velocities.  It includes the
pressure-projected meridional component and acts in the full physical
$L^2$ space.  The corresponding evolution about the exact reference
$\bar u=Q+e$ will be compared with it in
Section~\ref{sec:history}.

\begin{lemma}[Linear evolution about the localized array]\label{bg:propagator}
Set
\begin{equation}\label{bg:metric-constants}
 \zeta_Q=\frac{A}{\log R},\qquad
 K_Q=\exp\left(12A+\zeta_Q\log_+\frac{4\zeta_Q}{e}\right),\qquad
 d_Q=2\log K_Q,\qquad
 \epsilon_Q=\frac{20K_Q(d_Q+4)A}{3\sqrt\pi\,b}.
\end{equation}
The linearized equation about $Q=Q_{\mathscr A}$ has an evolution
$\mathcal U_Q(t,s)$ on axisymmetric solenoidal fields, for $0<s<t$, with
\begin{equation}\label{bg:linear-l2}
 \norm{\mathcal U_Q(t,s)}_{2\to2}
 \leq K_Q(t/s)^{\epsilon_Q}.
\end{equation}
For $j=1,2$, it also satisfies
\begin{equation}\label{bg:linear-smoothing}
 \begin{split}
 \norm{\nabla^j\mathcal U_Q(t,s)f}_2
 &\leq C_j(1+A)^{j/2}K_Q(t/s)^{\epsilon_Q}
                         (t-s)^{-j/2}\norm f_2,\\
 \norm{\mathcal U_Q(t,s)f}_\infty
 &\leq C(1+A)^{3/4}K_Q(t/s)^{\epsilon_Q}
                         (t-s)^{-3/4}\norm f_2.
 \end{split}
\end{equation}
If $\epsilon_Q\leq1/8$ and $F$ is an axisymmetric solenoidal source with $\norm{F(t)}_2\leq Mt^{-3/4}$, the
zero-trace response $\mathcal T_QF(t)=\int_0^t\mathcal U_Q(t,s)F(s)\idiff s$
is well defined, has square-integrable gradient on finite time intervals,
and satisfies
\begin{equation}\label{bg:source-map}
 \sup_{t>0}\left(t^{-1/4}\norm{\mathcal T_QF(t)}_2
       +t^{1/4}\norm{\nabla\mathcal T_QF(t)}_2
       +t^{1/2}\norm{\mathcal T_QF(t)}_\infty\right)
 \leq L_{\mathrm{ref}}(A) M,
 \quad L_{\mathrm{ref}}(A)\leq C(1+A)^{3/4}K_Q.
\end{equation}
\end{lemma}

\begin{proof}
We construct a bounded axial multiplier which absorbs the potentially
large invariant coupling.  Define
\begin{equation}\label{bg:multiplier}
 \Phi_Q(t,k)=\int_t^\infty
                   \left(2b_0(u,k)-\frac{k^2}{2}\right)_+\idiff u,
 \qquad m_Q(t,k)=e^{\Phi_Q(t,k)/2},\qquad m_Q(t,0)=1.
\end{equation}
The active time set in this integral is an interval.  The second estimate
in \eqref{bg:clock-bounds}, followed by maximization in $|k|\sqrt h$,
gives
\begin{equation}\label{bg:multiplier-bound}
 0\leq\Phi_Q\leq2\log K_Q,\qquad
 1\leq m_Q\leq K_Q,\qquad
 \partial_t\log m_Q=-\frac12
                       \left(2b_0-\frac{k^2}{2}\right)_+.
\end{equation}
For $k\ne0$, the active set is contained in
$0<u<64A/(3k^2)$.  Its moving endpoint contributes zero when the
integral is differentiated.  Since $0\leq k\partial_kb_0\leq b_0$,
one obtains
\begin{equation}\label{bg:multiplier-derivative}
 k\partial_k\Phi_Q
 =2\int_{\mathrm{active}}k\partial_kb_0\idiff u
       -k^2|\mathrm{active}|,
 \qquad |k\partial_k\Phi_Q|\leq d_Q.
\end{equation}
These statements hold almost everywhere and uniformly on compact
frequency intervals not containing zero.

For clarity, the two sides of the derivative bound follow differently.
On the active set,
\[
 k\partial_k\Phi_Q
 \leq \Phi_Q-\tfrac12k^2|\mathrm{active}|\leq d_Q,
 \qquad
 k\partial_k\Phi_Q\geq-k^2|\mathrm{active}|
 \geq-64A/3\geq-d_Q,
\]
since $d_Q=2\log K_Q\geq24A$.  This also shows why the estimate remains
uniform as the retained array is enlarged.

Put $B_{\ell}(k,k')=\min(c_{\ell,\mathrm{coup}},L_{\ell}\max(|k|,|k'|))$.  For comparable frequencies
of the same sign, integrate
$|\partial_km_Q|\leq K_Qd_Q/(2|k|)$.  Otherwise use
$\max(|k|,|k'|)\leq2|k-k'|$ and \eqref{bg:multiplier-bound}.  Both cases
give
\begin{equation}\label{bg:multiplier-commutator}
 B_{\ell}(k,k')\left|\frac{m_Q(t,k)}{m_Q(t,k')}-1\right|
 \leq K_Q(d_Q+2)L_{\ell}|k-k'|.
\end{equation}
This estimate includes opposite-sign and zero frequencies.

Multiplication by $\eta_{\ell}$ becomes convolution with the nonnegative
even kernel
\begin{equation}\label{bg:axial-kernel}
 G_{\ell}(q)=\frac{b\ell}{\sqrt{2\pi}}e^{-D_{\ell}q^2/2},\qquad
 \int G_{\ell}\leq1,\qquad
 \mu_{\ell}:=\int|q|G_{\ell}(q)\idiff q
       =\sqrt{2/\pi}\,\frac{b\ell}{D_{\ell}}.
\end{equation}
In the two nonshear terms, the radial pressure primitive uses the
output frequency $k$, whereas \eqref{bg:radial-from-vertical} uses the
input frequency $k'$.  Their combined norm is bounded by $B_{\ell}(k,k')$.
The remaining term is the axial shear
$-A\ell^2\partial_z\eta_{\ell}\,W_{T_{\ell}}v_ze_\theta$, where
$W_T$ is the profile in \eqref{bg:unit-profiles}; its radial multiplier norm is
$L_{\ell}/\sqrt e$, with convolution kernel bounded by $|q|G_{\ell}(q)$.

For $E=\norm{m_Qv}_2^2$ and $D=\norm{m_Q\nabla v}_2^2$, split the
multiplier ratio into one and its difference.  Symmetry of $G_{\ell}$ and
$B_{\ell}$, the estimate
$B_{\ell}(k,k')\leq\min(c_{\ell,\mathrm{coup}},L_{\ell}|k|)+L_{\ell}|k-k'|$, and
\eqref{bg:multiplier-commutator} imply
\begin{equation}\label{bg:coupling-energy}
 \operatorname{Re}\langle m_Qv,m_QB_Qv\rangle
 \leq\int b_0(t,k)\norm{m_Q\widehat v^{\,z}(k)}_{L^2_{x_h}}^2\idiff k
       +K_Q(d_Q+4)\sum_{\ell}L_{\ell}\mu_{\ell} E,
\end{equation}
where $B_Qv=-\PP\div(Q\otimes v+v\otimes Q)$.  The full scale sum is
bounded by
\begin{equation}\label{bg:width-sum}
 \sum_{\ell}L_{\ell}\mu_{\ell}
 \leq\frac{2A}{\sqrt\pi b}
       \sum_{\ell}\frac{\ell}{(\ell^2+t)^{3/2}}
 \leq\frac{20A}{3\sqrt\pi bt}.
\end{equation}
The multiplier commutes with $\PP$ and Cartesian derivatives.  Its
time derivative in \eqref{bg:multiplier-bound}, together with
\eqref{bg:coupling-energy}--\eqref{bg:width-sum}, therefore gives
\begin{equation}\label{bg:metric-energy}
 \frac12E'+\frac34D
 \leq\frac{\epsilon_Q}{t}E
       +\operatorname{Re}\langle m_Qv,m_QF\rangle.
\end{equation}
The inequality
$b_0-\frac12(2b_0-k^2/2)_+\leq k^2/4$
retains the stated fraction of the full Cartesian dissipation.

We next explain existence and smoothing for this estimate.  If
$Q=\Omega_Q(t,r,z)x_h^\perp$, where
\begin{equation}\label{bg:array-angular-coefficient}
 \Omega_Q(t,r,z)=\sum_{\ell\in\mathscr A}
       A\ell^2\eta_\ell(t,z)\Omega_{T_\ell(t)}(r),
\end{equation}
its axisymmetric coupling has the smooth Cartesian extension
\begin{equation}\label{bg:cartesian-matrix}
 B_Qv=\PP M_Qv,\qquad
 M_Qv=-2\Omega_Q v_h^\perp
       -\frac{\partial_r\Omega_Q}{r}(x_h\cdot v_h)x_h^\perp
       -\partial_z\Omega_Q\,v_zx_h^\perp.
\end{equation}
On every positive time interval,
$\norm{\nabla^jM_Q(t)}_\infty\leq C_jAt^{-1-j/2}$.
The heat Volterra equation therefore constructs $\mathcal U_Q$ there.
Omitted radii $\ell<\ell_-$ or $\ell>\ell_+$ contribute at most
$CA(\ell_-^2/s^2+\ell_+^{-2})$ on a time triangle starting at $s>0$.
Finite arrays converge in operator norm, and their common estimate
\eqref{bg:metric-energy} gives \eqref{bg:linear-l2} for the limiting
array.

For smoothing, take the last interval of length
$h_0=\min((t-s)/2,t/(8(1+A)))$.  One has
$h_0\geq(t-s)/(8(1+A))$.  Rescaling time by $h_0$ and space by
$\sqrt{h_0}$ bounds the matrix derivatives through the required orders.
For the Cartesian derivative energies $E_j=\sum_{|\alpha|=j}
\norm{\partial^\alpha v}_2^2$, differentiation gives
$E_j'+2E_{j+1}\leq C\sum_{i\leq j}E_i$ in these coordinates.
A positive combination of $E_0$, $\tau E_1$, and $\tau^2E_2$
absorbs the derivatives of the time weights.  This proves the first
line of \eqref{bg:linear-smoothing};
$\norm v_\infty\leq C\norm v_2^{1/4}\norm{\nabla^2v}_2^{3/4}$
proves the second line.

Finally, set $\kappa_Q=1/4-\epsilon_Q>0$ and initially start the
source response at a positive time.  The square root of
\eqref{bg:metric-energy} satisfies
$(E^{1/2})'\leq\epsilon_Qt^{-1}E^{1/2}+K_QMt^{-3/4}$.
Consequently,
\begin{equation}\label{bg:zero-trace-energy}
 E(t)^{1/2}\leq K_QM\kappa_Q^{-1}t^{1/4},\qquad
 \int_0^tD(u)\idiff u
 \leq\frac{2K_Q^2M^2}{3\kappa_Q^2}\sqrt t.
\end{equation}
Responses starting at two times below $\tau$ differ at $t\geq\tau$
by at most $K_QM\kappa_Q^{-1}t^{\epsilon_Q}\tau^{\kappa_Q}$ in
$L^2$.  The homogeneous energy estimate gives convergence of their
gradients as well.  Integrating \eqref{bg:linear-smoothing} uses the
finite beta integrals $\mathrm B(\kappa_Q,1/2)$ and
$\mathrm B(\kappa_Q,1/4)$, proving \eqref{bg:source-map}.
Since $\norm{B_Q\mathcal T_QF(t)}_2\leq CAt^{-1}
\norm{\mathcal T_QF(t)}_2$ is integrable at zero, the equation also
holds with the stated zero initial trace.
\end{proof}

\subsection{Construction of the exact reference solution}

The heat evolution of each azimuthal Gaussian is explicit, but the
sum is not itself a solution of the unforced Navier--Stokes equations.
Its nonlinear radial force involves products of \emph{all} pairs of
swirls, not just their self-interactions.  The pressure used below
cancels every ordered pair in the radial direction; what remains is
an axial derivative caused by localization.
The axial width makes the residual small.

There are three distinct tasks below.  First, compute a pressure that
cancels the radial force of \emph{every ordered pair} of Gaussian
columns.  Second, show that the uncancelled axial residual has an
$L^2$ bound proportional to $b^{-1/2}t^{-3/4}$; its power of $t$ is
integrable after the source estimate of Lemma~\ref{bg:propagator}.
Third, solve a zero-initial-trace fixed-point equation for $e$ in a
space controlling $t^{-1/4}\|e\|_2$, $t^{1/4}\|\nabla e\|_2$ and
$t^{1/2}\|e\|_\infty$.  That space is chosen to make the nonlinear
product $e\cdot\nabla e$ have precisely the same $t^{-3/4}$ source
size.  Lemma~\ref{bg:propagator} then maps this source class into the
correction space, where the quadratic term is controlled by
\eqref{bg:product}.

\begin{proposition}[The reference array]\label{bg:reference}
Fix $\varrho>0$ and set $R=e^{\varrho A}\geq2$.  There is a constant
$B_\varrho$ such that, if $b\geq e^{B_\varrho A}$, every nonempty radius set
$\mathscr A\subset R^\Z$ bounded above admits an unforced solution
$\bar u_{\mathscr A}=Q_{\mathscr A}+e_{\mathscr A}$ on
$0<t\leq \ell_{\max}^2$, where $\ell_{\max}=\max\mathscr A$.  It satisfies
\begin{equation}\label{bg:correction}
 \sup_t\left(t^{-1/4}\norm{e_{\mathscr A}(t)}_2
       +t^{1/4}\norm{\nabla e_{\mathscr A}(t)}_2
       +t^{1/2}\norm{e_{\mathscr A}(t)}_\infty\right)
 \leq C_\varrho A^P e^{C_\varrho A}b^{-1/2}.
\end{equation}
For every fixed integer $j\geq0$, one also has
\begin{equation}\label{bg:derivatives}
 \begin{split}
 \norm{\nabla^je_{\mathscr A}(t)}_2
   &\leq\Sigma_jt^{1/4-j/2},\\
 \norm{\nabla^je_{\mathscr A}(t)}_\infty
   &\leq\Sigma_jt^{-1/2-j/2},\qquad
 \Sigma_j\leq C_jA^{P_j}e^{C_jA}b^{-1/2}.
 \end{split}
\end{equation}
The constants are independent of the retained radii.  The solution is
centrally odd, axisymmetric, smooth for $t>0$, and strongly continuous
in $L^2$ at zero with datum $Q_{\mathscr A}(0)$.  It satisfies global
and local energy equality through zero.

For $\mathscr A=\{R^j:-I\leq j\leq J\}$, write
$\bar u^{I,J}$ and allow $I=\infty$.  The reference solutions satisfy
the exact covariance
\begin{equation}\label{bg:covariance}
 R\bar u^{I,J}(R^2t,Rx)=\bar u^{I+1,J-1}(t,x),\qquad J\geq1.
\end{equation}
For $\ell_n=R^{-n}$, $n\geq0$, their limiting datum is
\begin{equation}\label{bg:datum}
 u_0(x)=\sum_{n=0}^\infty
    \frac{A}{\ell_n^2}
    \exp\left(-\frac{r^2}{4\ell_n^2}
              -\frac{z^2}{2b^2\ell_n^2}\right)x_h^\perp.
\end{equation}
This datum is nonzero, divergence free, smooth away from the origin,
and belongs to $L^{\mathfrak p}(\R^3)$ for every $1\leq \mathfrak p<3$.
\end{proposition}

\begin{proof}
The pressure cancellation must retain every ordered pair of scales.
For radii $\ell,\tilde\ell$, define
\begin{equation}\label{bg:pair-pressure}
 P_{\ell,\tilde\ell}=-\frac{2A^2\ell^2\tilde\ell^2}{T_{\ell}T_{\tilde\ell}(T_{\ell}+T_{\tilde\ell})}\eta_{\ell}\eta_{\tilde\ell}
       \exp\left(-\frac{r^2(T_{\ell}+T_{\tilde\ell})}{4T_{\ell}T_{\tilde\ell}}\right).
\end{equation}
Direct differentiation yields
$\partial_rP_{\ell,\tilde\ell}=Q_{\ell,\theta}Q_{\tilde\ell,\theta}/r$.  Thus, with
$P_Q=\sum_{\ell,\tilde\ell}P_{\ell,\tilde\ell}$, the full heat array obeys
\begin{equation}\label{bg:residual}
 Q_t-\Delta Q+(Q\cdot\nabla)Q+\nabla P_Q
       =(\partial_zP_Q)e_z.
\end{equation}
For $\ell\leq \tilde\ell$, the effective axial variance is
$D_{\ell}D_{\tilde\ell}/(D_{\ell}+D_{\tilde\ell})\geq D_{\ell}/2$.  Its differentiated Gaussian has
$L_z^2$ norm at most $Cb^{-1/2}\ell^{-1/2}$, and the horizontal Gaussian
has norm $(2\pi T_{\ell}T_{\tilde\ell}/(T_{\ell}+T_{\tilde\ell}))^{1/2}$.  Consequently,
\begin{equation}\label{bg:pair-residual}
 \norm{\partial_zP_{\ell,\tilde\ell}(t)}_2
 \leq CA^2b^{-1/2}
       \frac{\ell^{3/2}\tilde\ell^2}{\sqrt{T_{\ell}T_{\tilde\ell}}(T_{\ell}+T_{\tilde\ell})^{3/2}}.
\end{equation}
After $\ell=\sqrt t\,u$ and $\tilde\ell=\sqrt t\,v$, the last factor equals
$t^{-3/4}f(u,v)$, where for $u\leq v$ it is bounded by
\begin{equation}\label{bg:residual-regions}
 f(u,v)\leq C
 \begin{cases}
 u^{3/2}v^2,&v\leq1,\\
 u^{3/2}v^{-2},&u\leq1\leq v,\\
 u^{1/2}v^{-2},&1\leq u\leq v.
 \end{cases}
\end{equation}
In the first region the inner sum leaves a summable $v^{7/2}$ tail;
in the second region the two geometric sums separate; in the third,
the inner sum leaves $v^{-3/2}$.  Doubling the half-sum therefore gives
\begin{equation}\label{bg:full-residual}
 \sum_{\ell,\tilde\ell}\norm{\partial_zP_{\ell,\tilde\ell}(t)}_2
 \leq C_0A^2b^{-1/2}t^{-3/4}.
\end{equation}
This proves absolute convergence of its axial derivative. Differentiating
\eqref{bg:pair-pressure} in $r$ gives the radial bound
$\sum_{\ell,\tilde\ell}\|Q_{\ell,\theta}Q_{\tilde\ell,\theta}/r\|_2
\leq(\sum_\ell\|Q_{\ell,\theta}/r\|_\infty)
(\sum_{\tilde\ell}\|Q_{\tilde\ell,\theta}\|_2)<\infty$
at every positive time, since the radii are bounded above and geometrically
separated. Fixed spatial derivatives of the axial residual add only $C_jt^{-j/2}$, since both effective
variances are bounded below by a fixed multiple of $t$.

Set $F_Q=-\PP((\partial_zP_Q)e_z)$.  The correction equation is
\begin{equation}\label{bg:correction-equation}
 e_t-\Delta e=B_Qe+F_Q-\PP(e\cdot\nabla e),\qquad e(0)=0.
\end{equation}

\textit{The fixed-point space.}
Let $\mathcal X$ be the Banach space of real, centrally odd,
axisymmetric, divergence-free fields on $(0,\ell_{\max}^2]$, continuous
into $L^2$ on positive compact time intervals, with weak spatial gradient
and finite norm given by the left side of \eqref{bg:correction}.
The gradient and supremum terms may be taken as essential time suprema;
the solution map supplies their positive-time representatives.
Completeness follows from the weighted supremum norms and closedness of
the distributional gradient, divergence and symmetry conditions.
The residual and the fixed-point map preserve this subspace, since the
symmetries commute with the equation and projection, and the forcing is
projected onto solenoidal fields.  Finite norm gives the zero $L^2$ trace.

Its product estimate is
\begin{equation}\label{bg:product}
 \norm{\PP(v\cdot\nabla w)(t)}_2
 \leq\norm{v(t)}_\infty\norm{\nabla w(t)}_2
 \leq\norm v_{\mathcal X}\norm w_{\mathcal X}t^{-3/4}.
\end{equation}
With $M_A=C_0A^2b^{-1/2}$, Lemma~\ref{bg:propagator} shows that
$e\mapsto\mathcal T_QF_Q-\mathcal T_Q\PP(e\cdot\nabla e)$
is a contraction on $\norm e_{\mathcal X}\leq2L_{\mathrm{ref}}(A)M_A$, provided
$\epsilon_Q\leq1/8$ and $8L_{\mathrm{ref}}(A)^2M_A\leq1$.  For $R=e^{\varrho A}$,
the constant $K_Q$ is $e^{12A+c_Q(\varrho)}$, where
$c_Q(\varrho):=\varrho^{-1}\log_+(4/(e\varrho))$.
A sufficiently large fixed $B_\varrho$ therefore ensures both
conditions and proves \eqref{bg:correction}.

Deleting radii decreases all absolute majorants. Scaling the largest
radius to one gives the asserted uniform lifespan. The norm
$\mathcal X$ is invariant under Navier--Stokes scaling.
To apply fixed-point uniqueness to a scaled solution, write its Duhamel
formula first from a positive time $\delta$. By \eqref{bg:linear-l2},
its homogeneous initial term satisfies
\begin{equation}\label{bg:vanishing-homogeneous-trace}
 \|\mathcal U_Q(t,\delta)e(\delta)\|_2
 \leq K_Qt^{\epsilon_Q}\delta^{1/4-\epsilon_Q}
                 \|e\|_{\mathcal X}\longrightarrow0
                 \qquad(\delta\downarrow0).
\end{equation}
The source integral converges by \eqref{bg:product} and
\eqref{bg:source-map}. Thus the scaled correction satisfies the same
zero-trace mild equation and belongs to the same contraction ball.
Uniqueness in that ball gives \eqref{bg:covariance} for the entire
exact solution, including its correction.

For derivative estimates, rescale $[t/2,t]$ by time $t$, length
$\sqrt t$, and velocity $\sqrt t$.  The matrix derivatives are
$O_j(A)$, the residual derivatives are $O_j(A^2b^{-1/2})$, and the
correction supremum is bounded by
$\sigma_0=CA^Pe^{CA}b^{-1/2}\leq1$.  Leibniz's rule and Sobolev
interpolation give
$\norm{\nabla^j(e\otimes e)}_2\leq C_j\norm e_\infty
\norm{\nabla^je}_2$.  The differentiated energy inequalities are
therefore
\begin{equation}\label{bg:derivative-energies}
 E_j'+E_{j+1}
 \leq C_m(A+1+\sigma_0^2)\sum_{i\leq j}E_i
          +C_mA^4b^{-1},\qquad 0\leq j\leq m+2.
\end{equation}
Successive time weights on this fixed interval cost at most
$A^{P_m}e^{C_mA}$.  Spatial interpolation with two additional
derivatives, followed by rescaling, proves \eqref{bg:derivatives}.
Only finitely many derivative orders are needed in choosing $b$;
ordinary positive-time parabolic regularity then gives smoothness.

The single-scale datum has the exact norm scaling
\begin{equation}\label{bg:datum-norm}
 \norm{Q_{\ell}(0)}_{\mathfrak p}=C_{\mathfrak p}Ab^{1/\mathfrak p}\ell^{3/\mathfrak p-1}.
\end{equation}
Summing this expression proves the integrability assertion in
\eqref{bg:datum}.  The summands have the same azimuthal orientation
and are individually divergence free.  On a compact set away from
zero, at least one of $r/\ell$ and $|z|/(b\ell)$ tends to infinity as
$\ell\to0$, which gives local convergence of every differentiated
series and proves the remaining assertions about the datum.

The heat field is strongly continuous in $L^2$ at zero and has finite
integrated gradient energy.  Equations~\eqref{bg:correction} and
\eqref{bg:zero-trace-energy} give the same properties for the
correction.  The tensors $Qe$, $eQ$, and $e^2$ are in $L_t^1L_x^2$
near zero, since their norms are $O(t^{-1/4})$; together with
\eqref{bg:full-residual}, this passes \eqref{bg:correction-equation}
to zero with initial datum $Q(0)$.

The pressure may be chosen as
$p=\sum_{i,j}\mathcal R_i\mathcal R_j(\bar u_i\bar u_j)$.
The energy-class bounds already obtained, without using energy equality,
give $\bar u\in L^{10/3}_{t,x}$ and $p\in L^{5/3}_{t,x}$ on every
finite time interval in the stated lifespan.  To remove a spatial cutoff,
we also need a \emph{global} flux estimate.  Interpolation and H\"older give
\begin{equation}\label{bg:global-flux}
 \begin{aligned}
 \int_0^T\|\bar u(t)\|_3^3\idiff t
 &\leq C\|\bar u\|_{L^\infty(0,T;L^2)}^{3/2}
       \left(\int_0^T\|\nabla\bar u(t)\|_2^2\idiff t\right)^{3/4}T^{1/4},\\
 \|p(t)\|_{3/2}&\leq C\|\bar u(t)\|_3^2.
 \end{aligned}
\end{equation}
Consequently $|\bar u|^3+|p\bar u|$ is globally integrable in space and
time.  On a positive compact time interval, the solution is smooth,
bounded, and finite energy.  In the cutoff energy identity, the quadratic
cutoff term tends to zero by the $L^\infty_tL^2_x$ bound and the cubic
and pressure fluxes by \eqref{bg:global-flux}.  This gives the global
energy equality between positive times.  Strong $L^2$ convergence at zero
and finite dissipation extend it through zero.

For a compactly supported space-time test, $\bar u\in L^{10/3}_{t,x}$ and
$p\in L^{5/3}_{t,x}$ also give local integrability of the fluxes.  The
same initial-trace passage therefore proves local energy equality through
zero.  The global cutoff argument uses \eqref{bg:global-flux}, rather
than an $L^{10/3}$-to-$L^3$ embedding on the infinite-measure space.

\textit{Convergence of the exact reference as the inner cutoff is removed.}
The convergence of the heat arrays alone does not automatically imply
convergence of their nonlinear corrections.  Here is the needed argument.
Fix the outer cutoff and a positive compact time interval $[\tau,T]$.
For the arrays with inner cutoff $I\to\infty$, the coefficient matrices
converge with all fixed derivatives uniformly on every positive time
triangle.  Their residuals $F_Q^I(t)$ converge in $L^2$ at each $t>0$
by absolute convergence of the ordered-pair sum in
\eqref{bg:full-residual}.  They have the common majorant
$C A^2b^{-1/2}t^{-3/4}$.

Use the Picard iteration
\[
 e^{I,(0)}=0,\qquad
 e^{I,(n+1)}=\mathcal T_{Q^I}
       \bigl(F_Q^I-\PP(e^{I,(n)}\cdot\nabla e^{I,(n)})\bigr).
\]
The contraction factor and the ball in $\mathcal X$ are uniform in $I$.
For a fixed Picard index, induction gives convergence at every positive
time in $L^2$, $\dot H^1$, and $L^\infty$.
Indeed, the linear propagators converge on positive time triangles;
their bounds in Lemma~\ref{bg:propagator} dominate the time integral.
At its initial endpoint the majorant is
$t^{\epsilon_Q}s^{-3/4-\epsilon_Q}$, which is integrable because
$\epsilon_Q<1/4$.  At the terminal endpoint the gradient and supremum
kernels have the integrable singularities $(t-s)^{-1/2}$ and
$(t-s)^{-3/4}$.  The nonlinear source converges in $L^2$ by
\[
 \|v_I\cdot\nabla v_I-v\cdot\nabla v\|_2
 \leq\|v_I-v\|_\infty\|\nabla v_I\|_2
       +\|v\|_\infty\|\nabla(v_I-v)\|_2,
\]
and has the same $s^{-3/4}$ majorant.
Convergence of the smoothed propagators used here follows from the
heat Volterra equation on a positive time triangle, followed by the
last-interval derivative estimate in Lemma~\ref{bg:propagator}.

First let $I\to\infty$ at a fixed Picard index, then let the index
tend to infinity.  The common geometric bound on the Picard remainder
identifies the limit with the correction constructed for the infinite
inner array.  Splitting off short endpoint intervals in the same
integrals gives uniform convergence on $[\tau,T]$.
The derivative estimates \eqref{bg:derivatives}, applied
at finitely many higher orders, and Sobolev interpolation then give
convergence of every fixed spatial derivative on positive compact
time intervals in the norms needed below.
This convergence also holds in $C([0,T];L^2)$ for the complete
reference velocities.  For the heat parts, the $L^2$ contraction gives
\[
 \sup_{0\leq t\leq T}\|Q^I(t)-Q^\infty(t)\|_2
       \leq\|Q^I(0)-Q^\infty(0)\|_2\longrightarrow0.
\]
For the corrections, \eqref{bg:correction} gives the uniform bound
$\|e^I(t)-e^\infty(t)\|_2\leq2\Sigma_0t^{1/4}$ near zero.
Combine it with uniform convergence on $[\tau,T]$, first letting
$I\to\infty$ and then $\tau\downarrow0$.
In particular the actual reference coefficients, not just their heat
parts, converge on positive time intervals.  The linearized propagators about these
references converge strongly on such intervals by their difference
equations and the ordinary energy estimate.  This justifies the inner
cutoff limit used later in the weighted and transfer estimates.

\end{proof}

\section{An unstable viscous response and its signed quadratic stress}
\label{sec:seed}

This section constructs the prescribed linear response that will drive
one transition and determines the sign of its complete quadratic stress.
The comparison background is an axially invariant Gaussian swirl.  It is
an exact viscous solution but has infinite three-dimensional energy;
Section~\ref{sec:history} will compare its response with the finite-energy
reference of Proposition~\ref{bg:reference}.

There are three conclusions to keep in view: a sharp physical energy upper
bound and actual growth of the prescribed point input, weighted regularity
of its full response, and a definite sign for the complete stress history.
A frozen spectral calculation alone supplies neither the time-dependent
viscous growth nor the stress sign.  The frequencies outside the growing
band and the azimuthal mean of the second variation are therefore retained
through Lemma~\ref{seed:signed-history}.

Using the profiles in \eqref{bg:unit-profiles}, we write the heat age
and invariant Gaussian swirl as
\begin{equation}\label{seed:background}
 \tau_t=1+t,\qquad
 U_{\mathrm d}(t,x)=A W_{\tau_t}(r)e_\theta,
 \qquad c_\star=\sqrt{2}\,\mathrm e^{-3/2}.
\end{equation}
Here $A\geq 1$, $x=(x_h,z)$, and $r=|x_h|$. The field $U_{\mathrm d}$ solves
the Navier--Stokes equations with viscosity one, because
\begin{equation}\label{seed:swirl-heat}
 \partial_\tau W_\tau
 =\bigl(\partial_r^2+r^{-1}\partial_r-r^{-2}\bigr)W_\tau,
 \qquad
 \partial_r p_{\mathrm d}=\frac{A^2W_{\tau_t}^2}{r}.
\end{equation}
We denote its homogeneous linearized evolution by
$\mathcal U_{\mathrm d}(t,s)$, always with the full Cartesian Leray projection.

\subsection{The physical Fourier representation and a sharp energy estimate}

We keep the physical velocity energy when reducing the problem to a
radial equation.  In particular, the vertical component forced by
incompressibility is part of the norm.  The bounded azimuthal weight is
chosen to cancel the coupling in the stable inner region and to expose
the optimal leading growth constant in the remaining region.

The unit angular velocity and azimuthal coupling coefficient, evaluated
at heat age $\tau_t=1+t$, are
\begin{equation}\label{seed:coefficients}
 \Omega_{\tau_t}(r)=\tau_t^{-2}\exp\Bigl(-\frac{r^2}{4\tau_t}\Bigr),
 \qquad
 B_{\tau_t}(r)=2\tau_t^{-2}\Bigl(1-\frac{r^2}{4\tau_t}\Bigr)
                  \exp\Bigl(-\frac{r^2}{4\tau_t}\Bigr).
\end{equation}
For an axisymmetric solenoidal perturbation $v$, the linearized system is
\begin{equation}\label{seed:linearized-components}
 \begin{aligned}
 \partial_t v_r-(\Delta-r^{-2})v_r&=2A\Omega_{\tau_t} v_\theta-\partial_r p,\\
 \partial_t v_\theta-(\Delta-r^{-2})v_\theta&=-AB_{\tau_t} v_r,\\
 \partial_t v_z-\Delta v_z&=-\partial_z p,\\
 \partial_r v_r+r^{-1}v_r+\partial_z v_z&=0,
 \end{aligned}
\end{equation}
where $\Delta=\partial_r^2+r^{-1}\partial_r+\partial_z^2$ in these scalar
component equations. We use the unitary Fourier transform in $z$. At frequency
$k>0$, a solenoidal axisymmetric vector can be represented as
\begin{equation}\label{seed:fiber-coordinates}
 \widehat v^{\,z}(r,k)
   =\bigl(kf_r(r),f_\theta(r),\mathrm iD_rf_r(r)\bigr),
 \qquad D_rf_r=f_r'+r^{-1}f_r.
\end{equation}
The third component makes the divergence vanish. Set
$H_{\mathrm{rad}}=L^2((0,\infty),r\idiff r)$ and let
\begin{equation}\label{seed:radial-operator}
 L_k=-\partial_r^2-r^{-1}\partial_r+r^{-2}+k^2
\end{equation}
have its positive Friedrichs realization in $H_{\mathrm{rad}}$. Its form space, denoted by
$H_{L_k}$, has norm
\begin{equation}\label{seed:form-norm}
 \|f_r\|_{H_{L_k}}^2
 =\|D_rf_r\|_{H_{\mathrm{rad}}}^2+k^2\|f_r\|_{H_{\mathrm{rad}}}^2
 =\int_0^\infty
       \bigl(|f_r'|^2+r^{-2}|f_r|^2+k^2|f_r|^2\bigr)r\idiff r.
\end{equation}
Thus, on $X_k=H_{L_k}\oplus H_{\mathrm{rad}}$, the physical energy of the Fourier component
is $2\pi\|V\|_{X_k}^2$, where $V=(f_r,f_\theta)$. Inner products are linear
in the first argument.

At fixed $k>0$, incompressibility reads
$D_r\widehat v_r+ik\widehat v_z=0$.  Writing $\widehat v_r=kf_r$
therefore forces $\widehat v_z=iD_rf_r$.  In particular the radial
coordinate is not the radial velocity itself, and its derivative is part
of the physical energy:
\[
 \frac1{2\pi}\|\widehat v^{\,z}(k)\|_{L^2_{x_h}}^2
   =k^2\|f_r\|_{H_{\mathrm{rad}}}^2+\|D_rf_r\|_{H_{\mathrm{rad}}}^2+\|f_\theta\|_{H_{\mathrm{rad}}}^2.
\]
For a smooth regular profile, the cross term in
$\int|f_r'+f_r/r|^2r\,\idiff r$ is
$\int\partial_r|f_r|^2\,\idiff r=0$.  This proves
\eqref{seed:form-norm} first on the smooth core and then by closure.
The Friedrichs realization means precisely the positive operator associated
with the closure of this quadratic form; no singular $r^{-1}$ mode is
being admitted at the axis.

For completeness, the pressure-free radial weak equation is
\[
 \langle\partial_t f_r,L_kh\rangle_{H_{\mathrm{rad}}}
  =-\langle L_k f_r,L_kh\rangle_{H_{\mathrm{rad}}}
       +2Ak\langle\Omega_{\tau_t} f_\theta,h\rangle_{H_{\mathrm{rad}}}.
\]
The test velocity $(kh,0,iD_rh)$ is solenoidal, so its pairing with the
full pressure gradient vanishes.  This identity gives
$\partial_t f_r=-L_kf_r+2AkL_k^{-1}(\Omega_{\tau_t}f_\theta)$, whereas the
azimuthal equation gives
$\partial_t f_\theta=-L_kf_\theta-AkB_{\tau_t}f_r$.
These are the two off-diagonal blocks of
\eqref{seed:fiber-evolution}, acting in the physical energy space $X_k$.

Testing the meridional equations against
$(kh,0,\mathrm iD_rh)$ eliminates the pressure and gives
\begin{equation}\label{seed:fiber-evolution}
 V'=A\mathcal C_k(t)V-\mathcal D_kV,
 \qquad
 \mathcal C_k(t)=
 \begin{pmatrix}
 0&2kL_k^{-1}\Omega_{\tau_t}\\
 -kB_{\tau_t}&0
 \end{pmatrix},
 \qquad
 \mathcal D_k=\operatorname{diag}(L_k,L_k).
\end{equation}
The operator $\mathcal D_k$ is positive and self-adjoint in $X_k$. For a
nonnegative integer $m$, write
$X_k^m=D(\mathcal D_k^{m/2})$, with norm $\|\mathcal D_k^{m/2}V\|_{X_k}$. On compact subsets of
$0<k<\infty$, this norm is uniformly equivalent to the horizontal Sobolev norm
of order $m$ of the physical vector in \eqref{seed:fiber-coordinates}. The
identity
\begin{equation}\label{seed:radial-intertwining}
 D_rL_k=\bigl(-\partial_r^2-r^{-1}\partial_r+k^2\bigr)D_r
\end{equation}
also verifies the correspondence for the vertical component.

Set $\zeta_r(t,r)=r^2/(4\tau_t)$.  The radial and azimuthal coupling coefficients are $2\Omega_{\tau_t}$ and
$-B_{\tau_t}$.  If the azimuthal energy has weight $F$, their combined contribution
is proportional to $(2\Omega_{\tau_t}-FB_{\tau_t})\operatorname{Re}(v_r\overline{v_\theta})$.
In the inner region, choosing
$F=2\Omega_{\tau_t}/B_{\tau_t}=(1-\zeta_r)^{-1}$ cancels this term.  This choice cannot
be used globally: it becomes singular and then negative.  We stop it when
it reaches two, at $\zeta_r=1/2$, and keep it equal to two afterwards.
Thus the resulting energy is uniformly equivalent to the physical one.

For any positive $F$,
$2\sqrt F\,|v_rv_\theta|\leq |v_r|^2+F|v_\theta|^2$.
The remaining scalar maximization is therefore exactly the one in
\eqref{seed:sharp-constant}.  The corner in the weight is harmless for a
specific reason: its radial derivative jumps down, making the interface
contribution to the diffusion nonpositive.

The sharp coefficient is needed in the comparison of upper and lower
growth bounds in Lemma~\ref{seed:signed-history}.  Upper estimates will carry $e^{c_\star A}$, whereas a
lower bound carries $e^{GA}$ with $G<c_\star$ close to $c_\star$.
Their ratio costs $e^{(c_\star-G)A}$.  Replacing $c_\star$ by an unspecified
larger constant would leave a loss that the subsequent scale choices need
not absorb.

\begin{lemma}[A sharp physical energy estimate]\label{seed:physical-energy}
There is a time-dependent quadratic energy $E_k$, between the physical
Fourier energy and twice that energy, such that solutions of
\eqref{seed:fiber-evolution} satisfy
\begin{equation}\label{seed:sharp-energy}
 E_k'(t)\leq
 \Bigl(\frac{2c_\star A}{(1+t)^2}+\frac{5}{1+t}-2k^2\Bigr)E_k(t).
\end{equation}
After integration in frequency, the same estimate holds without the last term.
Moreover, the projected coupling $\mathcal B_A(t,k)$ in physical velocity
coordinates obeys
\begin{equation}\label{seed:small-frequency-coupling}
 \|\mathcal B_A(t,k)\|
 \leq\min\Bigl(\frac{2A}{\tau_t^2},
                  \frac{\sqrt{2}A|k|}{\tau_t^{3/2}}\Bigr),
 \qquad
 \int_0^\infty\|\mathcal B_A(t,k)\|\idiff t
 \leq 2\sqrt{2}A|k|.
\end{equation}
\end{lemma}

\begin{proof}
For $\zeta_r=r^2/(4\tau_t)$, use the weight
\begin{equation}\label{seed:energy-weight}
 F_\star(\zeta_r)=
 \begin{cases}
 (1-\zeta_r)^{-1},&0\leq\zeta_r\leq\tfrac12,\\
 2,&\zeta_r\geq\tfrac12,
 \end{cases}
\end{equation}
on the azimuthal component, and weight one on both meridional components.
Since $1\leq F_\star\leq 2$, the resulting energy has the stated equivalence.
The pressure cancels in the unweighted solenoidal meridional energy.
After division by the angular factor $2\pi$, the coupling contribution
to one half of the energy derivative is
\begin{equation}\label{seed:weighted-coupling}
 A\int_0^\infty
       (2\Omega_{\tau_t}-F_\star B_{\tau_t})
       \operatorname{Re}(v_r\overline{v_\theta})r\idiff r.
\end{equation}
Its coefficient vanishes for $\zeta_r<1/2$. In the remaining region,
\begin{equation}\label{seed:sharp-constant}
 \frac{|2\Omega_{\tau_t}-F_\star B_{\tau_t}|}{2\sqrt{F_\star}}
 =\tau_t^{-2}\frac{(2\zeta_r-1)\mathrm e^{-\zeta_r}}{\sqrt2}
 \leq c_\star\tau_t^{-2},
\end{equation}
with equality in the last scalar maximization at $\zeta_r=3/2$.

Integration by parts in the weighted azimuthal diffusion contributes the
negative gradient form and the scalar radial Laplacian of the weight. In the
inner region, direct differentiation gives
\begin{equation}\label{seed:weight-derivative}
 (\partial_t+\partial_r^2+r^{-1}\partial_r)F_\star(\zeta_r)
 =\tau_t^{-1}\bigl(\zeta_r F_\star''+(1-\zeta_r)F_\star'\bigr)
 \leq 5\tau_t^{-1}F_\star.
\end{equation}
At $\zeta_r=1/2$, the radial derivative jumps downward; the distributional
Laplacian therefore has an additional nonpositive measure. The time derivative
has no interface measure, because the weight is continuous. Combining
\eqref{seed:sharp-constant} and \eqref{seed:weight-derivative} proves
\eqref{seed:sharp-energy}. Axial diffusion contributes exactly $-2k^2E_k$,
since the weight is independent of $z$.

For \eqref{seed:small-frequency-coupling}, the physical coupling is
\begin{equation}\label{seed:physical-coupling}
 \mathcal B_A(t,k)v
 =A\mathbb P_k(2\Omega_{\tau_t} v_\theta e_r)-AB_{\tau_t} v_re_\theta.
\end{equation}
The two outputs are orthogonal, and their multiplier bounds give
$2A\tau_t^{-2}$. To exploit small $|k|$, define
\begin{equation}\label{seed:radial-potential}
 p_v(r)=-\int_r^\infty 2\Omega_{\tau_t}(\rho)v_\theta(\rho)\idiff\rho.
\end{equation}
Then $\mathbb P_k(\nabla_hp_v,0)=-\mathrm ik\mathbb P_k(p_ve_z)$.
The squared Hilbert--Schmidt norm of $v_\theta\mapsto p_v$ on $H_{\rm rad}$ is
\begin{equation}\label{seed:hilbert-schmidt}
 \int_0^\infty\frac{4\Omega_{\tau_t}(\rho)^2}{\rho}
       \int_0^\rho u\idiff u\idiff\rho=2\tau_t^{-3}.
\end{equation}
For the second output, incompressibility and regularity at the axis yield
\begin{equation}\label{seed:small-frequency-radial}
 v_r(r)=-\frac{\mathrm ik}{r}\int_0^r\rho v_z(\rho)\idiff\rho,
 \qquad |v_r(r)|\leq\frac{|k|}{\sqrt2}\|v_z\|_{H_{\mathrm{rad}}}.
\end{equation}
Since $\|B_{\tau_t}\|_{H_{\mathrm{rad}}}=\sqrt2\tau_t^{-3/2}$, these estimates give the second
bound inside the minimum in \eqref{seed:small-frequency-coupling}. Integration
of $\tau_t^{-3/2}$ proves its final assertion.
\end{proof}

\subsection{The growing eigenspace and the viscous complementary estimate}

The variational problem identifies a simple growing direction for the
frozen coupling.  The next lemma has a different task: to follow this
direction as the Gaussian diffuses, without losing control of high heat
frequencies on its complement.  Smooth finite-rank changes of metric
make these two requirements compatible.

We next identify the growing eigenspace of the frozen inviscid operator. The
viscous growth estimate will then follow from a metric on its complement;
no time-dependent change of spatial variables is applied to an arbitrary
evolving perturbation.

At time zero, an eigenvector with eigenvalue $\sigma\ne0$ satisfies
$\sigma f_\theta=-kB_1f_r$.  Substitution into the radial equation gives
\[
 \sigma^2L_kf_r=-2k^2\Omega_1B_1f_r=-k^2\Phi f_r.
\]
Although the original two-component operator is not self-adjoint, this
squared-eigenvalue equation is self-adjoint in the form inner product.
Its Rayleigh quotient can be written as
\[
 \frac{\int_0^\infty(-\Phi)|h|^2r\,\idiff r}
      {\|h\|_{H_{\mathrm{rad}}}^2+k^{-2}\|D_rh\|_{H_{\mathrm{rad}}}^2}.
\]
The numerator is positive for suitable functions supported outside $r=2$.
As $k$ increases, the derivative penalty decreases.  Concentrating a trial
function near the maximum of $-\Phi$ therefore explains why the leading
rate approaches $c_\star$; the proof below makes this observation precise.

There are two distinct choices of size.  We first choose a sufficiently
high but fixed compact band of axial frequencies.  Only afterwards do we
increase $A$.  On this band, viscous and eigenspace-motion errors have
fixed bounds, whereas the leading growth is multiplied by $A$.

\begin{lemma}[The radial variational problem]\label{seed:rayleigh}
For each $k>0$, there are a unique positive function $h_k$, normalized by
$\|h_k\|_{H_{L_k}}=1$, and a number $\sigma_k>0$ such that
\begin{equation}\label{seed:rayleigh-equation}
 \sigma_k^2L_kh_k=-k^2\Phi h_k,
 \qquad \Phi(r)=(4-r^2)\mathrm e^{-r^2/2},
\end{equation}
and $\sigma_k^2$ is the largest positive eigenvalue of this self-adjoint
form-space problem. The numbers $\sigma_k$ satisfy
\begin{equation}\label{seed:rate-limit}
 0<\sigma_k<c_\star,\qquad
 \sigma_k\text{ is strictly increasing in }k,\qquad
 \lim_{k\to\infty}\sigma_k=c_\star.
\end{equation}
The largest real part in the spectrum of $\mathcal C_k(0)$ is the algebraically
simple eigenvalue $\sigma_k$. Its eigenfunction and spectral projection are
smooth in $k>0$, with uniform bounds in each fixed graph norm on compact
frequency intervals.
\end{lemma}

\begin{proof}
The relevant variational formula is
\begin{equation}\label{seed:rayleigh-quotient}
 \sigma_k^2=\sup_{0\ne h\in H_{L_k}}
 \frac{-k^2\int_0^\infty\Phi|h|^2r\idiff r}
      {\|D_rh\|_{H_{\mathrm{rad}}}^2+k^2\|h\|_{H_{\mathrm{rad}}}^2}.
\end{equation}
Multiplication by $\Phi$ is compact from the form space to its dual. Indeed,
on a bounded horizontal disk this follows from Rellich compactness in the
regular angular-order-one sector; outside that disk the operator norm is
bounded by a constant times the supremum of $|\Phi|$. Therefore
$-k^2L_k^{-1}\Phi$ is compact and self-adjoint in the form inner product.
A smooth test supported in $r>2$ makes the numerator positive, so the positive supremum is attained. The forms have real
coefficients, so a maximizer can be chosen real: its Rayleigh quotient is
a denominator-weighted average of those of its real and imaginary parts.
For a real form-space function, replacing it by its absolute value preserves
the potential and derivative forms. The radial differential equation then excludes an interior zero.
At the axis its regular expansion is
\begin{equation}\label{seed:axis-expansion}
 h_k(r)=h_k'(0)r+O(r^3),\qquad h_k'(0)>0;
\end{equation}
the other local solution has order $r^{-1}$ and is excluded by the form norm.
For two regular solutions, the constant Wronskian
$r(h_1h_2'-h_1'h_2)$ vanishes at the axis. This proves simplicity.

The exact maximum of $-\Phi$ is $c_\star^2=2\mathrm e^{-3}$. The derivative
form in the denominator of \eqref{seed:rayleigh-quotient} is positive for any
nonzero admissible $h$, giving $\sigma_k<c_\star$. For a fixed test with
positive numerator, its quotient is strictly increasing in $k$. Testing with
the maximizer at the smaller frequency proves strict monotonicity.

For a quantitative lower bound, take the compactly supported function
\begin{equation}\label{seed:rayleigh-test}
 h_\delta(r)=\Bigl(1-\frac{(r^2-6)^2}{\delta^2}\Bigr)_+,
 \qquad 0<\delta\leq1.
\end{equation}

The trial function is centered at the most unstable radius $r^2=6$.
Its width is measured in the variable $x=r^2-6$, not in $r$ itself.
Since $r\,\idiff r=\idiff x/2$, its norm is
\[
 \|h_\delta\|_{H_{\mathrm{rad}}}^2
  =\frac12\int_{-\delta}^{\delta}(1-x^2/\delta^2)^2\,\idiff x
  =\frac{8\delta}{15}.
\]
In the support interior $h_\delta'=-4rx/\delta^2$, and therefore
\[
 \int_0^\infty|h_\delta'|^2r\,\idiff r
  =\frac8{\delta^4}\int_{-\delta}^{\delta}(6+x)x^2\,\idiff x
  =\frac{32}{\delta}.
\]
The function is a valid form-space test despite the corner at the ends
of its support; smooth approximations give the same quotient in the limit.
The potential loses $O(\delta^2)$ from its maximum, while the derivative
penalty relative to $k^2\|h_\delta\|_{H_{\mathrm{rad}}}^2$ is
$O((k\delta)^{-2})$.  Balancing these two losses suggests
$\delta$ of order $k^{-1/2}$ and an $O(k^{-1})$ error in the squared
rate.  The explicit choice below realizes this balance with fixed
constants.

Its squared $H_{\rm rad}$ norm and derivative energy are $8\delta/15$ and $32/\delta$,
respectively, while its inverse-square energy is at most
$(6-\delta)^{-1}\|h_\delta\|_{H_{\mathrm{rad}}}^2$. For $|x|\leq\delta$,
\begin{equation}\label{seed:rayleigh-test-potential}
 \frac{-\Phi(\sqrt{6+x})}{c_\star^2}
 =\Bigl(1+\frac{x}{2}\Bigr)\mathrm e^{-x/2}
 \geq1-\frac{x^2}{4}.
\end{equation}
To obtain this inequality, integrate the derivative
$-x\mathrm e^{-x/2}/4$ from zero and use $\mathrm e^{1/2}<2$.
With $\delta=4/\sqrt{k}$ and $k\geq16$, substitution into
\eqref{seed:rayleigh-quotient} yields
\begin{equation}\label{seed:quantitative-rate-limit}
 \frac{\sigma_k^2}{c_\star^2}
 \geq\frac{1-4/k}{1+15/(4k)+1/(5k^2)}.
\end{equation}
This proves the limit in \eqref{seed:rate-limit}.

Conjugation by $(f_r,f_\theta)\mapsto(L_k^{1/2}f_r,f_\theta)$ gives
$\mathcal C_k(0)$ the two off-diagonal blocks
$2kL_k^{-1/2}\Omega_1$ and $-kB_1L_k^{-1/2}$. They are compact, by local
Rellich compactness and Gaussian localization. Every nonzero spectral value
$\lambda$ therefore has an eigenvector whose radial coordinate satisfies
$\lambda^2L_kf_r=-k^2\Phi f_r$. Thus $\lambda$ is real or purely imaginary,
and its largest real part is $\sigma_k$. A generalized eigenvector at
$\sigma_k$ would give
\begin{equation}\label{seed:no-jordan-chain}
 \bigl(-k^2L_k^{-1}\Phi-\sigma_k^2\bigr)f_r=2\sigma_kh_k,
\end{equation}
contradicting self-adjointness after pairing with $h_k$. There is consequently
a strict spectral gap between $\sigma_k$ and the remaining spectrum.

Elliptic regularity in \eqref{seed:rayleigh-equation} places the eigenvector
and its adjoint in every fixed heat graph domain. On a compact frequency
interval, the exterior potential is uniformly small relative to $k^2$.
An exponential supersolution $\mathrm e^{-kr/2}$, followed by interior
regularity, bounds all fixed radial derivatives uniformly and exponentially
at infinity. A resolvent contour surrounding the simple eigenvalue gives
smooth dependence on $k$; differentiation of the radial equation and use of
the reduced resolvent give the same dependence in graph norms. These bounds
also control the factors of $r$ that arise on differentiating a dilation.
\end{proof}

Write $j_k=-kB_1h_k/\sigma_k$. At positive heat age, define
\begin{equation}\label{seed:growing-rate}
 \kappa=k\sqrt{\tau_t},\qquad
 \lambda_k(t)=\tau_t^{-2}\sigma_\kappa,\qquad
 \mathscr G_k(t)=\int_0^t\lambda_k(s)\idiff s.
\end{equation}
The corresponding eigenvector and adjoint vector are
\begin{equation}\label{seed:eigenvectors}
 \varphi_k(t)=
 \begin{pmatrix}
 h_\kappa(r/\sqrt{\tau_t})\\
 \tau_t^{-1/2}j_\kappa(r/\sqrt{\tau_t})
 \end{pmatrix},
 \qquad
 \varphi_k^\dagger(t)=\frac12
 \begin{pmatrix}
 h_\kappa(r/\sqrt{\tau_t})\\
 2k\Omega_{\tau_t}(r)h_\kappa(r/\sqrt{\tau_t})/\lambda_k(t)
 \end{pmatrix}.
\end{equation}
Indeed, the map
$(f_r,f_\theta)\mapsto(f_r(r/\sqrt{\tau_t}),
\tau_t^{-1/2}f_\theta(r/\sqrt{\tau_t}))$ is an isometry from $X_\kappa$
to $X_k$ and conjugates the frozen generator to
$\tau_t^{-2}\mathcal C_\kappa(0)$. In particular,
$\langle\varphi_k,\varphi_k^\dagger\rangle_{X_k}=1$. We use the coordinate
and projections
\begin{equation}\label{seed:spectral-projections}
 \ell_k(t)V=\langle V,\varphi_k^\dagger(t)\rangle_{X_k},\qquad
 \Pi_k(t)V=\varphi_k(t)\ell_k(t)V,\qquad Q_k(t)=1-\Pi_k(t).
\end{equation}
These are spectral projections; $Q_k$ need not be orthogonal.

For an integer $m\geq0$, let
$X_k^m=D(\mathcal D_k^{m/2})$, with the domain taken in $X_k$.
This graph space records the corresponding physical regularity, with norm
\begin{equation}\label{seed:graph-norm}
 \|V\|_{X_k^m}^2
  =\|L_k^{(m+1)/2}f_r\|_{H_{\mathrm{rad}}}^2+\|L_k^{m/2}f_\theta\|_{H_{\mathrm{rad}}}^2,
 \qquad V=(f_r,f_\theta).
\end{equation}
These graph norms therefore measure the Sobolev regularity of the
physical velocity.

For a bounded operator $C$ with an isolated simple eigenvalue $\lambda$,
its spectral projection is
\[
 \Pi=\frac1{2\pi i}\int_{\mathcal C}(z-C)^{-1}\,\idiff z,
\]
where the positively oriented contour surrounds that eigenvalue and no
other spectrum.  This projection satisfies $\Pi^2=\Pi$ and $C\Pi=\Pi C$,
but need not be orthogonal.  It extracts the growing coordinate using the
adjoint eigenvector, not an orthogonal projection onto the eigenvector.
Smooth parameter dependence follows locally by differentiating the
resolvent identity; in particular
$\partial(z-C)^{-1}=(z-C)^{-1}(\partial C)(z-C)^{-1}$.
The contours are chosen locally and then on a finite cover of the compact
parameter set, rather than assumed to be one fixed contour for all parameters.

Finally, a positive metric $R$ means the energy $\langle Rq,q\rangle$.
The desired estimate on the complement is a \emph{relative} growth gap:
its rate is strictly below the leading rate, not necessarily below zero.
A spectral gap alone does not give that estimate in the original energy
for a non-self-adjoint operator.  The Lyapunov integral in the proof
constructs an energy in which the gap becomes an energy inequality.

An arbitrary equivalent energy can interact badly with the unbounded
operator $\mathcal D_k$.  In contrast, for one smooth rank-one term
$Kq=u\langle q,v\rangle$, with $v\in D(\mathcal D_k)$,
\[
 K\mathcal D_kq=u\langle q,\mathcal D_kv\rangle,
 \qquad
 \|K\mathcal D_kq\|\leq\|u\|\,\|\mathcal D_kv\|\,\|q\|.
\]
Only a bounded error remains after the heat operator is moved to the
adjoint vector.  A finite sum has the same property.  This is why the
proof first constructs a coercive metric, then approximates its compact
correction by smooth finite rank, and only then handles diffusion.
The same argument must be made in each fixed graph norm; the adjoint
vector in that norm is not automatically the adjoint in $X_k$.

\begin{lemma}[A complementary metric compatible with diffusion]\label{seed:metric}
Fix a compact interval $J\subset(0,\infty)$, a time $t_0>0$, and an integer
$m\geq0$. There are smooth positive metrics
$\mathcal R_k^m(t)=1+\mathcal F_k^m(t)$ on $X_k^m$, uniformly equivalent to
its original metric, where $\mathcal F_k^m$ is self-adjoint and finite rank,
and a number $g_m>0$ such that
\begin{equation}\label{seed:metric-inequalities}
 \begin{aligned}
 \operatorname{Re}\langle\mathcal R_k^m\mathcal C_kq,q\rangle_{X_k^m}
 &\leq(\lambda_k-g_m)
       \langle\mathcal R_k^mq,q\rangle_{X_k^m},\\
 \operatorname{Re}\langle\mathcal R_k^mQ_k\mathcal D_kq,q\rangle_{X_k^m}
 &\geq-C_m\|q\|_{X_k^m}^2,
 \end{aligned}
\end{equation}
where the first inequality holds for $q\in Q_k(t)X_k^m$ and the
second is initially an operator pairing for $q\in Q_k(t)X_k^{m+2}$,
with a quadratic-form extension to $Q_k(t)X_k^{m+1}$.
All fixed metric derivatives and the required graph norms of the finite-rank
ranges and adjoint ranges are bounded independently of $A$. If
$a=\ell_kV$, $q=Q_kV$, and
$q_m=\langle\mathcal R_k^mq,q\rangle_{X_k^m}^{1/2}$, then
\begin{equation}\label{seed:moving-coordinate-estimates}
 \begin{aligned}
 |a'-A\lambda_ka|&\leq C(|a|+q_0),\\
 D^+q_m&\leq\bigl(A(\lambda_k-g_m)+C_m\bigr)q_m+C_m|a|,
 \end{aligned}
\end{equation}
where $D^+$ denotes the upper right Dini derivative. In particular, for
arbitrary physical input,
\begin{equation}\label{seed:band-propagator}
 \|\mathcal U_{\mathrm d,k}(t,s)\|_{2\to2}
 \leq C\exp\Bigl(A\int_s^t\lambda_k(u)\idiff u+C(t-s)\Bigr),
 \qquad k\in J,\quad0\leq s\leq t\leq t_0.
\end{equation}
\end{lemma}

\begin{proof}
\textit{The spectral splitting on the diffusion graph spaces.}
The diffusion graph spaces are defined by \eqref{seed:graph-norm}.
We first show that the spectral projections of the bounded coupling
$\mathcal C_k(t)$ extend consistently to these spaces.

Fix the graph order $m$ before choosing any of the large
parameters, and write $H_m=X_k^m$ temporarily.
The operator $\mathcal D_k$ is positive and self-adjoint
on $H_m$, with operator domain $X_k^{m+2}$ and form
domain $X_k^{m+1}$. On a compact set of positive $k$,
the graph norms are uniformly equivalent to their
inhomogeneous versions.

The two off-diagonal blocks of $\mathcal C_k(t)$ together
gain one derivative in the graph scale: the upper block
gains two radial derivatives, whereas the first coordinate
already carries one more derivative than the second.
Multiplication by the Gaussian coefficients, local Rellich
compactness, and the exponentially decreasing kernel of
$L_k^{-1}$ show that $\mathcal C_k(t)$ is compact on every
$X_k^m$. To justify this compactness, first truncate both
the coefficient and the output to a fixed radial disk.
The one-derivative gain makes this truncated map compact.
The Gaussian coefficient tail tends to zero in every
required multiplier norm; the resolvent kernel across the
remaining radial support gap makes the output tail tend
to zero in operator norm. Both errors are uniform on
the compact parameter set.

Nonzero generalized eigenvectors bootstrap through all graph orders.
For a complex resolvent parameter $\zeta\ne0$, a solution of
$(\zeta-\mathcal C_k)V=F$ satisfies
\begin{equation}\label{seed:resolvent-bootstrap}
 V=\zeta^{-1}F+\zeta^{-1}\mathcal C_kV.
\end{equation}
If $F\in X_k^m$, successive application of this identity upgrades a
solution in $X_k$ to $X_k^m$.  Compactness and the Fredholm alternative
then identify the nonzero spectrum in the different graph spaces; the
bounded inverse on a graph space follows from the closed graph theorem.
Thus the simple growing eigenvalue and its Riesz projection are the
same operators on their common dense smooth domain.  The adjoint of
the coordinate functional in $X_k^m$ need not equal its adjoint in
$X_k$, but it is smooth: it is obtained from the latter by applying
$\mathcal D_k^{-m}$.  This observation is important when an unbounded
heat operator is moved to an adjoint range.

\textit{A coercive metric at one parameter value.}
First fix $(k,t)$ and let $Q=Q_k(t)$ and
$\mathcal C_Q=\mathcal C_k(t)|_{QH_m}$.  The range $QH_m$ is a
closed Hilbert subspace with the inherited inner product.  Choose
\begin{equation}\label{seed:metric-gap-choice}
 \max\bigl(0,\sup\operatorname{Re}\operatorname{spec}\mathcal C_Q\bigr)
       <b_\star<\lambda_k(t).
\end{equation}

The semigroup bound used here concerns the bounded operator
$\mathcal C_Q$, not the full unbounded viscous generator.
For a bounded operator, spectral mapping gives
$r(e^{\mathcal C_Q})=e^{\sup\operatorname{Re}\operatorname{spec}\mathcal C_Q}$.
Choose $b'$ strictly between that spectral abscissa and $b_\star$.
The spectral-radius formula implies that, for all sufficiently large
integers $n$, $\|(e^{\mathcal C_Q})^n\|\leq e^{b'n}$.
The finitely many smaller integers are included in one multiplicative
constant.  Writing $u=n+\theta$, $0\leq\theta<1$, and bounding
$e^{\theta\mathcal C_Q}$ yields the continuous-time estimate below.
This argument permits transient amplification in a non-self-adjoint
system.  Diffusion is handled afterwards by the smooth finite-rank
structure of the metric.

We have therefore obtained
$\|e^{u\mathcal C_Q}\|\leq C e^{b'u}$ for $u\geq0$.
Define, on $QH_m$,
\begin{equation}\label{seed:lyapunov-metric}
 \mathcal S=2b_\star\int_0^\infty e^{-2b_\star u}
                e^{u\mathcal C_Q^*}e^{u\mathcal C_Q}\idiff u.
\end{equation}
The integral converges in operator norm.  Since a bounded generator has
an invertible exponential and
$\|e^{u\mathcal C_Q}q\|\geq e^{-u\|\mathcal C_Q\|}\|q\|$,
\[
 \langle\mathcal S q,q\rangle
 \geq \frac{b_\star}{b_\star+\|\mathcal C_Q\|}\|q\|^2.
\]
Together with the norm convergence of \eqref{seed:lyapunov-metric},
this proves equivalence with the original norm.  Integration of the derivative of the integrand gives
\begin{equation}\label{seed:lyapunov-identity}
 \mathcal C_Q^*\mathcal S+\mathcal S\mathcal C_Q
       =2b_\star\mathcal S-2b_\star I.
\end{equation}
There is no boundary term at infinity because $b'<b_\star$.

Since $\mathcal C_Q$ is compact, $e^{u\mathcal C_Q}-I$ is compact.
Subtracting the integral of $e^{-2b_\star u}I$ shows that
$\mathcal S-I$ is compact and self-adjoint.  Choose finite-dimensional
orthogonal projections $P_n$ in $QH_m$, with smooth ranges, converging
strongly to the identity of that subspace.  Such ranges exist because
$Q$ maps the dense smooth graph core into itself.  Compactness implies
\[
 \|P_n(\mathcal S-I)P_n-(\mathcal S-I)\|\longrightarrow0.
\]
Set $\mathcal F=P_n(\mathcal S-I)P_n$ and $\mathcal R_0=I+\mathcal F$.
If the approximation error is $\varepsilon$, then
\begin{equation}\label{seed:finite-rank-margin}
 \operatorname{Re}\langle\mathcal R_0\mathcal C_Qq,q\rangle
       -b_\star\langle\mathcal R_0q,q\rangle
 \leq\bigl(-b_\star+
       \varepsilon(\|\mathcal C_Q\|+b_\star)\bigr)\,\|q\|^2.
\end{equation}
Choose $\varepsilon$ smaller than half the coercivity constant and
than $\frac{b_\star}{2(\|\mathcal C_Q\|+b_\star)}$.
The resulting metric is positive, equivalent to the original norm,
and has a strict complementary growth bound below $\lambda_k(t)$.
Extend $\mathcal F$ by zero on the orthogonal complement of $QH_m$.
This extension is self-adjoint and finite rank on all of $H_m$.

\textit{Smoothness and uniformity in the parameters.}
The spaces at different frequencies are first identified by the isometry
\begin{equation}\label{seed:graph-identification}
 (f_r,f_\theta)\longmapsto
 \bigl(L_{k_\circ}^{-(m+1)/2}L_k^{(m+1)/2}f_r,
       L_{k_\circ}^{-m/2}L_k^{m/2}f_\theta\bigr).
\end{equation}
Here $L_k=L_{k_\circ}+(k^2-k_\circ^2)I$, so the operators commute.
The spectral calculus shows that the identification and its inverse
are smooth, with bounded derivatives on compact positive frequency
intervals, in every fixed graph norm.  It conjugates $\mathcal D_k$
to $\mathcal D_{k_\circ}+(k^2-k_\circ^2)I$.
All projections in the next formula are understood after this identification.
Near a fixed pair $(k_\circ,t_\circ)$, put
\begin{equation}\label{seed:projection-transport}
 \mathcal T=\Pi_k(t)\Pi_{k_\circ}(t_\circ)
                  +Q_k(t)Q_{k_\circ}(t_\circ).
\end{equation}
At the central pair $\mathcal T=I$; nearby it is invertible and
$\Pi_k(t)\mathcal T=\mathcal T\Pi_{k_\circ}(t_\circ)$.
Moreover, $\mathcal T-I$ and $\mathcal T^{-1}-I$ are smooth finite-rank
operators with smooth ranges and adjoint ranges.  Pull back the fixed
metric by $\mathcal T^{-1}$, that is, use
$\mathcal T^{-*}\mathcal R_0\mathcal T^{-1}$.
The difference of this metric from the identity is still finite rank.
The transported generator differs continuously in operator norm from
the generator at the central pair. The strict margin in
\eqref{seed:finite-rank-margin} therefore persists on a neighborhood, with a smaller positive
gap.  A finite cover of $J\times[0,t_0]$ and a smooth nonnegative
partition of unity give a global metric by convex combination.
The inequalities are linear in the metric, so the minimum of the
finitely many gaps is a valid $g_m>0$.  Positivity, upper bounds, rank
bounds, and every fixed graph bound on its ranges and derivatives are
uniform.  None of these choices involves $A$.

\textit{Compatibility with the unbounded diffusion.}
This step is initially performed for
$q\in Q_k(t)X_k^{m+2}$, not for an arbitrary vector of $X_k^m$.
The identity part of the metric gives
\begin{equation}\label{seed:heat-metric-cost}
 \langle\mathcal D_kq,q\rangle_{X_k^m}
       -\langle\Pi_k\mathcal D_kq,q\rangle_{X_k^m}.
\end{equation}
A convenient way to include the remaining terms is to set
$\mathcal K=\mathcal R_k^mQ_k-I$.
This is finite rank, and its adjoint range is in $X_k^{m+2}$ with
uniform graph bounds.  For example, writing
$\mathcal Kq=\sum_{j=1}^N u_j\langle q,v_j\rangle_{X_k^m}$
with smooth $u_j,v_j$ gives
\[
 \mathcal K\mathcal D_kq
       =\sum_{j=1}^N u_j\langle q,\mathcal D_kv_j\rangle_{X_k^m}.
\]
Consequently $\mathcal K\mathcal D_k$ extends to a bounded operator
on $X_k^m$, and
\[
 \operatorname{Re}\langle\mathcal R_k^mQ_k\mathcal D_kq,q\rangle_{X_k^m}
 \geq\|\mathcal D_k^{1/2}q\|_{X_k^m}^2-C_m\|q\|_{X_k^m}^2.
\]
This also defines the pairing as a closed quadratic form on
$Q_k(t)X_k^{m+1}$.  It proves the second inequality in
\eqref{seed:metric-inequalities} and explains why large heat frequencies
remain dissipative.

\textit{The evolving decomposition.}
For a graph-smooth solution write $V=a\varphi_k+q$, where
$a=\ell_kV$ and $q=Q_kV$.  Differentiation gives exactly
\begin{equation}\label{seed:exact-moving-equations}
 \begin{aligned}
 a'&=A\lambda_ka+\ell_k'V-\ell_k\mathcal D_kV,\\
 q'&=A\mathcal C_kq-Q_k\mathcal D_kq
       +a(Q_k'\varphi_k-Q_k\mathcal D_k\varphi_k)+Q_k'q.
 \end{aligned}
\end{equation}
The coordinate $\ell_k\mathcal D_k$ is bounded on $X_k$ because
its adjoint vector is $\mathcal D_k\varphi_k^\dagger$.
All remaining finite-rank and metric-derivative terms are bounded
in the indicated graph space.  Differentiating $q_m^2$, retaining
the favorable heat form, and using norm equivalence yields
\eqref{seed:moving-coordinate-estimates} when $q_m>0$.
At a zero of $q_m$, the same statement follows by regularizing its
square root and taking the upper right Dini derivative.
There is no term of order $Aa$ in the complementary equation, since
$Q_k\mathcal C_k\varphi_k=0$.
For $m=0$ the combined norm $(|a|^2+q_0^2)^{1/2}$ therefore has
upper logarithmic growth rate $A\lambda_k+C$.
Integration and equivalence with the physical norm give
\eqref{seed:band-propagator}.
Finally, approximate arbitrary physical input by graph-domain input.
The bounded-coefficient parabolic equation is strongly continuous
and smooths at positive elapsed time; its energy estimate passes to
the limit. This extends the propagator bound to the full physical fiber.
\end{proof}

\subsection{Growth of the prescribed point response}

The incoming field is determined by the quadratic output of the preceding
scale.  We compute its projection onto the adjoint growing vector,
retaining the distributional contribution at the symmetry axis.
The subsequent short-time smoothing and complementary estimate then
give positivity and alignment for this prescribed input.

The input that enters the transfer calculation is the distribution
\begin{equation}\label{seed:point-input}
 \mathcal I=\PP\partial_z(e_z\delta_0),\qquad
 v_A(t)=\mathcal U_{\mathrm d}(t,0)\mathcal I.
\end{equation}
Its Fourier symbol is bounded by $C|\xi|$, so $\mathcal I\in H^{-3}(\R^3)$.
On axisymmetric vectors, the linearized equation has the Cartesian form
\begin{equation}\label{seed:cartesian-linearization}
 (\partial_t-\Delta)v_A=-\PP\mathsf M_A(t)v_A,
 \qquad
 \mathsf M_A(t)v=A\Bigl(2\Omega_{\tau_t} v_h^\perp
       -\frac{\Omega_{\tau_t}}{2\tau_t}x_h^\perp(x_h\cdot v_h)\Bigr).
\end{equation}
The sign agrees with \eqref{bg:cartesian-matrix}: for the invariant
column, $M_{U_{\mathrm d}}=-\mathsf M_A$ and hence
$\mathcal L_{U_{\mathrm d}}=\Delta+\PP M_{U_{\mathrm d}}
=\Delta-\PP\mathsf M_A$.
The vertical component of $\mathsf M_Av$ is zero. This matrix is smooth across
the axis, and its derivatives of every fixed order are bounded by $C A$ on a
fixed time interval. Bounded perturbation of the heat semigroup in $H^{-3}$
therefore constructs \eqref{seed:point-input} uniquely. At positive time it is
smooth, finite-energy, real, centrally odd, axisymmetric, and solenoidal.

\begin{lemma}[The growing coordinate of the point response]\label{seed:point-growth}
Fix $J\Subset(0,\infty)$ and $t_0>0$. The coordinate
$a_{A,k}(t)=\ell_k(t)V_{A,k}(t)$ of \eqref{seed:point-input} has initial
distributional value
\begin{equation}\label{seed:positive-point-coordinate}
 a_0(k)=\frac{k h_k'(0)}{2\pi\sqrt{2\pi}}>0.
\end{equation}
For sufficiently large $A$, uniformly in $k\in J$,
\begin{equation}\label{seed:point-growth-bound}
 a_{A,k}(t)>0,\qquad
 \left|\log\frac{a_{A,k}(t)}{a_0(k)}-A\mathscr G_k(t)\right|
 \leq Ct+CA^{-1/2},\qquad A^{-1}\leq t\leq t_0.
\end{equation}
For each fixed graph order $m$, there is $\Theta_m>0$ such that
\begin{equation}\label{seed:point-alignment}
 \|V_{A,k}(t)-a_{A,k}(t)\varphi_k(t)\|_{X_k^m}
 \leq C_m A^{-1}a_{A,k}(t),
 \qquad\frac{1+\Theta_m\log A}{A}\leq t\leq t_0.
\end{equation}
\end{lemma}

\begin{proof}
The unitary axial transform of the impulse is
$\widehat{\mathcal I}^{\,z}(k)=\mathrm ik\mathbb P_k(e_z\delta_{0,h})/
\sqrt{2\pi}$. Let $G_k(r)=(2\pi)^{-1}K_0(kr)$ be the fundamental solution
of $-\Delta_h+k^2$, where $K_0,K_1$ are the modified Bessel functions. Directly
applying the Leray projection gives
\begin{equation}\label{seed:point-fiber}
 \mathbb P_k(e_z\delta_{0,h})
 =\mathrm ikG_k'e_r+(\delta_{0,h}-k^2G_k)e_z,
 \qquad
 f_{r,0,k}(r)=\frac{k^2}{2\pi\sqrt{2\pi}}K_1(kr),\quad f_{\theta,0,k}=0.
\end{equation}

The fundamental solution satisfies
$(-\Delta_h+k^2)G_k=\delta_{0,h}$.  The vertical component in
\eqref{seed:point-fiber} therefore has a Dirac contribution even though
the Bessel profile is written as an ordinary function for $r>0$.
Indeed $K_1(kr)\sim(kr)^{-1}$ near the axis, and the radial vector
$r^{-1}e_r$ has distributional divergence $2\pi\delta_{0,h}$ in two
horizontal dimensions.  Differentiating only its expression away from
the axis would lose this mass.

To retain this contribution, pair the full projected distribution with the
smooth solenoidal adjoint.  Self-adjointness of $\PP_k$ moves the projection
to that test, where it is the identity.  Its vertical value is
$iD_rh_k(0)/2=ih_k'(0)$, since
$h_k(r)=h_k'(0)r+O(r^3)$.  With the convention that the inner product is
linear in its first argument, the product of the $ik$ from the input and
the conjugate of this value is positive.  The next display keeps also the
angular factor $2\pi$ and the unitary Fourier factor $\sqrt{2\pi}$.

The distribution $D_rf_{r,0,k}$ includes a horizontal Dirac mass.
The physical vector associated with $\varphi_k^\dagger(0)$ is
\begin{equation}\label{seed:physical-adjoint}
 \Psi_k^\dagger
 =\bigl(kh_k/2,k\Omega_1h_k/\sigma_k,\mathrm iD_rh_k/2\bigr).
\end{equation}
It is smooth and solenoidal. Pairing the full impulse with this vector yields
\begin{equation}\label{seed:point-pairing}
 2\pi a_0(k)
 =\left\langle\frac{\mathrm ik}{\sqrt{2\pi}}
             \mathbb P_k(e_z\delta_{0,h}),\Psi_k^\dagger\right\rangle
 =\frac{kD_rh_k(0)}{2\sqrt{2\pi}}
 =\frac{k h_k'(0)}{\sqrt{2\pi}}.
\end{equation}
Equivalently, the radial Green identity is
\begin{equation}\label{seed:green-pairing}
 \int_0^\infty K_1(kr)L_kh_k(r)r\idiff r=\frac{2h_k'(0)}{k}.
\end{equation}
Lemma~\ref{seed:rayleigh} now gives a positive minimum for $a_0$ on $J$.
The free heat response also has positive coordinate, since
\begin{equation}\label{seed:free-heat-coordinate}
 \begin{aligned}
 \left\langle\mathrm e^{-t\mathcal D_k}V_{0,k},
                        \varphi_k^\dagger(0)\right\rangle_{X_k}
 &=\frac{k}{2\pi\sqrt{2\pi}}(\mathrm e^{-tL_k}h_k)'(0),\\
 (\mathrm e^{-tL_k}h_k)'(0)
 &=\frac{\mathrm e^{-k^2t}}{8t^2}
       \int_0^\infty\mathrm e^{-r^2/(4t)}h_k(r)r^2\idiff r>0.
 \end{aligned}
\end{equation}
We estimate the actual evolution next.

On $J$, the horizontal Fourier symbol of the initial physical vector is
uniformly bounded. Hence
$\|\mathrm e^{-t\mathcal D_k}V_{0,k}\|_{X_k}\leq C_Jt^{-1/2}$.
The bounded coupling in \eqref{seed:fiber-evolution} and its Duhamel formula
give
\begin{equation}\label{seed:initial-regularization}
 \|V_{A,k}(t)\|_{X_k}\leq C_Jt^{-1/2}\mathrm e^{C_JAt},
 \qquad 0<t\leq A^{-1}.
\end{equation}
Regularizing on the final half of this interval yields, at $t_1=A^{-1}$,
\begin{equation}\label{seed:short-time-graphs}
 \|V_{A,k}(t_1)\|_{X_k^m}\leq C_{J,m}A^{(m+1)/2}.
\end{equation}
For completeness, rescale space by $A^{-1/2}$ and time by $A^{-1}$ in
\eqref{seed:cartesian-linearization}. The rescaled zero-order matrix and its
required derivatives are uniformly bounded. Differentiated heat energies,
multiplied by successive powers of elapsed time, prove
\eqref{seed:short-time-graphs}.

The initial physical vector is uniformly bounded in $H_h^{-2}$ on $J$, so
\eqref{seed:cartesian-linearization} also gives
$\|\widehat v_A^{\,z}(t,k)\|_{H_h^{-2}}\leq C_J\mathrm e^{C_JAt}$.
The exact coordinate equation is
\begin{equation}\label{seed:point-coordinate-equation}
 a'=A\lambda_ka+
       \langle V_{A,k},\partial_t\varphi_k^\dagger
                              -\mathcal D_k\varphi_k^\dagger\rangle.
\end{equation}
The last pairing is distributional duality against a smooth vector with
uniformly bounded dual graph norm. It is $O(1)$ for $0\leq t\leq t_1$.
Consequently,
\begin{equation}\label{seed:coordinate-after-regularization}
 a(t_1)=\mathrm e^{A\mathscr G_k(t_1)}a_0(k)+O(A^{-1})
       =\mathrm e^{\sigma_k}a_0(k)+O(A^{-1}).
\end{equation}
In particular $a(t_1)\geq c_J>0$, whereas
$q_0(t_1)/a(t_1)\leq C_J\sqrt A$.

At $t_1=A^{-1}$ the growing coordinate is bounded below, but the
complement can still be of order $\sqrt A$.
Dividing the complementary estimate by the growing
coordinate exposes the rate difference $Ag_0$.  The resulting quadratic
ratio inequality has an invariant barrier of order $A$, well above the
initial order-$\sqrt A$ ratio.  Below that barrier the negative linear
term dominates the quadratic one, giving exponential relaxation to
order $A^{-1}$.

This bootstrap works first on
the maximal interval where $a>0$ and the ratio stays below the barrier.
There
\[
 \frac{a'}a\geq A\lambda_k-C(1+r_q).
\]
The bound on the time integral of $r_q$ prevents $\log a$ from tending to
$-\infty$ at a finite endpoint, while the differential inequality prevents
escape through the ratio barrier.  Both properties extend the interval.
For the higher graph ratio, an initial polynomial $A^{(m+1)/2}$ is reduced
to order $A^{-1}$ in a time of size $(\log A)/A$.  This explains the
preparation time and why only a finite list of graph orders is fixed
before choosing its constant $\Theta$.

The radial and azimuthal coordinates of the point input in
\eqref{seed:point-fiber} are real.  Equation~\eqref{seed:fiber-evolution}
preserves this real subspace, so $a(t)$ is real.
While $a>0$, set $r_q=q_0/a$. Lemma~\ref{seed:metric} gives
\begin{equation}\label{seed:ratio-inequality}
 D^+r_q\leq-Ag_0r_q+C(1+2r_q+r_q^2).
\end{equation}
For large $A$, its initial value lies below $g_0A/(4C)$, and the right-hand
side is negative at that barrier. Within the barrier, after increasing $A$,
\begin{equation}\label{seed:ratio-decay}
 D^+r_q\leq-\tfrac12Ag_0r_q+C,
 \qquad
 r_q(t)\leq C_J\sqrt A\,
     \mathrm e^{-Ag_0(t-t_1)/2}+C/A,
 \qquad
 \int_{t_1}^{t_0}r_q(t)\idiff t\leq C_JA^{-1/2}.
\end{equation}
The logarithmic form of \eqref{seed:point-coordinate-equation} prevents
$a$ from reaching zero while this finite ratio persists. Integrating that
equation and using \eqref{seed:coordinate-after-regularization} proves
\eqref{seed:point-growth-bound}.

For higher graph norms, the scalar coordinate is already controlled. The
additional logarithmic derivative costs at most $C(1+C_J\sqrt A)$, which
is absorbed by the fixed gap times $A$. Thus we obtain
\begin{equation}\label{seed:higher-ratio-decay}
 \frac{q_m(t)}{a(t)}
 \leq C_m A^{(m+1)/2}\mathrm e^{-Ag_m(t-t_1)/2}+C_m/A.
\end{equation}
Choosing $\Theta_m>(m+3)/g_m$ proves \eqref{seed:point-alignment}.
\end{proof}

\subsection{Weighted regularity of the full prepared response}

A growing band is useful for a lower bound, but it is not the actual
incoming field.  This subsection constructs and estimates the full
prepared response.  Bessel smoothing handles its initial distribution,
ordinary moments control its pressure tail, and a separate Gaussian
azimuthal estimate supplies the positive energy used in the sign proof.

The construction uses all axial frequencies of the point response. Fix
$\Theta>0$ and define
\begin{equation}\label{seed:prepared-definition}
 \tp=\frac{2+\Theta\log A}{A},\qquad
 S_A=v_A(\tp),\qquad n_A=\|S_A\|_2,\qquad g=\frac{S_A}{n_A}.
\end{equation}
When alignment in a finite list of graph norms is needed, choose $\Theta$
larger than the corresponding constants in Lemma~\ref{seed:point-growth}.
This choice is made before increasing $A$.

The point input is singular but belongs to a negative Sobolev space:
its symbol is of order $|\xi|$, and
$\int_{\R^3}(1+|\xi|^2)^{-3}|\xi|^2\,\idiff\xi<\infty$.
Heat smoothing therefore gives a smooth state at a positive time.
It does not make the full projected velocity Gaussian at infinity.
The pressure part has an algebraic tail; after Bessel smoothing its order
$|x|^{-4}$ gives the moment condition
\[
 \int_1^\infty r^{2\mu}|r^{-4}|^2r^2\,\idiff r
     =\int_1^\infty r^{2\mu-6}\,\idiff r<\infty
       \quad\Longleftrightarrow\quad \mu<\frac52.
\]
The separate Gaussian estimate is imposed only on the azimuthal component.
It is generated from zero azimuthal input by a Gaussian-localized source;
there is no azimuthal pressure derivative in the axisymmetric equation.
The proof below uses this structure rather than assuming Gaussian decay
of all velocity components.

At preparation, apparently large short-time estimates have only polynomial
cost.  For fixed $j$ and $C$,
\[
 \tp^{-j/2}e^{CA\tp}
 =\left(\frac{A}{2+\Theta\log A}\right)^{j/2}e^{2C}A^{C\Theta}.
\]
This is compatible with exponential amplification later on a fixed
source-time interval.  In particular, a polynomial norm for $S_A$ does not
contradict the later lower bound $e^{GA}$ for its evolved response.

\begin{lemma}[The prepared point response]\label{seed:prepared}
For every fixed $0\leq\mu<5/2$ and integer $m\geq0$, there are constants
$C_{\mu,m},P_{\mu,m}$, independent of $A$, such that
\begin{equation}\label{seed:weighted-prepared-bounds}
 c_0\leq n_A\leq CA^P,\qquad
 \sum_{|\nu|\leq m}M_\mu(\partial^\nu S_A)
  +\sum_{|\nu|\leq m}M_\mu(\partial^\nu g)
 \leq C_{\mu,m}A^{P_{\mu,m}}.
\end{equation}
Both $S_A$ and $g$ are real, centrally odd, axisymmetric, and solenoidal.
The Gaussian azimuthal energy defined in \eqref{seed:gaussian-energy} also
satisfies $\mathcal Q_A(\tp)\leq CA^P$.
\end{lemma}

\begin{proof}
Let $\Lambda=(1-\Delta)^{1/2}$ and put
$v_{A,-4}=\Lambda^{-4}v_A$. We first verify that
$M_\mu(\Lambda^{-4}\mathcal I)<\infty$ for $\mu<5/2$.
The kernel of $(1-\Delta)^{-2}$ in three dimensions is
$\mathrm e^{-|x|}/(8\pi)$. In
$\PP\partial_z(e_z\mathrm e^{-|x|}/(8\pi))$, the local part is bounded
and exponentially decreasing. The remaining part is a third derivative of
the Newton potential of this radial kernel. Its first radial derivative
is a constant multiple of
\begin{equation}\label{seed:bessel-newton-tail}
 |x|^{-2}\int_0^{|x|}u^2\mathrm e^{-u}\idiff u.
\end{equation}
Near zero, the expression in \eqref{seed:bessel-newton-tail} is
$r/3-r^2/4+O(r^3)$. Thus the potential is a constant plus a multiple of
$r^2$ and an $O(r^3)$ term whose third Cartesian derivatives are bounded.
At infinity, its third derivatives have order $|x|^{-4}$.
The squared weighted exterior
integral is therefore bounded by a constant times
$\int_1^\infty\rho^{2\mu-6}\idiff\rho$, which converges in the asserted range.

For each fixed integer $j$, conjugation of the smooth matrix in
\eqref{seed:cartesian-linearization} by Bessel powers satisfies
\begin{equation}\label{seed:bessel-conjugation}
 M_\mu\bigl(\Lambda^j\mathsf M_A(t)\Lambda^{-j}h\bigr)
 \leq C_{\mu,j}A M_\mu(h),\qquad0\leq t\leq1.
\end{equation}
Here is a direct justification of the weighted estimate. For even integer
orders, the product rule expands the differential factor, and commutators
have one lower differential order. The remaining negative-order Bessel
kernels are locally integrable and decrease exponentially at infinity.
On unit cubes, ordinary Sobolev multiplication applies and the polynomial
weight is uniformly comparable on neighboring cubes. Interactions between
distant cubes are summable by exponential kernel decay. To justify
interpolation, define $\Lambda^{\mathrm i\vartheta}$ by the Fourier
multiplier $(1+|\xi|^2)^{\mathrm i\vartheta/2}$, where
$\vartheta\in\R$. The same localization gives
\begin{equation}\label{seed:imaginary-bessel-bound}
 M_\mu(\Lambda^{\mathrm i\vartheta}h)
 \leq C_\mu e^{C_\mu|\vartheta|}M_\mu(h).
\end{equation}
Indeed, the unweighted multiplier has norm one, which controls the
interactions between neighboring unit cubes. Off the diagonal, the
heat-kernel representation of the Bessel multiplier has exponential
spatial decay, with bound $Ce^{C|\vartheta|}$ outside the unit ball.
The ratio of the polynomial weights grows only polynomially in the
separation of the cubes, so the remaining interactions are absolutely
summable. The three-lines
theorem, with a Gaussian factor in the imaginary interpolation parameter,
now gives the intermediate integer orders from the consecutive even orders.
The same localization applies to $\Lambda^{-1}$ and
$\partial_i\Lambda^{-1}$, which are bounded on unweighted $L^2$
and have exponentially decreasing kernels away from the origin. The
identities
\[
 \Lambda\psi=\Lambda^{-1}\psi
      -\sum_{i=1}^3(\partial_i\Lambda^{-1})\partial_i\psi,
 \qquad
 \partial_i\psi=(\partial_i\Lambda^{-1})\Lambda\psi
\]
therefore show that $M_\mu(\Lambda\psi)$ is equivalent to
$M_\mu(\psi)+\sum_{i=1}^3M_\mu(\partial_i\psi)$, with constants
depending only on $\mu$.
Iteration gives the corresponding equivalence at every integer graph
order. Only finitely many derivatives of $\mathsf M_A$ enter each estimate.

The case $j=-4$ of \eqref{seed:bessel-conjugation} is needed to start
from the singular point input.  We give its kernel proof explicitly,
so that no weighted input derivative is assumed in this step.
Let $\Lambda=(1-\Delta)^{1/2}$ and let $M(x)$ be a smooth matrix
whose derivatives through order seven are bounded.  For every fixed
$\mu\geq0$,
\begin{equation}
 \|\langle x\rangle^\mu\Lambda^{-4}M\Lambda^4h\|_2
 \leq C_\mu
       \left(\max_{|\alpha|\leq7}\|\partial^\alpha M\|_\infty\right)
       \|\langle x\rangle^\mu h\|_2 .
 \label{seed:bessel-minus-four-explicit}
\end{equation}
Initially take $h$ smooth and compactly supported.  Since
$\Lambda^4=(1-\Delta)^2$ is a differential operator,
\[
 \Lambda^{-4}M\Lambda^4
     =M+\Lambda^{-4}[M,\Lambda^4],
 \qquad
 [M,\Lambda^4]=\sum_{|\alpha|\leq3}a_\alpha(x)\partial^\alpha.
\]
Every $a_\alpha$ is a finite linear combination of derivatives of
$M$ of orders at most four.  The kernel of $\Lambda^{-4}$ in
three dimensions is $K(x)=e^{-|x|}/(8\pi)$.  Its weak derivatives
through order three are locally integrable functions, and satisfy
\[
 |\partial^\beta K(x)|
 \leq C_\beta e^{-c|x|}
             \bigl(1+|x|^{1-|\beta|}\bigr),
 \qquad 0<|x|,\quad |\beta|\leq3.
\]
For $\beta=0$ a harmless enlargement of the right-hand side is
understood.  In particular, the worst local singularity is
$|x|^{-2}$, which is integrable in dimension three.  No Dirac term
occurs in these derivatives: integration by parts on $|x|>\varepsilon$
has boundary terms tending to zero through the third derivative.

For one commutator summand, integrate by parts to obtain
\[
 \begin{split}
 \Lambda^{-4}(a_\alpha\partial^\alpha h)(x)
  &=(-1)^{|\alpha|}\int
       \partial_y^\alpha\bigl(K(x-y)a_\alpha(y)\bigr)h(y)\idiff y.
 \end{split}
\]
Leibniz's rule expresses its kernel as a finite sum of
$\partial^\beta K(x-y)\partial^{\alpha-\beta}a_\alpha(y)$,
with $|\beta|\leq3$.  All coefficient factors are bounded by the
derivative maximum in \eqref{seed:bessel-minus-four-explicit}.
The elementary moderate-weight inequality
\[
 \langle x\rangle^\mu
       \leq C_\mu\langle x-y\rangle^\mu\langle y\rangle^\mu
\]
and the finiteness of
$\|\langle\cdot\rangle^\mu\partial^\beta K\|_{L^1}$ therefore
give the asserted bound by Young's convolution inequality.  Density
extends the conjugated operator to the entire weighted $L^2$ space.
For $M=\mathsf M_A(t)$ on $0\leq t\leq1$, each required coefficient
derivative is $O(A)$, which gives exactly the $C_\mu A$ bound needed
in the negative-space equation below.  This argument uses no global
$A_2$ assumption on the polynomial weight.

The extension is consistent with distributional conjugation: weighted
$L^2$ convergence implies $L^2$ convergence, $\Lambda^4$ maps $L^2$
to $H^{-4}$, multiplication by $M$ is bounded on $H^{-4}$ by duality
from its product-rule bound on $H^4$, and $\Lambda^{-4}$ maps
$H^{-4}$ to $L^2$.  These maps identify the limit of the displayed
identity.  In the present application
$M$ has zero vertical output and is axisymmetric.  Scalar Bessel powers
preserve these two properties, so the horizontal axisymmetric weighted
Leray estimate of Lemma~\ref{bg:weighted} applies to the resulting source.
The restriction $\mu<5/2$ enters at that Leray step and in the initial
point-response tail, not in the Bessel conjugation proved here.

With this negative-order estimate and the projection bound, the exact equation is
\begin{equation}\label{seed:negative-space-equation}
 (\partial_t-\Delta)v_{A,-4}
 =-\PP\Lambda^{-4}\mathsf M_A(t)\Lambda^4v_{A,-4}.
\end{equation}
At graph order $j\geq0$, the argument of the Leray projection in the
source is
$\Lambda^{j-4}\mathsf M_A\Lambda^{4-j}
(\Lambda^jv_{A,-4})$. It is horizontal and axisymmetric, because scalar
Bessel powers preserve these properties. Thus
\eqref{seed:bessel-conjugation} at order $j-4$ and
Lemma~\ref{bg:weighted} apply at every graph order, including odd $j$.
For $E_j=M_\mu(\Lambda^jv_{A,-4})^2$, weighted integration by parts gives
\begin{equation}\label{seed:weighted-heat-energies}
 E_j'+c_\mu E_{j+1}\leq C_{\mu,j}(1+A)E_j.
\end{equation}
In detail, for a smooth vector $\psi$ the diffusion form is
\begin{equation}\label{seed:weighted-diffusion-form}
 \begin{aligned}
 -\operatorname{Re}\int \langle x\rangle^{2\mu}
                  \overline\psi\cdot\Delta\psi\idiff x
 &=\int\langle x\rangle^{2\mu}|\nabla\psi|^2\idiff x
       -\frac12\int\Delta\langle x\rangle^{2\mu}|\psi|^2\idiff x\\
 &\geq c_\mu M_\mu(\Lambda\psi)^2-C_\mu M_\mu(\psi)^2,
 \end{aligned}
\end{equation}
where the last inequality uses
$|\Delta\langle x\rangle^{2\mu}|
\leq C_\mu\langle x\rangle^{2\mu}$ and the first-order weighted
Bessel norm equivalence. Taking $\psi=\Lambda^jv_{A,-4}$ proves
\eqref{seed:weighted-heat-energies}.

The time-weight device can be checked at a general fixed order $N$.
Suppose $E_j'+cE_{j+1}\leq M_jE_j$ for $0\leq j\leq N$.
Choose positive numbers $\beta_0=1,\beta_1,\ldots,\beta_N$ so that
$j\beta_j\leq c\beta_{j-1}/2$ and set
$\mathcal E_N=\sum_{j=0}^N\beta_jt^jE_j$.
Differentiation gives
\[
 \mathcal E_N'
 \leq \sum_{j=1}^N j\beta_jt^{j-1}E_j
      -c\sum_{j=0}^N\beta_jt^jE_{j+1}
      +(\max_{j\leq N}M_j)\mathcal E_N.
\]
Each positive time-weight derivative is absorbed by the preceding
negative heat term.  For smooth initial data,
$\mathcal E_N(0)=E_0(0)$, so Gronwall gives
\[
 E_j(t)\leq\beta_j^{-1}t^{-j}
                   e^{(\max M_j)t}E_0(0),\qquad 0\leq j\leq N.
\]
Apply this to smooth approximations whose $E_0(0)$ is bounded and then
pass to the limit.  No initial $E_j(0)$, $j>0$, is required.
In the present application $M_j\leq C_{\mu,j}(1+A)$; taking square roots
gives the factor $t^{-j/2}e^{C_{\mu,j}At}$ below.  This is also the
finite-order smoothing mechanism used later on rescaled last intervals.

Applied to \eqref{seed:weighted-heat-energies}, this time-weighted
estimate yields
\begin{equation}\label{seed:negative-space-smoothing}
 M_\mu(\Lambda^jv_{A,-4}(t))
 \leq C_{\mu,j}t^{-j/2}\mathrm e^{C_{\mu,j}At}
             M_\mu(\Lambda^{-4}\mathcal I),\qquad0<t\leq1.
\end{equation}
For a precise approximation, take the negative-space initial data
$e^{\varepsilon\Delta}\Lambda^{-4}\mathcal I$, with
$\varepsilon>0$. They are smooth, centrally odd and axisymmetric, have
uniformly bounded $E_0(0)$, and converge in the weighted $L^2$ norm as
$\varepsilon\downarrow0$. The last assertion follows by density from
\eqref{bg:odd-operators}. After applying $\Lambda^4$, the resulting
solutions converge in $C([0,t_0];H^{-4})$ on each fixed interval. The smooth matrix and the
projection are bounded in $H^{-4}$, so uniqueness for the bounded
perturbation of heat identifies their limit with $v_A$. Taking $j=m+4$ in
\eqref{seed:negative-space-smoothing}, and using
$A\tp=2+\Theta\log A$, proves the upper bounds for $S_A$ in
\eqref{seed:weighted-prepared-bounds}. On one fixed compact positive
frequency interval, Lemma~\ref{seed:point-growth} gives
$a_{A,k}(\tp)\geq c_J$. The adjoint has bounded norm there, so Plancherel's
identity gives $n_A\geq c_0>0$. Division by $n_A$ proves the estimates for $g$.

We now establish the Gaussian azimuthal estimate. Define
\begin{equation}\label{seed:gaussian-weight}
 \chi_0=\mathrm e-2,\qquad
 \mathfrak m_{\mathrm G}(\zeta_r)=\int_0^1(\mathrm e^{1-y}-1)\mathrm e^{y\zeta_r}\idiff y,
 \qquad m_{\mathrm G}(t,r)=\tau_t\mathfrak m_{\mathrm G}\Bigl(\frac{r^2}{4\tau_t}\Bigr),
\end{equation}
and, for the point response, set
\begin{equation}\label{seed:gaussian-energy}
 \mathcal Q_A(t)
 =\int m_{\mathrm G}(t,r)|v_{A,\theta}|^2r\idiff r\idiff z
   +\chi_0\tau_t\int
       \bigl(|v_{A,r}|^2+|v_{A,z}|^2\bigr)r\idiff r\idiff z.
\end{equation}
The integral representation implies
$\chi_0\leq\mathfrak m_{\mathrm G}(\zeta_r)\leq\chi_0\mathrm e^{\zeta_r}$.
The initial impulse has zero azimuthal component, and its azimuthal equation is
\begin{equation}\label{seed:azimuthal-source}
 (\partial_t-\Delta_1)v_{A,\theta}=-AB_{\tau_t}v_{A,r},
 \qquad \Delta_1=\Delta-r^{-2}.
\end{equation}
Set $\lambda_{\mathrm G}(t)=1/(8\tau_t)$. The conjugated Cartesian vector
heat kernel is
\begin{equation}\label{seed:gaussian-conjugated-kernel}
 K(t,s;x,y)=
 \mathrm e^{\lambda_{\mathrm G}(t)|x_h|^2}
 G_{t-s}(x-y)\mathrm e^{-\lambda_{\mathrm G}(s)|y_h|^2}.
\end{equation}
For $0\leq s<t\leq1$, write $h=t-s$ and $d=x_h-y_h$. Its horizontal
exponent is
\begin{equation}\label{seed:gaussian-quadratic-form}
 -\Bigl(\frac1{4h}-\frac1{8\tau_t}\Bigr)|d|^2
 +\frac{d\cdot y_h}{4\tau_t}
 -\frac{h|y_h|^2}{8\tau_s\tau_t}.
\end{equation}
The corresponding positive quadratic form has determinant
$1/(64\tau_s\tau_t)$. Under $d=\sqrt h\,d'$ and
$y_h=h^{-1/2}y'$, it belongs to a compact family of positive forms.
Consequently \eqref{seed:gaussian-quadratic-form} is bounded above by
$-c|d|^2/h-ch|y_h|^2$. Input derivatives of order $j$ cost at most
$C_jh^{-j/2}$ after Gaussian polynomials are absorbed. Schur's test and
duality, together with the axial heat bounds, give
\begin{equation}\label{seed:gaussian-kernel-bounds}
 \|K(t,s)\|_{2\to2}\leq C,\qquad
 \|K(t,s)\|_{H^{-3}\to2}\leq Ch^{-3/2}.
\end{equation}

To estimate the source in Cartesian negative norms, let
$S=\operatorname{diag}(1,-1,1)$ and define
$\mathscr R f(x)=Sf(Sx)$. On axisymmetric fields,
$(1+\mathscr R)/2$ selects the meridional component. A constant horizontal
rotation maps its horizontal component to $v_re_\theta$. These maps are
bounded in every ordinary Sobolev space. The multiplier
$AB_{\tau_t}\mathrm e^{\lambda_{\mathrm G}(t)r^2}$ and each of its required
derivatives are bounded by $C_jA$. From
\eqref{seed:cartesian-linearization} and ordinary heat smoothing,
\begin{equation}\label{seed:full-negative-bounds}
 \|v_A(s)\|_{H^{-3}}\leq C\mathrm e^{CAs},\qquad
 \|v_A(s)\|_2\leq Cs^{-3/2}\mathrm e^{CAs},\qquad0<s\leq1.
\end{equation}
In the Duhamel formula for \eqref{seed:azimuthal-source}, use the first
bound for $s<t/2$ and the second for $s>t/2$. Its zero initial azimuthal
component then gives
\begin{equation}\label{seed:gaussian-swirl-bound}
 \bigl\|\mathrm e^{r^2/(8\tau_t)}v_{A,\theta}(t)\bigr\|_2
 \leq CA\mathrm e^{CAt}
    \Bigl(\int_0^{t/2}t^{-3/2}\idiff s
                 +\int_{t/2}^t s^{-3/2}\idiff s\Bigr)
 \leq CAt^{-1/2}\mathrm e^{CAt}.
\end{equation}
For $A^{-1}\leq t\leq1$, \eqref{seed:full-negative-bounds} and
\eqref{seed:gaussian-swirl-bound} imply
$\mathcal Q_A(t)\leq CA^2t^{-1}\mathrm e^{CAt}$, after changing the
constant. This is polynomial in $A$ at $t=\tp$.
\end{proof}

\subsection{A strict growth deficit at low axial frequencies}

The later sign estimate applies to the high-frequency portion.  The
missing frequencies are controlled here rather than discarded.  Very
small frequencies exploit the explicit factor $|k|$ in the point input;
a compact intermediate band uses the strict spectral gap below
$c_\star$.  Both bounds concern the entire lifetime energy.

The actual prepared velocity $S_A$ retains all axial Fourier
frequencies and has the ordinary spatial moments needed later.
A sharp multiplier $\mathbf1_{\{|k|\geq2\}}$ does \emph{not} generally
preserve those moments: the inverse Fourier transform can acquire a
slow axial tail.  Below, high- and low-frequency pieces are used to
compare scalar energy histories; when a weighted integration is
required, the proof first uses smooth frequency truncations and then
passes only the quantities for which the cutoff limit is justified.
The low-frequency bound is thus an estimate of a part of the genuine
response, not a proposal to delete that part from the construction.

The sign estimate below will control the part of the response with
$|k|\geq2$. To retain the full point response, we need to show that its lower
frequencies contribute an exponentially smaller lifetime energy.

\begin{lemma}[Low-frequency lifetime estimate]\label{seed:low-frequency}
There are fixed constants $\eta_{\mathrm{low}}>0$, $C$, and $N$, independent of $A$, such that
\begin{equation}\label{seed:low-frequency-lifetime}
 \int_\tp^\infty\|\mathbf1_{\{|D_z|<2\}}v_A(t)\|_2^2\idiff t
 \leq CA^N\mathrm e^{2(c_\star-\eta_{\mathrm{low}})A}.
\end{equation}
On a compact frequency interval away from zero, the propagator estimate used
in this proof applies to arbitrary physical input.
\end{lemma}

\begin{proof}
Choose
\begin{equation}\label{seed:low-frequency-constants}
 k_0=\frac{c_\star}{16\sqrt2},\qquad
 \sigma_\star=\sigma_{2\sqrt2}<c_\star,\qquad
 \eta_{\mathrm{low}}=\frac{c_\star-\sigma_\star}{2}>0.
\end{equation}
For any real $k$, the free heat norm of the point Fourier component is at
most $C|k|t^{-1/2}\mathrm e^{-k^2t}$. The Cartesian matrix in
\eqref{seed:cartesian-linearization} has norm at most $CA$, uniformly in
$k$. Remove the scalar axial heat factor in Duhamel's formula and use the
integrability of $t^{-1/2}$. At $t_1=A^{-1}$ this gives
\begin{equation}\label{seed:point-frequency-small-time}
 \|\widehat v_A^{\,z}(t_1,k)\|_{L_h^2}
 \leq C\sqrt A\,|k|\mathrm e^{-k^2/A}.
\end{equation}
For arbitrary input at time $s$, \eqref{seed:small-frequency-coupling} gives
\begin{equation}\label{seed:small-band-propagator}
 \|\mathcal U_{\mathrm d,k}(t,s)\|
 \leq\exp\Bigl(-k^2(t-s)
       +\sqrt2 A|k|\int_s^t(1+u)^{-3/2}\idiff u\Bigr).
\end{equation}
Therefore, for $0<|k|\leq k_0$,
\begin{equation}\label{seed:small-band-solution}
 \|\widehat v_A^{\,z}(t,k)\|_2
 \leq C\sqrt A\,|k|\mathrm e^{2\sqrt2 A|k|-k^2t}.
\end{equation}
After integration in time, the factor $k^2$ cancels the factor
$1/(2k^2)$ from the heat integral. The remaining frequency integral is
bounded by $CA\mathrm e^{4\sqrt2 Ak_0}=CA\mathrm e^{c_\star A/4}$.
The zero-frequency component is zero by continuity. This calculation uses
the factor $|k|$ in the actual point input.

For $k_0\leq|k|\leq2$ and $0\leq t\leq1$, the rescaled frequency is at
most $2\sqrt2$. Lemma~\ref{seed:metric} gives the squared propagator bound
\begin{equation}\label{seed:middle-band-propagator}
 C\exp\Bigl(2A\int_s^t(1+u)^{-2}
                       \sigma_{|k|\sqrt{1+u}}\idiff u+C(t-s)\Bigr)
 \leq C\mathrm e^{\sigma_\star A}.
\end{equation}
Together with \eqref{seed:point-frequency-small-time}, this controls the
energy at time one and its integral from $\tp$ to one by
$CA^N\mathrm e^{\sigma_\star A}$. For $t\geq1$, the sharp metric
\eqref{seed:sharp-energy} retains axial diffusion and yields
\begin{equation}\label{seed:middle-band-late-decay}
 \|\widehat v_A^{\,z}(t,k)\|_2^2
 \leq C(1+t)^5\mathrm e^{c_\star A-2k_0^2(t-1)}
                       \|\widehat v_A^{\,z}(1,k)\|_2^2.
\end{equation}
Its time integral has an $A$-independent constant. Since
$c_\star+\sigma_\star=2(c_\star-\eta_{\mathrm{low}})$, combining the small and middle
frequency regions proves \eqref{seed:low-frequency-lifetime}.
\end{proof}

\subsection{The complete quadratic stress and the retained azimuthal mean}

This is the point at which a linear instability becomes a signed
nonlinear interaction.  The second variation creates an axial mean in
its azimuthal component, and its pairing with the reference contributes
to the same order as $v\otimes v$.  We calculate that contribution before
taking any long-time or large-scale limit.

The second variation begins at the positive preparation time, with zero
initial value. There is no factorial in its definition: for the point
response it solves
\begin{equation}\label{seed:second-variation}
 (\partial_t-\mathcal L_{U_{\mathrm d}})w_A
      =-\PP\operatorname{div}(v_A\otimes v_A),\qquad w_A(\tp)=0,
\end{equation}
where $\mathcal L_{U_{\mathrm d}}$ is \eqref{strategy:linearized-generator}
with $V=U_{\mathrm d}$.
The complete stress is
\begin{equation}\label{seed:complete-stress-definition}
 \mathsf S_A=v_A\otimes v_A
             +U_{\mathrm d}\otimes w_A+w_A\otimes U_{\mathrm d},
 \qquad
 \int_{\R^3}\mathsf S_A\idiff x
       =\operatorname{diag}(s_{h,A},s_{h,A},s_{z,A}).
\end{equation}
The adjective ``complete'' refers to the inclusion of both mixed terms.

For a general real first variation $v$, use the cylindrical integrals
\begin{equation}\label{seed:energy-correlation}
 E_v(t)=\int|v|^2r\idiff r\idiff z,\qquad
 Z_v(t)=\int|v_z|^2r\idiff r\idiff z,\qquad
 F_v(t,r)=\int_\R v_r(t,r,z)v_\theta(t,r,z)\idiff z.
\end{equation}
Thus the physical energy is $2\pi E_v$. Let $s_{h,v},s_{z,v}$ denote the integrated
diagonal entries of its complete stress, with the corresponding second
variation initialized at zero at $\tp$.

For an axisymmetric vector,
$v_x=v_r\cos\theta-v_\theta\sin\theta$ and
$v_y=v_r\sin\theta+v_\theta\cos\theta$.  Angular integration of
$v\otimes v$ gives
\[
 s_h^{(vv)}=\pi\int(v_r^2+v_\theta^2)r\,\idiff r\idiff z,
 \qquad s_z^{(vv)}=2\pi\int v_z^2r\,\idiff r\idiff z.
\]
Thus $s_h^{(vv)}-s_z^{(vv)}=\pi(E_v-3Z_v)$.
For $U_{\mathrm d}=AW_{1+t}e_\theta$, the two mixed tensors have zero
vertical diagonal entry and contribute
$2\pi A\int W_{1+t}\mu_\theta[w]\,r\,\idiff r$ to each horizontal diagonal
entry, where we define the azimuthal mean by
\begin{equation}\label{seed:azimuthal-mean}
 \mu_\theta[w](t,r):=\int_\R w_\theta(t,r,z)\,\idiff z.
\end{equation}
This accounts for the factors in the next identity.

The correlation $F_v=\int v_rv_\theta\,\idiff z$ need not have a
pointwise sign.  Accordingly, the sign of the complete history will not
be inferred merely from the negative sign of $C_\tau$.  We first derive
an exact mean equation and an exact history formula.  A different,
positive weighted-energy argument then controls that correlation after
integration in time.

\begin{lemma}[Exact stress and mean identities]\label{seed:stress-history}
Suppose that $v$ is a smooth real axisymmetric solenoidal solution of the homogeneous
linearized equation on $[\tp,t_e]$, with
$v(\tp)\in H^4(\R^3)$ and $M_1(v(\tp))<\infty$.
Let its second variation solve \eqref{seed:second-variation} with $v_A$
replaced by $v$.  The complete stress is integrable, and its diagonal
difference satisfies
\begin{equation}\label{seed:exact-stress-formula}
 (s_{h,v}-s_{z,v})(t)
 =\pi(E_v-3Z_v)(t)
   +2\pi A\int_\tp^t\int_0^\infty
       C_{1+2t-s}(r)F_v(s,r)r\idiff r\idiff s,
\end{equation}
where
\begin{equation}\label{seed:stress-kernels}
 C_\tau(r)=W_\tau'(r)-r^{-1}W_\tau(r)
          =-\frac{r^2}{2\tau^3}\mathrm e^{-r^2/(4\tau)},
 \qquad
 K_{\mathrm{hist}}(\tau,r)=\frac8{r^2}
       \Bigl(1-\Bigl(1+\frac{r^2}{4\tau}\Bigr)
                    \mathrm e^{-r^2/(4\tau)}\Bigr).
\end{equation}
The continuous value is used at $r=0$, and $0\leq K_{\mathrm{hist}}(\tau,r)\leq2/\tau$.
For each terminal time $t_e$, the integrated identity is
\begin{equation}\label{seed:finite-stress-history}
 \begin{aligned}
 \int_\tp^{t_e}(s_{h,v}-s_{z,v})(t)\idiff t
 &=\pi\int_\tp^{t_e}(E_v-3Z_v)(t)\idiff t
   -\pi A\int_\tp^{t_e}\int_0^\infty
        K_{\mathrm{hist}}(1+s,r)F_v(s,r)r\idiff r\idiff s+\mathcal R_v(t_e),\\
 \mathcal R_v(t_e)
 &=\pi A\int_\tp^{t_e}\int_0^\infty
        K_{\mathrm{hist}}(1+2t_e-s,r)F_v(s,r)r\idiff r\idiff s,\\
 |\mathcal R_v(t_e)|
 &\leq\frac{\pi A}{1+t_e}\int_\tp^{t_e}E_v(s)\idiff s.
 \end{aligned}
\end{equation}
If $\mu_{\theta,e}(r):=\mu_\theta[w](t_e,r)$ is the actual azimuthal
mean at the terminal time, then its subsequent mixed-stress pairing under
free heat is
\begin{equation}\label{seed:retained-mean-pairing}
 \mathcal T_v(t_e)
 =4\pi A\int_0^\infty
       \Bigl(1-\mathrm e^{-r^2/(4(1+t_e))}\Bigr)\mu_{\theta,e}(r)\idiff r
 =-\mathcal R_v(t_e).
\end{equation}
This last quantity evolves the already existing mean; it includes no
quadratic forcing after $t_e$.
\end{lemma}

\begin{proof}
The weighted projection and tensor-heat estimates of Lemma~\ref{bg:weighted},
with exponent one, propagate $v,zv,w,zw\in L^2$ on each finite positive
time interval. The invariant linear matrix has horizontal axisymmetric
output, and the quadratic source obeys
$M_1(v\otimes v)\leq C\|v\|_\infty M_1(v)$, while the differentiated heat
kernel has an integrable time singularity. Positive-time smoothing supplies
the fixed derivatives used below. In particular,
\begin{equation}\label{seed:stress-integrability}
 \int_{\R^3}|U_{\mathrm d}w|\idiff x
 \leq\|U_{\mathrm d}(1+z^2)^{-1/2}\|_2
                   \|(1+z^2)^{1/2}w\|_2<\infty.
\end{equation}
The tensor $v\otimes v$ is also integrable. Axisymmetry makes the integrated
symmetric tensor diagonal with equal horizontal entries.

Axial integration of incompressibility gives
$\partial_r(r\int_\R w_r\idiff z)=0$ in distributions. The axial mean
belongs to $L^2((0,\infty),r\idiff r)$ by Cauchy--Schwarz with the
weight $1+z^2$. A nonzero constant divided by $r$ does not belong to this
space, so $\int_\R w_r\idiff z=0$. Therefore the exact mean
$\mu_\theta[w](t,r)=\int_\R w_\theta(t,r,z)\idiff z$ solves
\begin{equation}\label{seed:mean-equation}
 (\partial_t-\Delta_h^{(1)})\mu_\theta[w]
       =-(\partial_r+2/r)F_v,\qquad \mu_\theta[w](\tp)=0,
 \qquad
 \Delta_h^{(j)}=\partial_r^2+r^{-1}\partial_r-j^2r^{-2},\quad j=0,1.
\end{equation}

The coefficient $2/r$ in this mean equation includes a geometric term.
For a solenoidal axisymmetric $v$, the azimuthal component of
$(v\cdot\nabla)v$ is
\[
 v_r\partial_rv_\theta+v_z\partial_zv_\theta+\frac{v_rv_\theta}{r}
 =\partial_r(v_rv_\theta)+\partial_z(v_zv_\theta)
                           +\frac{2v_rv_\theta}{r}.
\]
Here incompressibility was used in writing the right side.  Axial
integration removes the $z$ derivative and gives
$(\partial_r+2/r)F_v$.

There are also two heat evolutions in the stress pairing.  A mean created
at time $s$ is evolved for $t-s$, while the background at time $t$ already
has heat age $1+t$.  Moving the radial heat to the background gives
$W_{1+t+(t-s)}=W_{1+2t-s}$.  This explains the heat age in
\eqref{seed:exact-stress-formula} and, after time integration, the
factor two in \eqref{seed:kernel-time-integral}.

The linear term $-AB_{\tau_t}w_r$ has zero axial mean. The azimuthal quadratic
divergence gives exactly the operator in \eqref{seed:mean-equation}, and
pressure has no azimuthal component.

Angular integration of $v\otimes v$ gives $\pi(E_v-3Z_v)$. The contribution
of the two mixed stresses is
\begin{equation}\label{seed:mixed-stress-contribution}
 2\pi A\int_0^\infty W_{1+t}\mu_\theta[w](t,r)r\idiff r.
\end{equation}
Use the self-adjoint radial heat semigroup, the identity
$\mathrm e^{h\Delta_h^{(1)}}W_\tau=W_{\tau+h}$, and the integration by parts
formula
\begin{equation}\label{seed:mean-integration-by-parts}
 \int_0^\infty W_\tau\bigl(-\partial_r-2/r\bigr)F_v\,r\idiff r
 =\int_0^\infty C_\tau F_v\,r\idiff r.
\end{equation}
This proves \eqref{seed:exact-stress-formula}.

The heat-age substitution gives the exact kernel identity
\begin{equation}\label{seed:kernel-time-integral}
 -2\int_0^\infty C_{\tau+2h}(r)\idiff h=K_{\mathrm{hist}}(\tau,r).
\end{equation}
For every $q\geq0$, the identity
$1-(1+q)\mathrm e^{-q}=\int_0^q y\mathrm e^{-y}\idiff y\leq q$,
applied with $q=r^2/(4\tau)$, gives $K_{\mathrm{hist}}(\tau,r)\leq2/\tau$.
Applying
\eqref{seed:kernel-time-integral} to the finite history gives
\eqref{seed:finite-stress-history}; the remainder estimate follows from
$1+2t_e-s\geq1+t_e$ and
$\int_0^\infty|F_v(s,r)|r\idiff r\leq E_v(s)/2$. These estimates also
justify the time interchanges by absolute integrability.

The terminal mean has enough radial integrability for its entire future
pairing. Indeed, Cauchy--Schwarz in $z$, using
$\int_\R(1+r^2+z^2)^{-1}\idiff z=\pi(1+r^2)^{-1/2}$, gives
\begin{equation}\label{seed:terminal-mean-moment}
 \int_0^\infty (1+r^2)^{1/2}|\mu_{\theta,e}(r)|^2r\idiff r
 \leq C M_1(w(t_e))^2.
\end{equation}
In particular, the product of $\mu_{\theta,e}$ and
$1-e^{-r^2/(4(1+t_e))}$ belongs to $L^1(\idiff r)$: the square of
this factor, divided by $r(1+r^2)^{1/2}$, is integrable both at zero and
at infinity.

After $t_e$, evolve only $\mu_{\theta,e}$ by
$\mathrm e^{(t-t_e)\Delta_h^{(1)}}\mu_{\theta,e}$. Its subsequent mixed-stress pairing is
\begin{equation}\label{seed:mean-future-history}
 2\pi A\int_{t_e}^\infty\int_0^\infty
             W_{1+2t-t_e}(r)\mu_{\theta,e}(r)r\idiff r\idiff t.
\end{equation}
Since
\begin{equation}\label{seed:swirl-time-integral}
 \int_{t_e}^\infty W_{1+2t-t_e}(r)\idiff t
 =\frac2r\Bigl(1-\mathrm e^{-r^2/(4(1+t_e))}\Bigr),
\end{equation}
the absolute integrability supplied by \eqref{seed:terminal-mean-moment}
justifies integration in $t$ and identifies the result with
$\mathcal T_v(t_e)$. Substituting the finite source history from
\eqref{seed:mean-equation} into \eqref{seed:mean-future-history} instead gives
$-\mathcal R_v(t_e)$ by \eqref{seed:kernel-time-integral}. Its absolute
integral is controlled by the bounded kernel $K_{\mathrm{hist}}(1+2t_e-s,r)$ and the finite
energy history. This proves \eqref{seed:retained-mean-pairing}.
\end{proof}

\paragraph{Construction of the sign weight.}
The history formula contains $A\int K_{\mathrm{hist}}(1+t,r)F_vr\,\idiff r$.
Testing the azimuthal equation with $m_{\mathrm G}v_\theta$ produces
$A\int(-B_{\tau_t})m_{\mathrm G}F_vr\,\idiff r$, while the meridional energy produces
$2A\int\Omega_{\tau_t}F_vr\,\idiff r$.  We therefore seek
\[
 K_{\mathrm{hist}}(\tau_t,r)=(-B_{\tau_t})m_{\mathrm G}+2\chi_0\tau_t\Omega_{\tau_t}.
\]
This condition gives the identity in
\eqref{seed:gaussian-weight-identity}.  Writing $m_{\mathrm G}=\tau_t\mathfrak m_{\mathrm G}(\zeta_r)$ gives
\[
 (\zeta_r-1)\mathfrak m_{\mathrm G}(\zeta_r)
    =\frac{e^{\zeta_r}-1-\zeta_r}{\zeta_r}-\chi_0.
\]
To avoid a pole at $\zeta_r=1$, one must choose $\chi_0=e-2$.
The integral representation in \eqref{seed:gaussian-weight} proves that
the resulting function is positive, including across that removable
singularity.  It also gives the derivative inequality needed for heat.

The high-frequency assumption now has a transparent role: axial diffusion
bounds each weighted dissipation below by four times the corresponding
energy.  After the time and radial derivatives of the weight have been
subtracted, a positive multiple of the energy remains.  Inserting that
lower bound for the correlation into the exact history formula gives
\emph{negative} integrated $s_{h,v}-s_{z,v}$.  The azimuthal mean of the second variation supplies the term controlled
by this multiplier and is therefore part of the sign mechanism.

\begin{lemma}[A sign estimate for the high-frequency history]\label{seed:high-frequency}
Initially let $v$ be a real axisymmetric solution of the
homogeneous linearized equation about $U_{\mathrm d}$ satisfying the
regularity and moment hypotheses of Lemma~\ref{seed:stress-history}.
Suppose in addition that $v$ has axial Fourier support in $|k|\geq2$ and
finite initial Gaussian azimuthal energy $\mathcal Q_v(\tp)$, defined as
in \eqref{seed:gaussian-energy}. Its complete stress, with zero initial azimuthal mean of the second
variation, satisfies
\begin{equation}\label{seed:high-frequency-coercivity}
 \begin{aligned}
 \int_\tp^{t_e}(s_{h,v}-s_{z,v})(t)\idiff t
 \leq{}&\frac\pi2\mathcal Q_v(\tp)-\frac\pi2\mathcal Q_v(t_e)\\
 &-\pi\Bigl(\frac12-\frac{A}{1+t_e}\Bigr)
                    \int_\tp^{t_e}E_v(t)\idiff t
   -3\pi\int_\tp^{t_e}Z_v(t)\idiff t.
 \end{aligned}
\end{equation}
The assertion allows unbounded high-frequency support.
For the sharply high-frequency-truncated prepared point response,
the left side is understood as the scalar bounded-kernel history on
the right side of \eqref{seed:finite-stress-history}, obtained by the
smooth-cutoff limit below.  This extension does not assert preservation
of an ordinary axial moment by a sharp Fourier cutoff.
\end{lemma}

\begin{proof}
The function $\mathfrak m_{\mathrm G}$ in \eqref{seed:gaussian-weight} is positive and
increasing, with $\mathfrak m_{\mathrm G}\geq\chi_0=\mathrm e-2>1/2$. Its integral
representation gives
\begin{equation}\label{seed:gaussian-weight-identity}
 (\zeta_r-1)\mathfrak m_{\mathrm G}(\zeta_r)
       =\frac{\mathrm e^{\zeta_r}-1-\zeta_r}{\zeta_r}-\chi_0,
 \qquad
 K_{\mathrm{hist}}(\tau_t,r)=(-B_{\tau_t})m_{\mathrm G}+2\chi_0\tau_t\Omega_{\tau_t},
\end{equation}
where the expressions have removable values at zero and one. Direct
differentiation also gives
\begin{equation}\label{seed:gaussian-weight-differential}
 (\partial_t+\partial_r^2+r^{-1}\partial_r)m_{\mathrm G}
 =\mathfrak m_{\mathrm G}+(1-\zeta_r)\mathfrak m_{\mathrm G}'+\zeta_r\mathfrak m_{\mathrm G}''
 \leq2\mathfrak m_{\mathrm G}=2m_{\mathrm G}/\tau_t,
 \qquad \sup_r m_{\mathrm G}B_{\tau_t}^2\leq C\tau_t^{-3}.
\end{equation}
Indeed, the integrand in the first differentiated expression is
$(\mathrm e^{1-y}-1)\mathrm e^{y\zeta_r}
(1+y-\zeta_r y(1-y))$, at most twice the positive integrand of $\mathfrak m_{\mathrm G}$.
The second estimate follows from $\mathfrak m_{\mathrm G}\leq\chi_0\mathrm e^{\zeta_r}$.

Test the azimuthal equation first with $m_{\mathrm G,N}v_\theta$, where
$m_{\mathrm G,N}=\min(m_{\mathrm G},N)$. The radial derivative jumps downward at the cap, producing
a nonpositive distributional Laplacian, while the time derivative follows
the almost-everywhere chain rule. Thus
$(\partial_t+\Delta_h)m_{\mathrm G,N}\leq2m_{\mathrm G,N}/\tau_t$. Its forcing is bounded by
\begin{equation}\label{seed:capped-forcing}
 A\left|\int B_{\tau_t}m_{\mathrm G,N}v_rv_\theta\,r\idiff r\idiff z\right|
 \leq CA\tau_t^{-3/2}\|v_r\|_{L^2(r\idiff r\idiff z)}
       \Bigl(\int m_{\mathrm G,N}|v_\theta|^2r\idiff r\idiff z\Bigr)^{1/2}.
\end{equation}
This gives bounds uniform in $N$ on each finite interval. Monotone
convergence and lower semicontinuity justify passage to $m_{\mathrm G}$.

Let $E_{\theta,m_{\mathrm G}}=\int m_{\mathrm G}|v_\theta|^2r\idiff r\idiff z$ and let
$D_{\theta,m_{\mathrm G}}$ be its weighted gradient dissipation, including
$r^{-2}|v_\theta|^2$. Let $E_{\mathrm{mer}}$ and $D_{\mathrm{mer}}$ be
the unweighted meridional energy and dissipation. The two equations give
\begin{equation}\label{seed:component-energy-identities}
 \begin{aligned}
 A\int_0^\infty(-B_{\tau_t})m_{\mathrm G}F_v\,r\idiff r
       &\geq\tfrac12 E_{\theta,m_{\mathrm G}}'+D_{\theta,m_{\mathrm G}}-E_{\theta,m_{\mathrm G}}/\tau_t,\\
 \tfrac12 E_{\mathrm{mer}}'+D_{\mathrm{mer}}
       &=2A\int_0^\infty\Omega_{\tau_t}F_v\,r\idiff r.
 \end{aligned}
\end{equation}
Multiplying the second identity by $\chi_0\tau_t$ and adding, using
\eqref{seed:gaussian-weight-identity}, yields
\begin{equation}\label{seed:gaussian-history-energy}
 A\int_0^\infty K_{\mathrm{hist}}(1+t,r)F_v\,r\idiff r
 \geq\tfrac12\mathcal Q_v'
    +D_{\theta,m_{\mathrm G}}+\chi_0\tau_tD_{\mathrm{mer}}
    -E_{\theta,m_{\mathrm G}}/\tau_t-\tfrac12\chi_0E_{\mathrm{mer}}.
\end{equation}
All weights are independent of $z$. The frequency assumption gives
$D_{\theta,m_{\mathrm G}}\geq4E_{\theta,m_{\mathrm G}}$ and
$D_{\mathrm{mer}}\geq4E_{\mathrm{mer}}$. Since $\tau_t\geq1$, the four
terms after $\mathcal Q_v'/2$ are at least $3\mathcal Q_v$. Hence we obtain
\begin{equation}\label{seed:positive-history-bound}
 A\int_\tp^{t_e}\int_0^\infty K_{\mathrm{hist}}(1+t,r)F_v\,r\idiff r\idiff t
 \geq\tfrac12\bigl(\mathcal Q_v(t_e)-\mathcal Q_v(\tp)\bigr)
          +3\int_\tp^{t_e}\mathcal Q_v(t)\idiff t.
\end{equation}
Because $\mathcal Q_v\geq\chi_0E_v\geq E_v/2$, substitution into
\eqref{seed:finite-stress-history} proves
\eqref{seed:high-frequency-coercivity}.

To include unbounded frequency support, first use even smooth cutoffs
$0\leq\chi_N\leq1$, compactly supported in $|k|>2$, with
$\chi_N^2\uparrow\mathbf1_{\{|k|>2\}}$ almost everywhere. For each fixed $N$, the smooth compact frequency multiplier
is convolution with a Schwartz kernel in $z$ and is bounded on every
ordinary polynomial moment space used here.  It commutes with the
invariant linearized evolution.  Thus its prepared input satisfies the
hypotheses of Lemma~\ref{seed:stress-history}, and that lemma applies
to the smooth-cutoff history.  The moment constants may depend on $N$;
no uniform moment estimate is used in passing to the sharp cutoff.
Parseval's identity and monotone convergence pass the nonnegative energies
and dissipations to the limit. The scalar history passes by the finite-time
energy bound and the bounded kernel $K_{\mathrm{hist}}(\tau,r)$. No preservation of an axial
moment by a sharp frequency cutoff is required.
\end{proof}

\subsection{The positive stress mass of the full response}

The high-frequency sign and the low-frequency deficit now fit together.
The growing band gives a positive lower bound for the energy history;
the sign identity converts it into a lower bound for the complete mass.
Normalization by the prepared seed costs only polynomial factors, so the
exponential gap needed by the cascade is retained.

The growing compact band supplies a lower bound on the time-integrated
energy.  The high-frequency inequality converts that energy into a sign
for the complete stress.  The low-frequency lifetime estimate prevents
the omitted frequencies from overturning this sign.  Schematically,
before normalization by $n_A^2$, the competing sizes are
\[
 \underbrace{A^{-1}e^{2GA}}_{\text{growing-band energy}},\qquad
 \underbrace{CA^Pe^{2(c_\star-\eta_{\mathrm{low}})A}}_{\text{low-frequency error}},\qquad
 \underbrace{CA^P}_{\text{initial Gaussian energy}}.
\]
Choose the fixed band and observation time so that $G>c_\star-\eta_{\mathrm{low}}$,
then increase $A$.  The first term dominates the other two.  Taking
$t_e=\tf=\tp+16A^2$ also makes the terminal mean remainder at most
$O(A^{-1})$ of the energy history.

The splitting is made at the level of the scalar history in
\eqref{seed:finite-stress-history}: Plancherel eliminates products
of disjoint axial frequencies after integration in $z$.  The complete
second variation remains in the stress tensor throughout the argument.

\begin{lemma}[The complete signed history]\label{seed:signed-history}
Given $G_0<c_\star$, one can choose a compact interval $J\subset(2,\infty)$,
a finite time $t_{\mathrm d}>1$, and a number
\begin{equation}\label{seed:growth-choice}
 G=\min_{k\in J}\mathscr G_k(t_{\mathrm d})>G_0,
 \qquad G>c_\star-\eta_{\mathrm{low}},
\end{equation}
where $\eta_{\mathrm{low}}$ is from Lemma~\ref{seed:low-frequency}. Fix $\Theta$ and take
$A$ sufficiently large. For the unnormalized first and second variations
initialized by $S_A$ and zero at $\tp$, define
\begin{equation}\label{seed:unnormalized-signed-mass}
 \widetilde{\mathcal M}_A=-\int_\tp^\tf(s_{h,A}-s_{z,A})(t)\idiff t,
 \qquad \tf=\tp+16A^2.
\end{equation}
Then
\begin{equation}\label{seed:unnormalized-mass-bounds}
 \frac{c_J}{A}\mathrm e^{2GA}\leq\widetilde{\mathcal M}_A
       \leq CA^P\mathrm e^{2c_\star A},\qquad
 \int_\tp^\tf\|v_A(t)\|_2^2\idiff t\leq C\widetilde{\mathcal M}_A,\qquad
 |\mathcal T_{v_A}(\tf)|\leq\frac{CA}{1+\tf}\widetilde{\mathcal M}_A.
\end{equation}
The same estimates hold with $\tf$ replaced by any endpoint $\tp+L$ with
$L\geq4A$.

For the normalized seed $g=S_A/n_A$, let $v_{\mathrm d},w_{\mathrm d}$ have
initial values $g,0$ at $\tp$, respectively. Write
\begin{equation}\label{seed:normalized-complete-stress}
 \mathsf S_{\mathrm d}=v_{\mathrm d}\otimes v_{\mathrm d}
       +U_{\mathrm d}\otimes w_{\mathrm d}+w_{\mathrm d}\otimes U_{\mathrm d},
 \qquad
 \mathcal M_{\mathrm d}=-\int_\tp^\tf(s_{h,\mathrm d}-s_{z,\mathrm d})(t)\idiff t.
\end{equation}
Their bounds are
\begin{equation}\label{seed:normalized-mass-bounds}
 A^{-P}\mathrm e^{2GA}\leq \mathcal M_{\mathrm d}
       \leq CA^P\mathrm e^{2c_\star A},\qquad
 \int_\tp^\tf\|v_{\mathrm d}(t)\|_2^2\idiff t\leq C\mathcal M_{\mathrm d},\qquad
 |\mathcal T_{v_{\mathrm d}}(\tf)|\leq\frac{CA}{1+\tf}\mathcal M_{\mathrm d}.
\end{equation}
In particular, with $\varrho=c_\star/8$, the parameters can be chosen so that
$0<2(c_\star-G)<\varrho/32$. The frequency interval, observation time, and
preparation constants are fixed before $A$.
\end{lemma}

\begin{proof}
Strict monotonicity of $\sigma_k$ gives
\begin{equation}\label{seed:integrated-growth-lower}
 \mathscr G_k(t)=\int_0^t(1+s)^{-2}\sigma_{k\sqrt{1+s}}\idiff s
       \geq\sigma_k\frac{t}{1+t}.
\end{equation}
Lemma~\ref{seed:rayleigh} therefore permits a fixed sufficiently high
interval $J$, followed by a fixed sufficiently late time $t_{\mathrm d}$,
for which \eqref{seed:growth-choice} holds. The positive coordinate
$a_0(k)$ has a positive minimum on this interval. By
Lemma~\ref{seed:point-growth} and the bounded physical adjoint norm,
\begin{equation}\label{seed:observed-band-energy}
 \int_{t_{\mathrm d}-A^{-1}}^{t_{\mathrm d}}
       \|\mathbf1_J(D_z)v_A(t)\|_2^2\idiff t
 \geq\frac{c_J}{A}\mathrm e^{2GA}.
\end{equation}
Indeed, $\mathscr G_k(t)\geq \mathscr G_k(t_{\mathrm d})-C/A$ on this interval. Its reflected
negative frequencies are present because $v_A$ is real. For large $A$,
the observation interval is contained in $[\tp,\tp+4A]$.

Split the first variation into its disjoint frequency parts $v_<,v_>$ below
and above two. Define their scalar stress histories by
\eqref{seed:exact-stress-formula}, each with zero initial azimuthal mean
of the second variation.
Parseval's identity makes every axial cross product vanish at each radius.
Thus the exact integrated quantities satisfy
\begin{equation}\label{seed:frequency-history-splitting}
 s_{h,A}-s_{z,A}=(s_{h,<}-s_{z,<})+(s_{h,>}-s_{z,>}),\qquad
 E_{v_A}=E_<+E_>,\qquad Z_{v_A}=Z_<+Z_>.
\end{equation}
This is a decomposition of the integrated complete history. No part of the
physical second variation has been removed. From
\eqref{seed:finite-stress-history}, $K_{\mathrm{hist}}(\tau,r)\leq2/\tau$, and
Lemma~\ref{seed:low-frequency},
\begin{equation}\label{seed:low-history-error}
 \left|\int_\tp^{t_e}(s_{h,<}-s_{z,<})(t)\idiff t\right|
 \leq C(1+A)\int_\tp^\infty\|v_<(t)\|_2^2\idiff t
 \leq CA^P\mathrm e^{2(c_\star-\eta_{\mathrm{low}})A}.
\end{equation}
For the high part, $\mathcal Q_>(\tp)\leq\mathcal Q_A(\tp)\leq CA^P$.
If $t_e=\tp+L$ with $L\geq4A$, then
$1/2-A/(1+t_e)\geq1/4$. Lemma~\ref{seed:high-frequency} therefore gives
\begin{equation}\label{seed:full-history-coercivity}
 -\int_\tp^{t_e}(s_{h,A}-s_{z,A})(t)\idiff t
 \geq\frac\pi4\int_\tp^{t_e}E_>(t)\idiff t
             -CA^P-CA^P\mathrm e^{2(c_\star-\eta_{\mathrm{low}})A}.
\end{equation}
Since $G>c_\star-\eta_{\mathrm{low}}$, the lower bound
\eqref{seed:observed-band-energy} absorbs both error terms into half the
positive energy contribution. This proves the lower stress-mass bound and
also gives $\int E_>\leq C\widetilde{\mathcal M}_A$. The low-frequency energy is
exponentially smaller by \eqref{seed:low-frequency-lifetime}. Hence the full
energy is bounded by the positive stress mass.

For the upper bound, integrate the sharp energy estimate
\eqref{seed:sharp-energy} from the polynomially bounded prepared seed. On
$|k|\geq2$, the factor $\mathrm e^{-8(t-\tp)}$ remains, so
\begin{equation}\label{seed:high-frequency-lifetime}
 \int_\tp^\infty\|v_>(t)\|_2^2\idiff t
 \leq CA^P\mathrm e^{2c_\star A}.
\end{equation}
The polynomial time factor is integrable against that exponential uniformly
for $\tp\leq1$. The low-frequency lifetime bound is smaller. Taking
absolute values in \eqref{seed:exact-stress-formula}, and using
\eqref{seed:kernel-time-integral} before integrating the source time, yields
\begin{equation}\label{seed:absolute-stress-history}
 \int_\tp^\infty|s_{h,A}(t)-s_{z,A}(t)|\idiff t
 \leq C(1+A)\int_\tp^\infty\|v_A(t)\|_2^2\idiff t
 \leq CA^P\mathrm e^{2c_\star A}.
\end{equation}
Finally, \eqref{seed:retained-mean-pairing} and the remainder bound in
\eqref{seed:finite-stress-history} prove the estimate for
$\mathcal T_{v_A}$.

Linearity of the first variation and quadratic homogeneity of the second
variation give
\begin{equation}\label{seed:normalization-relations}
 v_{\mathrm d}=v_A/n_A,\qquad w_{\mathrm d}=w_A/n_A^2,\qquad
 \mathcal M_{\mathrm d}=\widetilde{\mathcal M}_A/n_A^2,\qquad
 \mathcal T_{v_{\mathrm d}}=\mathcal T_{v_A}/n_A^2.
\end{equation}
The positive lower bound and polynomial upper bound for $n_A$ in
Lemma~\ref{seed:prepared} prove \eqref{seed:normalized-mass-bounds}.
\end{proof}

\section{The complete quadratic history at the source scale}
\label{sec:history}

We now transfer the sign estimate of Lemma~\ref{seed:signed-history}
from the invariant column to the exact finite-energy reference solution
of Proposition~\ref{bg:reference}.  This requires control of axial
localization, the remaining columns, and perturbations of the prepared
input.  In each comparison we retain the pressure and the interaction
between the reference velocity and the second variation.

Proposition~\ref{hist:energy} first preserves the sharp growth exponent
in a weighted estimate.  Proposition~\ref{hist:baseline} then compares
the complete stress histories, and Proposition~\ref{hist:generic}
extends the comparison to the input class used in the nonlinear
iteration.  These results supply the positive stress mass, stress moments,
axial fibers, and azimuthal mean required in Section~\ref{sec:transfer}.

Fix the constants supplied by the seed construction and set
\begin{equation}
 c_\star=\sqrt2e^{-3/2},\qquad
 \varrho=\frac{c_\star}{8},\qquad
 G=c_\star-\omega,\qquad
 0<2\omega<\frac{\varrho}{32},\qquad
 R=e^{\varrho A},\qquad b=e^{\beta A}.
 \label{hist:parameters}
\end{equation}
Here $\omega=c_\star-G$ is the fixed deficit between the sharp energy
rate and the lower growth rate of Lemma~\ref{seed:signed-history}.
The axial-width exponent $\beta$ is fixed before $A$ is taken sufficiently
large.  We use the preparation time, final source time, and moment
exponents given by
\begin{equation}
 \tp=\frac{2+\Theta\log A}{A},\qquad
 \tf=\tp+16A^2,\qquad
 s_\star=\frac{19}{8},\qquad
 \alpha=\frac74,\qquad
 \rho_c=\frac{3\varrho}{8},\qquad
 \gamma=\frac{\varrho}{4}.
 \label{hist:times-exponents}
\end{equation}
Here and below, $CA^P$ denotes a constant times a fixed finite power of $A$.
These constants may depend on fixed exponents and derivative orders, but
are independent of the inner and outer cutoffs and of the admissible input
shape.  The constants determining a lower bound on $\beta$ are fixed before
that bound is imposed.

\subsection{The isolated history and the actual reference solution}

We use the invariant column to normalize the prepared input, and the
exact finite-energy array to evolve it.  The two histories begin with
identical first variations and zero second variations at $\tp$.
The estimates below are uniform in the number of smaller and larger
columns.

Let $\mathcal U_{\mathrm d}$ be the linearized evolution about
$U_{\mathrm d}=AW_{1+t}e_\theta$.  Recalling
\eqref{seed:point-input} and \eqref{seed:prepared-definition}, the prepared
field and its normalization are
\begin{equation}
 S_A=\mathcal U_{\mathrm d}(\tp,0)\mathcal I,\qquad
 n_A=\|S_A\|_2,\qquad g=S_A/n_A.
 \label{hist:prepared}
\end{equation}
Lemma~\ref{seed:prepared} gives
\begin{equation}
 c_0\leq n_A\leq CA^P,\qquad \|g\|_2=1,\qquad
 \|g\|_{H^m}+M_\mu(g)\leq C_{m,\mu}A^{P_{m,\mu}}
 \quad (m\geq0,\ 0\leq\mu<5/2).
 \label{hist:prepared-bounds}
\end{equation}
All fields under consideration are real, axisymmetric, and solenoidal.
The prepared field is also centrally odd.

Write $v_{\mathrm d}$ and $w_{\mathrm d}$ for the first and second variations
about $U_{\mathrm d}$ with initial values $g$ and zero at $\tp$.  Thus the
coefficient of the square of the perturbation parameter is $w_{\mathrm d}$,
without a factorial.  In terms of the point responses from
\eqref{seed:point-input} and \eqref{seed:second-variation},
$v_{\mathrm d}=v_A/n_A$ and $w_{\mathrm d}=w_A/n_A^2$.
Define its complete stress by
\begin{equation}
 \mathsf S_{\mathrm d}
   =v_{\mathrm d}\otimes v_{\mathrm d}
     +U_{\mathrm d}\otimes w_{\mathrm d}
     +w_{\mathrm d}\otimes U_{\mathrm d},\qquad
 \int_{\R^3}\mathsf S_{\mathrm d}\idiff x
      =\operatorname{diag}(s_{h,\mathrm d},s_{h,\mathrm d},s_{z,\mathrm d}).
 \label{hist:isolated-stress}
\end{equation}
The seed estimate in Lemma~\ref{seed:signed-history} states that
\begin{equation}
 \mathcal M_{\mathrm d}:=-\int_\tp^\tf
       (s_{h,\mathrm d}-s_{z,\mathrm d})\idiff t
       \geq A^{-P}e^{2GA},\qquad
 \int_\tp^\tf\|v_{\mathrm d}(t)\|_2^2\idiff t\leq C \mathcal M_{\mathrm d}.
 \label{hist:isolated-mass}
\end{equation}
Thus $\mathcal M_{\mathrm d}$ is a positive stress mass; the minus sign in
\eqref{hist:isolated-mass} records the negative integrated anisotropy.
It is distinct from the spatial moment norm $M_\mu$ in \eqref{bg:moment}.
The second estimate in \eqref{hist:isolated-mass} controls the full
prepared first variation by this positive mass.

The reference solution from Proposition~\ref{bg:reference} is
$\bar u^{I,J}=Q^{I,J}+e^{I,J}$, where
\begin{equation}
 Q^{I,J}=\sum_{j=-I}^J Q_{R^j},\qquad
 Q_\ell(t,x)=A\ell^2\eta_\ell(t,z)(\ell^2+t)^{-2}
       e^{-r^2/(4(\ell^2+t))}x_h^\perp,
 \label{hist:array}
\end{equation}
and the axial cutoff for radius $\ell$ is
\begin{equation}
 \eta_\ell(t,z)=\frac{b\ell}{\sqrt{b^2\ell^2+2t}}
       \exp\left(-\frac{z^2}{2(b^2\ell^2+2t)}\right).
 \label{hist:cutoff}
\end{equation}
We allow $I=\infty$ and take $J\geq1$.  For every fixed derivative order
$m$, the correction satisfies
\begin{equation}
 \|\nabla^m e^{I,J}(t)\|_2\leq\Sigma_m t^{1/4-m/2},\qquad
 \|\nabla^m e^{I,J}(t)\|_\infty\leq\Sigma_m t^{-1/2-m/2},\qquad
 \Sigma_m\leq C_mA^{P_m}e^{K_mA}b^{-1/2}.
 \label{hist:reference-error}
\end{equation}
The constants $K_m$ are independent of $\beta$.  A fixed sufficiently large
width exponent therefore makes any prescribed finite list of these
coefficients at most one.

\subsection{A weighted energy estimate with the sharp growth integral}

The comparison requires the coefficient $c_\star$ in the leading
exponential.  We therefore use the bounded physical metric of
Lemma~\ref{seed:physical-energy}, with a spatial weight whose derivatives
contribute only a polynomial factor over $[\tp,\tf]$.
Axisymmetry reduces the pressure term to a lower spatial moment.

A fixed weight $\langle x\rangle^{2\mu}$ is useful for membership in a
moment space, but its derivatives can give an unnecessary exponential
loss in a sharp growth estimate.  We instead use
$w_\mu^{\mathrm{red}}=1+A^{-2\mu}\langle x\rangle^{2\mu}$.
It still controls the same ordinary moment, with a polynomial conversion
factor, while
\[
 \frac{|\nabla w_\mu^{\mathrm{red}}|}{w_\mu^{\mathrm{red}}}=O(A^{-1}),\qquad
 \frac{|\Delta w_\mu^{\mathrm{red}}|}{w_\mu^{\mathrm{red}}}=O(A^{-2}).
\]
The transition in this weight occurs at distance of order $A$, the
same order as the diffusion length over the source interval of length
$16A^2$.  Its diffusion cost integrated over that interval is consequently
bounded, while conversion back to $M_\mu$ costs only a power of $A$.

There are two further structural cancellations.  The metric $K$ changes
only the azimuthal component.  For an axisymmetric field,
$\operatorname{div}(a(r,z)e_\theta)=0$, so
$\operatorname{div}(Kv)=\operatorname{div}v=0$.
Also $\nabla w_\mu^{\mathrm{red}}$ has no azimuthal component.  Hence integration by
parts moves the entire pressure contribution to
$2\int p\nabla w_\mu^{\mathrm{red}}\cdot v$, with one derivative on the weight.
This lowers the pressure moment from $\mu$ to $\mu-1$, exactly the
range where the unrestricted double-Riesz estimate is available.
The main coupling retains $c_\star A/(1+t)^2$.
The remaining coefficients will have an $O(\log A)$ time integral,
which contributes the polynomial factor $e^{O(\log A)}=A^{O(1)}$.

Set $T=1+t$ and recall the weight $F_\star$ from
\eqref{seed:energy-weight}.  The bounded physical metric is
\begin{equation}
 \kappa_{\mathrm{met}}(t,r):=F_\star\left(\frac{r^2}{4T}\right),\qquad
 K=I_3+(\kappa_{\mathrm{met}}-1)e_\theta\otimes e_\theta.
 \label{hist:metric}
\end{equation}
For $1\leq\mu<5/2$, define
\begin{align}
 w_\mu^{\mathrm{red}}(x)&=1+A^{-2\mu}\langle x\rangle^{2\mu},&
 N_\mu(t,v)^2&=\int w_\mu^{\mathrm{red}}\,v\cdot Kv\idiff x,
 \label{hist:weighted-energy}\\
 D_\mu(t,v)&=\int w_\mu^{\mathrm{red}}\sum_{j=1}^3
           \partial_jv\cdot K\partial_jv\idiff x,&
 G_\mu(F)^2&=\int w_\mu^{\mathrm{red}}|F|^2\idiff x.
 \label{hist:dissipation-source}
\end{align}
In particular, $N_\mu(t,v)$ is uniformly equivalent to
$\|v\|_2+A^{-\mu}M_\mu(v)$.  We also use the unweighted metric energy
$E_0(t,v)=\int v\cdot Kv\idiff x$.

\begin{lemma}[The pressure estimate at a lower moment]
\label{hist:pressure}
Let $0\leq\lambda<3/2$ and $i,j\in\{1,2,3\}$.  Every scalar
function or tensor component $F$ with $M_\lambda(F)<\infty$ satisfies
\begin{equation}
 M_\lambda(\mathcal R_i\mathcal R_jF)\leq C_\lambda M_\lambda(F).
 \label{hist:riesz-bound}
\end{equation}
The same bound, with a constant independent of $L$, holds after replacing
$\langle x\rangle^\lambda$ by
$\min(\langle x\rangle^\lambda,L^\lambda)$.
\end{lemma}

\begin{proof}
This is \eqref{bg:riesz}, including the cap-independent version proved
in Lemma~\ref{bg:weighted}.  We shall apply it at $\lambda=\mu-1<3/2$.
\end{proof}

\begin{proposition}[The actual source-scale energy estimate]
\label{hist:energy}
Fix $1\leq\mu<5/2$. Choose $\beta$ sufficiently large and then $A$
sufficiently large. Let $v$ be axisymmetric and solenoidal, with
$v\in L^\infty(\tp,\tf;L^2)\cap L^2(\tp,\tf;H^1)$ and
$N_\mu(\tp,v(\tp))<\infty$. Let $F$ be an axisymmetric Cartesian tensor
with $G_\mu(F)\in L^2(\tp,\tf)$. Suppose that
\begin{equation}
 \partial_tv-\Delta v+\bar u^{I,J}\cdot\nabla v
       +v\cdot\nabla\bar u^{I,J}+\nabla p=-\operatorname{div}F.
 \label{hist:linear-equation}
\end{equation}
There is a nonnegative integrable function $a_\mu$ such that the following
inequality holds on $[\tp,\tf]$ in the sense of distributions in time:
\begin{equation}
 \frac{\diff}{\diff t}N_\mu(t,v)^2+c_\mu D_\mu(t,v)
 \leq\bigl(2\Gamma_{\mathrm{src}}'(t)+a_\mu(t)\bigr)N_\mu(t,v)^2
       +C_\mu G_\mu(F)^2,\qquad
 \int_\tp^\tf a_\mu(t)\idiff t\leq C_\mu(1+\log A),
 \label{hist:sharp-energy}
\end{equation}
where
\begin{equation}
 \Gamma_{\mathrm{src}}(t)=c_\star A\bigl((1+\tp)^{-1}-(1+t)^{-1}\bigr).
 \label{hist:growth}
\end{equation}
The corresponding unweighted estimate holds for $E_0$.  If $v(\tp)=0$,
$m\geq1$, and
$G_\mu(F(t))+\|F(t)\|_{H^4}\leq B e^{m\Gamma_{\mathrm{src}}(t)}$, then
\begin{equation}
 N_\mu(t,v)+\|v(t)\|_{H^4}
       \leq C_\mu A^{P_\mu}B e^{m\Gamma_{\mathrm{src}}(t)}.
 \label{hist:source}
\end{equation}
Its constant is independent of $m$.  The bound for $N_\mu$ requires
only $G_\mu(F(t))\leq Be^{m\Gamma_{\mathrm{src}}(t)}$; the additional $H^4$ hypothesis
is used for the Sobolev bound.  Homogeneous solutions satisfy the
corresponding estimate with
$N_\mu(\tp,v(\tp))+\|v(\tp)\|_{H^4}$ and exponent $\Gamma_{\mathrm{src}}(t)$.
For a solenoidal vector source $H$ and zero initial value, one also has
\begin{equation}
 N_\mu(t,v)\leq CA^P\int_\tp^t
       e^{\Gamma_{\mathrm{src}}(t)-\Gamma_{\mathrm{src}}(u)}N_\mu(u,H(u))\idiff u.
 \label{hist:vector-source}
\end{equation}
\end{proposition}

\begin{proof}
Write the radius-one column as $Q_1=\Omega_{\mathrm{act}}x_h^\perp$, where
$\Omega_{\mathrm{act}}=A\eta_1\Omega_T(r)$ and $\zeta_r=r^2/(4T)$.
Apart from the axial shear, its radial and azimuthal couplings are
$(2\Omega_{\mathrm{act}} v_\theta,-2\Omega_{\mathrm{act}}(1-\zeta_r)v_r)$.  Their contribution to one half
of the metric energy derivative is
$2\Omega_{\mathrm{act}}(1-\kappa_{\mathrm{met}}(1-\zeta_r))v_rv_\theta$.
It vanishes for $\zeta_r\leq1/2$.  For $\zeta_r\geq1/2$, its coefficient
relative to $v_r^2+2v_\theta^2$ is bounded by
\begin{equation}
 \frac{A\eta_1}{\sqrt2T^2}(2\zeta_r-1)e^{-\zeta_r}
       \leq\frac{c_\star A}{T^2}.
 \label{hist:sharp-coupling}
\end{equation}
Indeed the maximum is attained at $\zeta=3/2$.  Multiplication by
$w_\mu^{\mathrm{red}}$ leaves this coefficient unchanged.  The axial shear contributes at
most $CAb^{-1}T^{-3/2}N_\mu(t,v)^2$.

The metric extends smoothly across the axis and obeys, in the sense of
matrix distributions,
\begin{equation}
 \|\nabla K\|_\infty\leq CT^{-1/2},\qquad
 \partial_tK+\Delta K\leq CT^{-1}I_3.
 \label{hist:metric-derivatives}
\end{equation}
The radial derivative of $\kappa_{\mathrm{met}}$ jumps downward at $r^2=2T$, so the
interface measure has the favorable sign.  Angular differentiation is
controlled by
$(\kappa_{\mathrm{met}}-1)/r^2\leq C/T$ and
$(\kappa_{\mathrm{met}}-1)/r\leq C/\sqrt T$.
Moreover,
\begin{equation}
 \frac{|\nabla w_\mu^{\mathrm{red}}|}{w_\mu^{\mathrm{red}}}\leq C_\mu A^{-1},\qquad
 \frac{|\Delta w_\mu^{\mathrm{red}}|}{w_\mu^{\mathrm{red}}}\leq C_\mu A^{-2}.
 \label{hist:weight-derivatives}
\end{equation}
Consequently the exact heat and metric contribution,
\begin{equation}
 -2D_\mu(t,v)+\int v\cdot
 \bigl(w_\mu^{\mathrm{red}}(\partial_tK+\Delta K)
          +2\nabla w_\mu^{\mathrm{red}}\cdot\nabla K+(\Delta w_\mu^{\mathrm{red}})K\bigr)v\idiff x,
 \label{hist:heat-identity}
\end{equation}
is at most
$-2D_\mu(t,v)+C_\mu(T^{-1}+A^{-2})N_\mu(t,v)^2$.
The divergence-free correction $e^{I,J}$ contributes at most
\begin{equation}
 C_\mu\bigl(\|\nabla e^{I,J}\|_\infty
       +(T^{-1/2}+A^{-1})\|e^{I,J}\|_\infty\bigr)N_\mu(t,v)^2.
 \label{hist:correction-energy}
\end{equation}

The pressure requires only the lower moment $\mu-1<3/2$.
Axisymmetry gives $\operatorname{div}(Kv)=0$ and
$\nabla w_\mu^{\mathrm{red}}\cdot Kv=\nabla w_\mu^{\mathrm{red}}\cdot v$, so its contribution is
$2\int p\nabla w_\mu^{\mathrm{red}}\cdot v\idiff x$, where
\begin{equation}
 p=\sum_{i,j}\mathcal R_i\mathcal R_j
   \bigl(\bar u^{I,J}_iv_j+v_i\bar u^{I,J}_j+F_{ij}\bigr).
 \label{hist:actual-pressure}
\end{equation}
Writing $N=N_\mu(t,v)$, the pointwise weight comparisons give
\[
 \frac{|\nabla w_\mu^{\mathrm{red}}|}{\sqrt{w_\mu^{\mathrm{red}}}}
       \leq C_\mu A^{-\mu}\langle x\rangle^{\mu-1},\qquad
 A^{-\mu}M_{\mu-1}(v)\leq C_\mu A^{-1}N,\qquad
 A^{-\mu}M_{\mu-1}(F)\leq C_\mu A^{-1}G_\mu(F).
\]
The latter two inequalities follow by comparing the weights pointwise,
separately on $\langle x\rangle\leq A$ and its complement.
Consequently, with $\lambda=\mu-1<3/2$,
\[
 \begin{split}
 \left|2\int p\nabla w_\mu^{\mathrm{red}}\cdot v\idiff x\right|
 &\leq C_\mu A^{-\mu}M_{\mu-1}(p)N\\
 &\leq C_\mu A^{-1}\|\bar u^{I,J}\|_\infty N^2
       +C_\mu A^{-1}G_\mu(F)N\\
 &\leq C_\mu\bigl(A^{-1}\|\bar u^{I,J}\|_\infty+A^{-2}\bigr)N^2
       +C_\mu G_\mu(F)^2.
 \end{split}
\]
The direct tensor pairing is estimated by Cauchy--Schwarz against the
weighted heat form.  Derivatives falling on $K$ or $w_\mu^{\mathrm{red}}$ cost only
$C_\mu(T^{-1}+A^{-2})N^2+C_\mu G_\mu(F)^2$.
For the cap $\min(\langle x\rangle^{2\mu},L^{2\mu})$, the derivative
vanishes outside the cap and the same pointwise comparisons hold
with the corresponding capped lower moment.  The interface measure
in its Laplacian is nonpositive.  This proves the energy estimate uniformly in the cap.
The direct tensor term is $2\int F:\nabla(w_\mu^{\mathrm{red}} Kv)\idiff x$;
its derivative of $v$ is absorbed into a fixed part of $D_\mu(t,v)$.
No derivative of the tensor source $F$ is used.

For an additional column of radius $\ell$, the complete axisymmetric
matrix, including its axial shear, has norm at most
$CA\ell^2(\ell^2+t)^{-2}$.  Its velocity norm is at most
$CA\ell^2(\ell^2+t)^{-3/2}$.
The inner and outer sums therefore satisfy
\begin{align}
 \sum_{j<0}\|M_{Q_{R^j}}\|_\infty&\leq CAR^{-2}t^{-2},&
 \sum_{j>0}\|M_{Q_{R^j}}\|_\infty&\leq CAR^{-2},
 \label{hist:matrix-tails}\\
 \sum_{j<0}\|Q_{R^j}\|_\infty&\leq CAR^{-2}t^{-3/2},&
 \sum_{j>0}\|Q_{R^j}\|_\infty&\leq CAR^{-1}.
 \label{hist:velocity-tails}
\end{align}
One may therefore take $a_\mu$ bounded by
\begin{align}
 C_\mu\bigl(&T^{-1}+A^{-2}+AR^{-2}(t^{-2}+1)
        +Ab^{-1}T^{-3/2}+T^{-3/2}+R^{-2}t^{-3/2}+R^{-1}
 \notag\\[-2pt]
 &+\|\nabla e^{I,J}\|_\infty
      +(T^{-1/2}+A^{-1})\|e^{I,J}\|_\infty\bigr).
 \label{hist:remainder-coefficient}
\end{align}
The integrals of its first two terms are $O(1+\log A)$ and $16$.
The remaining column terms have integrals bounded by
$CAR^{-2}(\tp^{-1}+16A^2)+C+CR^{-2}\tp^{-1/2}+16A^2R^{-1}+CA/b$.
By \eqref{hist:reference-error}, the correction terms have integrals bounded
by $C\Sigma_1\log(\tf/\tp)+C\Sigma_0(1+\log A)$.
This proves the asserted integral bound.

To justify the calculation for an ordinary moment, fix $L\geq A$ and cap
$\langle x\rangle^{2\mu}$ at $L^{2\mu}$.  Its interface measure is
nonpositive, and \eqref{hist:weight-derivatives} is uniform in $L$.
The pressure argument uses the corresponding capped $(\mu-1)$-moment.
Apply the bounded-weight estimate after heat regularization, remove the
regularization, and then let $L\to\infty$.

For zero data, the integrating factor in \eqref{hist:sharp-energy} gives
\begin{equation}
 N_\mu(t,v)^2\leq CA^PB^2e^{2\Gamma_{\mathrm{src}}(t)}
       \int_\tp^t e^{2(m-1)\Gamma_{\mathrm{src}}(u)}\idiff u
       \leq CA^P(16A^2)B^2e^{2m\Gamma_{\mathrm{src}}(t)}.
 \label{hist:source-integrating-factor}
\end{equation}

This estimate is uniform in $m\geq1$.  We next obtain its $H^4$
supplement without changing the growth exponent.  Let
$\mathcal U_{\bar u}(t,u)$ denote the homogeneous propagator of
\eqref{hist:linear-equation}.  Gaussian differentiation and
\eqref{hist:reference-error} give a polynomial bound $D_{\mathrm{coeff}}(A)\leq CA^{P_0}$
for the pure-swirl matrix derivatives through order four and the
correction derivatives through order five, after fixing $\beta$
sufficiently large.  The differentiated, unweighted energy inequalities
on a short interval are
\begin{equation}
 \frac{\diff}{\diff t}\|\nabla^jv\|_2^2
       +\|\nabla^{j+1}v\|_2^2
 \leq D_{\mathrm{coeff}}(A)\sum_{i=0}^j\|\nabla^iv\|_2^2,
 \qquad 0\leq j\leq4.
 \label{hist:short-time-energies}
\end{equation}
Here the undifferentiated correction transport integrates to zero,
and its commutators contain lower derivatives of $v$.
Multiplying the inequality of order $j$ by a sufficiently small fixed
multiple of $(t-u)^j$ and summing gives the usual time-weighted heat
energy.  Starting at orders zero and three, respectively, yields
\begin{equation}
 \|\mathcal U_{\bar u}(u+h,u)q\|_{H^4}
       \leq C h^{-2}\|q\|_2,\qquad
 \|\mathcal U_{\bar u}(u+h,u)q\|_{H^4}
       \leq C h^{-1/2}\|q\|_{H^3},
 \quad 0<h\leq A^{-D},
 \label{hist:short-time-smoothing}
\end{equation}
where the fixed $D$ is chosen so that $D_{\mathrm{coeff}}(A)A^{-D}\leq1$.
Since $\|\PP\div F\|_{H^3}\leq C\|F\|_{H^4}$, the second bound
is precisely the integrable tensor-source kernel needed here.

For zero data put $h=\min(A^{-D},t-\tp)$ and write
\[
 v(t)=\mathcal U_{\bar u}(t,t-h)v(t-h)
      -\int_{t-h}^t\mathcal U_{\bar u}(t,u)\PP\div F(u)\idiff u.
\]
If $t-h>\tp$, estimate the first term by the established $L^2$ energy
bound and final-interval smoothing; otherwise it is zero.
The recent source satisfies
\[
 CA^P\int_{t-h}^t(t-u)^{-1/2}B e^{m\Gamma_{\mathrm{src}}(u)}\idiff u
       \leq2CA^P\sqrt h\,B e^{m\Gamma_{\mathrm{src}}(t)}.
\]
Here only monotonicity of $\Gamma_{\mathrm{src}}$ is used, not a derivative of
$e^{m\Gamma_{\mathrm{src}}}$.  Thus neither the smoothing interval nor the constant
depends on the nonlinear degree $m$.  For homogeneous initial data
on a shorter initial interval, use the $H^4$-to-$H^4$ estimate instead
of singular smoothing.  This proves \eqref{hist:source}; the
homogeneous and vector-source versions follow from the same energy
estimate and Duhamel's formula.

\end{proof}

\subsection{Comparison with the isolated source}

The differences begin from zero, so \eqref{hist:source} applies
directly.  The main subtlety is global rather than local: a distant
column can have a large total energy.  Negative spatial weights exploit
its separation when it is paired with the localized second variation,
allowing the complete stress comparison to be uniform in the outer cutoff.

For the actual reference solution define
$\mathcal L_{\bar u}v=\Delta v-\PP\operatorname{div}(\bar u\otimes v+v\otimes\bar u)$.
Let $v_g(\tp)=g$ and $w_g(\tp)=0$, with
\begin{equation}
 (\partial_t-\mathcal L_{\bar u})v_g=0,\qquad
 (\partial_t-\mathcal L_{\bar u})w_g
        =-\PP\operatorname{div}(v_g\otimes v_g).
 \label{hist:actual-variations}
\end{equation}
Define
\begin{equation}
 \mathsf S_g=v_g\otimes v_g+\bar u\otimes w_g+w_g\otimes\bar u,
 \qquad
 \int\mathsf S_g\idiff x=\operatorname{diag}(s_{h,g},s_{h,g},s_{z,g}),\qquad
 \mathcal M_g=-\int_\tp^\tf(s_{h,g}-s_{z,g})\idiff t.
 \label{hist:actual-stress}
\end{equation}

\begin{proposition}[Comparison of complete histories]
\label{hist:baseline}
There is a fixed $K_\ast$, independent of $\beta$, such that, with
$\xi_A=R^{-2}+e^{K_\ast A}b^{-1/2}$, the following bounds hold uniformly
in both cutoffs:
\begin{equation}\label{hist:baseline-bounds}
 \begin{aligned}
 N_\mu(t,v_g(t))+\|v_g(t)\|_{H^4}&\leq CA^Pe^{\Gamma_{\mathrm{src}}(t)},\\
 N_\mu(t,w_g(t))+\|w_g(t)\|_{H^4}&\leq CA^Pe^{2\Gamma_{\mathrm{src}}(t)}.
 \end{aligned}
\end{equation}
The differences and masses satisfy
\begin{align}
 N_1(t,v_g-v_{\mathrm d})&\leq CA^P\xi_Ae^{\Gamma_{\mathrm{src}}},&
 N_1(t,w_g-w_{\mathrm d})&\leq CA^P\xi_Ae^{2\Gamma_{\mathrm{src}}},
 \label{hist:baseline-difference}\\
 |\mathcal M_g-\mathcal M_{\mathrm d}|&\leq CA^P
       \bigl(R^{-15/8}+e^{K_\ast A}b^{-1/2}\bigr)e^{2c_\star A}.
 \label{hist:mass-comparison}
\end{align}
The estimates in \eqref{hist:baseline-bounds} hold for each fixed $1\leq\mu<5/2$.
If $\beta/2>K_\ast+2\omega$, then the mass satisfies
\begin{equation}
 A^{-P}e^{2GA}\leq \mathcal M_g\leq CA^Pe^{2c_\star A},\qquad
 \mathcal M_g/\mathcal M_{\mathrm d}=1+o(1).
 \label{hist:baseline-mass}
\end{equation}
The complete stress also satisfies
\begin{equation}
 \int_\tp^\tf\int_{\R^3}
  \bigl(1+|x|^\alpha+(\tf-t)^{\alpha/2}\bigr)|\mathsf S_g(t,x)|
       \idiff x\idiff t\leq CA^Pe^{2c_\star A}.
 \label{hist:baseline-moment}
\end{equation}
\end{proposition}

The two product estimates used in this section have different meanings.
The ordinary moment is placed entirely on one factor:
\[
 M_\mu(v\otimes w)\leq\|v\|_\infty M_\mu(w).
\]
For an equal-factor Sobolev product one can use the tame bound
$\|v\otimes v\|_{H^4}\leq C\|v\|_\infty\|v\|_{H^4}$.
To check its highest derivatives, a typical term is
$\partial^jv\,\partial^{4-j}v$.  For $0<j<4$, H\"older and
Sobolev interpolation give
\[
 \begin{aligned}
 \|\partial^jv\,\partial^{4-j}v\|_2
 &\leq \|\partial^jv\|_{8/j}
                       \|\partial^{4-j}v\|_{8/(4-j)},\\
 \|\partial^jv\|_{8/j}
 &\leq C\|v\|_\infty^{1-j/4}\|\nabla^4v\|_2^{j/4}.
 \end{aligned}
\]
The complementary powers multiply to
$C\|v\|_\infty\|\nabla^4v\|_2$.  The endpoint products put one
undifferentiated factor in $L^\infty$, and lower derivative orders are
treated in the same way.  The same product-rule argument gives the
corresponding bilinear bounds.  Thus the first-variation estimate produces an $e^{2\Gamma_{\mathrm{src}}}$
quadratic source with only polynomial additional cost.  It does not
require adding a weighted spatial derivative to the data.

\begin{proof}
Proposition~\ref{hist:energy} and \eqref{hist:prepared-bounds} give the
first-variation estimate.  Since $H^4(\R^3)$ is an algebra and the whole
weight can be placed on one factor,
$G_\mu(v_g\otimes v_g)+\|v_g\otimes v_g\|_{H^4}\leq CA^Pe^{2\Gamma_{\mathrm{src}}}$.
The source estimate \eqref{hist:source} gives the second variation.
The same argument applies to the isolated variations.

Put $a=\bar u-U_{\mathrm d}$,
$v_\Delta=v_g-v_{\mathrm d}$, and
$w_\Delta=w_g-w_{\mathrm d}$.
Both differences start from zero and solve
\begin{align}
 (\partial_t-\mathcal L_{\bar u})v_\Delta
   &=-\PP\operatorname{div}(a\otimes v_{\mathrm d}+v_{\mathrm d}\otimes a),
 \label{hist:first-difference}\\
 (\partial_t-\mathcal L_{\bar u})w_\Delta
   &=-\PP\operatorname{div}
      \bigl(a\otimes w_{\mathrm d}+w_{\mathrm d}\otimes a
                 +v_\Delta\otimes v_g+v_{\mathrm d}\otimes v_\Delta\bigr).
 \label{hist:second-difference}
\end{align}
We estimate the moment-one differences directly from the weighted
energy, which uses no derivative of a tensor source.  In the next
three bounds, let $(q,n)$ be either $(v_{\mathrm d},1)$ or
$(w_{\mathrm d},2)$, so that the isolated variation estimates give
$M_2(q)+N_1(t,q)\leq CA^Pe^{n\Gamma_{\mathrm{src}}(t)}$.
The extra columns act through their axisymmetric matrices.  Since the
Leray projection is bounded at moment one by \eqref{bg:riesz}, their
vector source satisfies
\begin{equation}
 N_1\left(t,\PP\sum_{j\neq0}M_{Q_{R^j}}q\right)
 \leq CA^PR^{-2}e^{n\Gamma_{\mathrm{src}}(t)}.
 \label{hist:tail-vector-source}
\end{equation}
Here the factors $t^{-2}+1$ in \eqref{hist:matrix-tails} are polynomial
on $[\tp,\tf]$.
For the main cutoff, Gaussian differentiation gives
\begin{equation}
 |1-\eta_1(t,z)|\leq C(t/b^2+|z|/b),
 \label{hist:main-cutoff-error}
\end{equation}
and hence
\begin{align}
 G_1\bigl((Q_1-U_{\mathrm d})\otimes q\bigr)
 &\leq C\|U_{\mathrm d}\|_\infty
       \left(\frac{t}{b^2}N_1(t,q)
                      +\frac{M_1(q)+A^{-1}M_2(q)}{b}\right)
 \leq CA^Pb^{-1}e^{n\Gamma_{\mathrm{src}}(t)},\notag\\
 G_1(e^{I,J}\otimes q)
 &\leq C\|e^{I,J}\|_\infty N_1(t,q)
 \leq CA^Pe^{K_\ast A}b^{-1/2}e^{n\Gamma_{\mathrm{src}}(t)}.
 \label{hist:cutoff-correction-sources}
\end{align}
The transposed tensors obey the same estimates.
Equations \eqref{hist:source-integrating-factor} and
\eqref{hist:vector-source} therefore prove the first bound in
\eqref{hist:baseline-difference}.  For the second bound, the remaining
tensor in \eqref{hist:second-difference} satisfies
\begin{equation}
 G_1\bigl(v_\Delta\otimes v_g+v_{\mathrm d}\otimes v_\Delta\bigr)
 \leq C\bigl(\|v_g\|_\infty+\|v_{\mathrm d}\|_\infty\bigr)
                  N_1(t,v_\Delta)
 \leq CA^P\xi_Ae^{2\Gamma_{\mathrm{src}}(t)}.
 \label{hist:quadratic-difference-source}
\end{equation}
Applying the same two response estimates to
\eqref{hist:second-difference} proves the asserted second-variation
comparison.

The exact stress difference contains
\[
 v_g\otimes v_g-v_{\mathrm d}\otimes v_{\mathrm d}
 +U_{\mathrm d}\otimes(w_g-w_{\mathrm d})
 +(\bar u-U_{\mathrm d})\otimes w_g
\]
and the corresponding transposed mixed terms.  The first line of
variation bounds does not by itself estimate the integrals of the mixed
terms: $U_{\mathrm d}$ has infinite three-dimensional energy, and
large-radius columns need not have small global $L^2$ norm.
For the invariant column, use axial integrability from a moment; for a
large-radius column, place a negative spatial weight on the background.

The latter mechanism can be seen at one radius $D\geq R$ on the source
interval.  With $x=Dy$, Gaussian integration and $\eta_D\leq1$ give
\[
 \|\langle x\rangle^{-\mu}Q_D\|_2^2
 \leq C A^2D^{1-2\mu}
     \int_{\R^3}|y|^{-2\mu}|y_h|^2e^{-c|y_h|^2}\,\idiff y.
\]
After axial integration, the radial integrand is a constant times
$r^{4-2\mu}e^{-cr^2}$.  The axial tail requires $\mu>1/2$ and the
origin requires $\mu<5/2$.  Thus
$\|\langle x\rangle^{-\mu}Q_D\|_2\leq C_\mu A D^{1/2-\mu}$.
Pairing this with $M_\mu(w_g)$ produces a small contribution
from each larger column, with a geometric sum independent of the outer
cutoff.  For the $\alpha$-moment of the stress, the same calculation
uses $\mu=s_\star-\alpha=5/8>1/2$.
These different exponents explain why the unweighted mass comparison and
the fractional stress moment are proved separately.

For the stress comparison, the first-variation products are controlled in
$L^1$ by \eqref{hist:baseline-difference}.
The invariant main column obeys
$\|U_{\mathrm d}(t)\|_{L^2_{x_h}}\leq CA/(1+t)$, whereas
$\|h\|_{L^1_zL^2_{x_h}}\leq\sqrt{\pi A}\,N_1(t,h)$.
These bounds control $U_{\mathrm d}\otimes w_\Delta$.
The localized main column is treated using
\begin{equation}
 \|(1-\eta_1(t,z))(1+z^2)^{-1/2}\|_{L^2_z}
        \leq C(b^{-1/2}+t/b^2),
 \label{hist:cutoff-dual}
\end{equation}
obtained by splitting at $|z|=b$.
The correction product is bounded by
$\|e^{I,J}\|_2\|w_g\|_2$.

The remaining columns require negative weights.  For $\mu>1/2$,
\begin{equation}
 \left\|\langle x\rangle^{-\mu}\sum_{j<0}Q_{R^j}(t)\right\|_2
       \leq C_\mu AR^{-2}t^{-1},
 \label{hist:inner-negative}
\end{equation}
because $\eta_\ell\leq1$, the $z$-weight is integrable, and the horizontal
Gaussian norm is $CA\ell^2/(\ell^2+t)$.
For $D\geq R$ and $1/2<\mu<5/2$, isotropic rescaling gives
\begin{equation}
 \|\langle x\rangle^{-\mu}Q_D(t)\|_2\leq C_\mu AD^{1/2-\mu}.
 \label{hist:outer-negative}
\end{equation}
Indeed the scaled integral is bounded by the identities and integrability
conditions
\begin{equation}
 \int_\R(r^2+z^2)^{-\mu}\idiff z=C_\mu r^{1-2\mu},\qquad
 \int_0^\infty r^{4-2\mu}e^{-c_1r^2}\idiff r<\infty.
 \label{hist:negative-integrals}
\end{equation}
At $\mu=s_\star$, the outer sum is $O(AR^{-15/8})$.
Pairing with $M_{s_\star}(w_g)$ and integrating in time proves
\eqref{hist:mass-comparison}.  The error relative to
\eqref{hist:isolated-mass} tends to zero because
$(15/8)\varrho>2\omega$ and $\beta/2>K_\ast+2\omega$.
This proves \eqref{hist:baseline-mass}.

For \eqref{hist:baseline-moment}, use
$\mu=s_\star-\alpha=5/8$ in \eqref{hist:outer-negative}; the full outer
sum is then $CA\sum_{j\geq1}R^{-j/8}$.
The three stress terms satisfy
\begin{align}
 \int\langle x\rangle^\alpha|v_g\otimes v_g|\idiff x
       &\leq M_{\alpha/2}(v_g)^2,\notag\\
 \int\langle x\rangle^\alpha|Q^{I,J}\otimes w_g|\idiff x
       &\leq\|\langle x\rangle^{\alpha-s_\star}Q^{I,J}\|_2
                  M_{s_\star}(w_g),\notag\\
 \int\langle x\rangle^\alpha|e^{I,J}\otimes w_g|\idiff x
       &\leq\|e^{I,J}\|_2M_\alpha(w_g).
 \label{hist:fractional-products}
\end{align}
All moments are below $5/2$, and the time factor is at most
$(16A^2)^{\alpha/2}$.  Equations \eqref{hist:baseline-bounds} and
\eqref{hist:reference-error} finish the proof.
\end{proof}

\subsection{The retained azimuthal mean and the low frequencies}

Two pieces of information are needed beyond the integrated mass.  The
azimuthal mean changes the next leading amplitude, and the low-frequency
fibers can be large before smoothing in parent coordinates.  We bound
them separately.

The azimuthal mean persists during diffusion and contributes to the response
of the larger scale.  Its relevant scalar is the following explicit
functional:
\begin{equation}
 \mathcal T_{I,\nu}(h)=4\pi A\int_0^\infty\int_\R e^{-\iunit\nu z}
 \left(\sum_{j=-I}^0R^{2j}
       \left(1-e^{-r^2/(4(\tf+R^{2j}))}\right)\right)
 h_\theta(r,z)\idiff z\idiff r.
 \label{hist:memory-functional}
\end{equation}
We use the unitary axial transform
$\widehat h^{\,z}(x_h,\nu)=(2\pi)^{-1/2}\int e^{-\iunit\nu z}h(x_h,z)\idiff z$.
Thus \eqref{hist:memory-functional} equals $\sqrt{2\pi}$ times its
expression with $h_\theta$ replaced by $\widehat h_\theta^{\,z}$ and the
$z$-integral omitted.

In \eqref{hist:memory-functional}, $\nu$ is a \emph{physical
source-scale axial frequency}.  The parent-scale frequencies used in
Section~\ref{sec:transfer} will be denoted by $k$, with $\nu=k/R$.
The index $I$ counts only the old or inner columns in the sum.
The weight in parentheses is an explicit radial heat-history kernel;
$\mathcal T_{I,\nu}(h)$ is a scalar \emph{linear functional} of the
azimuthal component of $h$.
Its value at $\nu=0$ is understood from the continuous Fourier
representative established below, rather than from the value of an
arbitrary almost-everywhere representative.

\begin{lemma}[Bounds for the retained mean and the fibers]
\label{hist:mean-fiber}
The baseline second variation satisfies
\begin{equation}
 \frac{|\mathcal T_{I,0}(w_g(\tf))|}{\mathcal M_g}
       \leq\frac CA+CA^P\xi_Ae^{2\omega A}.
 \label{hist:baseline-memory}
\end{equation}
For every fixed $K_0\geq1$, its low-frequency fibers obey
\begin{equation}
 \sup_{|k|\leq K_0A}\frac{R^2}{\mathcal M_g}
       \|\widehat w_g^{\,z}(\tf,\cdot,k/R)\|_{L^2_{x_h}}
       \leq C_{K_0}A^2R^2.
 \label{hist:baseline-fiber}
\end{equation}
The bound includes $k=0$.
\end{lemma}

\begin{proof}
For the isolated second variation put
$\mu_{\theta,\mathrm d}:=\mu_\theta[w_{\mathrm d}]$, with the mean operator
from \eqref{seed:azimuthal-mean}, and
$F(t,r)=\int_\R v_{{\mathrm d},r}v_{{\mathrm d},\theta}\idiff z$.
The axial mean of the radial second variation vanishes: its radial divergence
is zero, and a nonzero multiple of $1/r$ is excluded by its physical
$L^2_{x_h}$ norm.  Therefore
\begin{equation}
 (\partial_t-\Delta_h^{(1)})\mu_{\theta,\mathrm d}=-(\partial_r+2/r)F,\qquad
 \mu_{\theta,\mathrm d}(\tp)=0.
 \label{hist:mean-equation}
\end{equation}
For a heat age $\ell>0$, define
\begin{equation}\label{hist:mean-kernels}
 K_\ell(r)=1-e^{-r^2/(4\ell)}.
\end{equation}
The integrated history kernel $K_{\mathrm{hist}}$ is the one in
\eqref{seed:stress-kernels}.  The radial test $\chi_\ell=K_\ell/r$
has a $1/r$ tail and is not in
$L^2((0,\infty),r\idiff r)$.  The heat identity below is understood in
weighted duality, not as a direct Hilbert-space pairing of two $L^2$
vectors.  Formally its two kernel relations are
\begin{equation}
 e^{u\Delta_h^{(1)}}\frac{K_\ell}{r}=\frac{K_{\ell+u}}r,\qquad
 (\partial_r-r^{-1})\frac{K_\ell}{r}=-\frac14K_{\mathrm{hist}}(\ell,r).
 \label{hist:kernel-identities}
\end{equation}
We justify the first relation and its use by a heat-age truncation.
With $W_\tau$ as in \eqref{bg:unit-profiles}, direct integration gives
\begin{equation}\label{hist:regularized-mean-test}
 \chi_\ell=\frac14\int_\ell^\infty W_\tau\idiff\tau,
 \qquad
 \chi_{\ell,L}=\frac14\int_\ell^L W_\tau\idiff\tau
 =\frac{e^{-r^2/(4L)}-e^{-r^2/(4\ell)}}r,
 \quad L>\ell.
\end{equation}
The truncated test is of order $r$ at the axis and decreases Gaussianly
at infinity, so it is an admissible smooth $L^2(r\idiff r)$ heat test.
Self-adjointness and $e^{u\Delta_h^{(1)}}W_\tau=W_{\tau+u}$ give
\[
 e^{u\Delta_h^{(1)}}\chi_{\ell,L}=\chi_{\ell+u,L+u},\qquad
 (\partial_r-r^{-1})\chi_{\ell,L}
      =-\tfrac14(K_{\mathrm{hist}}(\ell,r)-K_{\mathrm{hist}}(L,r)).
\]
Apply these identities first to the mean equation with
$\ell=1+\tf$.  Integration by parts yields
\begin{equation}\label{hist:regularized-mean-history}
 \begin{aligned}
 4\pi A\int_0^\infty\chi_{\ell,L}(r)\mu_{\theta,\mathrm d}(\tf,r)r\idiff r
 ={}&-\pi A\int_\tp^\tf\int_0^\infty
 \bigl(K_{\mathrm{hist}}(\ell+\tf-u,r)-K_{\mathrm{hist}}(L+\tf-u,r)\bigr)\\
 &\hspace{33mm}{}\cdot F(u,r)r\idiff r\idiff u.
 \end{aligned}
\end{equation}
The sign follows from the adjoint of $-(\partial_r+2/r)$ in
$L^2(r\idiff r)$, namely $\partial_r-r^{-1}$.
The extra $K_{\mathrm{hist}}(L+\tf-u,r)$ term tends to zero as $L\to\infty$, since
$0\leq K_{\mathrm{hist}}(L+\tf-u,r)\leq2/L$ and
$\int_\tp^\tf\int_0^\infty |F|r\idiff r\idiff u<\infty$.
On the left, $0\leq\chi_{\ell,L}\leq\chi_\ell$ and
$M_1(w_{\mathrm d}(\tf))<\infty$ give an integrable majorant:
\[
 \int_0^\infty\!\int_\R K_\ell(r)
       |w_{{\mathrm d},\theta}(\tf,r,z)|\idiff z\idiff r<\infty,
 \qquad
 \int_0^\infty\frac{K_\ell(r)^2}{r\sqrt{1+r^2}}\idiff r<\infty.
\]
The second integral is, up to a constant, the squared weighted dual norm
used in Cauchy--Schwarz.  Dominated convergence therefore proves
\begin{equation}
 4\pi A\int_0^\infty K_{1+\tf}(r)\mu_{\theta,\mathrm d}(\tf,r)\idiff r
 =-\pi A\int_\tp^\tf\int_0^\infty
      K_{\mathrm{hist}}(1+2\tf-u,r)F(u,r)r\idiff r\idiff u.
 \label{hist:mean-history}
\end{equation}
The calculation may first be made for smooth truncated sources.  The
same bounded-kernel estimate removes that regularization.  Since
$0\leq K_{\mathrm{hist}}(\ell,r)\leq2/\ell$ and
$\int|F|r\idiff r\leq\|v_{\mathrm d}\|_2^2/(4\pi)$,
\eqref{hist:isolated-mass} bounds the right side by
$CA\mathcal M_{\mathrm d}/(1+\tf)\leq C\mathcal M_{\mathrm d}/A$.

The comparison with the actual field uses the dual bounds
\begin{align}
 \left|\int e^{-\iunit\nu z}K_\ell(r)h_\theta\idiff z\idiff r\right|
       &\leq C\ell^{-1/4}M_1(h),
 \label{hist:memory-dual-one}\\
 \left|\int(\cos(\nu z)-1)K_\ell(r)h_\theta\idiff z\idiff r\right|
       &\leq C|\nu|^{3/2}\ell^{-1/8}M_{9/4}(h).
 \label{hist:memory-dual-fractional}
\end{align}
For the first, the squared dual norm is bounded by
$C\int_0^\infty K_\ell(r)^2/(r\sqrt{1+r^2})\idiff r\leq C\ell^{-1/2}$.
For the second, use $|\cos y-1|\leq2|y|^{3/2}$ and
\begin{equation}
 \int_\R |z|^3(a+z^2)^{-9/4}\idiff z=\frac{16}{5}a^{-1/4}.
 \label{hist:dual-axial-integral}
\end{equation}
The remaining radial integral is
$C\int_0^\infty K_\ell(r)^2/(r(1+r^2)^{1/4})\idiff r\leq C\ell^{-1/4}$.
Both estimates follow by splitting at $r^2=\ell$.
Apply \eqref{hist:memory-dual-one} to $w_g-w_{\mathrm d}$.
The additional old-column coefficients sum to $O(R^{-2})$ and their ages
are at least $\tf$.  Equations \eqref{hist:baseline-difference} and
\eqref{hist:mean-history} prove \eqref{hist:baseline-memory}.

For the fiber estimate, the isolated Cartesian coupling matrix has norm at
most $2A/(1+t)^2$.  Testing the first-variation equation against
$-\Delta v_{\mathrm d}$ gives
\begin{equation}
 \int_\tp^\tf\|\Delta v_{\mathrm d}\|_2^2\idiff t
       \leq\|\nabla g\|_2^2
              +4A^2\int_\tp^\tf\|v_{\mathrm d}\|_2^2\idiff t.
 \label{hist:isolated-H2-energy}
\end{equation}
The seed bounds and \eqref{hist:isolated-mass} imply
\begin{equation}
 \int_\tp^\tf\|v_{\mathrm d}(t)\|_{L^2_zH^2_{x_h}}^2\idiff t
        \leq CA^2\mathcal M_{\mathrm d}.
 \label{hist:horizontal-H2}
\end{equation}
At axial frequency $\nu$, the projected invariant coupling satisfies
\begin{equation}
 \|B_{\mathrm d}(t,\nu)\|_{2\to2}
       \leq\sqrt2A|\nu|(1+t)^{-3/2},\qquad
 \|\mathcal U_{\mathrm d}(t,u,\nu)\|_{2\to2}
       \leq e^{2\sqrt2A|\nu|}.
 \label{hist:low-frequency-evolution}
\end{equation}
For clarity, the first bound follows from the radial primitive
\[
 \phi_0(r)=-2\int_r^\infty (1+t)^{-2}e^{-a^2/(4(1+t))}v_\theta(a)\idiff a,
\]
whose Hilbert--Schmidt norm is $\sqrt2(1+t)^{-3/2}$, and the identities
\begin{equation}
 \PP_\nu(\nabla_h\phi_0,0)=-\PP_\nu(0,\iunit\nu\phi_0),\qquad
 v_r(r)=-\frac{\iunit\nu}{r}\int_0^r a v_z(a)\idiff a.
 \label{hist:fiber-cancellation}
\end{equation}
The latter gives $|v_r(r)|\leq |\nu|\|v_z\|_{L^2(r\idiff r)}/\sqrt2$.
The meridional and azimuthal outputs are orthogonal.  Integration in time
proves the propagator bound, including the vanishing coupling at $\nu=0$.

The two-dimensional product estimate and axial Fourier Cauchy--Schwarz give
\begin{equation}
 \|\widehat{v_{\mathrm d}\otimes v_{\mathrm d}}^{\,z}(\nu)\|_{H^1_{x_h}}
       \leq C\|v_{\mathrm d}\|_{L^2_zH^2_{x_h}}^2.
 \label{hist:fiber-product}
\end{equation}
The projected divergence costs $C(1+|\nu|)$ in $L^2_{x_h}$.
Duhamel's formula, \eqref{hist:horizontal-H2}, and
\eqref{hist:low-frequency-evolution} yield
\begin{equation}
 \|\widehat w_{\mathrm d}^{\,z}(\tf,\nu)\|_{L^2_{x_h}}
       \leq CA^2(1+|\nu|)e^{2\sqrt2A|\nu|}\mathcal M_{\mathrm d}.
 \label{hist:isolated-fiber}
\end{equation}
The actual difference satisfies
$\|\widehat{w_g-w_{\mathrm d}}^{\,z}(\tf,\nu)\|_{L^2_{x_h}}
\leq\sqrt{A/2}\,N_1(\tf,w_g(\tf)-w_{\mathrm d}(\tf))$.
After division by $\mathcal M_g$ this is $CA^P\xi_Ae^{2\omega A}=o(1)$.
For $|k|\leq K_0A$ and $\nu=k/R$, one has $A|\nu|\leq K_0A^2/R=o(1)$.
This proves \eqref{hist:baseline-fiber}.  The first moment makes the axial
Fourier transform continuous with values in $L^2_{x_h}$, so the estimate
holds at every specified frequency.
\end{proof}

\subsection{A neighborhood of the prepared shape}

The nonlinear iteration produces inputs close to $g$ in $H^4$, with a
separately controlled spatial moment.  These two bounds enter the
comparison in different ways.  In cross terms we place the moment on
the prepared factor; in the term quadratic in the input error we pair
its moment with its small unweighted norm.  This yields the factor
$\delta_A B_{\mathrm{mom}}$ below.  The outgoing field is normalized by the signed
mass of its own input, as required for the shape estimate in
Proposition~\ref{trans:quadratic}.

Consider an input
\begin{equation}
 f=g+\zeta_f,\qquad
 \|\zeta_f\|_{H^4}\leq\delta_A=e^{-\rho_cA},\qquad
 M_{s_\star}(\zeta_f)\leq B_{\mathrm{mom}}\leq e^{\gamma A},
 \label{hist:input-class}
\end{equation}
where $\zeta_f=f-g$ is real, centrally odd, axisymmetric, and solenoidal.
Every input in the later class \eqref{nl:shape-class} satisfies these
hypotheses; orthogonality to $g$ is not needed for the estimates here.
Let $v_f,w_f$ be the actual variations with initial values $f,0$.
Their complete stress and positive stress mass are
\begin{equation}
 \mathsf S_f=v_f\otimes v_f+\bar u\otimes w_f+w_f\otimes\bar u,
 \qquad
 \int\mathsf S_f\idiff x=\operatorname{diag}(s_{h,f},s_{h,f},s_{z,f}),\qquad
 \mathcal M_f=-\int_\tp^\tf(s_{h,f}-s_{z,f})\idiff t.
 \label{hist:generic-stress}
\end{equation}

\begin{proposition}[The complete input at the larger scale]
\label{hist:generic}
For every input in \eqref{hist:input-class}, the mass and second variation
satisfy the following bounds uniformly in both cutoffs:
\begin{align}
 |\mathcal M_f-\mathcal M_g|&\leq CA^P\delta_A(1+B_{\mathrm{mom}})e^{2c_\star A},&
 A^{-P}e^{2GA}&\leq \mathcal M_f\leq CA^Pe^{2c_\star A},
 \label{hist:generic-mass}\\
 M_{s_\star}(w_f(t))+\|w_f(t)\|_{H^4}&\leq CA^Pe^{2\Gamma_{\mathrm{src}}(t)},
 \label{hist:generic-second}\\
 \int_\tp^\tf\int
  \bigl(1+|x|^\alpha+(\tf-t)^{\alpha/2}\bigr)|\mathsf S_f|
       \idiff x\idiff t&\leq CA^Pe^{2c_\star A}.
 \label{hist:generic-moment}
\end{align}
Put $m_f=-\mathcal M_f$ and define
\begin{equation}
 F_{f,R}(y)=\frac{R^4}{m_f}w_f(\tf,Ry),\qquad
 F_{f,k}(y_h)=\frac{R^3}{m_f}
          \widehat w_f^{\,z}(\tf,Ry_h,k/R),\qquad
 \mathcal W_f=\frac{M_{9/4}(w_f(\tf))}{\mathcal M_f}.
 \label{hist:normalized-input}
\end{equation}
Then
\begin{equation}
 \mathcal W_f\leq CA^Pe^{2\omega A},\qquad
 \|F_{f,R}\|_2\leq R^{5/2}\mathcal W_f,\qquad
 \|y_3F_{f,R}\|_2\leq R^{3/2}\mathcal W_f,
 \label{hist:whole-input}
\end{equation}
and, for every fixed $K_0\geq1$,
\begin{equation}
 \sup_{|k|\leq K_0A}\|F_{f,k}\|_{L^2_{y_h}}\leq C_{K_0}A^2R^2.
 \label{hist:generic-fiber}
\end{equation}
Writing
$\mathcal I_k=(2\pi)^{-1/2}\iunit k\PP_k(e_z\delta_{0,h})$,
one has at every real $k$
\begin{equation}
 \|F_{f,k}-\mathcal I_k\|_{H^{-17/4}_{y_h}}
      \leq CA^P(1+|k|)^{11/4}R^{-7/4}e^{2\omega A}.
 \label{hist:point}
\end{equation}
Finally,
\begin{equation}
 c_{\mathrm{mean},f}(k):=
       \frac{\mathcal T_{I,k/R}(w_f(\tf))}{m_f},\qquad
 c_{\mathrm{mean},f}(0)\in\R,\qquad
 |c_{\mathrm{mean},f}(0)|\leq C/A.
 \label{hist:memory}
\end{equation}
For $|k|\leq K_0A$, the same coefficient satisfies
\begin{equation}\label{hist:memory-frequency}
 |c_{\mathrm{mean},f}(k)-c_{\mathrm{mean},f}(0)|
       \leq CA^PR^{-3/2}e^{2\omega A}.
\end{equation}
\end{proposition}

\begin{proof}
Set $d_v=v_f-v_g$ and $d_w=w_f-w_g$.
Apply the homogeneous estimate first without a weight and then with moment
$s_\star$.  The resulting bounds are
\begin{equation}
 \|d_v\|_{H^4}\leq CA^P\delta_Ae^{\Gamma_{\mathrm{src}}},
 \qquad
 M_{s_\star}(d_v)\leq CA^P(1+B_{\mathrm{mom}})e^{\Gamma_{\mathrm{src}}}.
 \label{hist:first-tame}
\end{equation}
The tensor source of $d_w$ is exactly
$F_\Delta=v_g\otimes d_v+d_v\otimes v_g+d_v\otimes d_v$.
Place the weight on $v_g$ in both cross terms and use the small unweighted
supremum norm in the last term.  Then
\begin{equation}
 \|F_\Delta\|_{H^4}\leq CA^P\delta_Ae^{2\Gamma_{\mathrm{src}}},\qquad
 M_{s_\star}(F_\Delta)\leq CA^P\delta_A(1+B_{\mathrm{mom}})e^{2\Gamma_{\mathrm{src}}}.
 \label{hist:tame-source}
\end{equation}
Proposition~\ref{hist:energy} yields
\begin{equation}
 \|d_w\|_{H^4}\leq CA^P\delta_Ae^{2\Gamma_{\mathrm{src}}},\qquad
 M_{s_\star}(d_w)\leq CA^P\delta_A(1+B_{\mathrm{mom}})e^{2\Gamma_{\mathrm{src}}}.
 \label{hist:second-tame}
\end{equation}
Since $\gamma<\rho_c$, this proves \eqref{hist:generic-second}.

The negative-weight bounds \eqref{hist:inner-negative} and
\eqref{hist:outer-negative}, together with
\eqref{hist:reference-error}, imply
$\sup_{[\tp,\tf]}\|\langle x\rangle^{-s_\star}\bar u\|_2\leq CA^P$.
The exact stress difference is
$\mathsf S_f-\mathsf S_g=F_\Delta+\bar u\otimes d_w+d_w\otimes\bar u$.
Cauchy--Schwarz and \eqref{hist:second-tame} give its time-integrated
$L^1$ bound in \eqref{hist:generic-mass}.
Its relative size is at most
\begin{equation}
 CA^P\delta_A(1+B_{\mathrm{mom}})e^{2\omega A}
       \leq CA^Pe^{-(\rho_c-\gamma-2\omega)A}=o(1),
 \qquad \rho_c-\gamma-2\omega>3\varrho/32.
 \label{hist:relative-tame}
\end{equation}
Equation \eqref{hist:baseline-mass} proves the lower and upper mass bounds.

The fractional stress moment can be small even when the full input
moment tolerance is large.  The elementary weighted interpolation used
here is worth recording.  If $0<\alpha<2s_\star$ and
$\theta=\alpha/(2s_\star)$, then H\"older applied to
\[
 \langle x\rangle^\alpha|d_v|^2
   =\bigl(\langle x\rangle^{2s_\star}|d_v|^2\bigr)^\theta
                         (|d_v|^2)^{1-\theta}
\]
gives
$M_{\alpha/2}(d_v)^2\leq
 M_{s_\star}(d_v)^{\alpha/s_\star}\|d_v\|_2^{2-\alpha/s_\star}$.
Thus a strong smallness factor is retained rather than discarded by
using only the large moment.  With the actual bounds and
$\alpha/s_\star=14/19$, its exponential contribution is
\[
 \delta_A^{24/19}(1+B_{\mathrm{mom}})^{14/19}
 \leq C\exp\bigg(
       \Big(-\frac{24}{19}\frac{3\varrho}{8}
                   +\frac{14}{19}\frac{\varrho}{4}\Big)A\bigg)
   =C e^{-11\varrho A/38}.
\]
This is why the fractional stress estimate does not require an extra
positive exponential in $A$.  It uses interpolation of two bounds for
an already controlled difference, unlike the operator interpolation for
arbitrary lower-moment initial data in Section~\ref{sec:decay}.

For the fractional first-variation moment, interpolate the separate bounds
in \eqref{hist:first-tame}:
\begin{equation}
 M_{\alpha/2}(d_v)^2
   \leq\|d_v\|_2^{2-\alpha/s_\star}M_{s_\star}(d_v)^{\alpha/s_\star}
   \leq CA^Pe^{2\Gamma_{\mathrm{src}}}
        \delta_A^{2-\alpha/s_\star}(1+B_{\mathrm{mom}})^{\alpha/s_\star}.
 \label{hist:fractional-tame}
\end{equation}
Here $\alpha/s_\star=14/19$ and
$-(2-\alpha/s_\star)\rho_c+(\alpha/s_\star)\gamma=-11\varrho/38<0$.
The background terms use \eqref{hist:generic-second} and the negative
weight $s_\star-\alpha=5/8$, exactly as in
\eqref{hist:fractional-products}.  This proves
\eqref{hist:generic-moment} without an extra exponential loss.

The second-variation and mass bounds give the first estimate in
\eqref{hist:whole-input}.  The other two follow from the exact scalings
\begin{equation}
 \|F_{f,R}\|_2=R^{5/2}\frac{\|w_f(\tf)\|_2}{\mathcal M_f},\qquad
 \|y_3F_{f,R}\|_2=R^{3/2}\frac{\|zw_f(\tf)\|_2}{\mathcal M_f}.
 \label{hist:input-scaling}
\end{equation}
Axial Cauchy--Schwarz gives
$\sup_\nu\|\widehat d_w^{\,z}(\tf,\nu)\|_{L^2_{x_h}}
\leq CM_1(d_w(\tf))$.
Its normalized contribution is bounded by
$CA^PR^2e^{-(\rho_c-\gamma-2\omega)A}$.
The baseline fiber bound \eqref{hist:baseline-fiber} and
$\mathcal M_g/\mathcal M_f=1+o(1)$ now prove \eqref{hist:generic-fiber}.

The functional $\mathcal T_{I,0}$ has polynomial norm on $M_{s_\star}L^2$.
Indeed, after division by $r$ to use the measure $r\idiff r\idiff z$,
its kernel is $O(Ar/\tf)$ at the axis and $O(A/r)$ at infinity.
The squared weighted dual tail is $O(r^{-2s_\star})$ after integration
in $z$, and $\sum_{j\leq0}R^{2j}\leq4/3$.
Apply this bound to \eqref{hist:second-tame}, and use
\eqref{hist:baseline-memory} and \eqref{hist:relative-tame}.
This proves \eqref{hist:memory}.
Central oddness and axisymmetry make $w_{f,\theta}(r,z)$ even in $z$.
Equation \eqref{hist:memory-dual-fractional}, summed over the old radii,
proves \eqref{hist:memory-frequency}.

For this calculation, write
$\widehat h(\xi)=\int e^{-\iunit x\cdot\xi}h(x)\idiff x$
for the unnormalized full spatial Fourier transform.
The negative-Sobolev estimate is a quantitative moment expansion. To see
the leading term before its error, denote the accumulated zero-frequency
stress by $cI_3+m_f\operatorname{diag}(1,1,0)$, where $c\in\R$ is its
isotropic coefficient and $m_f=-\mathcal M_f$.
For $\xi=(\xi_h,\xi_z)$, the isotropic part vanishes after projection,
and
\[
 -i\PP(\xi)\operatorname{diag}(1,1,0)\xi
   =-i\PP(\xi)(\xi-\xi_ze_z)
   =i\xi_z\PP(\xi)e_z.
\]
The sign is therefore $m_f\mathcal I$, not $-m_f\mathcal I$.
This calculation uses the complete stress; substituting only
$v_f\otimes v_f$ would change the coefficient.

For an even integrable tensor $S$, the sine part cancels and
\[
 \left|\int(e^{-ix\cdot\xi}-1)S(x)\,\idiff x\right|
 \leq C_\alpha|\xi|^\alpha\int |x|^\alpha|S(x)|\,\idiff x,
 \qquad 0<\alpha\leq2.
\]
The heat-age error has the same order because
$|1-e^{-h|\xi|^2}|\leq h^{\alpha/2}|\xi|^\alpha$.
A divergence adds one power of frequency.  Under the normalization
$F_R(y)=R^4w_f(\tf,Ry)/m_f$, its full Fourier transform is
$R\widehat w_f(\tf,\xi/R)/m_f$; the leading degree-one symbol is
unchanged, whereas the remainder gains $R^{-\alpha}$.
The following calculation places precisely this remainder in a negative
horizontal Sobolev space.  It is not a small-$L^2$ assertion about $F_R$.

It remains to establish \eqref{hist:point}.
The exact second variation has the free-heat representation
\begin{equation}
 w_f(\tf)=-\int_\tp^\tf
       e^{(\tf-u)\Delta}\PP\operatorname{div}\mathsf S_f(u)\idiff u.
 \label{hist:free-heat-history}
\end{equation}

The expansion is performed on the integrable stress before applying
the Leray multiplier.  The ordinary velocity moment available here does
not imply a finite first absolute $L^1$ moment of the velocity itself.
Keeping the projection outside the stress expansion accommodates its
algebraic velocity tail and preserves the correct point distribution.

The complete stress is symmetric and
componentwise even.  Therefore
\begin{equation}
 |\widehat{\mathsf S}_f(u,\xi)-\widehat{\mathsf S}_f(u,0)|
       \leq C_\alpha|\xi|^\alpha\int|x|^\alpha|\mathsf S_f(u,x)|\idiff x,
 \qquad
 |e^{-h|\xi|^2}-1|\leq h^{\alpha/2}|\xi|^\alpha.
 \label{hist:fourier-remainder}
\end{equation}
The accumulated zero-frequency stress is a scalar matrix plus
$m_f\operatorname{diag}(1,1,0)$.
The scalar matrix is removed by $\PP\operatorname{div}$; the remaining
part gives $m_f\mathcal I$ with the sign in
\eqref{hist:free-heat-history}.
Using \eqref{hist:generic-moment} and rescaling consequently gives
\begin{equation}
 |\widehat F_{f,R}(\xi)-\widehat{\mathcal I}(\xi)|
       \leq CA^PR^{-\alpha}e^{2\omega A}|\xi|^{1+\alpha}.
 \label{hist:full-fourier-error}
\end{equation}
If $\widehat h^{\,h}$ denotes the analogous unnormalized horizontal
transform, its relation to the unitary axial transform is
$\widehat{F_{f,k}}^{\,h}(\xi_h)
=(2\pi)^{-1/2}\widehat F_{f,R}(\xi_h,k)$.
With $\delta_{\rm Sob}=1/4$, the integral
\begin{equation}
 \int_{\R^2}\langle\xi_h\rangle^{-8-2\delta_{\rm Sob}}
       (1+|\xi_h|)^{2+2\alpha}\idiff\xi_h<\infty
 \label{hist:negative-sobolev-integral}
\end{equation}
is finite because $\alpha<2+\delta_{\rm Sob}$.
This proves \eqref{hist:point} at every real $k$; the first moment again
provides a continuous $L^2_{x_h}$-valued representative.
\end{proof}

For later use, the whole-input estimate also controls the full axial heat
tail:
\begin{equation}
 \|\boldsymbol1_{\{|D_z|>K_0A\}}
           e^{A^{-1}\partial_z^2}F_{f,R}\|_2
 \leq CA^P\exp\left(\left(\frac52\varrho+2\omega-K_0^2\right)A\right).
 \label{hist:axial-tail}
\end{equation}
The strict inequalities
\begin{equation}
 (\alpha-\delta_{\rm Sob})\varrho>2\omega,\qquad
 \frac32\varrho>2\omega,\qquad
 K_0^2>\frac52\varrho+2\omega+1
 \label{hist:frequency-margins}
\end{equation}
provide the margins needed in Section~\ref{sec:transfer}.
The age inequalities $A\leq\tf$, $\tf+1\leq R^2/(2A)$, and $R\geq A^2$
hold once $A$ is sufficiently large.

\section{Endpoint control of the first variation without moment loss}
\label{sec:decay}

The signed quadratic history from Section~\ref{sec:history} must be separated
from the first variation before it reaches the next larger scale.  Axial
diffusion provides an improved endpoint bound after a long interval; the endpoint estimate follows the transient growth phase and uses the same ordinary spatial moment at input and output.  Requiring a stronger moment at each
transition would prevent iteration.  The low-frequency argument below
is designed to avoid that loss.

Throughout this section, the retained radii are
\(\{R^j:-I\leq j\leq J\}\), where \(J\geq1\) and
\(I\in\N\cup\{0,\infty\}\). We use
\(R=e^{\varrho A}\), \(b=e^{\beta A}\), and
\(\tf=\tp+16A^2\). The constants may depend on a specified moment
exponent and on \(\varrho,\Theta\), but are uniform in the two cutoffs.
Every lower bound imposed on \(\beta\) is fixed before \(A\) is increased.

\begin{lemma}[The final-time moment estimate]\label{dec:retirement}
Let \(0\leq s<5/2\) and \(\lambda>0\). Choose a sufficiently large
fixed width exponent \(\beta\), whose threshold may depend on these
prescribed \(s,\lambda\) and the previously fixed parameters, and then
take \(A\) sufficiently large. The homogeneous linearized solution about
\(\bar u^{I,J}\), with real, centrally odd, axisymmetric, solenoidal
initial value \(v(\tp)=f\), satisfies
\begin{equation}\label{dec:moment-decay}
 M_s\bigl(v(\tf)\bigr)
 \leq C A^P e^{\lambda A}M_s(f).
\end{equation}
The hypothesis is \(M_s(f)<\infty\); no weighted derivative or stronger
moment of \(f\) is required.
\end{lemma}

\subsection{A fixed space for the axial Fourier variable}

Differentiating a solenoidal Fourier fiber in $k$ does not preserve its
solenoidal constraint.  We therefore extend the column operator to one
fixed physical Hilbert space before differentiating.  The second
frequency derivative has a mild singularity at zero; fractional
regularity strictly below $1/2$ is enough for every moment used here.

Let $F(k)=\widehat f^{\,z}(k)$ with values in $L^2_{x_h}$.
For an integer $j$, Plancherel gives
$\|\partial_k^jF\|_{L^2_kL^2_{x_h}}=\|z^jf\|_2$.
For $0<\vartheta<1$, the difference-quotient identity is
\[
 \int_{\R^2}\frac{\|F(k+h)-F(k)\|_{L^2_{x_h}}^2}
                         {|h|^{1+2\vartheta}}\,\idiff k\idiff h
 =c_\vartheta\int_{\R^3}|z|^{2\vartheta}|f(x_h,z)|^2\,\idiff x.
\]
Indeed, a translation in $k$ becomes multiplication by $e^{-ihz}$;
integration of $|e^{-ihz}-1|^2|h|^{-1-2\vartheta}$ and the substitution
$h|z|=u$ give the formula.  Applying it to $F''$ produces the moment
$|z|^{2+\vartheta}f$.  Thus the hierarchy at orders $0,1,2,2+\vartheta$
uses exactly the given positional moment, not weighted derivatives of $f$.

There is a constraint issue when those derivatives are taken.  If
$\operatorname{div}_hF_h+ikF_z=0$, differentiation gives
\[
 \operatorname{div}_h(\partial_kF_h)+ik\partial_kF_z=-iF_z.
\]
The differentiated field is not solenoidal at frequency $k$.
Consequently a bound proved only on the solenoidal fiber cannot simply be
applied to it.  The fixed-space extension below is introduced precisely
to make those differentiated equations legitimate.

Frequency differentiation is simplest before imposing the solenoidal
constraint on each fiber. Let \(\mathcal H_h\) be the Hilbert space of
axisymmetric Cartesian vector fields on \(\R^2\), with all three physical
components and the usual \(L^2(\R^2)\) norm.  On radial profiles this norm is
$\|F\|_{2,h}^2=2\pi\int_0^\infty|F(r)|^2r\,\idiff r$; all three
velocity components are included.  The projection \(\PP_k\) denotes
the three-dimensional Leray multiplier with axial frequency \(k\).
For a radius \(\ell\), write \(T_\ell=\ell^2+t\) and define
\begin{align}
 \Phi_{\ell,t}g(r)
 &=-2A\ell^2T_\ell^{-2}
       \int_r^\infty e^{-u^2/(4T_\ell)}g(u)\idiff u,
       \label{dec:primitive}\\
 K_{\ell,t}h(r)
 &=A\ell^2T_\ell^{-2}
       \left(2-\frac{r^2}{2T_\ell}\right)e^{-r^2/(4T_\ell)}
       \frac1r\int_0^r uh(u)\idiff u,\notag\\
 \mathcal B_\ell^{\mathrm{ext}}(t,k)F
 &=-ik\PP_k\bigl(e_z\Phi_{\ell,t}F_\theta\bigr)
       +ik e_\theta K_{\ell,t}F_z.\notag
\end{align}
These formulas define bounded operators on \(\mathcal H_h\).  For the
physical solenoidal fiber \(\mathcal H_k\) from Section~\ref{sec:background},
\begin{equation}\label{dec:extension-restriction}
 \mathcal B_\ell^{\mathrm{ext}}(t,k)|_{\mathcal H_k}
       =\mathcal B_\ell(t,k),
\end{equation}
where the right side is \eqref{bg:column-coupling}.  Frequency derivatives
below are taken on the fixed space \(\mathcal H_h\), not on the
varying subspaces \(\mathcal H_k\).

For an operator-valued function \(D(k)\) and \(0<\vartheta<1\), use the
norm
\begin{equation}\label{dec:operator-fractional-norm}
 \|D\|_{\mathfrak H^\vartheta}^2
 =\int_\R\|D(k)\|_{\mathrm{op}}^2\idiff k
  +\int_{\R^2}
   \frac{\|D(k+h)-D(k)\|_{\mathrm{op}}^2}{|h|^{1+2\vartheta}}
       \idiff k\idiff h.
\end{equation}
Here $\|\cdot\|_{\mathrm{op}}$ is the operator norm on
$\mathcal L(\mathcal H_h)$.  Thus \eqref{dec:operator-fractional-norm}
is the Banach-valued Slobodeckij norm for functions with values in this
operator space; it does not use a Hilbert--Schmidt norm.  In particular,
$L^2_k(\mathcal L(\mathcal H_h))$ denotes square integrability of
operator norms, not the operator space $\mathcal L(L^2_k(\mathcal H_h))$.

\begin{lemma}[Frequency derivatives of the extended column operator]\label{dec:fixed-space}
The operators in \eqref{dec:primitive} are twice differentiable in
operator norm for $k\ne0$, and their first derivatives extend
continuously to $k=0$.  Across $k=0$, the second derivative is understood
weakly, with the representatives in the spaces stated below; its value
at that single frequency is immaterial.  They satisfy
\begin{equation}\label{dec:first-frequency-bounds}
 \|\mathcal B_\ell^{\mathrm{ext}}(t,k)\|_{\mathrm{op}}
       \leq |k|b_\ell(t),\qquad
 \|\partial_k\mathcal B_\ell^{\mathrm{ext}}(t,k)\|_{\mathrm{op}}
       \leq b_\ell(t),\qquad
 b_\ell(t)=CA\ell^2T_\ell^{-3/2}.
\end{equation}
For \(0<\vartheta<1/2\), one also has
\begin{equation}\label{dec:second-frequency-bounds}
 \|\partial_k^2\mathcal B_\ell^{\mathrm{ext}}(t,\cdot)\|_{\mathfrak H^\vartheta}
 \leq C_\vartheta A\ell^2
       \left(T_\ell^{-5/4}+T_\ell^{-5/4+\vartheta/2}\right).
\end{equation}
The meridional-output block of \(\partial_k\mathcal B_\ell^{\mathrm{ext}}\) satisfies the
corresponding bound with powers \(-7/4\) and \(-7/4+\vartheta/2\).
Its azimuthal-output part is constant in \(k\). The sums over all
retained radii of the bounds in
\eqref{dec:first-frequency-bounds}--\eqref{dec:second-frequency-bounds}
have time integrals over \([\tp,\tf]\) bounded by \(C_\vartheta A\).
\end{lemma}

\begin{proof}
Hilbert--Schmidt integration with the radial measure \(r\idiff r\)
gives
\begin{equation}\label{dec:primitive-size}
 \|\Phi_{\ell,t}\|_{\mathrm{op}}+
 \|K_{\ell,t}\|_{\mathrm{op}}
 \leq CA\ell^2T_\ell^{-3/2}.
\end{equation}
For a solenoidal physical fiber, incompressibility and regularity at the
axis give
\begin{equation}\label{dec:physical-extension}
 F_r(r)=-\frac{ik}{r}\int_0^r uF_z(u)\idiff u,
 \qquad
 \PP_k(\nabla_h\Phi,0)=-ik\PP_k(e_z\Phi).
\end{equation}
Thus \eqref{dec:primitive} is the asserted extension of the physical
operator. It takes every input to a solenoidal output.

We give the frequency estimates at unit horizontal scale. Let
\(\Phi_\ast g(r)=-2\int_r^\infty e^{-u^2/4}g(u)\idiff u\) and
\(\mathcal R_\ast(k)g=-ik\PP_k(e_z\Phi_\ast g)\). Direct integration
gives
\begin{equation}\label{dec:primitive-unit}
 \|\Phi_\ast g\|_2+\|\Phi_\ast g\|_1+
 \|\nabla_h\Phi_\ast g\|_2\leq C\|g\|_2.
\end{equation}
The norms in \eqref{dec:primitive-unit} are over $\R^2$.
For example, the \(L^1\) bound follows from
\(\int_0^\infty u^2e^{-u^2/4}|g(u)|\idiff u
\leq C\|g\|_{L^2(r\idiff r)}\).

The projection extends holomorphically to
\(\{k:|\operatorname{Im}k|<|\operatorname{Re}k|/8,\ k\neq0\}\),
because its denominator has modulus bounded below by
\(c_0(|\xi_h|^2+|k|^2)\) there. Subtracting the linear term in
\(\mathcal R_\ast\) and using \eqref{dec:primitive-unit}, we obtain
\begin{align}
 \|\bigl(\mathcal R_\ast(k)+ik e_z\Phi_\ast\bigr)g\|_2^2
 &\leq C|k|^4\int_{\R^2}
    \frac{|\widehat{\Phi_\ast g}(\xi_h)|^2}
         {|\xi_h|^2+|k|^2}\idiff\xi_h\notag\\
 &\leq C|k|^4\bigl(1+|\log|k||\bigr)\|g\|_2^2,
       \qquad 0<|k|\leq\tfrac12.
       \label{dec:frequency-origin}
\end{align}

The logarithm in \eqref{dec:frequency-origin} has a concrete source.
The $L^1$ estimate for $\Phi_*g$ bounds its horizontal Fourier transform,
and on the unit disk
\[
 \int_{|\xi_h|\leq1}\frac{\idiff\xi_h}{|\xi_h|^2+k^2}
     =\pi\log\frac{1+k^2}{k^2}.
\]
Outside that disk the $L^2$ bound suffices.  The subtraction of the
linear term leaves an order-$k^2$ remainder with this logarithmic factor
in its norm.  Its second derivative can therefore be locally unbounded
without preventing an integrated estimate.

The translation bound proved below is of order
$|h|(1+|\log|h||)$.  Dividing by $|h|^{1+2\vartheta}$ leaves the
integrable expression $|h|^{-2\vartheta}(1+|\log|h||)$ exactly when
$\vartheta<1/2$.  This explains the upper moment threshold
$2+\vartheta<5/2$ in this part of the proof.  It comes from frequency
regularity of the projected operator, rather than from a pointwise bound
on its second derivative at $k=0$.

Cauchy's formula on disks of radius proportional to \(|k|\) now yields
\begin{equation}\label{dec:frequency-derivatives}
 \|\mathcal R_\ast''(k)\|_{\mathrm{op}}
       \leq C\bigl(1+|\log|k||\bigr)^{1/2},\qquad
 \|\mathcal R_\ast'''(k)\|_{\mathrm{op}}
       \leq C|k|^{-1}\bigl(1+|\log|k||\bigr)^{1/2}.
\end{equation}
Moreover \(\mathcal R_\ast'(k)\to-ie_z\Phi_\ast\) as \(k\to0\),
so its second weak derivative has no atom at zero. For \(|k|\geq1/2\),
the identity
\(\mathcal R_\ast(k)=\PP_k(\nabla_h\Phi_\ast,0)\) and Cauchy's
formula give
\(\|\mathcal R_\ast^{(j)}(k)\|_{\mathrm{op}}\leq C_j|k|^{-j}\)
for \(j=1,2,3\).

Put \(D=\mathcal R_\ast''\). Splitting a translation integral into
\(|k|\leq2|h|\), \(2|h|<|k|<2\), and their complement, we find
\begin{equation}\label{dec:translation-modulus}
 \int_\R\|D(k+h)-D(k)\|_{\mathrm{op}}^2\idiff k
 \leq C|h|\bigl(1+|\log|h||\bigr),
 \qquad 0<|h|<\tfrac14.
\end{equation}
The middle region uses
\(\int_{2|h|}^2 |h|^2u^{-2}(1+|\log u|)\idiff u\); the outer
region uses the large-frequency derivative estimate. Since
\(\int_0^{1/4}h^{-2\vartheta}(1+|\log h|)\idiff h<\infty\)
for \(\vartheta<1/2\), this proves the required fractional regularity.

Let \(\mathcal D_{\sqrt{T_\ell}}\) be horizontal unitary dilation.
The meridional-output block at general scale is
\begin{equation}\label{dec:operator-scaling}
 A\ell^2T_\ell^{-2}\,
 \mathcal D_{\sqrt{T_\ell}}
 \mathcal R_\ast(k\sqrt{T_\ell})
 \mathcal D_{\sqrt{T_\ell}}^{-1}.
\end{equation}
Two frequency derivatives contribute \(T_\ell\), the axial \(L^2\)
change of variables contributes \(T_\ell^{-1/4}\), and the fractional
seminorm contributes \(T_\ell^{\vartheta/2}\). These factors prove
\eqref{dec:second-frequency-bounds}; one derivative gives the asserted
bound for the meridional-output block. The azimuthal block is linear in \(k\).

For completeness, the slowest outer-radius sum after time integration is
bounded by
\begin{equation}\label{dec:operator-tail-sum}
 CA\tf\sum_{\ell\geq R}
       \left(\ell^{-1/2}+\ell^{-1/2+\vartheta}\right).
\end{equation}
Both geometric series converge. The inner radii contribute at most
\(CA R^{-2}\) times a fixed power of \(\tp^{-1}\). For \(\ell=1\),
the time integrals of the powers in
\eqref{dec:second-frequency-bounds} converge even at infinity.
The same decomposition proves
\(\int_\tp^\tf\sum_\ell b_\ell(t)\idiff t\leq CA\).
For sufficiently large \(A\), these estimates prove all the uniform
time-integral assertions.
\end{proof}

\subsection{The weighted decay estimate}

We first prove the endpoint estimate for $2<s<5/2$.
The axial Fourier argument controls the transform, its first two
frequency derivatives, and a fractional difference seminorm at the
highest order.  Plancherel converts this hierarchy into an axial
moment of order $s$.  Gaussian localization then recovers the
horizontal moment.  The estimates are arranged so that a prescribed
loss $e^{\lambda A}$ absorbs all exponential growth that remains at
$\tf$.

\paragraph{The axial heat filter.}
For an axially invariant background, the scalar axial heat factor commutes
with the horizontal fiber evolution.  A rough coupling bound with time
integral $CA|k|$ gives the schematic multiplier
$e^{-Lk^2+CA|k|}$.  If $L$ is of order $A^2$, then
\[
 -Lk^2+CA|k|
   =-L\left(|k|-\frac{CA}{2L}\right)^2+\frac{C^2A^2}{4L},
\]
whose maximum is bounded independently of $A$.  This explains why a
long diffusive interval can remove a large transient amplification.
For the localized array, frequencies are coupled by the axial cutoff,
so this factorization is no longer exact.  We replace it by a forward-heat
filter and explicitly estimate its commutator with the coefficients.

The filter is
$W(t)=e^{(\tf-t)\partial_z^2}v_Q(t)$.
Its elapsed heat time is nonnegative throughout $[\tp,\tf]$.  It has endpoints
\[
 W(\tp)=e^{16A^2\partial_z^2}f,
 \qquad W(\tf)=v_Q(\tf).
\]
Differentiating removes the axial Laplacian from the equation for $W$;
the price is the commutator displayed in \eqref{dec:filtered-equation}.
The high-frequency part is already exponentially small at the first
endpoint, while the low-frequency part evolves with a small rate.
The preliminary, larger-exponent estimate is used only for errors
multiplied by the small axial-width factors.

\begin{proof}[Proof of Lemma~\ref{dec:retirement}]
We first assume \(2<s<5/2\), and set \(\vartheta=s-2\). Let
\(v_Q\) be the linearized evolution about the actual localized heat
array \(Q^{I,J}\), with initial value \(f\). On axisymmetric fields
its equation has the Cartesian form
\(\partial_t v_Q-\Delta v_Q=\PP M_Qv_Q\), where
\(M_Q\) is the smooth even matrix obtained from the pure-swirl
linearization. Gaussian differentiation, including the axial shear,
gives
\begin{align}
 m_0(t):=\|M_Q(t)\|_\infty
 &\leq CA\sum_\ell\frac{\ell^2}{(\ell^2+t)^2},\notag\\
 d_1(t):=\|\partial_zM_Q(t)\|_\infty
 &\leq \frac{CA}{b}\sum_\ell\frac{\ell}{(\ell^2+t)^2},\qquad
 \int_\tp^\tf m_0(t)\idiff t\leq CA,\quad
 \int_\tp^\tf d_1(t)\idiff t\leq\frac{CA}{b}.
 \label{dec:matrix-integrals}
\end{align}
These bounds include radii with \(b\ell<1\). Indeed
\(\|\eta_\ell'\|_\infty\leq C/(b\ell)\) and
\(\|\eta_\ell''\|_\infty
\leq C(b\ell)^{-1}(b^2\ell^2+2t)^{-1/2}\).

We shall use a preliminary estimate with a larger exponential constant.
For a starting time \(u\), introduce the moving weight
\(w(t,x)=(1+|x|^2/(1+t-u))^{s/2}\). Its time derivative is
nonpositive, and
\(\Delta_\ast w^2\leq C_s(1+t-u)^{-1}w^2\) for either
\(\Delta_\ast=\Delta_h\) or \(\Delta_\ast=\Delta\).
The odd-field projection estimate in Lemma~\ref{bg:weighted} has the
same constant after isotropic dilation. Weighted heat energy therefore
shows that a centrally odd solution of
\(\partial_tY-\Delta_\ast Y=\PP M_QY+F\) satisfies
\begin{equation}\label{dec:crude-moment}
 M_s(Y(t))\leq C_s(1+\tf)^{N_s}
       \exp\left(C_s\int_u^t m_0(\tau)\idiff\tau\right)
       \left(M_s(Y(u))+\int_u^tM_s(F(\tau))\idiff\tau\right).
\end{equation}
In particular,
\begin{equation}\label{dec:crude-vq}
 \sup_{\tp\leq t\leq\tf}M_s(v_Q(t))
       \leq CA^Pe^{K_sA}M_s(f),
\end{equation}
where \(K_s\) is independent of \(\beta\).

To expose the axial diffusion already accumulated at \(\tf\), define
\(W(t)=e^{(\tf-t)\partial_z^2}v_Q(t)\). Its exact equation is
\begin{equation}\label{dec:filtered-equation}
 \partial_tW-\Delta_hW=\PP M_QW+\mathcal R,
 \qquad
 \mathcal R=\PP\bigl[e^{(\tf-t)\partial_z^2},M_Q\bigr]v_Q.
\end{equation}
For a Lipschitz matrix \(M\), the Gaussian kernel gives
\begin{equation}\label{dec:heat-commutator}
 M_s\bigl([e^{h\partial_z^2},M]F\bigr)
 \leq C_s\|\partial_zM\|_\infty
       \left(h^{1/2}+h^{(s+1)/2}\right)M_s(F).
\end{equation}
This follows from
\(|M(x_h,z-y)-M(x_h,z)|\leq|y|\|\partial_zM\|_\infty\)
and the weighted translation bound. The commutator in
\eqref{dec:filtered-equation} is centrally odd, so
\eqref{dec:matrix-integrals} and \eqref{dec:crude-vq} imply
\begin{equation}\label{dec:filtered-source}
 \int_\tp^\tf M_s(\mathcal R(t))\idiff t
       \leq CA^Pb^{-1}e^{K_sA}M_s(f).
\end{equation}

The one-dimensional heat kernel also gives
\[
 M_s(W(t))\leq C_s\bigl(1+(\tf-t)^{s/2}\bigr)M_s(v_Q(t))
       \leq C_sA^Pe^{K_sA}M_s(f).
\]
Its polynomial dependence on the elapsed time is sufficient for the
cutoff commutator estimate below.

Choose \(0<\kappa_{\mathrm{cut}}<1/8\), to be fixed later, and a smooth real even
function \(\chi_{\mathrm{cut}}\), equal to one on \((-1,1)\) and supported in
\((-2,2)\). Put
\(W_< =\chi_{\mathrm{cut}}(D_z/\kappa_{\mathrm{cut}})W\) and \(W_>=W-W_<\).
The Fourier cutoff commutator has Schwartz kernel and satisfies
\begin{equation}\label{dec:band-commutator}
 M_s\bigl(\PP[\chi_{\mathrm{cut}}(D_z/\kappa_{\mathrm{cut}}),M_Q]W\bigr)
       \leq C_{s,\kappa_{\mathrm{cut}}}d_1(t)M_s(W).
\end{equation}
Consequently the two fields satisfy the horizontal-heat equation in
\eqref{dec:filtered-equation}, with sources
\begin{equation}\label{dec:band-sources}
 \mathcal S_< =\chi_{\mathrm{cut}}(D_z/\kappa_{\mathrm{cut}})\mathcal R
       +\PP[\chi_{\mathrm{cut}}(D_z/\kappa_{\mathrm{cut}}),M_Q]W,
 \qquad \mathcal S_>=\mathcal R-\mathcal S_<,
\end{equation}
for which
\begin{equation}\label{dec:band-source-bound}
 \int_\tp^\tf\left(M_s(\mathcal S_<)+M_s(\mathcal S_>)\right)\idiff t
       \leq C_{s,\kappa_{\mathrm{cut}}}A^Pb^{-1}e^{K_sA}M_s(f).
\end{equation}
The sources need not be frequency supported; the field \(W_<\) is.

For the low frequencies we use the extension in
\eqref{dec:primitive}. If
\(U_\ell=A\ell^2T_\ell^{-2}e^{-r^2/(4T_\ell)}x_h^\perp\),
then the actual localized column operator has the exact extension
\begin{equation}\label{dec:localized-extension}
 B_{Q_\ell}^{\mathrm{ext}}F
 =\mathcal B_\ell^{\mathrm{ext}}(t,D_z)(\eta_\ell F)
   -\eta_\ell'e_\theta\bigl(K_{\ell,t}+U_{\ell,\theta}\bigr)F_z.
\end{equation}
In fact, the azimuthal part of the first term is
\(\partial_zK_{\ell,t}(\eta_\ell F_z)\); subtracting the
\(\eta_\ell'K_{\ell,t}F_z\) term leaves the radial-to-azimuthal
coupling in \eqref{dec:physical-extension}. The other subtracted term
is the axial shear.  Thus $B_{Q_\ell}^{\mathrm{ext}}$ agrees with the
physical $B_{Q_\ell}$ on solenoidal axisymmetric fields; away from that
subspace it is the extension in \eqref{dec:localized-extension}, not
necessarily the Cartesian matrix extension in \eqref{bg:cartesian-matrix}.
Write $B_Q^{\mathrm{ext}}=\sum_\ell B_{Q_\ell}^{\mathrm{ext}}$.

Let \(G_\ell\) be the convolution kernel of multiplication by
\(\eta_\ell\) in the unitary axial Fourier representation. Its mass
is at most one and its first absolute moment \(\mu_\ell\) is at most
\(C/(b\ell)\). Define the aggregate coefficients
\begin{equation}\label{dec:low-frequency-coefficients}
 b_{\Sigma}(t)=\sum_\ell b_\ell(t),\qquad
 g_{\Sigma}(t)=\sum_\ell b_\ell(t)\mu_\ell(t),\qquad
 \int_\tp^\tf b_{\Sigma}(t)\idiff t\leq CA,
 \quad\int_\tp^\tf g_{\Sigma}(t)\idiff t\leq\frac{CA}{b}.
\end{equation}
For \(F(k)\) supported in \(|k|\leq2\kappa_{\mathrm{cut}}\),
\eqref{dec:first-frequency-bounds} and
\(k(G_\ell*F)=G_\ell*(kF)+(\nu G_\ell)*F\) imply
\begin{equation}\label{dec:low-frequency-energy}
 \|\widehat{B_Q^{\mathrm{ext}}F}\|_{L^2_k(\mathcal H_h)}
       \leq C(\kappa_{\mathrm{cut}} b_{\Sigma}+g_{\Sigma})\|F\|_{L^2_k(\mathcal H_h)}.
\end{equation}
The shear term in \eqref{dec:localized-extension} has the
\(g_{\Sigma}\) bound as well.

Frequency differentiation commutes with \(G_\ell*\). The first two
differentiated equations therefore have the additional terms
\begin{equation}\label{dec:differentiated-sources}
 \sum_\ell\mathcal B_\ell^{\mathrm{ext}\prime}(G_\ell*F),\qquad
 \sum_\ell\left(2\mathcal B_\ell^{\mathrm{ext}\prime}(G_\ell*F')
                    +\mathcal B_\ell^{\mathrm{ext}\prime\prime}(G_\ell*F)\right).
\end{equation}
The differentiated inputs need not be solenoidal, which is the reason for
using \(\mathcal H_h\). The shear remains part of the principal
convolution operator. Primes on $\mathcal B_\ell^{\mathrm{ext}}$ below denote
$k$ derivatives. At integer order two,
\begin{equation}\label{dec:second-derivative-product}
 \|\mathcal B_\ell^{\mathrm{ext}\prime\prime}(G_\ell*F)\|_{L^2_k}
 \leq\|\mathcal B_\ell^{\mathrm{ext}\prime\prime}\|_{L^2_k(\mathrm{op})}
       \|G_\ell*F\|_{L^\infty_k}
 \leq C\|\mathcal B_\ell^{\mathrm{ext}\prime\prime}\|_{L^2_k(\mathrm{op})}\|F\|_{H^1_k}.
\end{equation}

\paragraph{Triangular estimates for frequency derivatives.}
The following elementary calculation is the organizing principle.
Suppose $x_0,x_1,x_2\geq0$ satisfy
\[
 x_0'\leq a x_0,\qquad
 x_1'\leq a x_1+b x_0,\qquad
 x_2'\leq a x_2+b x_1,
\]
where $a,b$ are nonnegative functions of time.  Set
$A_0(t)=\int_\tp^t a$ and $B_0(t)=\int_\tp^t b$.
Successive integration after multiplication by $e^{-A_0}$ gives
\[
 \begin{aligned}
 x_0(t)&\leq e^{A_0(t)}x_0(\tp),\\
 x_1(t)&\leq e^{A_0(t)}\bigl(x_1(\tp)+B_0(t)x_0(\tp)\bigr),\\
 x_2(t)&\leq e^{A_0(t)}\bigl(x_2(\tp)+B_0(t)x_1(\tp)
                      +\tfrac12B_0(t)^2x_0(\tp)\bigr).
 \end{aligned}
\]
A source at each level is treated by Duhamel with the same estimate on
each subinterval.  Additional lower-level terms change the polynomial,
not the leading exponential.  More generally, a strictly lower-triangular
four-by-four coupling has no nonzero product of four such matrices,
even at different times; its ordered expansion has only three
nontrivial coupling levels.

In the actual problem, $\int a=O(\kappa_{\mathrm{cut}} A)+O(A/b)$ and the off-diagonal
integrals are $O(A)$, with constants allowed to depend on the fixed
$\kappa_{\mathrm{cut}}$.  The resulting bound is a polynomial in $A$ times
$e^{C_s\kappa_{\mathrm{cut}} A}$.  Applying Gronwall instead to the sum of all levels
with coefficient $a+b$ would produce $e^{CA}$ and lose the purpose of
the argument.  The inequalities below are therefore kept separate until
they have been integrated in increasing derivative order.

\textit{A closed triangular hierarchy, with its constants.}
Let $F(t,k)=\widehat{W_<}^{\,z}(t,k)$, and let
$S(t,k)=\widehat{\mathcal S_<}^{\,z}(t,k)$.
The identity \eqref{dec:localized-extension} is first used on the
solenoidal field $W_<$.  It then defines the operator on all of
$\mathcal H_h$ before any frequency derivative is taken.
Denote this extension by $\mathscr B(t)$, so that
\[
 \partial_tF-\Delta_hF=\mathscr B(t)F+S.
\]
Its part containing $\eta_\ell'$ is a convolution in $k$ with
operator-valued kernel; its $L^1_k$ operator norm is at most $Cg_{\Sigma}$.
It commutes with frequency differentiation and translation.
The remaining part is the sum of
$\mathcal B_\ell^{\mathrm{ext}}(t,k)(G_\ell*F)(k)$.
Convolution by $G_\ell$ is a contraction in $L^2_k$, $L^\infty_k$,
and each $H^a_k(\mathcal H_h)$ used below.

Put
\[
 n_j(t)=\|\partial_k^jF(t)\|_{L^2_k(\mathcal H_h)},\quad j=0,1,2,
 \qquad \sigma_j(t)=\|\partial_k^jS(t)\|_{L^2_k(\mathcal H_h)}.
\]
All three fields $F,F',F''$ have frequency support in
$[-2\kappa_{\mathrm{cut}},2\kappa_{\mathrm{cut}}]$; their lack of solenoidality does not affect
\eqref{dec:low-frequency-energy}, which was proved for the extension.
There are nonnegative functions $c_2,c_\vartheta$ with
\[
 \int_\tp^\tf(b_{\Sigma}+c_2+c_\vartheta)\idiff t\leq C_sA
\]
such that the differentiated equations give
\[
 \begin{aligned}
 D^+n_0&\leq a_0n_0+\sigma_0,\\
 D^+n_1&\leq a_0n_1+b_{\Sigma}n_0+\sigma_1,\\
 D^+n_2&\leq a_0n_2+2b_{\Sigma}n_1+c_2(n_0+n_1)+\sigma_2,
 \qquad a_0=C_s(\kappa_{\mathrm{cut}} b_{\Sigma}+g_{\Sigma}).
 \end{aligned}
\]
Here $D^+$ denotes the upper right Dini derivative, and $c_2$ is a
constant multiple of
$\sum_\ell\|\mathcal B_\ell^{\mathrm{ext}\prime\prime}\|_{L^2_k(\mathrm{op})}$.
The second derivative is paired with $G_\ell*F$ in $L^\infty_k$,
which requires only the preceding $H^1_k$ level, exactly as in
\eqref{dec:second-derivative-product}.  No derivative of the initial
velocity in physical space occurs in these inequalities.

For clarity, the fractional product estimates used at the next level are
\begin{align}
 \|DF'\|_{H^\vartheta_k}
 &\leq C_\vartheta
   \left(\|D\|_\infty+\|D\|_{\mathfrak H^\vartheta}\right)
                 \|F\|_{H^2_k},\notag\\
 \|DF\|_{H^\vartheta_k}
 &\leq C_\vartheta\|D\|_{\mathfrak H^\vartheta}\|F\|_{H^2_k}.
 \label{dec:fractional-product}
\end{align}
They are Hilbert-valued estimates.  In a product difference, put the
difference on $D$ and use the $L^\infty_k$ norm of $F$ or $F'$.
In the other term use $\|D\|_\infty$ for $DF'$, and use
\[
 \|F(\cdot+h)-F\|_\infty
       \leq C\min(|h|,1)\|F\|_{H^2_k}
\]
together with $\|D\|_{L^2_k(\mathrm{op})}$ for $DF$.
The integrals of $\min(|h|,1)^2|h|^{-1-2\vartheta}$ are finite.
Thus the second estimate permits the logarithmic singularity of
$\mathcal B_\ell^{\mathrm{ext}\prime\prime}$ at zero.  The constant azimuthal block of
$\mathcal B_\ell^{\mathrm{ext}\prime}$ is estimated separately by
$\|F'\|_{H^\vartheta_k}\leq C\|F\|_{H^2_k}$; it is kept separate from the $L^2(\R_k)$ operator norm.
One may take $c_\vartheta$ to be a constant multiple of the sum of
$b_{\Sigma}$, the $\mathfrak H^\vartheta$ norms of the nonconstant blocks
of $\mathcal B_\ell^{\mathrm{ext}\prime}$, and those of $\mathcal B_\ell^{\mathrm{ext}\prime\prime}$.
Their asserted time integral follows from Lemma~\ref{dec:fixed-space}.

Define the truncated highest seminorm
\[
 q_\vartheta(t)^2=\int_{|h|\leq\kappa_{\mathrm{cut}}}
 \frac{\|F''(t,\cdot+h)-F''(t,\cdot)\|_2^2}
      {|h|^{1+2\vartheta}}\idiff h.
\]
The translated difference is supported in $[-3\kappa_{\mathrm{cut}},3\kappa_{\mathrm{cut}}]$.
The principal convolution operator therefore has the same bound
$C_s(\kappa_{\mathrm{cut}} b_{\Sigma}+g_{\Sigma})$, after increasing the absolute constant.
Its commutator with translation obeys
\begin{equation}\label{dec:principal-difference}
 \left(\int_{|h|\leq\kappa_{\mathrm{cut}}}
 \frac{\left\|\sum_\ell
       \bigl(\mathcal B_\ell^{\mathrm{ext}}(t,\cdot+h)-\mathcal B_\ell^{\mathrm{ext}}(t,\cdot)\bigr)
                         (G_\ell*F'')\right\|_2^2}
      {|h|^{1+2\vartheta}}\idiff h\right)^{1/2}
 \leq C_\vartheta b_{\Sigma}\kappa_{\mathrm{cut}}^{1-\vartheta}n_2.
\end{equation}
This follows directly from the global Lipschitz estimate
$\|\mathcal B_\ell^{\mathrm{ext}}(k+h)-\mathcal B_\ell^{\mathrm{ext}}(k)\|\leq b_\ell|h|$;
no division by $k$ is made, including when the two frequencies have
opposite signs.  The shear commutes with translation.
The remaining differentiated terms are the two terms in
\eqref{dec:differentiated-sources}; \eqref{dec:fractional-product}
bounds them using $n_0+n_1+n_2$.  Consequently
\[
 D^+q_\vartheta
 \leq a_0q_\vartheta+
    C_s\bigl(b_{\Sigma}\kappa_{\mathrm{cut}}^{1-\vartheta}+c_\vartheta\bigr)
                    (n_0+n_1+n_2)
       +C_s\|S\|_{H^{2+\vartheta}_k}.
\]
These norm inequalities can first be obtained with a positive
regularization of each norm; equivalently, apply Minkowski's inequality
to the Hilbert space in $(k,h)$ after the usual heat energy estimate.

The resulting system is triangular: the highest seminorm does not
appear on the right sides for $n_0,n_1,n_2$.  Divide all four inequalities
by $\exp(\int_\tp^t a_0)$ and integrate in increasing order.
There are only four levels, and every off-diagonal coefficient has a
time integral bounded by $C_{s,\kappa_{\mathrm{cut}}}A$.
Their iterated integrals therefore give a fixed polynomial in $A$.
The omitted translations satisfy
\[
 \int_{|h|>\kappa_{\mathrm{cut}}}
  \frac{\|F''(\cdot+h)-F''\|_2^2}{|h|^{1+2\vartheta}}\idiff h
       \leq C_\vartheta\kappa_{\mathrm{cut}}^{-2\vartheta}n_2^2.
\]
Thus $n_0+n_1+n_2+q_\vartheta$, with the factor
$\kappa_{\mathrm{cut}}^{-\vartheta}n_2$ included, controls the full $H^{2+\vartheta}_k$
norm.  In the equivalent combined-energy formulation,
$b_{\Sigma}\kappa_{\mathrm{cut}}^{1-\vartheta}n_2
 =\kappa_{\mathrm{cut}} b_{\Sigma}(\kappa_{\mathrm{cut}}^{-\vartheta}n_2)$.
Finally, $\int a_0\leq C_s\kappa_{\mathrm{cut}} A+C_sA/b$.
The last term is bounded for the chosen exponential width.

Under the unitary axial transform, $H^s_k(\mathcal H_h)$ is precisely
the axial moment space with weight $\langle z\rangle^s$.
The initial multiplier is $\chi_{\mathrm{cut}}(k/\kappa_{\mathrm{cut}})e^{-16A^2k^2}$.
For fixed $\kappa_{\mathrm{cut}}$, its derivatives through order three are polynomially
bounded in $A$; difference quotients give the same statement at order
$2+\vartheta$.  We conclude
\begin{equation}\label{dec:axial-low-moment}
 \sup_{\tp\leq t\leq\tf}\|\langle z\rangle^sW_<(t)\|_2
 \leq C_{s,\kappa_{\mathrm{cut}}}A^Pe^{C_s\kappa_{\mathrm{cut}} A}
       \left(M_s(f)+\int_\tp^\tf M_s(\mathcal S_<)\idiff t\right).
\end{equation}
The constant $C_s$ in the exponential is independent of $\kappa_{\mathrm{cut}}$.
All inverse powers of $\kappa_{\mathrm{cut}}$ enter only the prefactor.

For the complementary frequency range, this same Gaussian gives
\begin{equation}\label{dec:high-initial}
 M_s(W_>(\tp))
       \leq C_{s,\kappa_{\mathrm{cut}}}A^Pe^{-c_1\kappa_{\mathrm{cut}}^2A^2}M_s(f).
\end{equation}
The integer derivatives and the remaining fractional difference quotient
prove this estimate; horizontal weights commute with the axial
multiplier. Applying \eqref{dec:crude-moment} and
\eqref{dec:band-source-bound}, and then combining with
\eqref{dec:axial-low-moment}, gives
\begin{equation}\label{dec:axial-full-moment}
 \sup_{\tp\leq t\leq\tf}\|\langle z\rangle^sW(t)\|_2
 \leq C_{s,\kappa_{\mathrm{cut}}}A^P
       \left(e^{C_s\kappa_{\mathrm{cut}} A}+b^{-1}e^{K_{s,\kappa_{\mathrm{cut}}}A}\right)M_s(f).
\end{equation}

\paragraph{Recovery of the full spatial moment.}
For fixed $s\geq0$, the ordinary weight is equivalent to the sum of its
axial and horizontal parts:
\[
 \langle x\rangle^s\asymp_s
       \langle x_h\rangle^s+\langle z\rangle^s.
\]
The frequency
argument has controlled only the second part.  In the equation for $W$,
the smaller columns are Gaussian-localized in $x_h$, so their multiplication
operators can supply a horizontal weight at the cost of powers of their
heat ages.  These powers have only polynomial time integrals on the source
interval.  The larger columns are not localized near the active core,
but their complete matrix is small, of order $A/R^2$; their contribution
can be absorbed by a weighted Gronwall estimate.  The next calculation
combines these two estimates to recover the full spatial moment.

We next recover the horizontal part of the moment. For the columns of
radius at most one, Gaussian localization gives
\begin{equation}\label{dec:horizontal-recovery}
 M_s(M_{Q,\mathrm{in}}W)
 \leq C_sA\sum_{\ell\leq1}\ell^2T_\ell^{-2}
       \left((1+T_\ell^{s/2})\|W\|_2
                  +\|\langle z\rangle^sW\|_2\right).
\end{equation}
The coefficient has a polynomial time integral. The remaining matrix
has norm at most \(CA/R^2\). In the Duhamel formula for
\eqref{dec:filtered-equation}, use
\(M_s(e^{h\Delta_h}Y)\leq C_s(1+h^{s/2})M_s(Y)\) and the
odd-field projection bound from Lemma~\ref{bg:weighted}.
In detail, with $Y_* = \sup_{[\tp,\tf]}M_s(W)$, the outer
part of that Duhamel formula is bounded by $CA^PR^{-2}Y_*$.
The polynomial includes the partial-heat weight factor and the interval
length $16A^2$.  Since $R=e^{\varrho A}$, this coefficient is less
than $1/2$ for sufficiently large $A$ and can be absorbed.
The inner-column coefficient in \eqref{dec:horizontal-recovery} has
a polynomial time integral, so it multiplies the already established
axial moment bound, not $Y_*$.  Equations \eqref{dec:filtered-source},
\eqref{dec:axial-full-moment}, and \eqref{dec:horizontal-recovery}
therefore imply
\begin{equation}\label{dec:heat-array-decay}
 M_s(v_Q(\tf))
 \leq C_{s,\kappa_{\mathrm{cut}}}A^P
       \left(e^{C_s\kappa_{\mathrm{cut}} A}+b^{-1}e^{K_{s,\kappa_{\mathrm{cut}}}A}\right)M_s(f),
\end{equation}
because \(W(\tf)=v_Q(\tf)\).

It remains to include the correction \(e^{I,J}\) in the true reference
solution. The difference \(d=v-v_Q\) has zero initial data and satisfies
the heat-array equation with the additional tensor
\begin{equation}\label{dec:true-difference-tensor}
 -\PP\operatorname{div}
       \left(e^{I,J}\otimes(v_Q+d)+(v_Q+d)\otimes e^{I,J}\right).
\end{equation}
Proposition~\ref{bg:reference} gives
\(\sup_{\tp\leq t\leq\tf}\|e^{I,J}(t)\|_\infty
\leq CA^Pe^{K_0'A}b^{-1/2}=:\varepsilon_b\).
The tensor has moment at most
\(2\varepsilon_b(M_s(v_Q)+M_s(d))\).
Using the tensor-heat estimate from Lemma~\ref{bg:weighted} and
\eqref{dec:crude-moment}, we obtain, with
\(D_\ast=\sup_{\tp\leq t\leq\tf}M_s(d(t))\),
\begin{equation}\label{dec:true-difference-bound}
 D_\ast\leq CA^Pe^{K_sA}\varepsilon_b
       \left(A^Pe^{K_sA}M_s(f)+D_\ast\right).
\end{equation}
For completeness, this response estimate follows directly from the
full-heat mild equation.  Its matrix part has coefficient $C_sm_0(t)$
with integral $O_s(A)$; its tensor part has kernel
$C_s(1+(t-u)^{-1/2})\varepsilon_b$.
First apply the scalar Volterra Gronwall estimate to the matrix part,
then take a time supremum.  The time integral of the tensor kernel is
polynomial in $A$, giving exactly \eqref{dec:true-difference-bound}.
Thus no bounded weighted operator norm for $\PP\div$ at zero elapsed
time has been assumed.
For sufficiently large fixed \(\beta\), the last term is absorbed, and
\(D_\ast\leq CA^Pe^{K_{s,\kappa_{\mathrm{cut}}}A}b^{-1/2}M_s(f)\).
Choose first \(\kappa_{\mathrm{cut}}\) so that \(C_s\kappa_{\mathrm{cut}}<\lambda/2\), then
\(\beta>2K_{s,\kappa_{\mathrm{cut}}}+2\), and finally increase \(A\).
Together with \eqref{dec:heat-array-decay}, this proves
\eqref{dec:moment-decay} for \(2<s<5/2\).

The argument is first applied to smooth centrally odd solenoidal data
and finite arrays. Estimate \eqref{dec:crude-moment} provides weighted
existence before its improvement is used. Radial truncation, rotation
averaging, odd projection, Leray projection, and heat regularization
give dense admissible approximations in \(M_s\). The exact-reference cutoff convergence proved at the end of
Proposition~\ref{bg:reference}, followed by positive-time propagator
convergence and lower semicontinuity of $M_s$, removes the inner cutoff.
At fixed cutoff, the same argument removes the smooth data regularization;
unweighted uniqueness identifies the limit.  The constants are uniform
before either limit is taken. 

For the lower moments, it is useful to specify the interpolation principle.
If one linear operator $T$ is bounded on $L^2$ and on
$L^2(\langle x\rangle^{2s'}\idiff x)$, its interpolated bound at
$0<s<s'$ is
\[
 \|T\|_{L^2(\langle x\rangle^{2s})\to L^2(\langle x\rangle^{2s})}
 \leq \|T\|_{2\to2}^{1-s/s'}
       \|T\|_{L^2(\langle x\rangle^{2s'})\to
                        L^2(\langle x\rangle^{2s'})}^{s/s'}.
\]
The weighted spaces here are the interpolation spaces for multiplication
by a positive weight; equivalently one applies the three-lines inequality
to $\langle x\rangle^{zs'}T\langle x\rangle^{-zs'}$ on a dense core.
Both endpoint operator bounds are proved with the same prescribed loss,
after one width choice satisfying their two thresholds.
Interpolation defines the operator on the common core, and density extends it to the stated data class with finite $M_s$.

At \(s=0\), the same argument
uses \eqref{dec:low-frequency-energy} without frequency differentiation.
For \(0<s\leq2\), compose the propagator with the projection onto
centrally odd, axisymmetric, solenoidal fields. This projection is
bounded in every moment below \(5/2\), by Lemma~\ref{bg:weighted}.
Interpolation between its unweighted bound and the bound at any fixed
\(2<s'<5/2\), both proved with the prescribed exponent \(\lambda\),
establishes \eqref{dec:moment-decay} on the remaining weighted spaces.
The bounded symmetry and Leray projections extend the interpolated
operator from the common smooth core, so the resulting estimate requires
only $M_s(f)<\infty$.
\end{proof}

\subsection{An unweighted Sobolev estimate at the final time}

The moment estimate allows a small exponential loss.  For the continuing
quadratic source we also need a stronger unweighted endpoint bound whose
leading term is only polynomial in $A$.  The invariant evolution provides
this bound by balancing its $A|k|$ growth against axial heat; localization
and the reference correction are added as small errors.

The next estimate will control the quadratic forcing that continues
after \(\tf\). Its leading term is polynomial in \(A\).

\begin{lemma}[The final-time Sobolev estimate]\label{dec:endpoint}
For \(1\leq s<5/2\), the solution in Lemma~\ref{dec:retirement}
satisfies
\begin{equation}\label{dec:endpoint-bound}
 \|v(\tf)\|_{H^4}
 \leq CA^P\|f\|_2+CA^Pe^{K_sA}b^{-1/2}M_s(f),
\end{equation}
where \(K_s\) is fixed before the width exponent \(\beta\).
\end{lemma}

\begin{proof}
Let \(V\) evolve from \(f\) about the invariant sum
\(U=\sum_\ell U_\ell\). By \eqref{dec:first-frequency-bounds},
the horizontal-heat evolution at axial frequency \(k\) has norm at
most \(e^{CA|k|}\). Restoring the axial heat over a time interval of
length \(16A^2-O(1)\) and maximizing
\(-(16A^2-O(1))k^2+CA|k|\), we find
\begin{equation}\label{dec:invariant-terminal-l2}
 \|V(\tf-\tau)\|_2\leq C\|f\|_2,
 \qquad 0\leq\tau\leq2.
\end{equation}
Every fixed derivative of the invariant matrix is bounded on the final
two units of time. Cartesian parabolic smoothing gives
\(\|V(\tf)\|_{H^4}\leq CA^P\|f\|_2\), while
\eqref{dec:crude-moment} gives
\begin{equation}\label{dec:invariant-first-moment}
 \sup_{\tp\leq t\leq\tf}M_1(V(t))
       \leq CA^Pe^{KA}M_1(f).
\end{equation}

The complete cutoff difference, including the axial shear, obeys
\begin{equation}\label{dec:cutoff-difference}
 \| (M_Q-M_U)V\|_2
 \leq\frac{CA}{b}\sum_\ell\ell
       \left(T_\ell^{-2}
             \left(\sqrt t\|V\|_2+\|zV\|_2\right)
                    +T_\ell^{-3/2}\|V\|_2\right).
\end{equation}
Here
\(|1-\eta_\ell|\leq C(b\ell)^{-1}(\sqrt t+|z|)\).
The inner-radius sum uses \(\sum_{\ell<1}\ell\leq C/R\),
and the outer-radius sums use
\(\sum_{\ell\geq R}\ell^{-3}\) and
\(\sum_{\ell\geq R}\ell^{-2}\). Their time integrals are
polynomial in \(A\).
The exact difference \(v-V\) has this vector source and the full
reference-correction tensor from \eqref{dec:true-difference-tensor}.
The unweighted version of \eqref{dec:true-difference-bound} yields
\begin{equation}\label{dec:invariant-actual-difference}
 \sup_{\tp\leq t\leq\tf}\|v(t)-V(t)\|_2
       \leq CA^Pe^{KA}b^{-1/2}M_s(f).
\end{equation}

On \([\tf-2,\tf]\), all required matrix and correction derivatives
are bounded. For \(j\geq1\), the useful uniform cutoff estimate is
\begin{equation}\label{dec:cutoff-derivatives}
 \|\partial_z^j\eta_\ell\|_\infty
 \leq\frac{C_j}{b\ell}
       (b^2\ell^2+2t)^{-(j-1)/2}
 \leq\frac{C_j}{b\ell}T_\ell^{-(j-1)/2}.
\end{equation}
At moment exponent one the ordinary double-Riesz estimate applies to
arbitrary differentiated fields. The Cartesian derivative energies with
weight \(\langle x\rangle^2\), multiplied at derivative order \(j\)
by \((t-(\tf-2))^j\), give
\begin{equation}\label{dec:terminal-weighted-graphs}
 \sup_{\tf-1\leq t\leq\tf}
       \sum_{|\nu|\leq4}M_1(\partial^\nu V(t))
       \leq CA^Pe^{KA}M_1(f).
\end{equation}
Matrix commutators contain only lower derivatives, and derivatives of
the weight have lower order.

Leibniz's rule, \eqref{dec:cutoff-derivatives}, and
\eqref{dec:terminal-weighted-graphs} bound the cutoff source in
\(H^3\) by \(CA^Pe^{KA}b^{-1}M_s(f)\) on the last unit interval.
The correction source
\(-\PP\operatorname{div}(e^{I,J}\otimes V+V\otimes e^{I,J})\)
has \(H^3\) norm at most \(CA^Pe^{KA}b^{-1/2}M_s(f)\).
The true homogeneous evolution there satisfies
\begin{equation}\label{dec:terminal-smoothing}
 \|\mathcal U_{\bar u}(\tf,\tf-1)g\|_{H^4}\leq C\|g\|_2,
 \qquad
 \|\mathcal U_{\bar u}(\tf,t)g\|_{H^4}
       \leq C(\tf-t)^{-1/2}\|g\|_{H^3}.
\end{equation}
These estimates follow from the Cartesian derivative heat energies,
retaining the correction transport term. Applying
\eqref{dec:terminal-smoothing} to the difference equation, using
\eqref{dec:invariant-actual-difference} at \(\tf-1\), proves the
width-dependent term in \eqref{dec:endpoint-bound}.
The invariant estimate following
\eqref{dec:invariant-terminal-l2} supplies its leading term.
\end{proof}

\section{Transfer to the next larger scale}
\label{sec:transfer}

The quadratic history in Section~\ref{sec:history} produces a field whose
leading large-scale distribution is a multiple of
$\mathcal I=\PP\partial_z(e_z\delta_0)$.  Its physical norm can nevertheless
be large, and its subsequent evolution still interacts with every old,
smaller column.  Negative-Sobolev closeness to the point input alone is
therefore insufficient to discard those interactions.

The central step, Lemma~\ref{trans:invariant-entrance}, computes the
old-column azimuthal-mean contribution and bounds the other interactions.
Only after identifying that leading coefficient do we include the axial
cutoffs, the nonlinear reference correction, and the continuing quadratic
source.  Parent-coordinate smoothing and exact time alignment then return
the normalized output to the prepared shape $S_A$.

We retain \(c_\star=\sqrt2e^{-3/2}\),
\(\varrho=c_\star/8\), \(R=e^{\varrho A}\),
\(b=e^{\beta A}\), and the times \(\tp,\tf\) from
Section~\ref{sec:decay}. The background is
\(\bar u^{I,J}=Q^{I,J}+e^{I,J}\), where \(J\geq1\) and
\(0\leq I\leq\infty\). In the coordinates of the parent scale, set
\begin{equation}\label{trans:times}
 d=\frac{\tf}{R^2},\qquad t_1=A^{-1},\qquad
 s_\ast=\tp-d,\qquad \ts=\tf+R^2s_\ast=R^2\tp.
\end{equation}
For sufficiently large \(A\), one has
\(d<\tp/2\), \(s_\ast>3/A\), and \(\ts<R^2\), so the reference
solution in Proposition~\ref{bg:reference} exists on the entire interval.
All constants below are uniform in the cutoffs. A coefficient written
\(CA^P\) may depend on the fixed derivative orders and fixed exponents;
its exponent is independent of \(A\).

\subsection{Physical fibers and the point response}

These estimates fix the normalization used in transfer.  Working one
axial fiber at a time gives an integrable initial heat singularity for
the point input, and only then do we integrate in frequency.  The point
response and the smoothing estimates use physical velocity coordinates.

The source coordinates are normalized around a radius-one active
column.  During transfer the parent variable is $y=x/R$ and the
parent elapsed time is $s=(t-t_{\mathrm f})/R^2$.
The scalar $d=t_{\mathrm f}/R^2$ records the age of the parent column
at the beginning of transfer.  The terminal transfer time
$s_*=t_{\mathrm p}-d$ is chosen so that
$t_{\mathrm f}+R^2s_*=R^2t_{\mathrm p}$ exactly.
Consequently the parent heat age agrees with the prepared input at the
next scale.  Errors in this matching cannot be absorbed by merely
changing the scalar amplitude.

An invariant swirl of radius \(\ell\) and heat age \(T\) has velocity
\(A\ell^2T^{-2}e^{-r^2/(4T)}x_h^\perp\). Write
\(\Omega_T=T^{-2}e^{-r^2/(4T)}\) and
\(B_T=(2-r^2/(2T))\Omega_T\). Its physical linearization on the
solenoidal axisymmetric fiber of axial frequency \(k\) is
\begin{equation}\label{trans:physical-coupling}
 B_{\ell,T,k}v
 =A\ell^2\left(\PP_k(2\Omega_Tv_\theta e_r)
                       -B_Tv_re_\theta\right).
\end{equation}
Before projection, this is a smooth Cartesian matrix. In particular,
if \(J_h\) denotes horizontal rotation through a right angle, then
\(2\Omega_T(e_r\otimes e_\theta-e_\theta\otimes e_r)
=-2\Omega_TJ_h\) and
\(r^2e_\theta\otimes e_r=x_h^\perp\otimes x_h\).
Its norm is at most \(2A\ell^2T^{-2}\).

We shall also use the low-frequency estimate
\begin{equation}\label{trans:low-frequency-coupling}
 \|B_{\ell,T,k}\|_{2\to2}
       \leq\sqrt2A\ell^2|k|T^{-3/2}.
\end{equation}
To verify it, define
\(\Phi_v(r)=-2\int_r^\infty\Omega_T(u)v_\theta(u)\idiff u\).
The squared Hilbert--Schmidt norm of this radial integral operator on
\(L^2(r\idiff r)\) is \(2T^{-3}\), and full-gradient cancellation
gives \(\PP_k(\nabla_h\Phi_v,0)=-ik\PP_k(e_z\Phi_v)\).
Incompressibility gives
\begin{equation}\label{trans:radial-low-frequency}
 v_r(r)=-\frac{ik}{r}\int_0^r uv_z(u)\idiff u,
 \qquad |v_r(r)|\leq\frac{|k|}{\sqrt2}\|v_z\|_{L^2(r\idiff r)}.
\end{equation}
Since \(\int_0^\infty B_T^2r\idiff r=2T^{-3}\), the orthogonality
of the two outputs in \eqref{trans:physical-coupling} proves
\eqref{trans:low-frequency-coupling}. Both outputs vanish at \(k=0\).

With the unitary axial Fourier transform fixed in
Section~\ref{sec:background}, the point distribution has fiber
\begin{equation}\label{trans:point-fiber}
 \mathcal I_k=(2\pi)^{-1/2}ik\PP_k(e_z\delta_{0,h}).
\end{equation}
For \(1\leq T_{\mathrm{par}}\leq2\), let \(\mathcal V_{T_{\mathrm{par}}}(s)\) be its
full linearized response about the invariant radius-one swirl whose
initial heat age is \(T_{\mathrm{par}}\). Since
\(\mathcal I\in H^{-3}(\R^3)\), bounded perturbation of heat in
\(H^{-3}\), followed by positive-time smoothing, constructs this
response. Its fibers satisfy
\begin{equation}\label{trans:point-upper}
 \|\mathcal V_{T_{\mathrm{par}},k}(s)\|_{L^2_h}
       \leq C|k|s^{-1/2}e^{CAs-k^2s},\qquad
 \|\mathcal V_{T_{\mathrm{par}}}(s)\|_2\leq Cs^{-5/4}e^{CAs}.
\end{equation}
Indeed \(\|e^{s\Delta_h}\mathcal I_k\|_2\leq C|k|s^{-1/2}\),
whose singularity is integrable in the bounded-matrix Duhamel formula,
and axial heat factors out exactly. Plancherel gives the second bound.

For later use, suppose a Cartesian matrix and its fixed derivatives are
bounded by \(C_jA\). The associated horizontal fiber evolution, with
the axial heat removed and denoted by \(\widetilde{\mathcal U}_k\),
and its adjoint satisfy
\begin{equation}\label{trans:fiber-smoothing}
 \|\widetilde{\mathcal U}_k(s,u)\|_{H^j_h\to H^j_h}
       \leq e^{C_jA(s-u)},\qquad
 \|D_h^j\widetilde{\mathcal U}_k(s,u)\|_{2\to2}
       \leq C_jA^{j/2}\qquad\text{if }s-u\asymp A^{-1}.
\end{equation}
Here $D_h^j$ denotes the collection of Cartesian horizontal derivatives
of total order $j$.  These estimates are uniform in real \(k\). They
follow by differentiating the Cartesian equation, using the commutation of
\(\PP_k\) with horizontal derivatives, and summing derivative
energies with successive powers of \(s-u\). The adjoint extension is
\(\Delta_h+\PP_kM^\ast\PP_k\), to which the same argument applies.

Lemma~\ref{trans:point-size} is an independent normalization check:
the coupling does not cancel the leading size of the point response on
the short interval $[0,A^{-1}]$.  The subsequent cascade uses the prepared
normalization $n_A$ from Lemma~\ref{seed:prepared}; it does not require
this separate lower bound.  We retain it to make the short-time scale
of the point response explicit.

The incoming approximation is not small in the physical $L^2$ norm:
its limiting point distribution is not itself in that space.  What is
small is its action on sufficiently smooth tests.  For example, with
$\mathfrak q=4+\delta$ and a unit terminal test $\psi$,
\[
 \left|\left\langle
       \mathcal U_{\mathrm{out},k}(t_1,0)(F_k-\mathcal I_k),\psi
            \right\rangle\right|
 \leq \|F_k-\mathcal I_k\|_{H_h^{-\mathfrak q}}
       \|\mathcal U_{\mathrm{out},k}(t_1,0)^*\psi\|_{H_h^{\mathfrak q}}.
\]
Positive elapsed time smooths the adjoint test.  Over time $t_1=A^{-1}$,
that smoothing costs a fixed power of $A$, so the small negative-norm
error becomes a small $L^2$ output error.  This is the meaning of dual
smoothing in the transfer proof.

The outer evolution used in this explanation is defined in the next
subsection.  Its coefficients and fixed derivatives have polynomial
bounds.  The old columns have much shorter scales, so the same argument
cannot simply be applied to their full evolution with uniform derivative
bounds.  They are separated out in an interaction expansion.  This is
why the weak approximation alone is not the whole transfer theorem.

\begin{lemma}[Bounds for the point response]\label{trans:point-size}
Uniformly for \(1\leq T_{\mathrm{par}}\leq2\), the point response satisfies
\begin{equation}\label{trans:point-size-bound}
 c_1A^{5/4}\leq\|\mathcal V_{T_{\mathrm{par}}}(A^{-1})\|_2
                  \leq C_1A^{5/4}
\end{equation}
for sufficiently large \(A\).
\end{lemma}

\begin{proof}
The upper bound is \eqref{trans:point-upper}. The axis value of the
matrix is \(-2A(T_{\mathrm{par}}+s)^{-2}J_h\). On the solenoidal space,
\(\PP J_h\PP\) is skew-adjoint and commutes with heat. Its point
response \(Z(s)\) therefore obeys
\begin{equation}\label{trans:axis-response}
 \|Z(s)\|_2=\|e^{s\Delta}\mathcal I\|_2=c_{\mathcal I}s^{-5/4},
 \qquad \||x|^2Z(s)\|_2\leq Cs^{-1/4},
 \quad 0<s\leq A^{-1},
\end{equation}
where \(c_{\mathcal I}>0\). For the weighted bound, differentiate
twice the Fourier symbol of \(Z\). The rotation parameter is at most
two; its derivatives of order \(j\leq2\) are bounded by
\(C_j|\xi|^{-j}\). The worst resulting singularity is
\(|\xi|^{-1}\), which is square integrable near zero in three
dimensions, and the change of variables \(\xi=s^{-1/2}\zeta\)
gives the asserted power.

The true matrix differs from its axis value by at most \(CA|x_h|^2\)
globally. Duhamel's formula and the bounded \(L^2\) propagator on
\((0,A^{-1})\) therefore imply
\begin{equation}\label{trans:axis-error}
 \|\mathcal V_{T_{\mathrm{par}}}(A^{-1})-Z(A^{-1})\|_2
       \leq CA\int_0^{A^{-1}}s^{-1/4}\idiff s
       \leq CA^{1/4}.
\end{equation}
For large \(A\), this is smaller than half the leading norm in
\eqref{trans:axis-response}, proving the lower bound.
\end{proof}

\subsection{Passage through the invariant columns}

The normalized incoming field is close to a point derivative in a
negative Sobolev norm, but it can be large in $L^2$ after rescaling.
We estimate its pairing with a smooth terminal test at the parent scale.
Running that test backward through the adjoint evolution reveals its
vanishing at the symmetry axis and gives an additional scale factor.
The first interaction contains the retained azimuthal mean, whose
coefficient must be computed explicitly.  The subsequent estimates distinguish
that leading contribution from the small remainder.

Before introducing the entire array, consider an invariant old column of
radius $\ell$, at source time $\tf$, and put $L=\tf+\ell^2$.
Let $\varphi$ be the incoming field and
$\mu_\theta[\varphi](r)=\int_\R\varphi_\theta(r,z)\,\idiff z$ its azimuthal mean.
If this already formed mean is subsequently propagated by free heat, its
mixed-stress contribution against that old column is
\[
 \begin{aligned}
 \mathcal T_\ell(\varphi)
 &=2\pi A\ell^2\int_0^\infty\int_0^\infty
       W_{L+u}(r)(e^{u\Delta_h^{(1)}}\mu_\theta[\varphi])(r)r\,\idiff r\idiff u\\
 &=2\pi A\ell^2\int_0^\infty\int_0^\infty
       W_{L+2u}(r)\mu_\theta[\varphi](r)r\,\idiff r\idiff u\\
 &=4\pi A\ell^2\int_0^\infty
       (1-e^{-r^2/(4L)})\mu_\theta[\varphi](r)\,\idiff r.
 \end{aligned}
\]
Here $\Delta_h^{(1)}$ is the radial heat generator in
\eqref{seed:mean-equation}; the middle line follows from
self-adjointness of radial heat, and the last uses its explicit time
integral.  The weighted duality and heat-age regularization of
Section~\ref{sec:history} justify this calculation without treating a
non-$L^2$ radial test as an $L^2$ vector.

The calculation propagates an existing mean; it adds no quadratic forcing
after $\tf$.  Summing it over the old radii and dividing by the incoming
signed mass gives the zero-frequency scalar $c_{\mathrm{mean}}(0)$ below.  Thus a small
old-column effect need not vanish: its accumulated action changes the
coefficient of the next point response.

This calculation is only the leading part of the full transfer.  The
parent column is still evolving, nonzero axial frequencies carry heat,
and the incoming field can interact repeatedly with old columns.
The adjoint proof below compares the genuine evolution with the displayed
calculation and estimates each of these differences.  In particular, it
does not replace the actual incoming field by its axial mean.

The next lemma requires only
\(A\leq\tf\), \(R\geq A^2\), and
\(\tf+1\leq R^2/(2A)\), which hold for our parameters.
Let \(\ell_j=R^{-j}\), \(0\leq j\leq I\), be the old radii, and set
\(\rho_j=\ell_j/R\), \(\tau_j=(\tf+\ell_j^2)/R^2\), and
\(T_{\mathrm{par}}=1+d\). In parent coordinates \(y=(y_h,y_3)\), the invariant
old and outer fields are
\begin{align}
 U_{\mathrm{in}}(s,y)
 &=\sum_{j=0}^I A\rho_j^2(\tau_j+s)^{-2}
             e^{-|y_h|^2/(4(\tau_j+s))}y_h^\perp,
             \label{trans:invariant-mixture}\\
 U_{\mathrm{out}}(s,y)
 &=\sum_{j=0}^{J-1} AR^{2j}(R^{2j}+d+s)^{-2}
             e^{-|y_h|^2/(4(R^{2j}+d+s))}y_h^\perp.\notag
\end{align}
The \(j=0\) term in the second sum is the parent column. This
invariant sum is an exact unforced solution: its radial pressure
derivative is
\((U_{\mathrm{in},\theta}+U_{\mathrm{out},\theta})^2/|y_h|\).
Thus its linearization retains all ordered interactions. The old sums
converge with the majorant \(\sum_j\ell_j^2\leq4/3\), and the
farther outer sums with geometric derivative majorants. In this
subsection, the outer cutoff may also equal infinity.

Let \(\varphi\) be real, axisymmetric, and solenoidal, with
\(\varphi_r,\varphi_\theta\) even in \(z\) and \(\varphi_z\) odd.
For a nonzero real scalar \(m\), define the normalized incoming field
and its fibers by
\begin{equation}\label{trans:incoming-normalization}
 F_R(y)=\frac{R^4}{m}\varphi(Ry),\qquad
 F_k(y_h)=\frac{R^3}{m}\widehat\varphi^{\,z}(Ry_h,k/R),\qquad
 \mathcal W=\frac{M_{9/4}(\varphi)}{|m|}.
\end{equation}
Fix \(K\geq1\), \(0<\delta<1/2\), and \(P_{\mathrm{in}}\geq2\). Suppose that,
for \(|k|\leq KA\),
\begin{equation}\label{trans:incoming-bounds}
 \|F_k\|_2\leq C_0A^{P_{\mathrm{in}}}R^2,\qquad
 \|F_k-\mathcal I_k\|_{H^{-4-\delta}_h}\leq\epsilon_K.
\end{equation}
With $\varphi,m$ fixed, write $c_{\mathrm{mean}}(k)$ for
$c_{\mathrm{mean}}(k;\varphi,m)$, and define the retained-mean kernel
and coefficient by
\begin{align}
 K_{\mathrm{in}}(r)
 &=\sum_{j=0}^I\ell_j^2
       \left(1-e^{-r^2/(4(\tf+\ell_j^2))}\right),
       \label{trans:mean-scalar}\\
 c_{\mathrm{mean}}(k)&=\frac{4\pi A\sqrt{2\pi}}{m}
       \int_0^\infty K_{\mathrm{in}}(r)
                  \widehat\varphi_\theta^{\,z}(r,k/R)\idiff r,
 \qquad a_{\mathrm{in}}=1+c_{\mathrm{mean}}(0),\notag
\end{align}
and put
\begin{equation}\label{trans:entrance-error}
 E_K=(1+R^\delta)\epsilon_K+R^{-1}\log(2+R)
       +\mathcal W\left(R^{-7/4}+R^{-3/2}\log(2+R)\right).
\end{equation}
The factor \(\sqrt{2\pi}\) in \eqref{trans:mean-scalar} converts the
unitary axial transform to the physical axial integral. In particular,
\(c_{\mathrm{mean}}(0)\) is real.

\begin{lemma}[Entrance through the invariant array]\label{trans:invariant-entrance}
Let \(V\) be the exact linearized solution about
\(U_{\mathrm{in}}+U_{\mathrm{out}}\), with \(V(0)=F_R\). Then
\begin{align}
 \|V(t_1)-a_{\mathrm{in}}\mathcal V_{T_{\mathrm{par}}}(t_1)\|_2
 &\leq CA^{P_{\mathrm{in}}+10}E_K+C|a_{\mathrm{in}}|A^{3/2}R^{-2}\notag\\
 &\quad+C_Ke^{-K^2A}
       \left(R^{5/2}\mathcal W+|a_{\mathrm{in}}|A^{3/2}\right).
 \label{trans:invariant-entrance-bound}
\end{align}
The constants are independent of both cutoffs, including their
infinite limits.  The incoming field is assumed to satisfy
$M_{9/4}(\varphi)<\infty$; constants may depend on the fixed
$K,\delta,P_{\mathrm{in}},C_0$, but not on $A,R,I,J$.
\end{lemma}

\begin{proof}
Write \(\mathcal U_{\mathrm{out},k}\) for the outer propagator and
\(\mathcal V_{\mathrm{out},k}(s)
=\mathcal U_{\mathrm{out},k}(s,0)\mathcal I_k\).
The derivatives of the outer matrix are bounded by \(C_jA\), and
those of its tail after the parent column by \(C_jAR^{-2}\).
For a unit terminal test \(\psi\) in the physical fiber, set
\(\Psi_s=\mathcal U_{\mathrm{out},k}(t_1,s)^\ast\psi\).
On \(0\leq s\leq t_1/2\), horizontal interpolation and
\eqref{trans:fiber-smoothing} give
\begin{align}
 \|D_h^j\Psi_s\|_\infty&\leq C_jA^{(j+1)/2},\notag\\
 \Psi_{s,r}(r)&=c_r(s)r+O(A^2r^3),\qquad
 \Psi_{s,\theta}(r)=c_\theta(s)r+O(A^2r^3).
 \label{trans:adjoint-axis-expansion}
\end{align}
Rotational covariance makes the radial profiles odd on a diameter,
and Taylor's formula with the uniform third derivatives gives the
remainders for all \(r\geq0\). Incompressibility implies
\begin{equation}\label{trans:adjoint-slopes}
 2c_r(s)=-ik\Psi_{s,z}(0),\qquad
 |c_r(s)|+|c_\theta(s)|\leq CA,
 \qquad |c_r(0)|\leq C|k|\sqrt A\,e^{-k^2/A}.
\end{equation}
On \(|k|\leq KA\), the backward equation also gives
\begin{equation}\label{trans:slope-derivative}
 |c_r'(s)|\leq C_KA^3.
\end{equation}
The horizontal heat term costs \(CA^2\), the axial term costs
\(Ck^2A\), and the projected matrix term costs \(CA^2\).
For the last bound, interpolate
\(D_h\PP_kM_{\mathrm{out}}^\ast\Psi_s\), whose \(L^2\) and
second derivative norms are bounded by \(CA^{3/2}\) and
\(CA^{5/2}\), respectively. This calculation is uniform down to
\(k=0\).

The incoming fiber can have norm of order $R^2$.  Thus the small
integrated norm of one old-column interaction, by itself, is not enough
for a small output error.  The adjoint test's linear vanishing at the axis
provides an additional gain in the last such interaction.  With two or
more old interactions this gain closes the estimate below.  With exactly
one interaction, its leading mean contribution must instead be computed
and retained.

The Volterra expansion around the outer propagator is absolutely
convergent on each physical fiber.  Its term with $n$ old interactions
is integrated over $0<s_1<\cdots<s_n<t_1$.
We separate $n=0$, $n=1$, and $n\geq2$; interactions at different old
radii remain inside the sums.
The outer propagator has norm at most $e^{CA(t-u)}$, so the product
of all its factors on disjoint intervals costs at most $e^{CAt_1}$,
independently of $n$.  This prevents an artificial constant to the power
of the number of interactions.

We first sum the terms with at least two old-column interactions.
By \eqref{trans:low-frequency-coupling}, their integrated operator
size is
\begin{equation}\label{trans:old-integrated-operator}
 B_\ast:=\int_0^{t_1}\|B_{\mathrm{in}}(s,k)\|\idiff s
 \leq\frac{2\sqrt2A|k|}{R}
       \sum_{j=0}^I\frac{\ell_j^2}{\sqrt{\tf+\ell_j^2}}
 \leq\frac{CA|k|}{R\sqrt\tf}.
\end{equation}

\paragraph{Vanishing of the horizontal adjoint components at the axis.}
The gain can be seen with one old radius $\rho=\ell/R$ and heat age
$T=\tau+s$.  Its matrix is bounded pointwise by a constant times
$A\rho^2T^{-2}$ times a Gaussian polynomial.  A generic bounded test
would not exploit its concentration.  The horizontal adjoint components,
however, vanish linearly at the axis: $|\Psi_h(r)|\leq CA r$ on the
early half-interval.  Since
\[
 \|r e^{-r^2/(cT)}\|_{L^2(\R^2)}\leq C T,
\]
this gives
$\|B_{\mathrm{old}}^*\Psi_s\|_2\leq CA^2\rho^2T^{-1}$,
not the cruder scale $A\rho^2T^{-2}\|\Psi_s\|_2$.
Time integration costs a logarithm rather than an inverse age:
\[
 \int_0^{t_1/2}\frac{\idiff s}{d+s}
        =\log\left(1+\frac{t_1}{2d}\right)
        \leq C\log(2+R).
\]
The ensuing $R^{-2}$ adjoint gain cancels the allowed $R^2$ in the
incoming norm.  One additional old interaction contributes another small
factor of order $A^{3/2}/R$.  Thus terms with at least two interactions
are small.  With exactly one interaction that extra factor is absent,
which explains why its leading part must be computed rather than bounded
away.  The following estimates implement this reasoning with derivatives
and with every old radius included.

The last old interaction has a stronger adjoint bound. On the early
half-interval, the linear vanishing in
\eqref{trans:adjoint-axis-expansion} and Gaussian integration give
\begin{equation}\label{trans:old-adjoint-bound}
 \|B_{\mathrm{in}}(s,k)^\ast\Psi_s\|_{H^j_h}
 \leq C_jA^2R^{-2}\sum_{n=0}^I
       \ell_n^2(\tau_n+s)^{-1-j/2}
 \leq C_jA^2R^{-2}(d+s)^{-1-j/2}.
\end{equation}
With no derivative on \(\Psi_s\), the vanishing and horizontal
Gaussian volume each contribute a factor \((\tau_n+s)^{1/2}\).
If \(j_1\geq1\) derivatives fall on \(\Psi_s\), the remaining
coefficient is bounded by
\(A^{(j_1+3)/2}(\tau_n+s)^{(j_1-1)/2}\leq A^2\), since
\(\tau_n+s\leq1/A\) on this interval. Projection is a contraction
on horizontal Sobolev spaces. On the late half-interval use
\eqref{trans:low-frequency-coupling}. Thus we obtain
\begin{equation}\label{trans:last-interaction}
 Q_\psi:=\int_0^{t_1}\|B_{\mathrm{in}}(s,k)^\ast\Psi_s\|_2\idiff s
 \leq CR^{-2}\left(A^2\log(2+R)+A^{3/2}|k|\right)
 \leq C_KA^3R^{-2}\log(2+R).
\end{equation}
In the ordered Duhamel expansion, the earlier old interactions
contribute \(B_\ast^{n-1}/(n-1)!\). The outer energy exponentials
multiply over disjoint time intervals, so their total cost remains
bounded. Since \(\tf\geq A\) and \(|k|\leq KA\), one has
\(B_\ast\leq C_KA^{3/2}/R\). The entire adjoint sum with
at least two old interactions has norm at most
\begin{equation}\label{trans:multiple-interactions}
 CQ_\psi(e^{B_\ast}-1)
       \leq C_KA^5R^{-3}\log(2+R).
\end{equation}
Pairing with \eqref{trans:incoming-bounds} leaves
\(CA^{P_{\mathrm{in}}+5}R^{-1}\log(2+R)\). This includes interactions between
different old radii in every order.

The outer response to \(F_k-\mathcal I_k\) has norm at most
\(CA^{(4+\delta)/2}\epsilon_K\) after time \(t_1/2\), by dual
smoothing. Hence the late half of the single old interaction is
bounded by
\begin{equation}\label{trans:late-single-interaction}
 CA^{3/2}|k|R^{-2}
       \left(|k|\sqrt A+A^{(4+\delta)/2}\epsilon_K\right).
\end{equation}
For its early half, we must also retain the genuine outer evolution
before the old interaction. The adjoint error from replacing that
part by free heat is
\begin{equation}\label{trans:incoming-leg-error}
 E_\psi=\int_0^{t_1/2}
       \left(\mathcal U_{\mathrm{out},k}(s,0)^\ast-e^{s\Delta_k}\right)
       B_{\mathrm{in}}(s,k)^\ast\Psi_s\idiff s,
 \qquad \Delta_k=\Delta_h-k^2.
\end{equation}
The parenthesized operator has \(H^j_h\) norm at most \(C_jAs\).
Using \eqref{trans:old-adjoint-bound}, we obtain
\begin{equation}\label{trans:incoming-leg-graphs}
 \|E_\psi\|_{H^2_h}\leq CA^3R^{-2}\log(2+R),\quad
 \|E_\psi\|_{H^4_h}\leq CA^3\tf^{-1},\quad
 \|E_\psi\|_{H^5_h}\leq CA^3R\tf^{-3/2}.
\end{equation}
Interpolation gives
\(\|E_\psi\|_{H^{4+\delta}_h}
\leq CA^3R^\delta\tf^{-1-\delta/2}\).
Pairing separately with \(\mathcal I_k\) and
\(F_k-\mathcal I_k\), and using
\(\|\mathcal I_k\|_{H^{-2}_h}\leq C|k|\), bounds this error by
\begin{equation}\label{trans:incoming-leg-pairing}
 CA^3\left(|k|R^{-2}\log(2+R)+R^\delta\epsilon_K\right).
\end{equation}

\textit{The exact single-interaction pairing.}
It remains to calculate the early free-heat interaction.  Put
$H=R^2/(2A)$, $\nu=k/R$, $L_j=\tf+\ell_j^2$, and
$\Phi_u(r)=R\Psi_{u/R^2}(r/R)$.  Let
$w(u)=e^{u\Delta_h}\widehat\varphi^{\,z}(\cdot,\nu)$, with the
vector horizontal Laplacian.  In cylindrical components this is
$e^{u\Delta_h^{(1)}}$ on $w_r,w_\theta$ and $e^{u\Delta_h^{(0)}}$ on $w_z$.
The contribution of the $j$th old column is exactly
\[
 \frac{2\pi A\ell_j^2}{m}\int_0^H e^{-\nu^2u}
 \int_0^\infty
 \left(2\Omega_{L_j+u}w_\theta(u,r)\overline{\Phi_{u,r}(r)}
       -B_{L_j+u}w_r(u,r)\overline{\Phi_{u,\theta}(r)}\right)
                         r\idiff r\idiff u.
\]
The factors in this formula come from $\idiff s=R^{-2}\idiff u$,
$\idiff y_h=R^{-2}\idiff x_h$, the incoming factor $R^3/m$,
and the factor $R^2$ in the old-column matrix.  The remaining $R$
is included in $\Phi_u$.
{
Since the adjoint test is solenoidal, self-adjointness of $\mathbb P_k$ removes the projection from the pairing.}
The hypotheses give $1\leq L_j\leq H$.

The adjoint expansion becomes
\[
 \Phi_{u,r}(r)=c_r(u/R^2)r+\mathcal E_r(u,r),\qquad
 \Phi_{u,\theta}(r)=c_\theta(u/R^2)r+\mathcal E_\theta(u,r),
 \qquad |\mathcal E_r|+|\mathcal E_\theta|
       \leq CA^2R^{-2}r^3.
\]
This global bound is a Cartesian Taylor estimate on a diameter;
rotational covariance removes the constant and quadratic horizontal
terms.

\textit{The dual radial norms.}
For $0\leq\zeta<7/4$,
\begin{equation}\label{trans:joint-moment-integral}
 \int_\R\frac{|z|^{2\zeta}}{(1+r^2+z^2)^{9/4}}\idiff z
       =C_\zeta(1+r^2)^{\zeta-7/4}.
\end{equation}
Joint Cauchy--Schwarz in $(r,z)$ therefore bounds a pairing against
an arbitrary component of $\widehat\varphi^{\,z}$ by $C M_{9/4}(\varphi)$
times the test norm in $L^2((1+r^2)^{-7/4}r\idiff r)$.
For the odd vertical component, its Fourier transform is a sine
transform.  The inequality $|\sin(\nu z)|\leq|\nu z|^{1/2}$
instead gives the additional factor $|\nu|^{1/2}$ and the dual weight
$(1+r^2)^{-5/4}r\idiff r$.  These are estimates of joint integrals;
no pointwise-in-$r$ bound on the incoming velocity is assumed.

For $T\geq1$ and $0<p_0<2$, splitting at $r=1$ and $r=\sqrt T$
gives
\begin{align}
 \|W_T\|_{L^2((1+r^2)^{-p_0}r\idiff r)}
       &\leq C_{p_0}T^{-1-p_0/2},\notag\\
 \|r^3T^{-2}e^{-r^2/(4T)}\|_{L^2((1+r^2)^{-p_0}r\idiff r)}
       &\leq C_{p_0}T^{-p_0/2}.
       \label{trans:radial-gaussian-norms}
\end{align}
For example, the square of the first norm is bounded by
$CT^{-4}T^{2-p_0}$; the square of the second is bounded by
$CT^{-4}T^{4-p_0}$.
For the radial heat generators $\Delta_h^{(1)}$ and $\Delta_h^{(0)}$
in \eqref{seed:mean-equation}, differentiation of Gaussian
heat evolution gives
\begin{align}
 e^{u\Delta_h^{(1)}}W_T&=W_{T+u},\notag\\
 e^{u\Delta_h^{(1)}}\bigl(r^3T^{-2}e^{-r^2/(4T)}\bigr)
       &=\left(\frac{8Tu}{T+u}+\frac{T^2r^2}{(T+u)^2}\right)W_{T+u},
       \label{trans:radial-heat-identities}\\
 e^{u\Delta_h^{(0)}}\bigl(r^2T^{-2}e^{-r^2/(4T)}\bigr)
       &=\left(\frac{4u}{(T+u)^2}+\frac{Tr^2}{(T+u)^3}\right)
                         e^{-r^2/(4(T+u))}.\notag
\end{align}

These identities can be reproduced without integrating a Bessel kernel.
For the angular-order-one profile,
\[
 r^2W_T=4T^2\partial_TW_T+8TW_T.
\]
Apply $e^{u\Delta_h^{(1)}}$ and use $e^{u\Delta_h^{(1)}}W_T=W_{T+u}$; differentiation in
$T$ commutes with the fixed semigroup.  Since
$\partial_TW_{T+u}=(-2/(T+u)+r^2/(4(T+u)^2))W_{T+u}$,
this gives exactly the second line of
\eqref{trans:radial-heat-identities}.
For the scalar radial heat, put $P_T=T^{-1}e^{-r^2/(4T)}$.
Then $e^{u\Delta_h^{(0)}}P_T=P_{T+u}$ and
\[
 r^2T^{-2}e^{-r^2/(4T)}=4T\partial_TP_T+4P_T.
\]
The same differentiation gives the final line.  The formulas are useful
because the cubic Taylor error and the affine azimuthal term involve
precisely these Gaussian polynomials.  Positivity will permit their use
as majorants even when the actual remainder changes sign.

Positivity of the order-one radial kernel can also be read from ordinary
two-dimensional heat.  Integrating its Gaussian kernel in the angular
variable gives, for $u>0$,
\[
 (e^{u\Delta_h^{(1)}}f)(r)=\frac1{2u}\int_0^\infty
 e^{-(r^2+\rho^2)/(4u)}
 I_1\left(\frac{r\rho}{2u}\right)f(\rho)\rho\,\idiff\rho,
\]
where
$ I_1(z)=\sum_{n=0}^\infty (z/2)^{2n+1}/(n!(n+1)!)\geq0$
for $z\geq0$.  The order-zero kernel has the same expression with
$I_0(z)=\sum_{n=0}^\infty(z/2)^{2n}/(n!)^2$.
These formulas follow by expanding the angular Fourier modes of the
Gaussian and also prove positivity for the regular radial realizations
used here.  Consequently $|e^{u\Delta_h^{(j)}}f|\leq e^{u\Delta_h^{(j)}}|f|$ for $j=0,1$.
This is a scalar-sector assertion; it is not a componentwise positivity
claim for the Leray-projected vector evolution.

Both radial semigroups have nonnegative scalar kernels.  They can
therefore also be applied to absolute Gaussian majorants of the
Taylor remainders.  A fixed polynomial in $r^2/T$ is absorbed by
replacing $T$ with $2T$.  In particular, at $T=L_j+u$, the weighted
norm of the evolved cubic majorant is at most
$C(L_j+u)^{-7/8}$ for $p_0=7/4$.
For $p_0>1$, the final line has weighted norm at most
$C_{p_0}(T+u)^{-1}$: its constant Gaussian term uses the integrability
of $(1+r^2)^{-p_0}r\idiff r$, and its quadratic term is smaller by
\eqref{trans:radial-gaussian-norms} and direct Gaussian integration.

\textit{The leading coefficient.}
In the radial affine term, move $e^{u\Delta_h^{(1)}}$ onto the Gaussian test,
freeze $c_r(u/R^2)$ at $c_r(0)$, replace $e^{-\nu^2u}$ by one,
and extend the time integral to infinity.  Since
\begin{equation}\label{trans:mean-heat-integral}
 \int_0^\infty W_{L_j+2u}(r)\idiff u
       =\frac2r\left(1-e^{-r^2/(4L_j)}\right),
\end{equation}
the sum of these leading terms is
$8\pi A\overline{c_r(0)}m^{-1}
 \int K_{\mathrm{in}}(r)\widehat\varphi_\theta^{\,z}(r,\nu)\idiff r$.
The inner product is linear in its first argument, so
\[
 \langle\mathcal I_k,\Psi_0\rangle
       =\frac{ik}{\sqrt{2\pi}}\overline{\Psi_{0,z}(0)}
       =\frac{2\overline{c_r(0)}}{\sqrt{2\pi}}.
\]
Using the definition of $c_{\mathrm{mean}}(k)$ proves that the leading pairing is
\begin{equation}\label{trans:leading-mean-pairing}
 c_{\mathrm{mean}}(k)\langle\mathcal V_{\mathrm{out},k}(t_1),\psi\rangle.
\end{equation}
In particular, the factor $\sqrt{2\pi}$ and the complex conjugation
are both fixed by the stated unitary Fourier convention.

\paragraph{Moment exponents in the remainder.}
The exponents $7/8$ and $1/8$ below are consequences of the input moment.
The moment $M_{9/4}$ first gives the dual radial
weight $(1+r^2)^{-7/4}$ after axial integration.  The Gaussian test then
has norm of order $T^{-15/8}$.  Its omitted time tail is bounded by
\[
 \int_H^\infty(L+2u)^{-15/8}\,\idiff u\leq C H^{-7/8}.
\]
The time variation of the affine coefficient supplies one factor $u$,
whereas a cubic spatial remainder supplies two extra radial powers.
Their respective integrals are bounded by
$\int_0^H u(L+2u)^{-15/8}\,\idiff u$ and
$\int_0^H(L+2u)^{-7/8}\,\idiff u$, both of order $H^{1/8}$.
Since $H=R^2/(2A)$, either $H^{-7/8}$ or $R^{-2}H^{1/8}$ supplies
$R^{-7/4}$ up to powers of $A$.

The affine azimuthal term has a different mechanism.  Incompressibility
converts a radial derivative into one factor $\nu=k/R$, and oddness of
the vertical input gives an additional $|\nu|^{1/2}$ through a sine
transform estimate.  Its time integral is logarithmic.  It therefore
has size $R^{-3/2}\log(2+R)$, again up to polynomial factors and the
normalized input moment.  These are the two moment-dependent errors in
$E_K$.

\textit{Four error estimates.}
For the radial affine term, the omitted time tail has test norm
at most $CH^{-7/8}$, because
$\int_H^\infty(L_j+2u)^{-15/8}\idiff u\leq CH^{-7/8}$.
The coefficient error on $[0,H]$ is bounded by
\[
 |e^{-\nu^2u}c_r(u/R^2)-c_r(0)|
 \leq uR^{-2}\bigl(k^2\sup|c_r|+\sup|c_r'|\bigr).
\]
The cubic remainders use the evolved Gaussian majorant just computed.
Consequently, after division by $|m|$, the first three errors are
bounded respectively by
\begin{align}
 &CA\ell_j^2|c_r(0)|H^{-7/8}\mathcal W,\notag\\
 &CA\ell_j^2R^{-2}
       \bigl(k^2\sup|c_r|+\sup|c_r'|\bigr)H^{1/8}\mathcal W,\notag\\
 &CA^3\ell_j^2R^{-2}H^{1/8}\mathcal W.
       \label{trans:radial-remainders}
\end{align}
The elementary time integrals used here are
\begin{equation}\label{trans:radial-time-integrals}
 \int_0^H u(L_j+2u)^{-15/8}\idiff u\leq CH^{1/8},\qquad
 \int_0^H(L_j+2u)^{-7/8}\idiff u\leq CH^{1/8}.
\end{equation}

For the fourth error, the affine azimuthal term, use the exact identities
\[
 rB_T=\partial_r(r^2\Omega_T),\qquad
 \partial_r(rw_r)=-i\nu rw_z.
\]
An integration by parts transforms its radial integral to
\[
 \int_0^\infty r^2 B_T w_r\idiff r
       =i\nu\int_0^\infty r^3\Omega_T w_z\idiff r.
\]
Move the $\Delta_h^{(0)}$ heat semigroup onto $r^2\Omega_T$, and apply the
last line of \eqref{trans:radial-heat-identities} with the dual
weight $p_0=5/4$.  Odd vertical parity supplies the additional
$|\nu|^{1/2}$.  The resulting fourth error is bounded by
\[
 CA\ell_j^2\sup|c_\theta|\,|\nu|^{3/2}
                  \log(1+H/L_j)\mathcal W.
\]
Boundary terms in this calculation are first removed by smooth radial
cutoffs.  The displayed dual norms bound their limits and also justify
moving the heat semigroup and interchanging the finite time integrals.
No absolute first moment in $L^1$ of the velocity is used.
Substitute $H=R^2/(2A)$, the slope estimates, and
$\sum_j\ell_j^2\leq4/3$.
The sum of all four errors is at most
\[
 C_KA^4\mathcal W\bigl(R^{-7/4}+R^{-3/2}\log(2+R)\bigr).
\]

Combining \eqref{trans:multiple-interactions},
\eqref{trans:late-single-interaction},
\eqref{trans:incoming-leg-pairing}, and
\eqref{trans:leading-mean-pairing} proves, on \(|k|\leq KA\),
\begin{equation}\label{trans:fiber-entrance}
 \|V_k(t_1)-(1+c_{\mathrm{mean}}(k))\mathcal V_{\mathrm{out},k}(t_1)\|_2
       \leq CA^{P_{\mathrm{in}}+8}E_K.
\end{equation}
Even axial parity permits replacement of \(c_{\mathrm{mean}}(k)\) by \(c_{\mathrm{mean}}(0)\).
Indeed the exact joint integral
\begin{equation}\label{trans:mean-frequency-dual}
 \int_\R\frac{|z|^3}{(1+r^2+z^2)^{9/4}}\idiff z
       =\frac{16}{5}(1+r^2)^{-1/4}
\end{equation}
and \(|\cos(\nu z)-1|\leq2|\nu z|^{3/2}\) imply
\begin{equation}\label{trans:mean-frequency-error}
 |c_{\mathrm{mean}}(k)-c_{\mathrm{mean}}(0)|
 \leq CA|k/R|^{3/2}\mathcal W
       \sum_j\ell_j^2L_j^{-1/8}
 \leq C_KA^{5/2}R^{-3/2}\mathcal W.
\end{equation}
The radial dual integral is bounded by \(CL_j^{-1/4}\) before
taking its square root; split at \(r^2=L_j\) and use
\(1-e^{-r^2/(4L_j)}\leq\min(r^2/(4L_j),1)\).
Integrating \eqref{trans:fiber-entrance} on \(|k|\leq KA\)
therefore costs at most \(CA^{P_{\mathrm{in}}+10}E_K\).

On \(|k|>KA\), retain the exact axial heat. The remaining integrated
matrix norm is bounded by
\begin{equation}\label{trans:high-frequency-matrix}
 2A\sum_j\frac{\ell_j^2}{\tf+\ell_j^2}
       +2At_1\sum_{j\geq0}R^{-2j}\leq C.
\end{equation}
Consequently this part of the full incoming response is at most
\(Ce^{-K^2A}\|F_R\|_2
\leq Ce^{-K^2A}R^{5/2}\mathcal W\).
Equation \eqref{trans:point-upper} bounds the outer point tail by
\(C_KA^{3/2}e^{-K^2A}\).

Finally, the farther outer matrix has \(W^{3,\infty}\) norm at most
\(CAR^{-2}\). In the difference between its point response and
\(\mathcal V_{T_{\mathrm{par}}}\), use \(H^{-3}\) on \((0,t_1/2)\), where
the point norm is bounded and the remaining smoothing costs
\(CA^{3/2}\). On \((t_1/2,t_1)\), use
\(\|\mathcal V_{T_{\mathrm{par}}}(s)\|_2\leq CA^{5/4}\).
Duhamel's formula gives
\begin{equation}\label{trans:outer-point-comparison}
 \|\mathcal V_{\mathrm{out}}(t_1)-\mathcal V_{T_{\mathrm{par}}}(t_1)\|_2
       \leq CA^{3/2}R^{-2}.
\end{equation}
This comparison acts on the point response, where the singular initial
distribution is controlled by smoothing. The full incoming state has
already been treated in \eqref{trans:fiber-entrance}.

The preceding calculations are first made for smooth admissible
incoming fields and finite column sums.  Smooth centrally odd,
axisymmetric, solenoidal fields are dense in the $M_{9/4}$ norm:
truncate radially, average rotations, take the odd part, apply the
bounded odd-field Leray projection, and smooth.
For this approximation, axial Fourier transforms converge uniformly
on bounded frequency intervals in $L^2_h$, because $M_1$ controls
$L^1_zL^2_h$.  Thus the incoming fiber bounds hold with errors tending
to zero in their constants and in $\epsilon_K$.
The weighted dual estimates above also pass every retained-mean
functional to the limit.  Apply the proved inequality and then pass to
the limit in $L^2$ by continuity of the linear evolution.

The old and outer coefficient sums converge in the relevant
$L^1_t$ operator norms, dominated respectively by
$\sum_j\ell_j^2$ and $\sum_jR^{-2j}$.
The same is true for their displayed adjoint and derivative majorants.
Absolute convergence of the Volterra series therefore identifies the
infinite-cutoff limit with the stated evolution, and preserves its
bounds.  No division by $k$ occurred in the proof.  At $k=0$ both
invariant couplings vanish on the solenoidal physical fiber and
$\mathcal I_0=0$, so the estimates include that fiber as well.
Plancherel, together with the separately controlled high-frequency
tail, proves \eqref{trans:invariant-entrance-bound}.
\end{proof}

\subsection{The actual incoming field and the localized evolution}

We now verify that the second variation constructed in Section~\ref{sec:history}
satisfies the abstract incoming assumptions.  Evenness of its complete
stress supplies the negative-Sobolev approximation.  The retained-mean
factor is positive, and comparison with the exact localized background
uses its small correction and one axial moment of the incoming field.

We first express the incoming bounds in terms of the complete stress.
Suppose
\begin{equation}\label{trans:stress-representation}
 \varphi=-\int_\tp^\tf e^{(\tf-t)\Delta}
                    \PP\operatorname{div}\mathsf S(t)\idiff t,
 \qquad \mathsf S(t,-x)=\mathsf S(t,x),
\end{equation}
and suppose its accumulated spatial integral is
\(m\operatorname{diag}(1,1,0)+c_0I_3\), with \(m\neq0\).
For \(1<\alpha<2\), let
\begin{equation}\label{trans:stress-moment}
 \mathfrak M_\alpha=
 \int_\tp^\tf\int_{\R^3}
       \left(|x|^\alpha+(\tf-t)^{\alpha/2}\right)
       |\mathsf S(t,x)|\idiff x\idiff t.
\end{equation}
Evenness, \(|1-\cos y|\leq C_\alpha|y|^\alpha\), and
\(1-e^{-u}\leq u^{\alpha/2}\) give the full Fourier remainder
\(C_\alpha R^{-\alpha}|\xi|^{1+\alpha}
\mathfrak M_\alpha/|m|\) after normalization by
\eqref{trans:incoming-normalization}. The isotropic term \(c_0I_3\)
is annihilated by \(\PP\). Hence, at every real axial frequency,
\begin{equation}\label{trans:negative-sobolev-certificate}
 \|F_k-\mathcal I_k\|_{H^{-4-\delta}_h}
 \leq C_{\alpha,\delta}(1+|k|)^{1+\alpha}
       R^{-\alpha}\frac{\mathfrak M_\alpha}{|m|}.
\end{equation}
Indeed the squared horizontal integral is bounded by
\begin{equation}\label{trans:horizontal-integral}
 (1+k^2)^{1+\alpha}
       \int_{\R^2}(1+|\xi_h|^2)^{-3+\alpha-\delta}\idiff\xi_h
 =\frac{\pi(1+k^2)^{1+\alpha}}{2-\alpha+\delta}.
\end{equation}
This calculation retains all horizontal frequencies and uses only the
fractional moment in \eqref{trans:stress-moment}.

For the actual second variation from Proposition~\ref{hist:generic},
the incoming state is \(\varphi=w_f(\tf)\) and its signed scalar is
\(m_f=-\mathcal M_f<0\).  With this input, the notation in
\eqref{trans:incoming-normalization} and \eqref{trans:mean-scalar} specializes to
\begin{equation}\label{trans:actual-input-identification}
 F_R=F_{f,R},\qquad F_k=F_{f,k},\qquad
 c_{\mathrm{mean}}(k;w_f(\tf),m_f)=c_{\mathrm{mean},f}(k).
\end{equation}
The first two identities follow from \eqref{hist:normalized-input}; the
last follows from \eqref{hist:memory-functional} and \eqref{hist:memory}.
Propositions~\ref{hist:baseline} and~\ref{hist:generic} give
\begin{equation}\label{trans:actual-input-estimates}
 A^{-P}e^{2GA}\leq \mathcal M_f\leq CA^Pe^{2c_\star A},\qquad
 \mathcal W_f+\frac{\mathfrak M_\alpha}{\mathcal M_f}
       \leq CA^Pe^{2\omega A},\qquad
 \sup_{|k|\leq KA}\|F_k\|_2\leq C_KA^2R^2,
\end{equation}
where \(\mathcal W_f=M_{9/4}(w_f(\tf))/\mathcal M_f\) and
\(\omega=c_\star-G\). They also give, for some fixed \(c_{\mathrm{mean,err}}>0\),
\begin{equation}\label{trans:main-mean-bound}
 \left|\frac{4\pi A}{m_f}
       \int_0^\infty\left(1-e^{-r^2/(4(1+\tf))}\right)
       \int_\R w_{f,\theta}(\tf,r,z)\idiff z\idiff r\right|
 \leq\frac CA+CA^Pe^{-c_{\mathrm{mean,err}}A}.
\end{equation}

The additional old radii in \eqref{trans:mean-scalar} have an absolute
bound. For \(L\geq1\), joint Cauchy--Schwarz gives
\begin{equation}\label{trans:old-mean-dual}
 \left|\int_0^\infty\int_\R
       \left(1-e^{-r^2/(4L)}\right)h_\theta(r,z)
                         \idiff z\idiff r\right|
       \leq CL^{-1/4}M_1(h).
\end{equation}
The squared dual integral is
\(\pi\int_0^\infty(1-e^{-r^2/(4L)})^2
       (r\sqrt{1+r^2})^{-1}\idiff r\leq CL^{-1/2}\).
Since \(\sum_{j\geq1}\ell_j^2\leq CR^{-2}\), their contribution to
\(c_{\mathrm{mean},f}(0)\) is at most
\(CAR^{-2}\tf^{-1/4}\mathcal W_f\).
Equations \eqref{trans:actual-input-estimates} and
\eqref{trans:main-mean-bound} therefore imply
\begin{equation}\label{trans:positive-mean-factor}
 a_f:=1+c_{\mathrm{mean},f}(0)=1+O(A^{-1})>0,
\end{equation}
because \(2\omega<2\varrho\). Thus the retained mean changes the
leading coefficient by a controlled real factor.

We now compare the invariant evolution with the actual one. In parent
coordinates write
\begin{equation}\label{trans:actual-fields}
 U_{\mathrm{act}}(s,y)=R\bar u^{I,J}(\tf+R^2s,Ry),\qquad
 E(s,y)=Re^{I,J}(\tf+R^2s,Ry),
\end{equation}
and let \(U_{\mathrm{inv}}=U_{\mathrm{in}}+U_{\mathrm{out}}\) be
\eqref{trans:invariant-mixture}. Their homogeneous propagators are
\(\mathcal U_{\mathrm{act}}\) and
\(\mathcal U_{\mathrm{inv}}\), respectively. The actual heat columns
have radii \(R^{j-1}\), ages \(R^{2j-2}+d+s\), axial widths
\(bR^{j-1}\), and elapsed axial heat time \(d+s\).

\begin{lemma}[Comparison with the localized reference]\label{trans:actual-comparison}
For \(F_R\) as in \eqref{trans:incoming-normalization}, one has
\begin{equation}\label{trans:actual-comparison-bound}
 \sup_{0\leq s\leq s_\ast}
       \|\mathcal U_{\mathrm{act}}(s,0)F_R
                    -\mathcal U_{\mathrm{inv}}(s,0)F_R\|_2
 \leq CA^P\mathcal W
       \left(R^{7/2}b^{-1}+R^{5/2}e^{C_\ast A}b^{-1/2}\right),
\end{equation}
where \(C_\ast\) is fixed before \(\beta\). Moreover,
\begin{equation}\label{trans:homogeneous-polynomial}
 \|\mathcal U_{\mathrm{act}}(s,u)\|_{2\to2}
       +\|\mathcal U_{\mathrm{inv}}(s,u)\|_{2\to2}
       \leq CA^P,
 \qquad 0\leq u\leq s\leq s_\ast.
\end{equation}
\end{lemma}

\begin{proof}
The actual and invariant heat matrices, including the axial shear,
are bounded by
\begin{equation}\label{trans:comparison-matrix}
 a(s)=CA\sum_{j=-I}^J
          \frac{R^{2j-2}}{(R^{2j-2}+d+s)^2}.
\end{equation}
For the old radii, the time integral is at most
\(CA\sum_j\ell_j^2/(\tf+\ell_j^2)\leq CA/\tf\).
The parent column contributes \(CAs_\ast\), and the farther columns
contribute \(CAs_\ast\sum_{j\geq1}R^{-2j}\).
Consequently \(\int_0^{s_\ast}a(s)\idiff s\leq C(1+\log A)\).
The velocity has the separate bound
\begin{equation}\label{trans:invariant-velocity-integral}
 \int_0^{s_\ast}\|U_{\mathrm{inv}}(s)\|_\infty\idiff s
 \leq\frac{CA}{R\sqrt\tf}+CAs_\ast
                   +CAs_\ast\sum_{j\geq1}R^{-j}
 \leq C(1+\log A).
\end{equation}
Proposition~\ref{bg:reference} gives
\begin{equation}\label{trans:correction-integrals}
 \|E(s)\|_\infty\leq\sigma_0(d+s)^{-1/2},\qquad
 \|\nabla E(s)\|_\infty\leq\sigma_1(d+s)^{-1},\qquad
 \sigma_0+\sigma_1\leq CA^Pe^{C_\ast A}b^{-1/2}.
\end{equation}
Since \(\log(1+s_\ast/d)\leq C\varrho A\), its integrated
strain is bounded after the fixed width choice. Correction transport
cancels in the ordinary energy because \(\operatorname{div}E=0\).
These observations prove \eqref{trans:homogeneous-polynomial}.

Let \(V=\mathcal U_{\mathrm{inv}}(s,0)F_R\).
Its energy controls \(\sup_s\|V(s)\|_2\) and
\(\|\nabla V\|_{L^2_sL^2_y}\) by \(CA^P\|F_R\|_2\).
The exact axial moment equation is
\begin{align}
 (\partial_s-\Delta)(y_3V)
 &=\PP M_{\mathrm{inv}}(y_3V)-2\partial_{y_3}V\notag\\
 &\quad-[y_3,\PP]\operatorname{div}
       \left(U_{\mathrm{inv}}\otimes V+V\otimes U_{\mathrm{inv}}\right).
 \label{trans:axial-moment-equation}
\end{align}

The projection term in this moment equation should be regarded as one
composed operator.  In full Fourier variables $\xi$, multiplication by
$y_3$ becomes $i\partial_{\xi_3}$, so the symbol of
$[y_3,\PP]\operatorname{div}$ is assembled from
$(i\partial_{\xi_3}P_{ij}(\xi))(i\xi_\ell)$.
The derivative of the degree-zero projection has degree minus one;
the divergence restores degree zero.  The resulting matrix symbol is
bounded on the unit sphere and hence, by homogeneity and Plancherel,
\[
 \|[y_3,\PP]\operatorname{div}F\|_2\leq C\|F\|_2.
\]
Its value at the single point $\xi=0$ is immaterial.  The commutator
$[y_3,\PP]$ alone need not have this $L^2$ bound, so separating the
divergence from it would lose the estimate.  With
$F=U_{\mathrm{inv}}\otimes V+V\otimes U_{\mathrm{inv}}$, the source
is bounded by $C\|U_{\mathrm{inv}}\|_\infty\|V\|_2$.
The other new term, $-2\partial_{y_3}V$, comes directly from commuting
$y_3$ through the Laplacian and is integrable using the gradient energy.

The smooth matrix is defined on arbitrary axisymmetric Cartesian
fields, so \(y_3V\) need not be solenoidal. The composed multiplier
\([y_3,\PP]\operatorname{div}\) has order zero and bounded symbol.
Its source is bounded by
\(C\|U_{\mathrm{inv}}\|_\infty\|V\|_2\).
The heat commutator is integrated by Cauchy--Schwarz using the gradient
energy. Equations \eqref{trans:comparison-matrix} and
\eqref{trans:invariant-velocity-integral} therefore imply
\begin{equation}\label{trans:axial-moment-bound}
 \sup_s\|y_3V(s)\|_2
       \leq CA^P\left(\|F_R\|_2+\|y_3F_R\|_2\right).
\end{equation}
Axial difference quotients justify
\eqref{trans:axial-moment-equation} for this input class.

For every width \(w>0\), its Gaussian heat cutoff satisfies
\begin{equation}\label{trans:cutoff-elementary}
 |1-\eta_w(\tau,y_3)|\leq C(\sqrt\tau+|y_3|)/w,
 \qquad \|\partial_{y_3}\eta_w(\tau)\|_\infty\leq C/w.
\end{equation}
Separating amplitude loss and the Gaussian factor proves the first
bound wherever its right side is small; the trivial bound proves it
elsewhere. If \(M_Q\) is the actual heat matrix, these
inequalities give
\begin{align}
 \|(M_Q-M_{\mathrm{inv}})V\|_2
 &\leq\frac{CA}{b}\sum_{j=-I}^JR^{j-1}
       \left((R^{2j-2}+d+s)^{-3/2}\|V\|_2\right.\notag\\
 &\hspace{44mm}\left.
                    +(R^{2j-2}+d+s)^{-2}\|y_3V\|_2\right).
 \label{trans:complete-cutoff-source}
\end{align}
This includes the axial shear and is valid for arbitrarily small old
widths. The two old-radius time integrals are at most
\begin{equation}\label{trans:cutoff-old-integrals}
 2\sum_j\frac{\ell_j}{\sqrt{\tf+\ell_j^2}}
       \leq C\tf^{-1/2},\qquad
 R\sum_j\frac{\ell_j}{\tf+\ell_j^2}\leq CR/\tf.
\end{equation}
The parent and farther contributions have bounded geometric
integrals. Thus the vector source in
\eqref{trans:complete-cutoff-source} has \(L^1_sL^2_y\) norm at most
\(CA^P(R/b)(\|F_R\|_2+\|y_3F_R\|_2)\).

The reference-correction source consists of both
\(-\PP(E\cdot\nabla V)\) and \(-\PP(V\cdot\nabla E)\).
For transport, use the product of
\(\|E\|_{L^2_sL^\infty_y}\) and
\(\|\nabla V\|_{L^2_sL^2_y}\); the first squared norm is at most
\(\sigma_0^2\log(1+s_\ast/d)\).
For stretching, use
\(\sigma_1\log(1+s_\ast/d)\sup_s\|V(s)\|_2\).
Their sum is bounded by
\(CA^Pe^{C_\ast A}b^{-1/2}\|F_R\|_2\).
Apply \eqref{trans:homogeneous-polynomial} to the exact zero-data
difference equation. Finally the scaling in
\eqref{trans:incoming-normalization} gives
\(\|F_R\|_2\leq R^{5/2}\mathcal W\) and
\(\|y_3F_R\|_2\leq R^{3/2}\mathcal W\), proving
\eqref{trans:actual-comparison-bound}.
\end{proof}

\subsection{Late propagation and tensor sources}

When the perturbation is expanded in powers of its initial amplitude,
the $m$th coefficient is driven by products of lower coefficients.
A bound with a constant depending uncontrollably on $m$ is not enough
to sum this expansion.  The source estimate here retains a common
linear growth factor and integrates the resulting forcing without
differentiating the $m$-dependent exponential.  The resulting
polynomial-in-$A$ constant is independent of $m$, so the familiar
Catalan majorant can later be applied to the full nonlinear series.

The second variation continues to be forced after $\tf$.  That source
cannot be dropped merely because its earlier history dominates.  The
first variation has already decayed to a polynomial bound, and the
estimates below show that its entire later quadratic response is small
relative to the positive stress mass after parent-scale normalization.

On the physical interval \([\tf,\ts]\), use the growth integral
\begin{equation}\label{trans:late-growth}
 \Gamma_{\mathrm{tr}}(t)=2AR^2
       \left((R^2+\tf)^{-1}-(R^2+t)^{-1}\right).
\end{equation}
The coefficient two is the physical matrix bound in
\eqref{trans:physical-coupling}. Its larger value causes only a
polynomial cost, because
\begin{equation}\label{trans:late-growth-bounds}
 \Gamma_{\mathrm{tr}}'(t)\geq\frac{A}{2R^2},\qquad
 \Gamma_{\mathrm{tr}}(\ts)\leq2A\tp\leq C(1+\log A).
\end{equation}

\begin{lemma}[Late propagation and tensor sources]\label{trans:late}
Let \(\mathcal U(t,u)\) be the homogeneous evolution about
\(\bar u^{I,J}\) for \(\tf\leq u\leq t\leq\ts\).
For each fixed integer \(j\geq0\),
\begin{equation}\label{trans:late-homogeneous}
 \|\mathcal U(t,u)f\|_{H^j}
       \leq C_je^{\Gamma_{\mathrm{tr}}(t)-\Gamma_{\mathrm{tr}}(u)}
                    \|f\|_{H^j}.
\end{equation}
Suppose \(v_F(\tf)=0\) and
\(\partial_t v_F-\mathcal L_{\bar u^{I,J}}v_F=-\PP\operatorname{div}F\).
If \(n\geq2\) and \(\|F(t)\|_2\leq B_0e^{n\Gamma_{\mathrm{tr}}(t)}\),
then
\begin{equation}\label{trans:late-tensor-l2}
 \|v_F(t)\|_2\leq\frac{CR}{\sqrt{(n-1)A}}
                      B_0e^{n\Gamma_{\mathrm{tr}}(t)}.
\end{equation}
If also \(\|F(t)\|_{H^j}\leq B_je^{n\Gamma_{\mathrm{tr}}(t)}\),
then
\begin{equation}\label{trans:late-tensor-hj}
 \|v_F(t)\|_{H^j}\leq C_j\frac{R}{\sqrt A}
                  (B_0+B_j)e^{n\Gamma_{\mathrm{tr}}(t)}.
\end{equation}
The constants are independent of \(n,I,J,A\). The finite derivative
orders needed in these estimates are fixed before \(\beta\).
\end{lemma}

\begin{proof}
In the physical energy, the radius-\(R\) column is bounded by
\(\Gamma_{\mathrm{tr}}'\), except for its axial shear. All remaining
coefficients have the majorant
\begin{equation}\label{trans:late-remainder-coefficient}
 b_{\mathrm{rem}}(t)=C\left(At^{-2}+AR^{-4}
                         +\frac{A}{bR^2}+\frac{\sigma_1}{t}\right).
\end{equation}
The old columns use \(\sum_{\ell\leq1}\ell^2\leq C\), the farther
columns use \(\sum_{\ell\geq R^2}\ell^{-2}\leq CR^{-4}\), and
the parent shear is at most \(CA/(bR^2)\).
Correction transport cancels by incompressibility. Moreover,
\begin{equation}\label{trans:late-remainder-integral}
 \int_\tf^\ts b_{\mathrm{rem}}(t)\idiff t
 \leq C\left(A/\tf+AR^{-2}\tp+A\tp/b
                         +\sigma_1\log(\ts/\tf)\right)\leq C.
\end{equation}
Pair the tensor directly with \(\nabla v_F\), and absorb part of the
heat form. Pressure pairs to zero against the solenoidal velocity.
This gives
\begin{equation}\label{trans:late-energy}
 \frac{\diff}{\diff t}\|v_F\|_2^2+\|\nabla v_F\|_2^2
 \leq2\left(\Gamma_{\mathrm{tr}}'+b_{\mathrm{rem}}\right)\|v_F\|_2^2
                       +C\|F\|_2^2.
\end{equation}
For zero data, the integrating factor and
\eqref{trans:late-growth-bounds} yield
\begin{align}
 \|v_F(t)\|_2^2
 &\leq CB_0^2e^{2n\Gamma_{\mathrm{tr}}(t)}
       \int_\tf^t
        e^{-2(n-1)(\Gamma_{\mathrm{tr}}(t)-\Gamma_{\mathrm{tr}}(u))}
                                                     \idiff u\notag\\
 &\leq\frac{CR^2}{(n-1)A}B_0^2e^{2n\Gamma_{\mathrm{tr}}(t)}.
 \label{trans:late-integrating-factor}
\end{align}
This proves \eqref{trans:late-tensor-l2} and retains the corresponding
weighted gradient energy.

All fixed Cartesian matrix derivatives are uniformly bounded for
\(t\geq\tf\): old columns have size \(C_jAt^{-2-j/2}\), while the
parent and outer columns have size \(C_jAR^{-2-j}\).
The required correction derivatives are also bounded after the fixed
width choice. Cartesian derivative energies and successive time
weights therefore give uniform \(L^2\)-to-\(H^j\) smoothing on a
final interval of length one, and the tensor-source kernel
\(C_j(t-u)^{-1/2}\) in \(H^j\).
For correction transport, integration by parts cancels its highest
derivative term; the remaining commutators are controlled by heat and
the bounded coefficient derivatives.

For elapsed time at most one, these are local estimates. For longer
homogeneous evolution, apply the \(L^2\) energy up to the last unit
interval and then smooth. This proves
\eqref{trans:late-homogeneous}. For the source equation, smooth the
earlier state using \eqref{trans:late-tensor-l2} and integrate the
last-unit tensor kernel. Monotonicity of
\(\Gamma_{\mathrm{tr}}\) bounds each source factor by its value at
\(t\) without a constant depending on \(n\).
Since \(R/\sqrt A\geq1\), this gives
\eqref{trans:late-tensor-hj}.
\end{proof}

One consequence will be used repeatedly. Suppose \(v_f(\tp)=f\),
\(\|f\|_2\leq C\), and \(M_s(f)\leq Ce^{\gamma A}\) for a fixed
\(1<s<5/2\). Lemma~\ref{dec:endpoint} and
\eqref{trans:late-homogeneous} give, after increasing \(\beta\),
\begin{equation}\label{trans:continued-first-variation}
 \sup_{\tf\leq t\leq\ts}\|v_f(t)\|_{H^4}\leq CA^P.
\end{equation}
Thus the entire continuing quadratic source satisfies
\begin{equation}\label{trans:continuing-quadratic-source}
 \int_\tf^\ts
       \left(\|v_f\otimes v_f\|_1+\|v_f\otimes v_f\|_{H^4}
              +\|\PP\operatorname{div}(v_f\otimes v_f)\|_{H^3}\right)
                                                       \idiff t
 \leq CA^PR^2.
\end{equation}
The corresponding vector-source response has the same
\(CA^PR^2\) bound in physical \(L^2\), by
\eqref{trans:late-homogeneous}.

\subsection{Smoothing in parent coordinates and exact time alignment}

An $L^2$ transfer estimate is not automatically an $H^4$ transfer estimate:
rescaling derivatives costs additional powers of $R$.  We gain derivatives
only after evolving on a final interval measured in parent coordinates.
The second lemma then aligns the actual reference and the prepared point
response at exactly the time needed for the next transition.

The comparison must hold in \(H^4\) after rescaling, since that is the
norm used to define the shape class. Smoothing the actual difference
equation on its final interval preserves the already proved
\(L^2\) gain.

For a velocity rescaling $H(y)=Rh(Ry)$, direct change of variables gives
\[
 \|H\|_2=R^{-1/2}\|h\|_2,
 \qquad \|\partial_y^4H\|_2=R^{7/2}\|\partial_x^4h\|_2.
\]
A small physical $L^2$ error therefore says nothing by itself about the
rescaled fourth derivatives.  We instead evolve the \emph{difference}
equation on a final interval in the parent variables.  At that time all
relevant coefficient derivatives in those variables are polynomial in
$A$, rather than carrying uncontrolled powers of $R$.

For an interval of length $h=A^{-1}$ and an error satisfying
$\partial_sZ=\mathcal L_{\mathrm{act}}Z+H_0-\PP\operatorname{div}F_0$,
Duhamel separates the earlier state from recent sources:
\[
 \begin{aligned}
 Z(s)={}&\mathcal U_{\mathrm{act}}(s,s-h)Z(s-h)
       +\int_{s-h}^s\mathcal U_{\mathrm{act}}(s,u)H_0(u)\,\idiff u\\
 &-\int_{s-h}^s\mathcal U_{\mathrm{act}}(s,u)
                              \PP\operatorname{div}F_0(u)\,\idiff u.
 \end{aligned}
\]
The first term uses the already small $L^2$ error and positive-time
smoothing.  The vector source is estimated in $H^4$.  For an $H^4$
tensor source, the one derivative in its divergence is supplied by heat
with kernel $(s-u)^{-1/2}$, whose integral over the last interval is
$2\sqrt h$.  There is no need to demand five derivatives of that tensor.
This procedure upgrades the gain already proved; it does not attempt to
derive a fresh global derivative estimate over the entire long interval.

\begin{lemma}[Smoothing in parent coordinates]\label{trans:smoothing}
Set \(\mathcal J_2=(s_\ast-2/A,s_\ast)\) and
\(\mathcal J_1=(s_\ast-1/A,s_\ast)\). For \(s\in\mathcal J_1\),
the actual parent-coordinate evolution satisfies
\begin{equation}\label{trans:parent-smoothing}
 \|\mathcal U_{\mathrm{act}}(s,s-A^{-1})\|_{L^2\to H^4}\leq CA^P,
 \qquad
 \|\mathcal U_{\mathrm{act}}(s,u)\PP\operatorname{div}F\|_{H^4}
       \leq CA^P(s-u)^{-1/2}\|F\|_{H^4}
\end{equation}
for \(s-A^{-1}\leq u<s\). The aged point response
\(\Psi(s)=\mathcal V_{1+d}(s)\) satisfies, on \(\mathcal J_2\),
\begin{equation}\label{trans:aged-point-graphs}
 \|\Psi(s)\|_{H^j}+\|y_3\Psi(s)\|_{H^j}\leq C_jA^{P_j}
\end{equation}
at every fixed required order. If \(B_{\mathrm{main}}\) is its
single-column operator, then
\begin{equation}\label{trans:reference-coefficient-error}
 \|(B_{\mathrm{act}}-B_{\mathrm{main}})\Psi(s)\|_{H^4}
 \leq CA^P\left(R^{-2}+b^{-1}+e^{C_\ast A}b^{-1/2}\right).
\end{equation}
The estimates extend to the closed final intervals by continuity.
\end{lemma}

\begin{proof}
All fixed derivatives of the heat-part matrix on \(\mathcal J_2\)
are at most \(C_jA\). The parent age is comparable to one; the old
derivative sum is at most \(C_jAR^{-2}s^{-2-j/2}\), and the farther
sum at most \(C_jAR^{-2}\). Since \(s\asymp\tp\), both are bounded
by \(C_jA\). The correction derivatives satisfy
\begin{equation}\label{trans:parent-correction-derivatives}
 \|\nabla_y^jE(s)\|_\infty
       \leq\sigma_j(d+s)^{-1/2-j/2},\qquad
 \sigma_j\leq C_jA^{P_j}e^{C_jA}b^{-1/2}.
\end{equation}
They are small for the finite derivative list after choosing
\(\beta\). Rescale an interval of length \(A^{-1}\) by the spatial
factor \(A^{-1/2}\). The matrix, drift, and their required rescaled
derivatives are uniformly bounded. The Cartesian derivative heat
energies give \eqref{trans:parent-smoothing}. This argument
differentiates the valid Cartesian equation and retains both the
Leray projection and correction transport.

To prove the weighted point estimates, begin in negative Sobolev
spaces. The distribution \(\mathcal I\) belongs to \(H^{-3}\),
and \(y_3\mathcal I\) belongs to \(H^{-4}\), since its Fourier
symbol has degree zero. Equation
\eqref{trans:axial-moment-equation} applies in \(H^{-4}\).
The heat commutator is bounded by the \(H^{-3}\) norm, the projection
commutator has order zero, and multiplication by the invariant
Gaussian velocity costs \(CA\) in these negative spaces.
The two norms grow by at most \(e^{CA\tp}\), a polynomial in \(A\).
Smoothing on a positive final fraction of the interval proves
\eqref{trans:aged-point-graphs}.

The old and farther columns in
\eqref{trans:reference-coefficient-error} are controlled by their
complete derivative sums and \eqref{trans:aged-point-graphs}.
For the parent cutoff use
\begin{equation}\label{trans:parent-cutoff-bound}
 |1-\eta_b(d+s,y_3)|\leq C\left((d+s)/b^2+|y_3|/b\right)
\end{equation}
and its differentiated forms. The first weighted graph norm in
\eqref{trans:aged-point-graphs} controls these products.
The correction contributes both \(E\cdot\nabla\Psi\) and
\(\Psi\cdot\nabla E\). Equation
\eqref{trans:parent-correction-derivatives} and the \(H^5\) point
bound give the last term in
\eqref{trans:reference-coefficient-error}, proving the lemma.
\end{proof}

\begin{lemma}[Alignment of the heat ages]\label{trans:alignment}
At the physical endpoint \(\ts=R^2\tp\), the exact reference
solution satisfies
\begin{equation}\label{trans:reference-alignment}
 R\bar u^{I,J}(\ts,Ry)=\bar u^{I+1,J-1}(\tp,y).
\end{equation}
For the isolated propagator with continuing heat age \(1+t\), one
also has, for each fixed integer \(j\geq0\),
\begin{equation}\label{trans:point-alignment}
 \|\mathcal V_{1+d}(s_\ast)-S_A\|_{H^j}
       \leq C_jA^{P_j}d\leq C_jA^{P_j}R^{-2}.
\end{equation}
\end{lemma}

\begin{proof}
Equation \eqref{trans:reference-alignment} is the covariance
\eqref{bg:covariance} at the time specified in
\eqref{trans:times}. For the point response,
\(\mathcal V_{1+d}(s_\ast)=\mathcal U_{\mathrm d}(\tp,d)\mathcal I\).
Let \(L_{\mathrm d}(s)=\Delta-\PP\mathsf M_A(s)\), where
\(\mathsf M_A\) is the matrix in \eqref{seed:cartesian-linearization}.
The matrix derivatives are bounded by \(C_jA\) on \((0,\tp)\),
and \(\mathcal I\in H^{-4}\), \(\Delta\mathcal I\in H^{-6}\).
Consequently
\(\sup_{0\leq s\leq d}\|L_{\mathrm d}(s)\mathcal I\|_{H^{-6}}
\leq CA\).
The remaining evolution time is at least \(\tp/2\), so
\(\|\mathcal U_{\mathrm d}(\tp,s)\|_{H^{-6}\to H^j}
\leq C_jA^{P_j}\). Differentiation in the initial time gives
\begin{equation}\label{trans:initial-time-difference}
 \mathcal U_{\mathrm d}(\tp,d)\mathcal I
       -\mathcal U_{\mathrm d}(\tp,0)\mathcal I
 =-\int_0^d\mathcal U_{\mathrm d}(\tp,s)
                         L_{\mathrm d}(s)\mathcal I\idiff s.
\end{equation}
One first proves this identity for a smooth regularization; the
negative-space bound just stated permits its removal. Integration
proves \eqref{trans:point-alignment}, with the polynomial factor
\(\tf\) absorbed into \(A^{P_j}\).
\end{proof}

\subsection{The complete quadratic output}

We now compare the transfer errors using the same positive stress mass.
The leading coefficient includes the retained mean.  We first propagate
the already formed second variation and then add the continuing quadratic
source, preserving a common normalization.  The conclusion is a strong
shape estimate with an exponential margin, uniform over all admissible
inputs and both cutoffs.

We normalize every contribution by the positive stress mass $\mathcal M_f$ of the actual
input and by the velocity-transfer factor $R^{-3}$.  The retained mean
is included in the scalar coefficient $a_f=1+O(A^{-1})$; the remaining
error must be exponentially small in the parent $H^4$ norm.
This separation is used in the normalization of
Proposition~\ref{nl:transition}.

The following quantities summarize the contributions to the normalized
output.  Fixed powers of $A$ and constants depending only on the already
chosen auxiliary parameters are suppressed in this table, but not in the
proof.  The input moment and fractional-stress errors obey
$\mathcal W_f+\mathfrak M_\alpha/\mathcal M_f\leq CA^Pe^{2\omega A}$.
\begin{center}
\small
\renewcommand{\arraystretch}{1.15}
\begin{tabular}{@{}p{0.43\linewidth}p{0.48\linewidth}@{}}
\toprule
Origin of the error & Size or control after normalization\\\midrule
Point approximation and the early incoming leg &
 $R^{-(\alpha-\delta)}e^{2\omega A}$\\
At least two old-column interactions & $R^{-1}\log(2+R)$\\
Retained-mean frequency error &
 $R^{-3/2}\log(2+R)e^{2\omega A}$\\
Axial localization and reference correction &
 $R^{7/2}b^{-1}e^{2\omega A}
   +R^{5/2}e^{C_\ast A}b^{-1/2}e^{2\omega A}$\\
Continuing quadratic source &
 $R^{9/2}/\mathcal M_f$ in $L^2$; $R^{15/2}/\mathcal M_f$ for final $H^4$ control\\
Exact heat-age alignment & $R^{-2}$\\\bottomrule
\end{tabular}
\end{center}
The high-frequency tail is separately controlled by $e^{-K^2A}$ times
the incoming norm.  The coefficient $a_f=1+O(A^{-1})$ remains in the leading term.
The table distinguishes smallness from scale separation, from width, and
from division by the large positive stress mass.  These mechanisms are fixed in
that order before taking the cascade depth to infinity.

For $t\geq\tf$, denote the homogeneous continuation by
\begin{equation}\label{trans:homogeneous-second-variation}
 w_f^{\mathrm{hom}}(t):=\mathcal U_{\bar u}(t,\tf)w_f(\tf).
\end{equation}
It agrees with $w_f$ at $\tf$, but has no quadratic forcing afterward.
The difference $w_f-w_f^{\mathrm{hom}}$ is the zero-data response to
$-\PP\operatorname{div}(v_f\otimes v_f)$ on $[\tf,t]$, estimated in
\eqref{trans:continuing-quadratic-source}.

\begin{proposition}[Transfer of the complete quadratic response]\label{trans:quadratic}
Let \(s_\star=19/8\), \(\gamma=\varrho/4\), and
\(g=S_A/n_A\). Suppose \(f\) is real, centrally odd, axisymmetric,
and solenoidal, with
\begin{equation}\label{trans:generic-shape}
 \langle f,g\rangle=1,\qquad
 \|f-g\|_{H^4}\leq e^{-3\varrho A/8},\qquad
 M_{s_\star}(f-g)\leq e^{\gamma A}.
\end{equation}
Let \(v_f,w_f\) be the actual first and complete second variations
about \(\bar u^{I,J}\), with \(v_f(\tp)=f\) and \(w_f(\tp)=0\),
and let \(m_f=-\mathcal M_f<0\) be their signed complete-stress integral
through \(\tf\). Fix \(G\) so that
\(0<2(c_\star-G)<\varrho/32\).
For sufficiently large fixed \(\beta\), there is a real
\(a_f=1+O(A^{-1})>0\), uniform in \(f,I,J\), such that, for every
fixed \(0<\kappa<\varrho\),
\begin{equation}\label{trans:quadratic-aligned}
 \left\|\frac{R^4}{m_f}w_f(\ts,R\cdot)-a_fS_A\right\|_{H^4}
       \leq C_\kappa A^Pe^{-\kappa A}.
\end{equation}
Before replacing the aged point response by \(S_A\), the corresponding
\(L^2\) estimate holds uniformly for \(A^{-1}\leq s\leq s_\ast\),
and its \(H^4\) version holds uniformly for
\(s_\ast-A^{-1}\leq s\leq s_\ast\), with
\(w_f(\ts,R\cdot)\) replaced by
\(w_f(\tf+R^2s,R\cdot)\) and \(S_A\) by
\(\mathcal V_{1+d}(s)\).
The same conclusions, with the same \(m_f,a_f\), hold for the
homogeneous continuation \(w_f^{\mathrm{hom}}\) defined in
\eqref{trans:homogeneous-second-variation}.
\end{proposition}

\begin{proof}
Propositions~\ref{hist:energy}, \ref{hist:baseline}, and
\ref{hist:generic} apply to \eqref{trans:generic-shape}.
Lemma~\ref{seed:prepared} gives
\(M_{s_\star}(g)\leq CA^P\), so
\(M_{s_\star}(f)\leq CA^P+e^{\gamma A}\), as required in
\eqref{trans:continued-first-variation}.
Use \(\alpha=7/4\) and \(\delta=1/4\) in
\eqref{trans:negative-sobolev-certificate}.
The fractional moment is admissible because
\(\alpha<s_\star-1/2=15/8\).
Equations \eqref{trans:actual-input-estimates} and
\eqref{trans:positive-mean-factor} therefore apply to
\(\varphi=w_f(\tf)\).

Lemma~\ref{trans:invariant-entrance} gives the homogeneous invariant
output at \(s=A^{-1}\). Up to fixed polynomial factors, its error
rates are
\begin{equation}\label{trans:entrance-rates}
 \varrho,\qquad
 (\alpha-\delta)\varrho-2\omega,\qquad
 \tfrac32\varrho-2\omega,\qquad
 K^2-\tfrac52\varrho-2\omega,
 \qquad \omega=c_\star-G.
\end{equation}
The \(R^{-7/4}\mathcal W_f\) and \(R^{-2}\) terms have larger
rates. After time \(A^{-1}\), the integrated matrix norm of the
old and farther columns is at most \(CA^PR^{-2}\): integrate
\(CAR^{-2}s^{-2}\) for the old columns and
\(CAR^{-2}\) for the farther columns.
The point response has polynomial norm by
\eqref{trans:point-upper}. Its exact difference equation thus
continues the estimate to every \(s\leq s_\ast\).

Lemma~\ref{trans:actual-comparison} adds the error rates
\begin{equation}\label{trans:localization-rates}
 \beta-\tfrac72\varrho-2\omega,
 \qquad \beta/2-C_\ast-\tfrac52\varrho-2\omega.
\end{equation}
The complete second variation also includes its quadratic source
after \(\tf\). By \eqref{trans:continuing-quadratic-source}, its
physical \(L^2\) response is at most \(CA^PR^2\).
Normalization by \(R^4/m_f\) contributes the \(L^2\) factor
\(R^{5/2}/\mathcal M_f\), so this term is bounded by
\begin{equation}\label{trans:continuing-source-l2}
 CA^P\frac{R^{9/2}}{\mathcal M_f}
       \leq CA^P e^{-(2G-9\varrho/2)A}.
\end{equation}
These estimates are uniform for every earlier parent time in
\((A^{-1},s_\ast)\).
For our \(\alpha,\delta\), the two middle rates in
\eqref{trans:entrance-rates} exceed \(\varrho\).
Choose \(K\) so that its rate also exceeds \(\varrho\), then
choose \(\beta\) so that both rates in
\eqref{trans:localization-rates} do. The rate in
\eqref{trans:continuing-source-l2} is larger than \(\varrho\).
This proves the uniform aged \(L^2\) assertion.

For the \(H^4\) assertion, put
\begin{equation}\label{trans:rescaled-variations}
 V(s,y)=Rv_f(\tf+R^2s,Ry),\qquad
 W(s,y)=\frac{R^4}{m_f}w_f(\tf+R^2s,Ry),\qquad
 Z=W-a_f\mathcal V_{1+d}.
\end{equation}

The factor in the next tensor source can be checked directly.
Multiplying the $w_f$ equation by $R^6/m_f$ accounts for the $R^4/m_f$
normalization and the time factor $R^2$.  A physical divergence contributes
$R^{-1}\operatorname{div}_y$, and
$v_f\otimes v_f=R^{-2}V\otimes V$.  Hence the coefficient is
$R^6R^{-1}R^{-2}/m_f=R^3/m_f$.

Their exact difference equation is
\begin{equation}\label{trans:quadratic-difference-equation}
 \partial_sZ-\Delta Z
 =B_{\mathrm{act}}Z
   +a_f(B_{\mathrm{act}}-B_{\mathrm{main}})\mathcal V_{1+d}
   -\frac{R^3}{m_f}\PP\operatorname{div}(V\otimes V).
\end{equation}
By \eqref{trans:continued-first-variation} and exact scaling,
\begin{equation}\label{trans:continuing-source-h4}
 \|V\|_\infty\leq CA^PR,\qquad
 \|V\|_{H^4_y}\leq CA^PR^{7/2},\qquad
 \left\|\frac{R^3}{m_f}V\otimes V\right\|_{H^4_y}
       \leq CA^P\frac{R^{15/2}}{\mathcal M_f}.
\end{equation}
The last rate is \(2G-15\varrho/2>\varrho\).
For each time on the final interval, apply
Lemma~\ref{trans:smoothing} on its preceding interval of length
\(A^{-1}\). The initial \(L^2\) error is already bounded uniformly;
the reference vector source is controlled by
\eqref{trans:reference-coefficient-error}, and the tensor source by
\eqref{trans:continuing-source-h4} with its integrable half-order
kernel. This proves the aged \(H^4\) assertion.

For clarity, the normalized continuing-source response itself has
\(H^4\) norm at most \(CA^PR^{15/2}/\mathcal M_f\) on the final interval.
Its history through \(s-A^{-1}\) has the \(L^2\) bound
\eqref{trans:continuing-source-l2}; smooth that state on the last
interval and integrate the recent tensor using
\eqref{trans:continuing-source-h4}. Subtracting this controlled
response proves the homogeneous-continuation assertion with the same
mass and mean factor. Finally Lemma~\ref{trans:alignment} replaces
the aged point response by \(S_A\), at the additional rate
\(2\varrho\). This proves \eqref{trans:quadratic-aligned}.
\end{proof}

We record the corresponding smoothing statement for a true nonlinear
perturbation. It will be used to pass from a quadratic comparison to
the complete nonlinear output.

\begin{lemma}[Smoothing of a nonlinear perturbation]\label{trans:nonlinear-smoothing}
Let \(h\) solve the exact perturbation equation about
\(\bar u^{I,J}\). Suppose, on the enlarged final parent interval,
\begin{equation}\label{trans:nonlinear-physical-bound}
 \|h(t)\|_{H^4_x}\leq CA^Pa^2e^{2c_\star A},\qquad
 a=e^{-x_0A},\qquad
 12\varrho\leq x_0\leq13\varrho+\varrho/32.
\end{equation}
Set \(H(s,y)=Rh(\tf+R^2s,Ry)\), and let
\(\zeta\geq A^{-P}a^2e^{2GA}R^{-3}\) be a positive scalar.
If
\(\|H(s)+\zeta\mathcal V_{1+d}(s)\|_2
\leq CA^P\zeta e^{-\kappa A}\)
uniformly on \(\mathcal J_2\), then
\begin{equation}\label{trans:nonlinear-h4-upgrade}
 \|H(s)+\zeta\mathcal V_{1+d}(s)\|_{H^4}
       \leq CA^P\zeta e^{-\kappa A},\qquad s\in\mathcal J_1,
\end{equation}
provided \(0<\kappa<15\varrho/32-2\omega\) and \(\beta\) is
chosen as in Proposition~\ref{trans:quadratic}.
\end{lemma}

\begin{proof}
For \(Z=H+\zeta\mathcal V_{1+d}\), the exact equation is
\begin{equation}\label{trans:nonlinear-difference-equation}
 \partial_sZ-\Delta Z
 =B_{\mathrm{act}}Z
   -\zeta(B_{\mathrm{act}}-B_{\mathrm{main}})\mathcal V_{1+d}
   -\PP\operatorname{div}(H\otimes H).
\end{equation}
The tame product inequality and scaling give
\(\|H\otimes H\|_{H^4_y}
\leq C\|H\|_\infty\|H\|_{H^4_y}
\leq CA^PR^{9/2}a^4e^{4c_\star A}\).
After division by \(\zeta\), its exponent is at most
\begin{equation}\label{trans:nonlinear-source-rate}
 \tfrac{15}{2}\varrho-2x_0+4c_\star-2G
       \leq-\tfrac12\varrho+2\omega
       <-\tfrac{15}{32}\varrho.
\end{equation}
Apply \eqref{trans:parent-smoothing} to
\eqref{trans:nonlinear-difference-equation} on the last interval of
length \(A^{-1}\). The given \(L^2\) error supplies its initial
state, \eqref{trans:reference-coefficient-error} controls the vector
source, and \eqref{trans:nonlinear-source-rate} controls the tensor
source. Only fixed polynomial factors occur. The prescribed strict
inequality for \(\kappa\) leaves room to absorb them, proving
\eqref{trans:nonlinear-h4-upgrade}.
\end{proof}

\section{The nonlinear transition and finite cascades}
\label{sec:nonlinear}

We now turn the complete quadratic transfer into an exact nonlinear
evolution.  The amplitude and shape are normalized separately.  Higher
nonlinear terms must remain smaller than the second variation in parent
coordinates, while the normalized output must satisfy the same strong-norm
and ordinary-moment conditions as the input.  The moment is indispensable
because pressure creates algebraic tails.

After proving this shape invariance, we prepare the smallest-scale input
from smooth data at physical time zero.  A finite shooting argument then
selects a genuine orbit with a prescribed terminal coordinate.  It needs
neither monotonicity of the amplitude map nor compactness of the
infinite-dimensional shape class.

\subsection{Parameters and the normalized shape class}

The ordering of choices is part of the argument: spectral and regularity
parameters come first, then the axial-width exponent, then $A$, and only
then the cascade depth.  The shape normalization removes amplitude
errors along $g$.  Its strong norm and ordinary moment impose different
requirements and will be propagated by different estimates.

The fractions below provide one compatible set of strict margins; they
are not asserted to be optimal.  There are three kinds of requirements.
The ordinary moment must remain below $5/2$, while the complete stress
must have a fractional moment $\alpha$ satisfying
$1<\alpha<2$ and $s_\star-\alpha>1/2$.
The choices $s_\star=19/8$ and $\alpha=7/4$ leave
$s_\star-\alpha=5/8$.  The extra moment used only in smooth preparation
is $s_\star+1/16=39/16<5/2$; it is not imposed on arbitrary shapes.

Next, the exponential gap $\omega=c_\star-G$ must be small compared with
the scale-separation exponent $\varrho$.  The condition
$2\omega<\varrho/32$ makes the source-comparison and transfer errors
small relative to the lower bound on the positive stress mass.  Finally the allowed shape and
moment tolerances must survive normalization by that mass.  The proof
will obtain a strong error with exponent $13\varrho/32$, whereas it
allows the weaker tolerance $\rho_c=12\varrho/32$.  Its quadratic moment
cost is $\varrho/8+2\omega<5\varrho/32$, below the allowed
$\gamma=8\varrho/32$.

These are separate restrictions: increasing the axial width cannot
repair a bad mass-normalized exponent in $R$, and increasing the depth
cannot repair any of them.  Once the strict exponent margins are fixed,
$A^Pe^{-\epsilon A}\to0$ for each fixed $P$ and $\epsilon>0$ absorbs
the polynomial constants.  Only the finite list of width-dependent
errors is handled by increasing $\beta$.

We retain the constants
\begin{equation}\label{nl:parameters}
 \begin{gathered}
 c_\star=\sqrt2e^{-3/2},\qquad \varrho=c_\star/8,\qquad
 G=c_\star-\omega,\qquad 0<2\omega<\varrho/32,\\
 R=e^{\varrho A},\qquad b=e^{\beta A},\qquad
 \tp=(2+\Theta\log A)/A,\qquad
 \tf=\tp+16A^2,\qquad \ts=R^2\tp,\\
 s_\star=19/8,\qquad \rho_c=3\varrho/8,\qquad
 \gamma=\varrho/4,\qquad \delta_A=e^{-\rho_cA}.
 \end{gathered}
\end{equation}

The number $G$ and the fixed preparation constant $\Theta$ are supplied
by Lemma~\ref{seed:signed-history}.  Its frequency interval and every
required derivative order are fixed before $\beta$ and $A$.
Also prescribe the endpoint loss $\lambda_0=2\varrho$ in
Lemma~\ref{dec:retirement} at moment $s=s_\star$.  Include the width
thresholds for this choice and for Lemma~\ref{dec:endpoint} in the finite
list of conditions from Sections~\ref{sec:history}--\ref{sec:transfer}
and the preparation argument below.  Choose one $\beta$ satisfying that
list, and only then take $A$ sufficiently large.  All these parameters
remain fixed as the cascade depth tends to infinity.

Recall the point distribution, prepared response, and its normalization:
\begin{equation}\label{nl:seed}
 \mathcal I=\PP\partial_z(e_z\delta_0),\qquad
 S_A=\mathcal U_{\mathrm d}(\tp,0)\mathcal I,\qquad
 n_A=\norm{S_A}_2,\qquad g=S_A/n_A.
\end{equation}
By Lemma~\ref{seed:prepared}, one has $n_A\geq c_0>0$, $\norm g_2=1$,
and polynomial bounds for every fixed weighted graph norm used below.
All fields in this section are real, centrally odd, axisymmetric, and
solenoidal.  Define the shape class and amplitude interval by
\begin{equation}\label{nl:shape-class}
 \begin{split}
 \mathcal C_A&=\left\{g+\zeta:
     \langle\zeta,g\rangle=0,\quad
     \norm\zeta_{H^4}\leq\delta_A,\quad
     M_{s_\star}(\zeta)\leq e^{\gamma A}\right\},\\
 a_-&=e^{-(13\varrho+\varrho/32)A},\qquad
 a_+=e^{-12\varrho A}.
 \end{split}
\end{equation}
Thus $\langle f,g\rangle=1$ for every $f\in\mathcal C_A$.  The
pairing in this section is the real $L^2$ pairing.

The clocks and the distinct roles of the three norms are collected in the
\hyperref[nl:phase-guide]{phase guide in Section~\ref{sec:strategy}}.
Here we use that guide with the precise parameters just fixed.

For reference, we state the quantitative consequences of the linear
analysis in the form needed for the nonlinear argument.  Use
$\mathcal L_{\bar u}$ from \eqref{strategy:linearized-generator} with
$\bar u=\bar u^{I,J}$, and recall the source and transfer growth integrals
from \eqref{hist:growth} and \eqref{trans:late-growth}:
\begin{equation}\label{nl:growth}
 \Gamma_{\mathrm{src}}(t)=c_\star A
       \left((1+\tp)^{-1}-(1+t)^{-1}\right),\qquad
 \Gamma_{\mathrm{tr}}(t)=2AR^2
       \left((R^2+\tf)^{-1}-(R^2+t)^{-1}\right).
\end{equation}
The weighted energy of Proposition~\ref{hist:energy}, together with
its final-interval smoothing estimate, gives the following implication
for the zero-data equation $(\partial_t-\mathcal L_{\bar u})v=-\PP\div F$:
\begin{equation}\label{nl:early-source}
 \begin{split}
 \norm{F(t)}_{H^4}+M_\mu(F(t))&\leq B_0e^{m\Gamma_{\mathrm{src}}(t)},\quad
                  v(\tp)=0\\
 &\Longrightarrow\quad
 \norm{v(t)}_{H^4}+M_\mu(v(t))
       \leq CA^PB_0e^{m\Gamma_{\mathrm{src}}(t)},
 \end{split}
\end{equation}
for $1\leq\mu<5/2$ and $m\geq2$.  The polynomial is independent of
$m$.  Indeed, the integrating factor bounds the weighted squared
source by $CA^P(\tf-\tp)B_0^2e^{2m\Gamma_{\mathrm{src}}(t)}$;
smoothing only on a final interval of inverse-polynomial length gives
the same degree-independent estimate in $H^4$.

For $f\in\mathcal C_A$, write $\zeta_f=f-g$ and let $v_f,w_f$ be the actual first
and second variations about $\bar u$, with $v_f(\tp)=f$ and
$w_f(\tp)=0$.  Proposition~\ref{hist:generic} supplies their own
positive stress mass $\mathcal M_f$ and the bounds
\begin{equation}\label{nl:mass}
 A^{-P}e^{2GA}\leq \mathcal M_f\leq CA^Pe^{2c_\star A},\qquad
 \norm{w_f(t)}_{H^4}+M_{s_\star}(w_f(t))
       \leq CA^Pe^{2\Gamma_{\mathrm{src}}(t)}.
\end{equation}
The mass here is computed from the complete tensor
$\mathsf S_f=v_f\otimes v_f+\bar u\otimes w_f+w_f\otimes\bar u$,
as defined in Section~\ref{sec:history}.

For the moment and loss prescribed above, Lemmas~\ref{dec:retirement}
and \ref{dec:endpoint} give
\begin{equation}\label{nl:linear-decay}
 \begin{split}
 M_{s_\star}(v_f(\tf))&\leq CA^Pe^{2\varrho A}M_{s_\star}(f),\\
 \norm{v_f(\tf)}_{H^4}
 &\leq CA^P\norm f_2+CA^Pe^{K_{s_\star}A}b^{-1/2}M_{s_\star}(f).
 \end{split}
\end{equation}
The chosen width makes the second line polynomial uniformly on
$\mathcal C_A$.  Each other prescribed loss in Lemma~\ref{dec:retirement}
is available after its corresponding width choice; no simultaneous
assertion for every loss at this fixed width is needed here.

Proposition~\ref{trans:quadratic} also gives a
number $a_f=1+O(A^{-1})>0$ such that, for any fixed $0<\kappa_{\mathrm r}<\varrho$,
\begin{equation}\label{nl:quadratic-transfer}
 \norm{R\bigl(\mathcal U_{\bar u}(\ts,\tf)w_f(\tf)\bigr)(R\cdot)
       +a_f\mathcal M_fR^{-3}S_A}_{H^4}
 \leq CA^P\mathcal M_fR^{-3}e^{-\kappa_{\mathrm r}A}.
\end{equation}
The same comparison in $L^2$ holds uniformly on the enlarged final
interval in parent-scale coordinates.  Finally, Lemma~\ref{trans:late}
propagates $H^4$ homogeneously by
$Ce^{\Gamma_{\mathrm{tr}}(t)-\Gamma_{\mathrm{tr}}(u)}$, and sends a
zero-data tensor bounded by $B_0e^{m\Gamma_{\mathrm{tr}}}$ to an
$H^4$ response bounded by
$C(R/\sqrt A)B_0e^{m\Gamma_{\mathrm{tr}}}$, independently of $m\geq2$.
Its total growth is polynomial because
$\Gamma_{\mathrm{tr}}(\ts)\leq2A\tp$.

\subsection{An invariant class for the exact nonlinear map}

The proof has three parts.  A convergent expansion controls the source
interval; a second continuation controls the long transfer interval and
identifies the leading amplitude; a separate weighted argument returns
the output to the moment class.  The last part must avoid a quadratic
loss in the incoming moment tolerance.  These conclusions concern the
full nonlinear solution, including all higher-order terms.

Here $a$ is the positive \emph{perturbation amplitude}, and the active
Gaussian column has radius one in the local transition coordinates.
The normalized shape $f$ lies in a hyperplane
$\langle f,g\rangle_2=1$.  Write $h_{\mathrm{out}}$ for the final
perturbation in parent coordinates, as defined in \eqref{nl:map}.
Its amplitude and shape are
$\mathcal A=-\langle h_{\mathrm{out}},g\rangle_2$ and
$\mathcal F=-h_{\mathrm{out}}/\mathcal A$; since $\|g\|_2=1$, this gives
$\langle\mathcal F,g\rangle_2=1$ \emph{exactly} whenever
$\mathcal A>0$.

\begin{proposition}[One nonlinear transition]\label{nl:transition}
Let $I\geq0$, $J\geq1$, $f\in\mathcal C_A$, and $a\in[a_-,a_+]$.
The datum $u(\tp)=\bar u^{I,J}(\tp)-af$ has a strong finite-energy
solution through $\ts$.  Set $h=u-\bar u^{I,J}$ and define
\begin{equation}\label{nl:map}
 h_{\mathrm{out}}(y)=Rh(\ts,Ry),\qquad
 \mathcal A^{I,J}(a,f)=-\langle h_{\mathrm{out}},g\rangle,\qquad
 \mathcal F^{I,J}(a,f)=-h_{\mathrm{out}}/\mathcal A^{I,J}(a,f).
\end{equation}
The amplitude is positive, the shape belongs to $\mathcal C_A$, and
both maps are continuous in $a$ and in the $H^4$ topology of $f$.
Uniformly in $I,J,a,f$, one has
\begin{equation}\label{nl:amplitude-bounds}
 A^{-P}e^{(2G-3\varrho)A}a^2
 \leq\mathcal A^{I,J}(a,f)
 \leq A^Pe^{(2c_\star-3\varrho)A}a^2
\end{equation}
and
\begin{equation}\label{nl:crossings}
 \mathcal A^{I,J}(a_-,f)<a_-/2,\qquad
 \mathcal A^{I,J}(a_+,f)>2a_+.
\end{equation}
There is a finite constant $C(A)$, independent of the cutoffs and
of the input in this class, such that
\begin{equation}\label{nl:interstage}
 \sup_{\tp\leq t\leq\ts}
       \bigl(\norm{h(t)}_{H^4}+M_{s_\star}(h(t))\bigr)\leq C(A).
\end{equation}
The output reference is exactly $\bar u^{I+1,J-1}(\tp)$, and the
shape assertion holds even if the output amplitude leaves $[a_-,a_+]$.
\end{proposition}

\begin{proof}
\smallskip\noindent\textit{Construction on the source-scale interval.}\quad
We first construct the entire nonlinear solution on $[\tp,\tf]$.
Let $V_1=v_f$, and for $m\geq2$ solve
\begin{equation}\label{nl:series}
 (\partial_t-\mathcal L_{\bar u})V_m
       =-\PP\div\sum_{j=1}^{m-1}V_j\otimes V_{m-j},\qquad
 V_m(\tp)=0.
\end{equation}
At $m=2$, this is precisely the second-variation equation, so
$V_1=v_f$ and $V_2=w_f$; no factorial is included.

\textit{The unweighted source estimate for the nonlinear recurrence.}
Apply the unweighted version of Proposition~\ref{hist:energy}:
if $v(\tp)=0$ and $\|F(t)\|_{H^4}\leq B e^{m\Gamma_{\mathrm{src}}(t)}$,
then $\|v(t)\|_{H^4}\leq CA^P B e^{m\Gamma_{\mathrm{src}}(t)}$,
with $C,P$ independent of $m\geq2$.
The unweighted energy estimate and final-interval smoothing prove this
bound; monotonicity of $\Gamma_{\mathrm{src}}$ controls the source without
differentiating $e^{m\Gamma_{\mathrm{src}}}$.  No weighted norm is used
in this recurrence.  The ordinary moment is estimated separately below.
Combining this bound with the $H^4$ algebra estimate gives, with
$L_{\mathrm{src}}(A)=CA^P$ and
$b_m=\sup_{\tp\leq t\leq\tf}
e^{-m\Gamma_{\mathrm{src}}(t)}\norm{V_m(t)}_{H^4}$,

\begin{equation}\label{nl:catalan}
 b_1\leq L_{\mathrm{src}}(A),\qquad
 b_m\leq L_{\mathrm{src}}(A)\sum_{j=1}^{m-1}b_jb_{m-j},\qquad
 b_m\leq4^{m-1}L_{\mathrm{src}}(A)^{2m-1}.
\end{equation}

To see the recurrence explicitly, introduce nonnegative coefficients
$d_m$ by $d_1=L_{\mathrm{src}}(A)$ and
$d_m=L_{\mathrm{src}}(A)\sum_{j=1}^{m-1}d_jd_{m-j}$.
The generating series is determined near $z=0$ by
\[
 \sum_{m\geq1}d_mz^m
   =L_{\mathrm{src}}(A)z+L_{\mathrm{src}}(A)\left(\sum_{m\geq1}d_mz^m\right)^2
   =\frac{1-\sqrt{1-4L_{\mathrm{src}}(A)^2z}}{2L_{\mathrm{src}}(A)}.
\]
Its coefficients are the Catalan coefficients times $L_{\mathrm{src}}(A)^{2m-1}$,
so $d_m\leq4^{m-1}L_{\mathrm{src}}(A)^{2m-1}$ and induction gives $b_m\leq d_m$.
With $z=a e^{\Gamma_{\mathrm{src}}(t)}$, absolute convergence follows from
$4L_{\mathrm{src}}(A)^2a e^{\Gamma_{\mathrm{src}}(t)}<1$ uniformly on the interval.
The same bound with the first two terms removed gives the cubic remainder.

This argument controls every nonlinear order, not just a formal Taylor
expansion.  The series converges in $H^4$ and its quadratic products
converge absolutely in the source norm.  One may therefore pass the sum
through the mild equation using the integrable tensor-source kernel.
The resulting mild solution is the strong solution furnished by local
parabolic theory, and its uniform norm continues it across the interval.
The key prior fact is that the source estimate has a constant independent
of the degree $m$.  A bound that introduced an uncontrolled constant at
each degree would not imply this generating-function estimate.

Since $ae^{c_\star A}$ is exponentially small throughout the
amplitude interval, the series $h=\sum_{m\geq1}(-a)^mV_m$
and its nonlinear source converge in the strong solution space.
In particular,
\begin{equation}\label{nl:early-remainder}
 \norm{h(t)}_{H^4}\leq CA^Pae^{\Gamma_{\mathrm{src}}(t)},\qquad
 \norm{h(\tf)+av_f(\tf)-a^2w_f(\tf)}_{H^4}
       \leq CA^Pa^3e^{3c_\star A}.
\end{equation}
By \eqref{nl:linear-decay} and \eqref{nl:mass}, it follows that
\begin{equation}\label{nl:at-tf}
 \norm{h(\tf)}_{H^4}\leq CA^Pa^2e^{2c_\star A}.
\end{equation}
Here $a\leq a^2e^{2c_\star A}$ for sufficiently large $A$, because
$13\varrho+\varrho/32<2c_\star$.

\smallskip\noindent\textit{Continuation and the leading amplitude.}\quad
Apply the tensor-source estimate of Lemma~\ref{trans:late} to a
second series starting from $h(\tf)$.  Its convergence parameter is
bounded by $CA^P(R/\sqrt A)a^2e^{2c_\star A}$.
Writing $a=e^{-xA}$ gives $12\varrho\leq x\leq13\varrho+\varrho/32$
and
$\varrho-2x+2c_\star\leq-7\varrho$.
Thus the exact solution continues to $\ts$, and
\begin{equation}\label{nl:late-bound}
 \sup_{\tf\leq t\leq\ts}\norm{h(t)}_{H^4}
       \leq CA^Pa^2e^{2c_\star A}.
\end{equation}
The difference from the homogeneous propagation of $a^2w_f(\tf)$,
in source-coordinate $L^2$, is bounded by
\begin{equation}\label{nl:joining-error}
 CA^P\left(a+a^3e^{3c_\star A}+Ra^4e^{4c_\star A}\right).
\end{equation}
The three terms are, respectively, the decayed first variation,
the early terms of degree at least three, and the entire later
nonlinear response.  Proposition~\ref{trans:quadratic} estimates the continuing quadratic source, which remains part of the exact perturbation equation.

Define
\begin{equation}\label{nl:leading}
 A_{\mathrm{lead}}=a_fa^2\mathcal M_fR^{-3}n_A>0,
 \qquad h_{\mathrm{out}}=-A_{\mathrm{lead}}g+r_{\mathrm{out}}.
\end{equation}
Normalize each error by the outgoing amplitude $A_{\rm lead}$.
For example, the decayed first variation has physical $L^2$ size
at most $CA^Pa$.  Velocity rescaling multiplies that norm by $R^{-1/2}$.
Dividing by the lower bound for $A_{\mathrm{lead}}$ therefore gives
\[
 \frac{CA^PR^{-1/2}a}{A_{\mathrm{lead}}}
 \leq CA^{P'}R^{5/2}a^{-1}e^{-2GA}
 =CA^{P'}e^{(x+5\varrho/2-2G)A}.
\]
The factor $R^{5/2}$ is the combination of the $L^2$ rescaling and the
$R^{-3}$ in the leading amplitude.  The early cubic term and the late
nonlinear term must be divided by this same amplitude.  This produces
the next three exponents.
At the worst endpoint $x=13\varrho+\varrho/32$, the first exponent
is $-15\varrho/32+2\omega$.  The available decay is therefore stronger
than $13\varrho/32$, and that last value is stronger than the tolerance
$12\varrho/32$.  The two strict gaps have different jobs: one absorbs
constants in the error estimate, and the other puts the normalized
output back inside the prescribed shape class.

The $L^2$ scaling from source to parent coordinates is $R^{-1/2}$.
After division by $A_{\mathrm{lead}}$, the terms in
\eqref{nl:joining-error} have exponential rates
\begin{equation}\label{nl:three-rates}
 x+\tfrac52\varrho-2G,\qquad
 -x+3c_\star+\tfrac52\varrho-2G,\qquad
 -2x+4c_\star+\tfrac72\varrho-2G.
\end{equation}
Over the chosen amplitude interval, these are at most
$-15\varrho/32+2\omega$, $-3\varrho/2+2\omega$, and
$-9\varrho/2+2\omega$.  The transfer error in
\eqref{nl:quadratic-transfer} can be taken with
$\kappa_{\mathrm r}=3\varrho/4$.

The shape requires an $H^4$ estimate in parent coordinates; rescaling
the $L^2$ error alone does not provide it.  We now apply
Lemma~\ref{trans:nonlinear-smoothing} on the enlarged final parent
interval.  Put
\[
 h_{\mathrm{par}}(s,y)=Rh(\tf+R^2s,Ry),\qquad
 \zeta_{\mathrm{tr}}=a_fa^2\mathcal M_fR^{-3}=A_{\mathrm{lead}}/n_A,
 \qquad \kappa=13\varrho/32.
\]
Here $h_{\mathrm{par}}(s_\ast,\cdot)$ is the output $h_{\mathrm{out}}$ in \eqref{nl:leading},
and the scalar in Lemma~\ref{trans:nonlinear-smoothing} is chosen as
$\zeta_{\mathrm{tr}}$, distinct from the input shape error.  It has
the lower bound required in that lemma by
\eqref{nl:mass} and $a_f=1+O(A^{-1})>0$.
The physical $H^4$ hypothesis is \eqref{nl:late-bound}.
Moreover, \eqref{nl:joining-error} and the uniform aged version of
\eqref{nl:quadratic-transfer} give
\[
 \|h_{\mathrm{par}}(s)+\zeta_{\mathrm{tr}}\mathcal V_{1+d}(s)\|_2
       \leq CA^P\zeta_{\mathrm{tr}} e^{-\kappa A}
\]
on the enlarged final interval $\mathcal J_2$.
The rates in \eqref{nl:three-rates} are strictly stronger than $\kappa$, so the
polynomial bounds for $n_A$ cause no loss of this exponent.
For reference, the nonlinear tensor in this application obeys
\begin{equation}\label{nl:rescaled-product}
 \norm{Rh(t,R\cdot)\otimes Rh(t,R\cdot)}_{H^4}
       \leq CA^PR^{9/2}a^4e^{4c_\star A}.
\end{equation}
After division by $\zeta_{\mathrm{tr}}$, its exponent is at most
$15\varrho/2-2x+4c_\star-2G\leq-\varrho/2+2\omega$,
exactly the source margin used in Lemma~\ref{trans:nonlinear-smoothing}.
Since $2\omega<\varrho/32$, the choice
$\kappa=13\varrho/32$ satisfies its strict requirement
$\kappa<15\varrho/32-2\omega$.
That lemma upgrades the comparison to $H^4$ on $\mathcal J_1$.
At the endpoint, Lemma~\ref{trans:alignment} replaces
$\mathcal V_{1+d}(s_\ast)$ by $S_A=n_Ag$ with the stronger error
$CA^PR^{-2}$.  Consequently
\begin{equation}\label{nl:shape-error}
 \norm{r_{\mathrm{out}}}_{H^4}\leq CA^Pe^{-\kappa A}A_{\mathrm{lead}},
 \qquad \kappa=13\varrho/32>\rho_c.
\end{equation}

Write $r_{\mathrm{out}}=\langle r_{\mathrm{out}},g\rangle g+r_{\mathrm{out},\perp}$, with
$\langle r_{\mathrm{out},\perp},g\rangle=0$ and $\|g\|_2=1$.
Then
\[
 h_{\mathrm{out}}=-\bigl(A_{\mathrm{lead}}-\langle r_{\mathrm{out}},g\rangle\bigr)g+r_{\mathrm{out},\perp}.
\]
The parallel correction changes the exact amplitude; only $r_{\mathrm{out},\perp}$
changes the normalized shape.  If $\|r_{\mathrm{out}}\|_2\leq A_{\mathrm{lead}}/2$,
then $\mathcal A\geq A_{\mathrm{lead}}/2$, and
\[
 \|\mathcal F-g\|_{H^4}
 \leq \frac{2(1+\|g\|_{H^4})}{A_{\mathrm{lead}}}\|r_{\mathrm{out}}\|_{H^4}.
\]
The factor involving $g$ is polynomial in $A$, so the strict difference
between $13\varrho/32$ and $\rho_c$ absorbs it.
This is also why the input-dependent mass and the $O(A^{-1})$ mean
correction were included in $A_{\mathrm{lead}}$: they need not be as
small as the permitted shape error, because they act in the amplitude
direction.  The ordinary moment still requires a separate estimate;
smallness in $H^4$ alone cannot control a perturbation spread far away.

It follows directly from \eqref{nl:map} and \eqref{nl:leading} that
\begin{equation}\label{nl:normalization}
 \mathcal A=A_{\mathrm{lead}}-\langle r_{\mathrm{out}},g\rangle>0,\qquad
 \mathcal F-g=\frac{\langle r_{\mathrm{out}},g\rangle g-r_{\mathrm{out}}}{\mathcal A}.
\end{equation}
The polynomial graph bounds for $g$ and the strict gap
$\kappa>\rho_c$ give the $H^4$ and hyperplane conditions in
\eqref{nl:shape-class}, as well as \eqref{nl:amplitude-bounds}.

\smallskip\noindent\textit{Propagation of the ordinary moment.}\quad
The large permitted moment of $\zeta_f$ must not enter quadratically.
Let $h_g$ be the exact perturbation from $-ag$ on the same reference,
and put $\delta h=h-h_g$.  Telescoping one input at a time in
\eqref{nl:series}, the extra factor $m$ is summable against the
Catalan majorant.  Therefore
\begin{equation}\label{nl:tame-h4}
 \norm{\delta h(t)}_{H^4}\leq CA^Pa\delta_Ae^{\Gamma_{\mathrm{src}}(t)},\qquad
 \norm{h_g(t)}_{H^4}\leq CA^Pae^{\Gamma_{\mathrm{src}}(t)}.
\end{equation}
Apply the capped weighted energy of Proposition~\ref{hist:energy}
to $h_g$.  The additional coefficient from its own tensor has
integral at most $CA^Pa^2e^{2c_\star A}=o(1)$, so
$M_{s_\star}(h_g(t))\leq CA^Pae^{\Gamma_{\mathrm{src}}(t)}$.
For the difference, the exact tensor is
$h_g\otimes \delta h+\delta h\otimes h_g+\delta h\otimes \delta h$.
Put the weight on $h_g$ in each cross term.  If
$E_{\delta h}(t)=N_{s_\star}(t,\delta h(t))^2$ denotes the normalized weighted energy from
Proposition~\ref{hist:energy}, its inequality becomes
\begin{equation}\label{nl:moment-difference}
 E_{\delta h}'+cD_{s_\star}(t,\delta h(t))
 \leq\left(2\Gamma_{\mathrm{src}}'+a_{s_\star}
                          +C\norm{\delta h}_\infty^2\right)E_{\delta h}
        +C\norm{\delta h}_\infty^2N_{s_\star}(t,h_g(t))^2.
\end{equation}
Here $N_{s_\star}$, $D_{s_\star}$, and $a_{s_\star}$ are exactly
the energy, dissipation, and integrable remainder coefficient of
that proposition.  If $B_{\mathrm{mom}}=M_{s_\star}(\zeta_f)\leq e^{\gamma A}$,
the initial square-root energy is bounded by $CA^Pa(1+B_{\mathrm{mom}})$.
Equation~\eqref{nl:tame-h4} gives
\begin{equation}\label{nl:moment-difference-integrals}
 \int_\tp^\tf\norm{\delta h}_\infty^2\idiff t
       \leq CA^Pa^2\delta_A^2e^{2c_\star A}=o(1),
 \qquad
 \int_\tp^\tf e^{-2\Gamma_{\mathrm{src}}}
       \norm{\delta h}_\infty^2N_{s_\star}(t,h_g(t))^2\idiff t
       \leq CA^Pa^4\delta_A^2e^{2c_\star A}.
\end{equation}
Gronwall's inequality in \eqref{nl:moment-difference} therefore
gives the bound
\begin{equation}\label{nl:moment-difference-bound}
 M_{s_\star}(\delta h(t))\leq CA^Pa(1+B_{\mathrm{mom}})e^{\Gamma_{\mathrm{src}}(t)}.
\end{equation}
All estimates are uniform in the cap, so monotone convergence
proves the moment membership of the constructed strong solution.

Expand the tensor with $h=h_g+\delta h$ and put the weight on only one factor:
\[
 \begin{aligned}
 M_{s_\star}(h_g\otimes h_g)&\leq
       \|h_g\|_\infty M_{s_\star}(h_g),\\
 M_{s_\star}(h_g\otimes \delta h)+M_{s_\star}(\delta h\otimes h_g)&\leq
       2\|\delta h\|_\infty M_{s_\star}(h_g),\\
 M_{s_\star}(\delta h\otimes \delta h)&\leq
       \|\delta h\|_\infty M_{s_\star}(\delta h).
 \end{aligned}
\]
Only the last weighted factor carries $1+B_{\mathrm{mom}}$; it is multiplied by the
small unweighted norm of $\delta h$.  Equations \eqref{nl:tame-h4} and
\eqref{nl:moment-difference-bound} give
\begin{equation}\label{nl:tame-tensor}
 M_{s_\star}(h\otimes h)+\norm{h\otimes h}_2
       \leq CA^Pa^2e^{2\Gamma_{\mathrm{src}}(t)}(1+\delta_A B_{\mathrm{mom}}).
\end{equation}
Since $h+av_f$ starts at zero, \eqref{nl:early-source},
\eqref{nl:tame-tensor}, and the prescribed-loss estimate
\eqref{nl:linear-decay} imply
\begin{equation}\label{nl:moment-tf}
 Y_f:=M_{s_\star}(h(\tf))
 \leq CA^P\left(ae^{2\varrho A}(1+B_{\mathrm{mom}})
                     +a^2e^{2c_\star A}(1+\delta_A B_{\mathrm{mom}})\right).
\end{equation}

On $[\tf,\ts]$, the time integral of the complete heat matrix
$M_Q$ is $O(1+\log A)$.  The radii at most one contribute
$CA/t^2$, the radius $R$ contributes $CAR^2/(R^2+t)^2$, and
the remaining radii contribute $CA/R^4$.  The remaining tensor
is exactly
$e^{I,J}\otimes h+h\otimes e^{I,J}+h\otimes h$.
Its moment is bounded by
$C(\norm{e^{I,J}}_\infty+\norm h_\infty)M_{s_\star}(h)$.
Using Lemma~\ref{bg:weighted}, \eqref{bg:derivatives}, and
\eqref{nl:late-bound}, the weighted mild inequality gives
\begin{equation}\label{nl:moment-late}
 Y_\ast\leq CA^PY_f
       +CA^PR^2\left(e^{KA}b^{-1/2}+a^2e^{2c_\star A}\right)Y_\ast,
 \qquad Y_\ast=\sup_{\tf\leq t\leq\ts}M_{s_\star}(h(t)).
\end{equation}
The nonlinear exponent in the last coefficient is at most
$2\varrho-24\varrho+16\varrho=-6\varrho$.
Choose $\beta/2>K+2\varrho$ to absorb the correction term as well.
Then $Y_\ast\leq CA^PY_f$.  This argument may first be made on
short weighted existence intervals or with capped weights;
unweighted strong uniqueness identifies these with the solution
already constructed.

The homogeneous part of the parent moment gains
$R^{-s_\star-1/2}$ from velocity rescaling.  The leading amplitude is
bounded below by a polynomial times $a^2e^{2GA}R^{-3}$.
For the quadratic part of $Y_f$, division therefore costs, up to powers
of $A$,
\[
 R^{-s_\star-1/2}
 \frac{a^2e^{2c_\star A}}{a^2e^{2GA}R^{-3}}
   =e^{\bigl((5/2-s_\star)\varrho+2\omega\bigr)A}
   =e^{(\varrho/8+2\omega)A}.
\]
For the linear part, substitute $a=e^{-xA}$ and the endpoint loss
$e^{2\varrho A}$; its exponent is
$(5/2-s_\star)\varrho+x-2G+2\varrho$.
The large input moment $B_{\mathrm{mom}}$ enters this linear contribution as $1+B_{\mathrm{mom}}$,
but in the quadratic contribution only as $1+\delta_A B_{\mathrm{mom}}$, by the tame
tensor estimate.  These are exactly the two terms that must fit below
the allowed output tolerance.  The constant part of the inhomogeneous
moment is controlled instead by the normalized $L^2$ estimate.
Treating both parts as if they had the same scaling would lose this
useful distinction.

The constant part and homogeneous part of the output weight scale
differently.  Their exact relations are
\begin{equation}\label{nl:moment-scaling}
 \norm{h_{\mathrm{out}}}_2=R^{-1/2}\norm{h(\ts)}_2,\qquad
 \norm{|y|^{s_\star}h_{\mathrm{out}}}_2
       =R^{-s_\star-1/2}\norm{|x|^{s_\star}h(\ts)}_2.
\end{equation}
Combining \eqref{nl:moment-tf}--\eqref{nl:moment-scaling} with
\eqref{nl:normalization} gives
\begin{equation}\label{nl:output-moment}
 \begin{split}
 M_{s_\star}(\mathcal F)
 &\leq CA^P\left(1+e^{\gamma_0A}(1+\delta_A B_{\mathrm{mom}})
                       +e^{\xi_{\mathrm{lin}}A}(1+B_{\mathrm{mom}})\right),\\
 \gamma_0&=\varrho/8+2\omega,\qquad
 \xi_{\mathrm{lin}}=(5/2-s_\star)\varrho+x-2G+2\varrho.
 \end{split}
\end{equation}

At the worst endpoint $x=13\varrho+\varrho/32$, the explicit choices give
\[
 \gamma_0<5\varrho/32<\gamma=8\varrho/32,\qquad
 \delta_A B_{\mathrm{mom}}\leq e^{-\varrho A/8},\qquad
 \xi_{\mathrm{lin}}+\gamma\leq-19\varrho/32+2\omega<0.
\]
These strict inequalities hold before $A$ is increased, so every fixed
polynomial factor in \eqref{nl:output-moment} is absorbed independently
of the cascade depth.  Subtracting the polynomial moment of $g$ proves
$M_{s_\star}(\mathcal F-g)<e^{\gamma A}$ and completes the shape invariance.

  Equations~\eqref{nl:early-remainder},
\eqref{nl:late-bound}, and \eqref{nl:moment-late} also prove
\eqref{nl:interstage}.

At $a_-$ the upper bound in \eqref{nl:amplitude-bounds}, divided by
$a_-$, is at most $A^Pe^{-\varrho A/32}$; at $a_+$ the lower
bound, divided by $a_+$, is at least
$A^{-P}e^{(\varrho-2\omega)A}$.  This proves
\eqref{nl:crossings}.  The two convergent nonlinear constructions
depend continuously on their $H^4$ initial data, and
\eqref{nl:normalization} divides by a positive scalar.  The asserted
continuity follows.  Finally, \eqref{bg:covariance} identifies the
output reference, including its nonlinear correction.
\end{proof}

\subsection{Preparation from smooth initial data}

We prepare the first admissible state from smooth data at the smallest
retained scale.  Thereafter Proposition~\ref{nl:transition} supplies
successive states of the same solution.  The regularization parameter
is chosen before the depth, and the estimates must place the prepared
state in $\mathcal C_A$ uniformly in the number of larger columns.

The regularized impulse obeys the exact scaling
$\mathcal I_{\eta_{\mathrm{reg}}}(x)=\eta_{\mathrm{reg}}^{-2}\mathcal I_1(x/\sqrt{\eta_{\mathrm{reg}}})$.
Consequently its $j$th derivative has $L^2$ scale
$\eta_{\mathrm{reg}}^{-5/4-j/2}$; in particular the $H^4$ cost is at most
$C\eta_{\mathrm{reg}}^{-13/4}$ for $0<\eta_{\mathrm{reg}}\leq1$.
With $\eta_{\mathrm{reg}}=e^{-\varrho A}$ and $|\varepsilon|\leq2a_+$,
\[
 |\varepsilon|\eta_{\mathrm{reg}}^{-13/4}\leq2e^{-35\varrho A/4},\qquad
 |\varepsilon|\eta_{\mathrm{reg}}^{-13/2}\leq2e^{-11\varrho A/2}.
\]
The first is a smallness parameter for the nonlinear preparation; the
second controls its normalized quadratic remainder.  These facts alone
do not prove that the prepared shape is near $g$: the linear endpoint
must also be compared to the unregularized point response.  That
comparison is made in negative norms followed by positive-time smoothing,
which avoids a negative power of $\eta_{\mathrm{reg}}$ in its leading endpoint error.
Finally $\eta_{\mathrm{reg}}$ remains fixed as the depth tends to infinity.  The physical
initial perturbation vanishes because of $\ell_N^{1/2}$ scaling, not
because $\eta_{\mathrm{reg}}$ is sent to zero at each stage.

To realize a finite cascade as an actual solution starting at physical
time zero, we prepare its smallest scale from a regularized point
distribution.  Only arrays whose smallest radius is one occur in this
step.  Define
\begin{equation}\label{nl:regularization}
 \eta_{\mathrm{reg}}=e^{-\varrho A},\qquad \chi_{\mathrm{seed}}=1/16,\qquad
 \mathcal I_{\eta_{\mathrm{reg}}}=e^{\eta_{\mathrm{reg}}\Delta}\mathcal I,
 \qquad \chi_{\mathrm{seed}}\beta>6\varrho.
\end{equation}
The regularization parameter is fixed before the depth limit.  The
field $\mathcal I_{\eta_{\mathrm{reg}}}$ is smooth and retains the full axial Fourier
spectrum and all prescribed symmetries.

\begin{proposition}[Initial preparation]\label{nl:preparation}
For every finite $J\geq0$ and every real $\varepsilon$ with
$a_-/2\leq|\varepsilon|\leq2a_+$, the datum
\begin{equation}\label{nl:prepared-datum}
 u^{0,J}_\varepsilon(0)=Q^{0,J}(0)
                    +(\varepsilon/n_A)\mathcal I_{\eta_{\mathrm{reg}}}
\end{equation}
has a strong finite-energy solution on $[0,\tp]$.  Its prepared
perturbation has the form
\begin{equation}\label{nl:preparation-error}
 \begin{split}
 \frac{u^{0,J}_\varepsilon(\tp)-\bar u^{0,J}(\tp)}{\varepsilon}
       &=g+d_{J,\varepsilon},\\
 \norm{d_{J,\varepsilon}}_{H^4}
       &\leq CA^P\left(\eta_{\mathrm{reg}}+R^{-2}
              +\eta_{\mathrm{reg}}^{-4}(b^{-\chi_{\mathrm{seed}}}+b^{-1/2})
              +|\varepsilon|\eta_{\mathrm{reg}}^{-13/2}\right),\\
 M_{s_\star}(g+d_{J,\varepsilon})&\leq CA^P.
 \end{split}
\end{equation}
The solution and its prepared state depend continuously on
$\varepsilon$ in $H^4$, uniformly in $J$.

For $\varepsilon=-\lambda$, $a_-/2\leq\lambda\leq2a_+$, define
\begin{equation}\label{nl:preparation-map}
 B_J(\lambda)=-\langle u^{0,J}_{-\lambda}(\tp)
                               -\bar u^{0,J}(\tp),g\rangle,
 \qquad
 F_J(\lambda)=-\frac{u^{0,J}_{-\lambda}(\tp)
                               -\bar u^{0,J}(\tp)}{B_J(\lambda)}.
\end{equation}
Then $B_J>0$, $F_J\in\mathcal C_A$, both maps are continuous,
and their endpoint amplitudes satisfy
\begin{equation}\label{nl:preparation-crossings}
 B_J(a_-/2)<a_-,\qquad B_J(2a_+)>a_+.
\end{equation}
\end{proposition}

\begin{proof}
For the radii $1,R,\ldots,R^J$, the Cartesian matrix in
\eqref{bg:cartesian-matrix} and the velocity have, uniformly on
$[0,\tp]$, the derivative bounds
\begin{equation}\label{nl:smooth-coefficients}
 \norm{\nabla^mM_{Q^{0,J}}}_\infty
       +\norm{\nabla^mQ^{0,J}}_\infty\leq C_mA.
\end{equation}
The matrix and velocity sums use $\sum_{j\geq0}R^{-(2+m)j}$
and $\sum_{j\geq0}R^{-(1+m)j}$, respectively.
For $1\leq\ell\leq\tilde\ell$ and $0\leq t\leq1$, direct Gaussian
integration of \eqref{bg:pair-pressure} gives
\begin{equation}\label{nl:smooth-residual}
 \norm{\nabla^m\partial_zP_{\ell,\tilde\ell}(t)}_2
       \leq C_mA^2b^{-1/2}\ell^{1/2}\tilde\ell^{-2},\qquad
 \sup_{[0,\tp]}\norm{F_Q(t)}_{H^m}\leq C_mA^2b^{-1/2}.
\end{equation}
For $m=0$, the pressure amplitude is $O(A^2\tilde\ell^{-2})$, its
horizontal $L^2$ factor is $O(\ell)$, and the axial derivative
factor is $O((b\ell)^{-1/2})$. Each extra derivative costs at most
$C_m/\ell$. Summation first in $\ell\leq\tilde\ell$ leaves
$C\tilde\ell^{-3/2}$, which is summable over $\tilde\ell$.

For the correction starting at zero, ordinary Sobolev energy with
$m\geq3$ yields
\begin{equation}\label{nl:regular-correction-energy}
 \frac{\diff}{\diff t}\norm e_{H^m}^2+\norm{\nabla e}_{H^m}^2
 \leq C_m(A+1)\norm e_{H^m}^2+C_m\norm e_{H^m}^4+C_mA^4b^{-1}.
\end{equation}
Since $A\tp=2+\Theta\log A$, a bootstrap from zero gives
\begin{equation}\label{nl:regular-correction}
 \sup_{[0,\tp]}\norm{e^{0,J}(t)}_{H^m}\leq C_mA^{P_m}b^{-1/2}.
\end{equation}
Applying this at finitely many higher orders gives the needed
supremum bounds.  Weak--strong uniqueness, using the energy identity of
Proposition~\ref{bg:reference}, identifies this smooth correction with
the reference solution on the initial interval.

Let $V_{\eta_{\mathrm{reg}}}^J$ be the linear response with initial value $\mathcal I_{\eta_{\mathrm{reg}}}$
about the axially invariant field obtained from $Q^{0,J}$ by replacing
each cutoff $\eta_{R^j}(t,z)$ by one. Because $\mathcal I\in H^{-4}$ and
$\Delta\mathcal I\in H^{-6}$,
\begin{equation}\label{nl:impulse-approximation}
 \norm{\mathcal I_{\eta_{\mathrm{reg}}}-\mathcal I}_{H^{-6}}
 \leq\int_0^{\eta_{\mathrm{reg}}}\norm{e^{u\Delta}\Delta\mathcal I}_{H^{-6}}\idiff u
 \leq C\eta_{\mathrm{reg}}.
\end{equation}
The coefficient estimates \eqref{nl:smooth-coefficients} make
negative-to-positive smoothing over $\tp$ polynomial in $A$.
The extra invariant matrices have all fixed derivative norms
$O_m(AR^{-2})$.  Split their Duhamel difference from the isolated
point response at $\tp/2$.  On the first half use the negative
source norm and the remaining positive smoothing time; on the
second use the already regularized response.  This proves
\begin{equation}\label{nl:invariant-preparation}
 \norm{V_{\eta_{\mathrm{reg}}}^J(\tp)-S_A}_{H^4}\leq CA^P(\eta_{\mathrm{reg}}+R^{-2}).
\end{equation}
In particular, no negative power of $\eta_{\mathrm{reg}}$ is needed for this
endpoint comparison.

The endpoint moment is also uniform as $\eta_{\mathrm{reg}}\to0$.  Write
$\Lambda=(1-\Delta)^{1/2}$ and $Z_{\eta_{\mathrm{reg}}}=\Lambda^{-4}V_{\eta_{\mathrm{reg}}}^J$.
The Bessel-smoothed distribution $\Lambda^{-4}\mathcal I$
is locally square integrable and has an $O(|x|^{-4})$ tail;
thus $M_\mu(Z_{\eta_{\mathrm{reg}}}(0))\leq C_\mu$ for $\mu<5/2$.
The weighted graph argument of Lemma~\ref{seed:prepared} applies
to this axially invariant comparison field using
\eqref{nl:smooth-coefficients}.
For every fixed $0\leq\mu<5/2$ and nonnegative integers $j,k$, it gives
\begin{equation}\label{nl:uniform-seed-moment}
 M_\mu(\Lambda^kZ_{\eta_{\mathrm{reg}}}(t))
       \leq C_{\mu,k}t^{-k/2}e^{C_{\mu,k}At}M_\mu(Z_{\eta_{\mathrm{reg}}}(0)),
 \qquad
 M_\mu(\Lambda^jV_{\eta_{\mathrm{reg}}}^J(\tp))\leq C_{\mu,j}A^{P_{\mu,j}}.
\end{equation}
The projection here is applied to a horizontal axisymmetric matrix
output, so its weighted bound is precisely Lemma~\ref{bg:weighted}.

For intermediate comparisons one may use the rough scaling
$\mathcal I_{\eta_{\mathrm{reg}}}(x)=\eta_{\mathrm{reg}}^{-2}\mathcal I_1(x/\sqrt{\eta_{\mathrm{reg}}})$.
The differentiated kernel satisfies
$|\nabla^j\mathcal I_1(x)|\leq C_j(1+|x|)^{-4-j}$, and hence
\begin{equation}\label{nl:rough-seed-moment}
 M_\mu(\Lambda^j\mathcal I_{\eta_{\mathrm{reg}}})
       \leq C_{\mu,j}\eta_{\mathrm{reg}}^{-5/4-j/2},\qquad
 0\leq\mu<5/2,\quad j=0,1,\ldots,\quad0<\eta_{\mathrm{reg}}\leq1.
\end{equation}
The weighted graph energy used in \eqref{nl:uniform-seed-moment}, now
started from \eqref{nl:rough-seed-moment}, also gives the intermediate-time
bound
\begin{equation}\label{nl:rough-response-moment}
 \sup_{0\leq t\leq\tp}M_\mu(\Lambda^jV_{\eta_{\mathrm{reg}}}^J(t))
       \leq C_{\mu,j}A^{P_{\mu,j}}\eta_{\mathrm{reg}}^{-5/4-j/2}.
\end{equation}
Let $v_{\eta_{\mathrm{reg}}}^J$ be the actual linearization about $\bar u^{0,J}$
with the same initial value.  In the space
$\norm v_{\mathcal Z}=\norm v_{H^4}+M_{s_\star}(v)$,
the propagator and the vector and tensor source maps on $[0,\tp]$
have polynomial constants.  To see this directly, the heat
Duhamel equation contains the matrix term of norm $CA$ and a
reference-correction tensor kernel bounded by
$C\epsilon_A(1+(t-u)^{-1/2})$, with
$\epsilon_A\leq CA^Pb^{-1/2}$ from \eqref{nl:regular-correction}.
After weighting the time supremum by $e^{-K(A+1)t}$, the two
kernel integrals are bounded by
\begin{equation}\label{nl:volterra-kernels}
 \frac{CA}{K(A+1)},\qquad
 C\epsilon_A\left(\frac1{K(A+1)}
                    +\sqrt{\frac\pi{K(A+1)}}\right).
\end{equation}
A fixed large $K$ absorbs both, while removing the weight costs
only $e^{K(A+1)\tp}=CA^P$.

The cutoff comparison uses the global bounds
\begin{equation}\label{nl:cutoff-fractional}
 \begin{aligned}
 |1-\eta_\ell(t,z)|&\leq C_{\chi_{\mathrm{seed}}}(b\ell)^{-\chi_{\mathrm{seed}}}
                   \langle z\rangle^{\chi_{\mathrm{seed}}},\\
 |\partial_z\eta_\ell(t,z)|&\leq C_{\chi_{\mathrm{seed}}}(b\ell)^{-1-\chi_{\mathrm{seed}}}
                   |z|^{\chi_{\mathrm{seed}}},
 \qquad\ell\geq1,\quad0\leq t\leq1.
 \end{aligned}
\end{equation}
Their differentiated forms follow directly from the Gaussian.
Since $s_\star+\chi_{\mathrm{seed}}=39/16<5/2$, the full matrix cutoff
difference applied to $V_{\eta_{\mathrm{reg}}}^J$ costs
$CAb^{-\chi_{\mathrm{seed}}}M_{s_\star+\chi_{\mathrm{seed}}}(V_{\eta_{\mathrm{reg}}}^J)$, with a summable
scale coefficient $\ell^{-2-\chi_{\mathrm{seed}}}$.  The actual reference-correction difference is
$-\PP(e^{0,J}\cdot\nabla V_{\eta_{\mathrm{reg}}}^J+V_{\eta_{\mathrm{reg}}}^J\cdot\nabla e^{0,J})$.
Its $H^4$ norm uses at most five derivatives of $V_{\eta_{\mathrm{reg}}}^J$;
its moment uses the same first derivative and ordinary moment.
The largest loss in \eqref{nl:rough-response-moment} is
$\eta_{\mathrm{reg}}^{-15/4}\leq\eta_{\mathrm{reg}}^{-4}$.  The zero-data difference equation
and \eqref{nl:volterra-kernels} consequently give
\begin{equation}\label{nl:actual-preparation}
 \sup_{[0,\tp]}\norm{v_{\eta_{\mathrm{reg}}}^J-V_{\eta_{\mathrm{reg}}}^J}_{\mathcal Z}
       \leq CA^P\eta_{\mathrm{reg}}^{-4}(b^{-\chi_{\mathrm{seed}}}+b^{-1/2}).
\end{equation}
The extra moment used here belongs to the prescribed regularized
seed and was proved in \eqref{nl:uniform-seed-moment}; it is not
assumed for an arbitrary member of $\mathcal C_A$.

It remains to solve the nonlinear preparation.  Set
$H_\varepsilon=u^{0,J}_\varepsilon-\bar u^{0,J}$; its initial
value is $(\varepsilon/n_A)\mathcal I_{\eta_{\mathrm{reg}}}$ and its exact tensor
source is $-\PP\div(H_\varepsilon\otimes H_\varepsilon)$.
The algebra estimate
$\norm{v\otimes w}_{H^4}+M_{s_\star}(v\otimes w)
\leq C\norm v_{\mathcal Z}\norm w_{\mathcal Z}$
and the polynomial source map give a contraction, since
\begin{equation}\label{nl:initial-smallness}
 \sup_{[0,\tp]}\norm{(\varepsilon/n_A)v_{\eta_{\mathrm{reg}}}^J}_{\mathcal Z}
       \leq CA^P|\varepsilon|\eta_{\mathrm{reg}}^{-13/4},\qquad
 |\varepsilon|\eta_{\mathrm{reg}}^{-13/4}\leq2e^{-35\varrho A/4}.
\end{equation}
The entire nonlinear remainder obeys
\begin{equation}\label{nl:initial-nonlinear-error}
 \sup_{[0,\tp]}\norm{H_\varepsilon
                       -(\varepsilon/n_A)v_{\eta_{\mathrm{reg}}}^J}_{\mathcal Z}
       \leq CA^P\varepsilon^2\eta_{\mathrm{reg}}^{-13/2}.
\end{equation}
Combining \eqref{nl:invariant-preparation},
\eqref{nl:uniform-seed-moment}, \eqref{nl:actual-preparation}, and
\eqref{nl:initial-nonlinear-error} proves
\eqref{nl:preparation-error}.  The error rates are explicit:
$\eta_{\mathrm{reg}}=e^{-\varrho A}$, $R^{-2}=e^{-2\varrho A}$,
$\eta_{\mathrm{reg}}^{-4}(b^{-\chi_{\mathrm{seed}}}+b^{-1/2})\leq2e^{-2\varrho A}$, and
$|\varepsilon|\eta_{\mathrm{reg}}^{-13/2}\leq2e^{-11\varrho A/2}$.

For the negative path, direct normalization gives
\begin{equation}\label{nl:initial-normalization}
 B_J(\lambda)=\lambda\bigl(1+\langle d_{J,-\lambda},g\rangle\bigr),
 \qquad
 F_J(\lambda)-g
 =\frac{d_{J,-\lambda}-\langle d_{J,-\lambda},g\rangle g}
        {1+\langle d_{J,-\lambda},g\rangle}.
\end{equation}
The denominator is positive for large $A$.  The strict inequality
$\varrho>\rho_c$ and the polynomial moment bound place this shape
in $\mathcal C_A$.  Moreover
$B_J(\lambda)=\lambda(1+O(A^Pe^{-\varrho A}))$, which proves
\eqref{nl:preparation-crossings}.  Continuous dependence follows
from the same contraction and the positive normalization.
\end{proof}

\subsection{Selection of an exact finite cascade}

We now select the initial scalar parameter so that the exact finite
cascade has a prescribed terminal amplitude.  The transition map
need not be monotone in the amplitude or the shape.  The argument uses
only the endpoint inequalities \eqref{nl:crossings} and
\eqref{nl:preparation-crossings}, together with continuity of the maps.
An auxiliary clipping map produces a continuous scalar path; an
interior terminal value will imply that clipping was inactive along
the selected orbit.

\begin{lemma}[Finite shooting]\label{nl:shooting}
Fix $a_\ast\in(a_-,a_+)$ and an integer $N\geq0$.
There are $\lambda_N\in(a_-/2,2a_+)$ and a genuine finite orbit
$(a_j^{(N)},f_j^{(N)})_{j=0}^N$ of the nonlinear transition maps, with
$a_j^{(N)}\in(a_-,a_+)$ and $f_j^{(N)}\in\mathcal C_A$, such that
\begin{equation}\label{nl:orbit}
 \begin{gathered}
 a_0^{(N)}=B_N(\lambda_N),\qquad f_0^{(N)}=F_N(\lambda_N),\\
 a_{j+1}^{(N)}=\mathcal A^{j,N-j}(a_j^{(N)},f_j^{(N)}),\qquad
 f_{j+1}^{(N)}=\mathcal F^{j,N-j}(a_j^{(N)},f_j^{(N)}),\quad0\leq j<N,\\
 a_N^{(N)}=a_\ast.
 \end{gathered}
\end{equation}
\end{lemma}

\begin{proof}
Fix $N$ and suppress the superscript $(N)$ within this proof.
Let $\pi_{[a_-,a_+]}$ be the continuous clipping map from $\R$
to $[a_-,a_+]$.  On the entire interval
$a_-/2\leq\lambda\leq2a_+$, define the continuous initial path
$a_0(\lambda)=\pi_{[a_-,a_+]}B_N(\lambda)$,
$f_0(\lambda)=F_N(\lambda)$, and then recursively set
\begin{equation}\label{nl:clipped-orbit}
 f_{j+1}(\lambda)=\mathcal F^{j,N-j}(a_j(\lambda),f_j(\lambda)),
 \qquad
 a_{j+1}(\lambda)=\pi_{[a_-,a_+]}
          \mathcal A^{j,N-j}(a_j(\lambda),f_j(\lambda)).
\end{equation}
Proposition~\ref{nl:transition} makes this recursion well defined
and continuous, including when the unmodified output amplitude
leaves the interval.  By \eqref{nl:crossings}, either boundary
amplitude is absorbing under the clipped recursion, regardless
of the shape.  Equation~\eqref{nl:preparation-crossings} therefore
implies that the terminal amplitude is $a_-$ at the left endpoint
and $a_+$ at the right endpoint of the path.

The intermediate value theorem gives a $\lambda_N$ whose terminal
amplitude is $a_\ast$.  No earlier amplitude on that path can be
a boundary amplitude, since the terminal amplitude would then be
the same boundary value.  The initial and all subsequent clippings
are consequently inactive on this selected orbit.  The recursion
is precisely \eqref{nl:orbit}, and each stage is the actual
nonlinear solution of Proposition~\ref{nl:transition}.
\end{proof}

\section{The common initial datum and the suitable solutions}
\label{sec:limit}

At finite depth we have exact smooth solutions whose initial difference is
concentrated at the smallest retained scale.  That difference tends to zero
in physical $L^2$, while an observation at the largest scale remains
prescribed.  We pass these two properties to the limit in the same sequence
of finite constructions.

The first task is energy-class compactness, including spatial tightness and
the common strong initial trace. The bounds between successive scales then
give positive-time regularity, which together with the strong trace extends
the energy equalities through time zero. Finally, fixed physical tests
preserve the selected observations and distinguish the branches.

\subsection{Physical solutions at finite depth}

All rescaled transitions must refer to one physical time and one physical
velocity.  The covariance of the exact reference performs this matching,
including its correction.  The formulas below also exhibit the small
factor $\ell_N^{1/2}$ in the initial energy norm, while leaving the
terminal amplitude independent of depth.

Fix a terminal amplitude $a_\ast\in(a_-,a_+)$ and let $\lambda_N$ be the
parameter selected by Lemma~\ref{nl:shooting}.  Define
\begin{equation}
 \ell_n=R^{-n},\qquad T_n=\ell_n^2\tp,\qquad T_0=\tp,
 \label{limit:scales}
\end{equation}
and the physical reference solution by
\begin{equation}
 \bar u_N(t,x)=\ell_N^{-1}
    \bar u^{0,N}(t/\ell_N^2,x/\ell_N),\qquad
 \bar u_N(0,x)=\sum_{n=0}^NQ_{\ell_n}(0,x).
 \label{limit:reference}
\end{equation}
Its physical correction is
\begin{equation}\label{limit:physical-correction}
 e_N(t,x):=\ell_N^{-1}e^{0,N}(t/\ell_N^2,x/\ell_N)
          =e^{N,0}(t,x),
\end{equation}
where the second equality is the exact covariance \eqref{bg:covariance}.
Thus $\bar u_N=\sum_{n=0}^NQ_{\ell_n}+e_N$.
Start the perturbed solution from
\begin{equation}
 U_N(0,x)=\sum_{n=0}^NQ_{\ell_n}(0,x)
       -\frac{\lambda_N}{n_A}\ell_N^{-1}\mathcal I_{\eta_{\mathrm{reg}}}(x/\ell_N).
 \label{limit:finite-data}
\end{equation}
Proposition~\ref{nl:preparation} identifies its state at $T_N$ with the
first state of the selected finite orbit.  Proposition~\ref{nl:transition}
and the covariance \eqref{bg:covariance} then identify the successive
states in the same physical coordinates.  The strong solutions coincide
at the joining times and therefore concatenate into one exact solution
on $[0,T_0]$.  Every change of variables is the Navier--Stokes scaling,
so the viscosity remains one.

The chronological index $j$ in \eqref{nl:orbit} increases from the
smallest to the largest radius, whereas the spatial index $n$ decreases.
We therefore set
\begin{equation}\label{limit:orbit-indices}
 a_{N,n}:=a_{N-n}^{(N)},\qquad f_{N,n}:=f_{N-n}^{(N)},
 \qquad 0\leq n\leq N.
\end{equation}
Thus $n=N$ is the prepared smallest scale and $n=0$ is the terminal
largest scale.  The numbers $a_{N,n}\in(a_-,a_+)$ and shapes
$f_{N,n}\in\mathcal C_A$ satisfy
\begin{equation}
 U_N(T_n,x)-\bar u_N(T_n,x)
      =-\ell_n^{-1}a_{N,n}f_{N,n}(x/\ell_n),
 \qquad a_{N,0}=a_\ast,
 \qquad 0\leq n\leq N.
 \label{limit:physical-orbit}
\end{equation}
Both finite solutions are smooth and have finite energy.  Their pressures
are the full Riesz pressures
\begin{equation}
 p_N=\sum_{i,j}\mathcal R_i\mathcal R_j(U_{N,i}U_{N,j}),\qquad
 \bar p_N=\sum_{i,j}\mathcal R_i\mathcal R_j(\bar u_{N,i}\bar u_{N,j}).
 \label{limit:finite-pressure}
\end{equation}
They satisfy the global and local energy identities.

Recalling \eqref{bg:datum}, the common limiting datum is
\begin{equation}
 u_0(x)=\sum_{n=0}^\infty\frac{A}{\ell_n^2}
       \exp\left(-\frac{|x_h|^2}{4\ell_n^2}
                   -\frac{z^2}{2b^2\ell_n^2}\right)x_h^\perp.
 \label{limit:initial-data}
\end{equation}
Gaussian integration gives
$\|Q_{\ell_n}(0)\|_{\mathfrak q}=C_{\mathfrak q}Ab^{1/\mathfrak q}\ell_n^{3/\mathfrak q-1}$.
Thus the series converges in every $L^{\mathfrak q}$, $1\leq \mathfrak q<3$, and is divergence
free.  Every nonzero summand has the same azimuthal orientation, so
$u_0\neq0$.  On any compact set avoiding the origin, each differentiated
summand has a Gaussian factor bounded by
$\exp(-c_b|x|^2/\ell_n^2)$, and its derivative series converges uniformly.
Hence $u_0$ is smooth away from the origin.

To estimate its weak-$L^3$ norm, set
$\rho(x)=(|x_h|^2+z^2/b^2)^{1/2}$. Splitting the radii at
$\ell_n=\rho$ and using $e^{-s^2/4}\leq Cs^{-4}$ for $s>0$ gives
\[
 \begin{aligned}
 \sum_{n=0}^\infty\ell_n^{-2}
          e^{-\rho^2/(4\ell_n^2)}
 &\leq\sum_{\ell_n\geq\rho}\ell_n^{-2}
       +C\rho^{-4}\sum_{\ell_n<\rho}\ell_n^2\\
 &\leq C_R\rho^{-2},\qquad \rho>0.
 \end{aligned}
\]
Here both sums are geometric, and $C_R$ depends only on the fixed
ratio $R$. Since $|x_h|\leq\rho$, the datum in
\eqref{limit:initial-data} therefore obeys
$|u_0(x)|\leq C_RA/\rho(x)$ for $x\ne0$.
The ellipsoid $\{x:\rho(x)<a\}$ has volume $4\pi ba^3/3$ for
$a>0$. The distribution function and the weak-$L^3$ quasinorm thus satisfy
\begin{equation}\label{limit:datum-weak-three}
 \bigl|\{x:|u_0(x)|>\lambda\}\bigr|
       \leq C_RbA^3\lambda^{-3},\qquad \lambda>0,
 \qquad
 \|u_0\|_{L^{3,\infty}}\leq C_RAb^{1/3}.
\end{equation}

The geometric sum and the scaling of the regularized input give
\begin{equation}
 \left\|u_0-\sum_{n=0}^NQ_{\ell_n}(0)\right\|_2
       \leq CA\sqrt b\,\frac{\ell_{N+1}^{1/2}}{1-R^{-1/2}},\qquad
 \left\|\frac{\lambda_N}{n_A}\ell_N^{-1}
                  \mathcal I_{\eta_{\mathrm{reg}}}(\cdot/\ell_N)\right\|_2
       \leq C(A)\ell_N^{1/2}.
 \label{limit:initial-convergence}
\end{equation}
Here $\eta_{\mathrm{reg}}$ is fixed before $N\to\infty$, and
$\|\mathcal I_{\eta_{\mathrm{reg}}}\|_2\leq C\eta_{\mathrm{reg}}^{-5/4}$.
Consequently $U_N(0)$ and $\bar u_N(0)$ both converge strongly to $u_0$
in $L^2$.

To see that $u_0\notin L^3$, let
$D_n=\{x:\ell_n\leq r\leq2\ell_n,\ |z|\leq\ell_n\}$.
These cylinders are pairwise disjoint because $R>4$.
On $D_n$, the $n$th summand in \eqref{limit:initial-data} has magnitude
at least $Ae^{-3/2}/\ell_n$, using $b\geq1$.
All summands have the same azimuthal orientation, and
$|D_n|=6\pi\ell_n^3$. Therefore
\begin{equation}\label{limit:datum-not-three}
 \int_{D_n}|u_0(x)|^3\idiff x
 \geq |D_n|\left(\frac{Ae^{-3/2}}{\ell_n}\right)^3
 =6\pi e^{-9/2}A^3,
\end{equation}
which proves the claimed failure of $L^3$ integrability after summing in $n$.

In fact, $u_0$ belongs to no $L^{3,\mathfrak r}$ with $1\leq \mathfrak r<\infty$.
In this assertion, $\mathfrak r$ is the second Lorentz exponent.
Let $c_0=e^{-3/2}$ and
$\mu_0(\lambda)=|\{x:|u_0(x)|>\lambda\}|$.  The same lower bound on
$D_n$ shows that, for
\[
 J_n=\left(\frac{c_0A}{2\ell_n},\frac{c_0A}{\ell_n}\right),
 \qquad \lambda\in J_n,
\]
we have $\mu_0(\lambda)\geq6\pi\ell_n^3$ and hence
$\lambda\mu_0(\lambda)^{1/3}\geq
 \frac12c_0A(6\pi)^{1/3}$.  Since $R>4$, the intervals $J_n$
are disjoint and each has logarithmic length $\log2$.  Thus,
by the Lorentz quasinorm convention \eqref{app:lorentz-norm},
for every $1\leq \mathfrak r<\infty$,
\begin{equation}\label{limit:datum-not-finite-lorentz}
 \begin{aligned}
 \|u_0\|_{L^{3,\mathfrak r}}^{\mathfrak r}
 &=3\int_0^\infty
       \bigl(\lambda\mu_0(\lambda)^{1/3}\bigr)^{\mathfrak r}
            \frac{\idiff\lambda}{\lambda}\\
 &\geq3\left(\frac{c_0A}{2}(6\pi)^{1/3}\right)^{\mathfrak r}
      \sum_{n=0}^\infty\int_{J_n}\frac{\idiff\lambda}{\lambda}
 =\infty.
 \end{aligned}
\end{equation}

\subsection{Global compactness and the initial energy inequality}

We first obtain convergence against fixed spatial tests uniformly in
time; this will preserve the terminal coordinate. A uniform spatial
tail estimate then upgrades local compactness to global strong
$L^2_{t,x}$ convergence. Interpolation gives the $L^3_{t,x}$ convergence
needed for the energy fluxes. Finally, the energy inequality and weak
continuity at zero recover the common strong initial trace.

The two finite energy identities and \eqref{limit:initial-convergence}
give a constant $E_A$, independent of $N$ and $a_\ast$, such that
\begin{equation}
 \sup_{0\leq t\leq T_0}
       \bigl(\|U_N(t)\|_2^2+\|\bar u_N(t)\|_2^2\bigr)
 +\int_0^{T_0}
       \bigl(\|\nabla U_N\|_2^2+\|\nabla\bar u_N\|_2^2\bigr)\idiff t
       \leq E_A.
 \label{limit:uniform-energy}
\end{equation}

The following argument is applied to either sequence, denoted by $v_N$,
and the extractions are made jointly. In this argument, $\pi_N$ denotes
the Riesz pressure associated with $v_N$, and $\pi$ denotes the pressure
of its limit $v$.

\begin{proposition}[Compactness of the finite solutions]
\label{limit:compactness}
Along a subsequence, $U_N$ and $\bar u_N$ converge to suitable Leray--Hopf
solutions $U$ and $\bar u$ with the same strong initial trace $u_0$.
The convergence is strong in $L^2((0,T_0)\times\R^3)$, weak in the gradient
energy space, and uniform in time against every fixed $L^2$ test.
\end{proposition}

\begin{proof}
Energy interpolation bounds $v_N$ in $L^{10/3}_{t,x}$ and its pressure
in $L^{5/3}_{t,x}$.  Moreover,
\begin{equation}
 \|\PP\operatorname{div}(v_N\otimes v_N)\|_{H^{-1}}
       \leq C\|v_N\|_4^2
       \leq C\|v_N\|_2^{1/2}\|\nabla v_N\|_2^{3/2}.
 \label{limit:time-derivative}
\end{equation}
Together with the heat term, this bounds $\partial_tv_N$ in
$L^{4/3}(0,T_0;H^{-1})$. Compactness on bounded spatial balls gives
strong local $L^2_{t,x}$ convergence. For every $\varphi\in H^1$,
H\"older's inequality also gives
\begin{equation}\label{limit:equicontinuity}
 \left|\langle v_N(t)-v_N(s),\varphi\rangle\right|
 \leq |t-s|^{1/4}\|\partial_tv_N\|_{L^{4/3}(0,T_0;H^{-1})}
                    \|\varphi\|_{H^1}.
\end{equation}
A diagonal extraction and density in $L^2$, using
\eqref{limit:uniform-energy}, therefore give
\begin{equation}
 \langle v_N(t),\varphi\rangle\longrightarrow\langle v(t),\varphi\rangle
 \quad\hbox{uniformly for }t\in[0,T_0],\qquad \varphi\in L^2.
 \label{limit:weak-time}
\end{equation}
This convergence includes every geometric time in \eqref{limit:scales}.

To specify the representative, take a countable $L^2$-dense set of
$H^1$ tests and apply Arzel\`a--Ascoli to their pairings.  The uniform
energy bound extends the limiting functionals to all of $L^2$,
uniformly in time.  Riesz representation defines $v(t)$ for every
$t$, and the same density argument gives weak continuity.  This
representative agrees almost everywhere with the local strong limit.
The values at observation times and energy endpoints are consequently
those of this same weakly continuous representative.

To obtain global strong convergence, choose an exterior cutoff $\chi_K$
that vanishes on $|x|\leq K$, equals one on $|x|\geq2K$, and satisfies
$|\nabla\chi_K|\leq C/K$ and $|\Delta\chi_K|\leq C/K^2$.
The smooth local energy identity gives
\begin{align}
 \sup_{0\leq t\leq T_0}\int_{|x|\geq2K}|v_N(t,x)|^2\idiff x
 &\leq\int_{|x|\geq K}|v_N(0,x)|^2\idiff x+\frac{CE_AT_0}{K^2}
 \notag\\
 &\quad+\frac C K\int_0^{T_0}
       \bigl(\|v_N(t)\|_3^3+
             \|\pi_N(t)\|_{3/2}\|v_N(t)\|_3\bigr)\idiff t.
 \label{limit:spatial-tail}
\end{align}
An additional exterior compact cutoff first makes the test admissible.
It can be removed because the Riesz estimate gives
$\|\pi_N(t)\|_{3/2}\leq C\|v_N(t)\|_3^2$ and
\begin{equation}
 \int_0^{T_0}\|v_N(t)\|_3^3\idiff t
 \leq C\left(\sup_t\|v_N(t)\|_2\right)^{3/2}
       \left(\int_0^{T_0}\|\nabla v_N(t)\|_2^2\idiff t\right)^{3/4}T_0^{1/4}.
 \label{limit:integrable-flux}
\end{equation}
The initial tails tend to zero uniformly in $N$ by strong initial
convergence.  Equation \eqref{limit:spatial-tail} therefore makes the
spatial energy tails uniformly small in both $N$ and $t$.

For the limit's tails, multiply by the exterior indicator and use weak
lower semicontinuity at each time; the same uniform tail bound passes
to $v$.  Split the spacetime norm of $v_N-v$ into a fixed ball and its
complement to obtain global, rather than merely local, strong $L^2$
convergence.  The additional extraction can be chosen so that
$\sum_N\|v_N-v\|_{L^2_{t,x}}^2<\infty$; Fubini then gives the
claimed strong convergence times.  Make this extraction jointly for
the two sequences.

Interpolation with the uniform $L^{10/3}_{t,x}$ bound gives strong
$L^3_{t,x}$ convergence through the estimate
\[
 \|v_N-v\|_{L^3_{t,x}}
 \leq \|v_N-v\|_{L^2_{t,x}}^{1/6}
       \|v_N-v\|_{L^{10/3}_{t,x}}^{5/6}\longrightarrow0.
\]
The quadratic tensors consequently converge strongly in
$L^{3/2}_{t,x}$, and boundedness of the Riesz transforms gives
$\pi_N\to\pi$ in the same space. H\"older's inequality then yields
$\pi_Nv_N\to\pi v$ in $L^1_{t,x}$, which is the pressure-flux
convergence needed for the local energy inequality.  The uniform
$L^{5/3}_{t,x}$ pressure bound identifies the same weak limit in that
space.  In particular, the limiting pressures are
\begin{equation}
 p=\sum_{i,j}\mathcal R_i\mathcal R_j(U_iU_j),\qquad
 \bar p=\sum_{i,j}\mathcal R_i\mathcal R_j(\bar u_i\bar u_j).
 \label{limit:pressures}
\end{equation}
These convergences pass the momentum equations, including the strongly
convergent initial terms.  They also pass the cubic and pressure fluxes
in the integrated local energy identities.  Lower semicontinuity of the
localized dissipation proves the local energy inequality with initial
term $\frac12\int|u_0|^2\phi(0)\idiff x$.

For every $t\in[0,T_0]$, weak convergence at that endpoint and weak
lower semicontinuity in the dissipation give
\begin{equation}
 \frac12\|U(t)\|_2^2+\int_0^t\|\nabla U(s)\|_2^2\idiff s
       \leq\frac12\|u_0\|_2^2.
 \label{limit:energy-zero}
\end{equation}
Weak continuity at zero identifies $U(0)=u_0$.
The inequality in \eqref{limit:energy-zero} and weak lower
semicontinuity give matching upper and lower bounds on
$\lim_{t\downarrow0}\|U(t)\|_2$; hence $U(t)\to u_0$ strongly in $L^2$.
The same argument applies to $\bar u$.

Let $\mathcal T\subset(0,T_0)$ be the common full-measure set of strong
$L^2_x$ convergence times for the two sequences.
For every $s\in\mathcal T$ and every $t\geq s$, passing the finite energy
identities gives
\begin{equation}
 \frac12\|U(t)\|_2^2+\int_s^t\|\nabla U(\tau)\|_2^2\idiff\tau
       \leq\frac12\|U(s)\|_2^2,
 \label{limit:energy-ae-start}
\end{equation}
and the same inequality for $\bar u$.
The right side converges strongly, whereas lower semicontinuity is used
only on the left.

The same subsequence works for every $t\geq s$: convergence of endpoint
pairings already holds for all times, and restriction of the global
weak gradient convergence to any fixed interval $[s,t]$ preserves
weak convergence.  The full-measure set of admissible starting times
is therefore chosen once, not anew for each terminal time.

This proves the full Leray--Hopf requirements.
\end{proof}

\subsection{Mixed norms, the Lorentz endpoint, and positive-time regularity}

The same scale bounds are used for every pair of exponents; no new datum
or subsequence is chosen for each mixed norm.  Their time integrability
explains the strict inequality at the Serrin boundary.  Away from zero
they give bounded velocities, hence smoothness.  Only after that step do
we recover energy equalities through the initial time.

The constants below are uniform in the depth $N$ and the terminal
coordinate $a_\ast$ once $A$, $R$, and $\tp$ have been fixed. They may
depend on a chosen spatial exponent $\mathfrak q$.

The estimate between successive scales in
Proposition~\ref{nl:transition} controls the physical perturbation
$h_N=U_N-\bar u_N$ throughout each interval $T_n\leq t\leq T_{n-1}$:
\begin{equation}
 \|h_N(t)\|_2\leq C(A)\ell_n^{1/2},\qquad
 \|h_N(t)\|_\infty\leq C(A)\ell_n^{-1}.
 \label{limit:interstage}
\end{equation}
Since $R$ and $\tp$ are fixed, these imply
\begin{equation}
 \|h_N(t)\|_2\leq C_At^{1/4},\qquad
 \|h_N(t)\|_\infty\leq C_At^{-1/2},\qquad T_N\leq t\leq T_0.
 \label{limit:perturbation-time}
\end{equation}

Indeed, $T_n\leq t\leq T_{n-1}=R^2T_n$ implies
\[
 \ell_n^{1/2}\leq\tp^{-1/4}t^{1/4},\qquad
 \ell_n^{-1}\leq R\sqrt{\tp}\,t^{-1/2}.
\]
The constants in these two inequalities are independent of $N$.

For each fixed $t>0$, the lower restriction holds for all sufficiently
large $N$.  Weak convergence at that time and distributional closedness
of a bounded $L^\infty$ ball pass both bounds to $h=U-\bar u$.
Interpolation then gives
\begin{equation}
 \|h(t)\|_{\mathfrak q}\leq C_{A,\mathfrak q}t^{(3/\mathfrak q-1)/2},\qquad 2\leq \mathfrak q\leq\infty.
 \label{limit:perturbation-q}
\end{equation}

For the reference solution, Gaussian integration in \eqref{hist:array}
gives
\begin{equation}
 \|Q_\ell(t)\|_{\mathfrak q}
  =C_{\mathfrak q} A b\,\ell^3(\ell^2+t)^{-3/2+1/\mathfrak q}
            (b^2\ell^2+2t)^{-1/2+1/(2\mathfrak q)}
  \leq C_{\mathfrak q} A b\,\ell^3(\ell^2+t)^{-2+3/(2\mathfrak q)}.
 \label{limit:gaussian-q}
\end{equation}
The supremum interpretation gives the same expression for $\mathfrak q=\infty$.
Split the physical radii at $\ell=\sqrt t$.
The smaller radii have a geometric sum of $\ell^3$ bounded by $Ct^{3/2}$.
For the larger radii, sum $C_{A,\mathfrak q}\ell^{3/\mathfrak q-1}$.
This sum is bounded for $\mathfrak q<3$, is $O(1+|\log t|)$ for $\mathfrak q=3$, and is
$O(t^{(3/\mathfrak q-1)/2})$ for $\mathfrak q>3$.
The correction in \eqref{hist:reference-error} satisfies
$\|e_N(t)\|_2\leq C_At^{1/4}$ and
$\|e_N(t)\|_\infty\leq C_At^{-1/2}$ uniformly in the physical cutoff.
Interpolation and \eqref{limit:perturbation-q} therefore yield
\begin{equation}
 \|U(t)\|_{\mathfrak q}+\|\bar u(t)\|_{\mathfrak q}\leq C_{A,\mathfrak q}
 \begin{cases}
 1,&2\leq \mathfrak q<3,\\
 1+|\log t|,&\mathfrak q=3,\\
 t^{(3/\mathfrak q-1)/2},&3<\mathfrak q\leq\infty,
 \end{cases}
 \qquad 0<t\leq T_0.
 \label{limit:mixed-bounds}
\end{equation}

\begin{proposition}[Regularity and both energy identities]
\label{limit:mixed}
Both $U$ and $\bar u$ belong simultaneously to every space
\begin{equation}
 L^{\mathfrak p}(0,T_0;L^{\mathfrak q}(\R^3)),\qquad
 1\leq \mathfrak p\leq\infty,\quad 2\leq \mathfrak q\leq\infty,\quad
 \frac2{\mathfrak p}+\frac3{\mathfrak q}>1.
 \label{limit:strict-range}
\end{equation}
They also satisfy
\begin{equation}\label{limit:lorentz-bound}
 \sup_{0<t\leq T_0}
 \left(\|U(t)\|_{L^{3,\infty}}
          +\|\bar u(t)\|_{L^{3,\infty}}\right)\leq C(A),
\end{equation}
where $C(A)$ is independent of the terminal coordinate $a_\ast$.
They are smooth for positive times, belong to $C([0,T_0];L^2)$, and
satisfy the global and local energy identities through time zero.
\end{proposition}

\begin{proof}
For $\mathfrak q<3$, the bound in \eqref{limit:mixed-bounds} lies in every time
$L^{\mathfrak p}$, including $\mathfrak p=\infty$.  For $\mathfrak q=3$, the logarithm lies in every
finite time $L^{\mathfrak p}$.  For $\mathfrak q>3$, the power is integrable precisely when
$\mathfrak p(1-3/\mathfrak q)/2<1$, which is \eqref{limit:strict-range}.
At $\mathfrak q=\infty$ this includes every $\mathfrak p<2$.

For the Lorentz endpoint, write $\bar u=Q+e$, where
$Q(t)=e^{t\Delta}u_0$, and retain $h=U-\bar u$.
The heat kernel $G_t$ is nonnegative and has integral one.
For every measurable set $E\subset\R^3$ of finite measure, Fubini's
theorem and \eqref{intro:weak-three} give
\[
 \begin{aligned}
 \int_E|Q(t,x)|\idiff x
 \leq\int_{\R^3}G_t(y)
           \int_{E-y}|u_0(x)|\idiff x\idiff y
 \leq\frac32\|u_0\|_{L^{3,\infty}}|E|^{2/3}.
 \end{aligned}
\]
The last inequality follows by integrating
$\min\bigl(|E|,\|u_0\|_{L^{3,\infty}}^3\lambda^{-3}\bigr)$
over $\lambda>0$.
Taking $E=\{x:|Q(t,x)|>\lambda\}$, which has finite measure because
$Q(t)\in L^2$, yields
\begin{equation}\label{limit:heat-weak-three}
 \sup_{t>0}\|Q(t)\|_{L^{3,\infty}}
       \leq\frac32\|u_0\|_{L^{3,\infty}}
       \leq C_RAb^{1/3},
\end{equation}
where the second inequality is \eqref{limit:datum-weak-three}.
At $\mathfrak q=3$, \eqref{limit:perturbation-q} bounds $h$ uniformly in
$L^3$. For the correction, \eqref{bg:correction} gives
\begin{equation}\label{limit:correction-three}
 \|e(t)\|_3
       \leq\|e(t)\|_2^{2/3}\|e(t)\|_\infty^{1/3}
       \leq C(A) t^{1/6}t^{-1/6}=C(A),
       \qquad 0<t\leq T_0.
\end{equation}
Since $L^3\subset L^{3,\infty}$, the quasitriangle inequality,
\eqref{limit:heat-weak-three}, \eqref{limit:correction-three}, and
\eqref{limit:perturbation-q} at $\mathfrak q=3$ prove
\eqref{limit:lorentz-bound}. Each constant is independent of
$a_\ast$, as are the bounds in \eqref{limit:perturbation-time}.

Equation \eqref{limit:mixed-bounds} bounds each limiting velocity in
$L^\infty_x$ on every closed positive-time interval.
Choose an earlier time $s\in\mathcal T$ at which the velocity belongs to
$H^1$.  The local strong solution from that value agrees with the
Leray--Hopf limit by weak--strong uniqueness.  While they agree,
\begin{equation}
 \frac{\diff}{\diff t}\|\nabla v(t)\|_2^2+\|\Delta v(t)\|_2^2
       \leq C\|v(t)\|_\infty^2\|\nabla v(t)\|_2^2,
       \qquad v=U,\bar u.
 \label{limit:H1-continuation}
\end{equation}
The coefficient is integrable on the chosen positive-time interval.
This estimate continues the strong solution throughout that interval,
and parabolic regularity gives smoothness for every $t>0$.

Between positive times, smooth finite-energy solutions satisfy the energy
identity.  To remove the spatial cutoff, use
\eqref{limit:integrable-flux}; on a positive-time interval the pressure
also belongs to $L^2_{t,x}$ and its product with the velocity is integrable.
Sending the starting time to zero uses the strong initial trace and finite
total dissipation.  Thus we obtain
\begin{equation}
 \frac12\|v(t)\|_2^2+\int_s^t\|\nabla v(\tau)\|_2^2\idiff\tau
       =\frac12\|v(s)\|_2^2,
 \qquad 0\leq s\leq t\leq T_0,\qquad v=U,\bar u.
 \label{limit:energy-equality}
\end{equation}
The norm is continuous by \eqref{limit:energy-equality}; weak continuity
then gives strong $L^2$ continuity.
Let $\pi$ denote the pressure corresponding to $v$ in
\eqref{limit:pressures}. Then $v\in L^{10/3}_{t,x}$ and
$\pi\in L^{5/3}_{t,x}$ give local integrability of $|v|^3$ and
$\pi v$ through zero.
The local energy identity between positive times therefore extends to
zero by the strong trace and finite dissipation.
\end{proof}

\subsection{Fixed observations and the family of branches}

Distinct approximating solutions can converge to the same limit, so a
quantitative observation must survive compactness.  The normalized shape
coordinate provides such an observation at each fixed geometric time.
One compactly supported test separates every branch from the reference;
different terminal coordinates separate the branches from one another.

To obtain a family of solutions rather than unrelated subsequential
limits, the reference solutions must converge to the \emph{same}
$\bar u$ independently of the chosen terminal coordinate.  This is
where convergence of the \emph{exact nonlinear reference correction}
from Proposition~\ref{bg:reference} is needed; convergence of the heat
array alone would not suffice.  The terminal scalar is detected by
one fixed $L^2$ test at $T_0$, and weak convergence uniformly in time
preserves that pairing.  Only the perturbed subsequence is allowed to
depend on the chosen terminal coordinate.

For each $n\geq0$, introduce the fixed physical test
\begin{equation}
 \psi_n(x)=-\ell_n^{-3/2}g(x/\ell_n),\qquad \|\psi_n\|_2=1.
 \label{limit:tests}
\end{equation}
The normalization $\langle f_{N,n},g\rangle=1$ and
\eqref{limit:physical-orbit} give
\begin{equation}
 \langle U_N(T_n)-\bar u_N(T_n),\psi_n\rangle
       =\ell_n^{1/2}a_{N,n}\geq a_-\ell_n^{1/2}.
 \label{limit:finite-observations}
\end{equation}
Uniform weak time convergence in \eqref{limit:weak-time} passes this
inequality to the limiting fields.  At the terminal time the observation
is exact:
\begin{equation}
 \langle U(T_0)-\bar u(T_0),-g\rangle=a_\ast.
 \label{limit:terminal-coordinate}
\end{equation}

There is a single compactly supported smooth test that distinguishes every
branch from the reference solution at the same interior time.
Choose $\phi_1\in C_c^\infty(\R^3;\R^3)$ with
\begin{equation}
 \|\phi_1-\psi_1\|_2<\frac{a_-\ell_1^{1/2}}{4\sqrt{E_A}}.
 \label{limit:compact-test}
\end{equation}
By \eqref{limit:uniform-energy} and \eqref{limit:finite-observations},
\begin{equation}
 \langle U(T_1)-\bar u(T_1),\phi_1\rangle
       \geq\frac12a_-\ell_1^{1/2}>0.
 \label{limit:observations}
\end{equation}
Both $T_1$ and $\phi_1$ are independent of the terminal amplitude and of
the subsequence used to obtain the branch.

The physical reference in \eqref{limit:reference} is exactly
$\bar u^{N,0}$ by the covariance \eqref{bg:covariance}.
The cutoff-convergence argument at the end of
Proposition~\ref{bg:reference} therefore gives
$\bar u_N\to\bar u^{\infty,0}$ in $C([0,T_0];L^2)$.
Set $\bar u=\bar u^{\infty,0}$ before choosing a terminal amplitude.
For each $a_\ast\in(a_-,a_+)$, carry out the finite shooting and
extract a convergent subsequence only of the perturbed solutions.
The reference limit is the same on every such subsequence, while
\eqref{limit:terminal-coordinate} fixes the chosen branch's terminal
coordinate.  Different terminal amplitudes therefore give different
solutions with the same initial datum.
  Setting
\begin{equation}
 a_\ast=a_-+\vartheta(a_+-a_-),\qquad 0<\vartheta<1,
 \label{limit:parameter}
\end{equation}
defines the family $u^\vartheta$.

This parameter labels the branches by their terminal coordinate.
No continuity of the chosen branch as a function of $\vartheta$ is
asserted or needed; the exact coordinate already makes the labeling
injective.

Let $u,v$ be two members of the family and put $w=u-v$.
For $0<\epsilon<T_0$, their smoothness and the difference equation give
\begin{equation}\label{limit:positive-time-difference}
 \begin{aligned}
 \frac12\frac{\diff}{\diff t}\|w\|_2^2+\|\nabla w\|_2^2
 &=-\int_{\R^3}(w\cdot\nabla)v\cdot w\idiff x\\
 &=\int_{\R^3}(w\cdot\nabla)w\cdot v\idiff x\\
 &\leq\frac12\|\nabla w\|_2^2+C\|v\|_\infty^2\|w\|_2^2,
 \qquad \epsilon\leq t\leq T_0.
 \end{aligned}
\end{equation}
The coefficient is integrable on $[\epsilon,T_0]$ by
\eqref{limit:mixed-bounds}. Gronwall's inequality therefore propagates
agreement at any positive time to $T_0$.

Any two different members, including the reference member, differ at every
time $0<t\leq T_0$.  If they agreed at one such time,
\eqref{limit:positive-time-difference} would imply agreement
through $T_0$, contradicting \eqref{limit:terminal-coordinate}.
Their strong time continuity consequently gives
\begin{equation}
 \int_0^\varepsilon\|u(t)-v(t)\|_2^2\idiff t>0
       \qquad (0<\varepsilon\leq T_0)
 \label{limit:early-separation}
\end{equation}
for every distinct pair.  The geometric observations in
\eqref{limit:finite-observations} also provide explicit witnesses at
times tending to zero for each nonreference branch.

\subsection{Completion and global continuation}

We now fix the parameters in their order of dependence and combine
Propositions~\ref{nl:transition}, \ref{nl:preparation},
\ref{limit:compactness}, and \ref{limit:mixed} with
Lemma~\ref{nl:shooting} to prove Theorem~\ref{intro:main}. The final
corollary extends each branch by the suitable energy construction.

\label{limit:main}

\begin{proof}[Proof of Theorem~\ref{intro:main}]
The seed construction fixes $G$ satisfying \eqref{hist:parameters} and the
preparation constant $\Theta$.  Fix the exponents in
\eqref{hist:times-exponents}.  Choose the frequency threshold in the
transfer estimate sufficiently large, and prescribe
$\lambda=2\varrho$ at $s=s_\star$ in Lemma~\ref{dec:retirement}
before choosing the width, as in Section~\ref{sec:nonlinear}.
All lower bounds on $\beta$ in the reference construction, source history,
decay, transfer, and preparation estimates involve constants fixed
independently of $\beta$.  Choose one width exponent satisfying these
finitely many bounds, including $\beta/16>6\varrho$.
Finally take $A$ sufficiently large for the stated estimates and the
strict exponential inequalities.  These choices are independent of the
finite depth.

Propositions~\ref{nl:transition} and \ref{nl:preparation} and
Lemma~\ref{nl:shooting} then give the finite solutions used in
\eqref{limit:finite-data} and \eqref{limit:physical-orbit}.
Equation \eqref{limit:initial-convergence} identifies their common strong
initial trace as \eqref{limit:initial-data}.
Proposition~\ref{limit:compactness} gives suitable Leray--Hopf limits,
and Proposition~\ref{limit:mixed} proves the strict mixed-norm range,
the Lorentz endpoint, positive-time smoothness, and both energy identities.
The Gaussian estimates following \eqref{limit:initial-data} give
its integrability, nonvanishing, and smoothness away from the origin,
while \eqref{limit:datum-not-finite-lorentz} excludes every finite-index
$L^{3,\mathfrak r}$, in particular $L^3$.
Its weak-$L^3$ membership is \eqref{limit:datum-weak-three}.
Finally, \eqref{limit:terminal-coordinate}--\eqref{limit:early-separation}
prove distinctness, the fixed compact observation, and arbitrarily early
separation for the family in \eqref{limit:parameter}.
\end{proof}

\begin{corollary}[Global suitable Leray--Hopf solutions]
\label{limit:global}
The common datum in Theorem~\ref{intro:main} admits a family of distinct
global suitable Leray--Hopf solutions.  Each solution retains the
properties of that theorem on $[0,T_0]$.
\end{corollary}

\begin{proof}
At $T_0$, each branch has a smooth finite-energy value.

Write this value as $v_0$.  For a smooth spatial approximate identity
$\rho_\varepsilon$, solve
\[
 \partial_t v_\varepsilon-\Delta v_\varepsilon
 +(\rho_\varepsilon*v_\varepsilon)\cdot\nabla v_\varepsilon
 +\nabla p_\varepsilon=0,\qquad
 \operatorname{div}v_\varepsilon=0,\qquad v_\varepsilon(0)=v_0.
\]
For fixed $\varepsilon$, derivatives of the advecting field are bounded
by $C_{\varepsilon,j}\|v_\varepsilon\|_2$.  Its energy identity and
differentiated parabolic energies therefore give continuation on every
finite interval.  Its local energy flux and pressure are
\[
 \frac{|v_\varepsilon|^2}{2}(\rho_\varepsilon*v_\varepsilon)
       +p_\varepsilon v_\varepsilon,
 \qquad
 p_\varepsilon=\sum_{i,j}\mathcal R_i\mathcal R_j
       ((\rho_\varepsilon*v_\varepsilon)_i v_{\varepsilon,j}).
\]
The energy, time-derivative and exterior-tail estimates in
Proposition~\ref{limit:compactness} hold uniformly in $\varepsilon$,
using the $L^{\mathfrak p}$ contraction of spatial mollification.  On every finite
time interval they give global strong $L^3_{t,x}$ convergence after
extraction.  The mollified velocities have the same limit, since
\[
 \|\rho_\varepsilon*v_\varepsilon-v\|_{L^3_{t,x}}
 \leq\|v_\varepsilon-v\|_{L^3_{t,x}}
       +\|\rho_\varepsilon*v-v\|_{L^3_{t,x}}\longrightarrow0.
\]
The pressures converge in $L^{3/2}$ and the fluxes in $L^1$.
Lower semicontinuity of localized dissipation proves suitability.
The initial energy inequality gives the strong initial trace, and the
almost-everywhere starting-time argument of
Proposition~\ref{limit:compactness} gives the full Leray--Hopf inequality.
A diagonal extraction over increasing time intervals produces a global
suitable solution from $v_0$; this is the suitable energy construction
used here \cite{Ler34,Hop51,CKN82}.  Concatenate that solution with
the constructed branch.
  Strong $L^2$ continuity at $T_0$ passes the
initial term of the later global and local energy inequalities across
the join.  The earlier branch satisfies equality by
\eqref{limit:energy-equality}, so the concatenation is suitable and
Leray--Hopf on every finite time interval.  Its distinction from the
other branches persists on the already constructed interval.
\end{proof}

\vspace{0.2in}
\noindent{\bf Acknowledgments.}
This work was supported in part by NSF Grant DMS-2508463,
the Choi Family Gift Fund, and the Mike Yan Gift Fund.

\vspace{0.2in}
\noindent{\bf Disclosure of AI assistance.} The cascade approach builds on the authors’ earlier work on two-dimensional Euler nonuniqueness \cite{HQPZ26b,HQPZ26a} and the related literature discussed in the Introduction. AI tools, including ChatGPT and Astra, were employed to accelerate mathematical exploration and the development of the analytical strategy. The authors directed this process, critically evaluated the AI-generated suggestions, and independently revised and verified the resulting mathematical arguments. The authors take full responsibility for the mathematical content, correctness, and presentation of the results.

\appendix

\section{Uniqueness in nonendpoint Lorentz spaces}
\label{app:lorentz}

This appendix proves the uniqueness criterion used to compare
Theorem~\ref{intro:main} with the finite-index Lorentz classes.
Phuc's regularity theorem applies to a Leray--Hopf solution in
$L^\infty_tL^{3,\mathfrak r}_x$, $3<\mathfrak r<\infty$; the uniqueness assertion in
\cite[Theorem~1.7]{Phu15} additionally assumes that its initial value
belongs to $L^3$.  The finite-index condition also yields uniqueness by
the classical local mild-solution construction and relative energy
method \cite{Kat84,Dub03}. We include the argument with the initial
trace and the energy inequality at every restarting time.

All spatial spaces in this appendix are over $\R^3$, and $\mathfrak r$ denotes
the second Lorentz exponent.
For a measurable scalar or vector field $f$, $1<\mathfrak p<\infty$, and
$1\leq \mathfrak q<\infty$, we use the Lorentz quasinorm
\begin{equation}\label{app:lorentz-norm}
 \|f\|_{L^{\mathfrak p,\mathfrak q}}
 =\left(\mathfrak p\int_0^\infty
       \Bigl(\lambda\bigl|\{x:|f(x)|>\lambda\}\bigr|^{1/\mathfrak p}\Bigr)^{\mathfrak q}
       \frac{\idiff\lambda}{\lambda}\right)^{1/\mathfrak q}.
\end{equation}
For $\mathfrak q=\infty$, the integral is replaced by the supremum of
$\lambda|\{|f|>\lambda\}|^{1/\mathfrak p}$, where $|E|$ denotes the Lebesgue measure of a measurable set $E$. In particular, the $L^{3,\infty}$ space is equipped with
\begin{equation}\label{intro:weak-three}
 \|f\|_{L^{3,\infty}(\R^3)}
 =\sup_{\lambda>0}\lambda
       \bigl|\{x\in\R^3:|f(x)|>\lambda\}\bigr|^{1/3}.
\end{equation}
The spaces $L^{\mathfrak p,\mathfrak q}$ increase with $\mathfrak q$,
and $L^{\mathfrak p,\mathfrak p}=L^{\mathfrak p}$.  When $\mathfrak q<\infty$, smooth compactly supported
functions are dense in $L^{\mathfrak p,\mathfrak q}$.

\begin{proposition}[Weak--strong uniqueness at finite Lorentz index]
\label{app:lorentz-uniqueness}
Let $T>0$, $1\leq \mathfrak r<\infty$, and $u_0\in L^2_\sigma$.
Suppose that $u$ is a Leray--Hopf solution of \eqref{intro:NS} on
$[0,T]$ with initial value $u_0$ and
\begin{equation}\label{app:lorentz-hypothesis}
 u\in L^\infty(0,T;L^{3,\mathfrak r}).
\end{equation}
Then $u_0\in L^{3,\mathfrak r}$, and every Leray--Hopf solution $v$ on $[0,T]$
with the same initial value satisfies
\begin{equation}\label{app:lorentz-uniqueness-conclusion}
 u(t)=v(t)\quad\text{in }L^2,\qquad 0\leq t\leq T.
\end{equation}
\end{proposition}

The proof uses the Lorentz H\"older and Sobolev inequalities
\begin{equation}\label{app:lorentz-products}
 \|fg\|_2\leq C\|f\|_{L^{3,\infty}}\|g\|_{L^{6,2}},
 \qquad
 \|g\|_{L^{6,2}}\leq C\|\nabla g\|_2,
 \qquad g\in H^1.
\end{equation}
The first inequality also applies to tensor products.  These are
standard consequences of the Lorentz product and convolution
inequalities \cite{ONe63}.  We first obtain an energy identity and
then prove the local comparison needed at each initial value.

\begin{lemma}[Energy identity under a weak-$L^3$ bound]
\label{app:lorentz-energy}
Let $T>0$ and $a\in L^2_\sigma$.  Suppose that $w$ is a divergence-free
distributional solution of \eqref{intro:NS} with initial value $a$
in the sense of \eqref{intro:weak-formulation}, and that
\begin{equation}\label{app:lorentz-energy-class}
 w\in L^\infty(0,T;L^2)\cap L^2(0,T;H^1)
                  \cap L^\infty(0,T;L^{3,\infty}).
\end{equation}
Then $w\in C([0,T];L^2_\sigma)$, $w(0)=a$, and
\begin{equation}\label{app:lorentz-energy-identity}
 \frac12\|w(t)\|_2^2+\int_s^t\|\nabla w(\tau)\|_2^2\idiff\tau
   =\frac12\|w(s)\|_2^2,
 \qquad 0\leq s\leq t\leq T.
\end{equation}
\end{lemma}

\begin{proof}
By \eqref{app:lorentz-products}, the nonlinear tensor satisfies
\begin{equation}\label{app:lorentz-tensor-energy}
 \|w\otimes w\|_{L^2(0,T;L^2)}
 \leq C\|w\|_{L^\infty(0,T;L^{3,\infty})}
          \|\nabla w\|_{L^2(0,T;L^2)}.
\end{equation}
Set $H^1_\sigma=H^1\cap L^2_\sigma$, and let $H^{-1}_\sigma$
denote its dual with respect to the $L^2$ pairing.  The projected equation therefore gives
$\partial_t w\in L^2(0,T;H^{-1}_\sigma)$.
Together with $w\in L^2(0,T;H^1_\sigma)$, the Hilbert-space energy
chain rule gives a representative in $C([0,T];L^2_\sigma)$.
Its initial value is $a$ by \eqref{intro:weak-formulation}.
The nonlinear pairing is integrable by
\eqref{app:lorentz-tensor-energy}, and, for almost every time,
\begin{equation}\label{app:lorentz-cancellation}
 \int_{\R^3}(w\otimes w):\nabla w\idiff x=0.
\end{equation}
Indeed, the cancellation holds for smooth compactly supported
solenoidal fields and passes by $H^1$ approximation, using
$H^1\hookrightarrow L^4$.  The chain rule and
\eqref{app:lorentz-cancellation} prove
\eqref{app:lorentz-energy-identity}.
\end{proof}

\begin{lemma}[Local comparison from a finite-index Lorentz datum]
\label{app:lorentz-local}
Let $1\leq \mathfrak r<\infty$ and $a\in L^2_\sigma\cap L^{3,\mathfrak r}$.
There is $\tau>0$ such that all Leray--Hopf solutions with initial
value $a$ coincide on $[0,\tau]$.  They agree there with a mild
solution $W$ satisfying
\begin{equation}\label{app:lorentz-local-solution}
 W\in C([0,\tau];L^2_\sigma)
       \cap C([0,\tau];L^{3,\infty})
       \cap L^2(0,\tau;H^1),\qquad W(0)=a.
\end{equation}
\end{lemma}

\begin{proof}
We first construct a mild solution by the heat-kernel iteration,
then compare it with an arbitrary energy solution by splitting the
mild solution into a bounded part and a small weak-$L^3$ remainder.
For $\tau>0$, let $\mathcal K_\tau$ be the Banach space of fields
$V\in C((0,\tau];L^4)$ with norm
\begin{equation}\label{app:lorentz-kato-norm}
 \|V\|_{\mathcal K_\tau}
   =\sup_{0<t\leq\tau}t^{1/8}\|V(t)\|_4<\infty,
 \qquad \lim_{t\downarrow0}t^{1/8}\|V(t)\|_4=0.
\end{equation}
The heat estimate
$\|e^{t\Delta}a\|_4\leq Ct^{-1/8}\|a\|_{L^{3,\infty}}$
follows by interpolating its weak-$L^3$ bound with
$\|e^{t\Delta}a\|_\infty\leq Ct^{-1/2}\|a\|_{L^{3,\infty}}$.
Density at finite Lorentz index therefore implies
\begin{equation}\label{app:lorentz-small-heat}
 \lim_{\tau\downarrow0}\|e^{t\Delta}a\|_{\mathcal K_\tau}=0.
\end{equation}
More explicitly, for any $b\in C_c^\infty(\R^3;\R^3)$ one has
\[
 \|e^{t\Delta}a\|_{\mathcal K_\tau}
 \leq C\|a-b\|_{L^{3,\infty}}+\tau^{1/8}\|b\|_4.
\]
First approximate $a$ by $b$ in $L^{3,\mathfrak r}\hookrightarrow L^{3,\infty}$,
then let $\tau\downarrow0$.  The same approximation, together with
boundedness of the heat semigroup on $L^{3,\infty}$, proves its
strong continuity at $a$ in that space.

For time-dependent fields $V,Z$, define the bilinear integral
\begin{equation}\label{app:lorentz-bilinear}
 \mathcal N(V,Z)(t)
   =\int_0^t e^{(t-s)\Delta}\PP\operatorname{div}
                   (V(s)\otimes Z(s))\idiff s.
\end{equation}
The heat kernel and boundedness of $\PP$ on $L^3$ and $L^4$ give
\begin{equation}\label{app:lorentz-bilinear-kernels}
 \begin{aligned}
 \|\mathcal N(V,Z)(t)\|_4
   &\leq C\int_0^t(t-s)^{-7/8}\|V(s)\|_4\|Z(s)\|_4\idiff s,\\
 \|\mathcal N(V,Z)(t)\|_3
   &\leq C\int_0^t(t-s)^{-3/4}\|V(s)\|_4\|Z(s)\|_4\idiff s.
 \end{aligned}
\end{equation}
For example, the change of variables $s=t\sigma$ gives
\begin{equation}\label{app:lorentz-time-kernel}
 \int_0^t(t-s)^{-7/8}s^{-1/4}\idiff s
 =t^{-1/8}\int_0^1(1-\sigma)^{-7/8}\sigma^{-1/4}\idiff\sigma.
\end{equation}
Both endpoint singularities in the last integral are integrable.
Replacing $7/8$ by $3/4$ gives a finite constant independent of $t$.
The two bounds in \eqref{app:lorentz-bilinear-kernels} therefore imply
\begin{equation}\label{app:lorentz-bilinear-bounds}
 \begin{aligned}
 \|\mathcal N(V,Z)\|_{\mathcal K_\tau}
   &\leq C\|V\|_{\mathcal K_\tau}\|Z\|_{\mathcal K_\tau},\\
 \sup_{0<s\leq t}\|\mathcal N(V,Z)(s)\|_3
   &\leq C\|V\|_{\mathcal K_t}\|Z\|_{\mathcal K_t},
             \qquad 0<t\leq\tau.
 \end{aligned}
\end{equation}
In particular, the second quantity tends to zero as $t\downarrow0$.
The kernel bounds also give continuity at positive times by splitting
the integral away from its endpoints.

Choose $\tau$ so that \eqref{app:lorentz-small-heat} makes the first
bound in \eqref{app:lorentz-bilinear-bounds} contractive on the ball
of radius $2\|e^{t\Delta}a\|_{\mathcal K_\tau}$.
The resulting fixed point solves the mild equation
\begin{equation}\label{app:lorentz-mild}
 W(t)=e^{t\Delta}a-\mathcal N(W,W)(t).
\end{equation}
It is solenoidal, smooth for positive times, and bounded on each
closed positive-time strip.  The usual positive-time heat smoothing
applies since $W$ is continuous into $L^4$ there.
The heat term is strongly continuous at $a$ in $L^{3,\infty}$,
and the second bound in \eqref{app:lorentz-bilinear-bounds} makes the
nonlinear term tend to zero in $L^3$. Hence
$W\in C([0,\tau];L^{3,\infty})$ with $W(0)=a$.

To obtain the energy bounds, use $a_j=e^{j^{-1}\Delta}a$ as initial
values in the same iteration.  These data converge to $a$ in both
$L^2$ and $L^{3,\infty}$, and heat contraction on $L^4$ gives
$\|e^{t\Delta}a_j\|_{\mathcal K_\tau}
\leq\|e^{t\Delta}a\|_{\mathcal K_\tau}$.
The corresponding fixed points $W_j$ converge to $W$ in
$\mathcal K_\tau$.  Each agrees initially with the classical solution
from $a_j$ and continues to $\tau$: on every closed positive-time
interval its $L^4$ norm is bounded, so the classical continuation
criterion applies \cite{Kat84}.  Their energy identities yield
\begin{equation}\label{app:lorentz-approximate-energy}
 \frac12\|W_j(t)\|_2^2+\int_0^t\|\nabla W_j(s)\|_2^2\idiff s
   =\frac12\|a_j\|_2^2\leq\frac12\|a\|_2^2,
 \qquad 0\leq t\leq\tau.
\end{equation}
Weak compactness and convergence in $\mathcal K_\tau$ place $W$ in
$L^\infty(0,\tau;L^2)\cap L^2(0,\tau;H^1)$.
The mild equation \eqref{app:lorentz-mild} includes the initial
value $a$ in the sense of \eqref{intro:weak-formulation}.
Lemma~\ref{app:lorentz-energy} now gives strong $L^2$ continuity and
energy equality for $W$, proving \eqref{app:lorentz-local-solution}.

It remains to compare $W$ with an arbitrary Leray--Hopf solution $z$
from $a$.  Given $\varepsilon>0$, density and the
$L^{3,\infty}$ continuity in \eqref{app:lorentz-local-solution}
provide a smooth bounded field $b$ and a smaller $\tau>0$ such that
\begin{equation}\label{app:lorentz-splitting}
 W(t)=b+\zeta(t),\qquad
 \sup_{0\leq t\leq\tau}\|\zeta(t)\|_{L^{3,\infty}}\leq\varepsilon.
\end{equation}
Write $h=z-W$.  The relative energy inequality is
\begin{equation}\label{app:lorentz-relative-energy}
 \frac12\|h(t)\|_2^2+\int_0^t\|\nabla h(s)\|_2^2\idiff s
 \leq\int_0^t\left|\int_{\R^3}W\cdot(h\cdot\nabla)h\idiff x\right|
                  \idiff s.
\end{equation}
To justify the initial endpoint, first use the smooth solution $W$
as a test on $[\delta,t]$, where $\delta$ is an admissible starting
time for the energy inequality of $z$.
Choose such times tending to zero.  The initial difference tends
to zero strongly in $L^2$, and the integrand is bounded by
$C\|W\|_{L^\infty(0,\tau;L^{3,\infty})}\|\nabla h(s)\|_2^2$
by \eqref{app:lorentz-products}.  This integrable majorant permits
passage to \eqref{app:lorentz-relative-energy}.

Using \eqref{app:lorentz-splitting} in the nonlinear term gives
\begin{equation}\label{app:lorentz-relative-estimate}
 \begin{aligned}
 \left|\int_{\R^3}\zeta\cdot(h\cdot\nabla)h\idiff x\right|
   &\leq C\varepsilon\|\nabla h\|_2^2,\\
 \left|\int_{\R^3}b\cdot(h\cdot\nabla)h\idiff x\right|
   &\leq\frac14\|\nabla h\|_2^2
                    +C\|b\|_\infty^2\|h\|_2^2.
 \end{aligned}
\end{equation}
Choose $\varepsilon$ so that $C\varepsilon\leq1/4$.
Equations~\eqref{app:lorentz-relative-energy} and
\eqref{app:lorentz-relative-estimate}, followed by Gronwall's
inequality, give $h=0$ on $[0,\tau]$.
\end{proof}

\begin{proof}[Proof of Proposition~\ref{app:lorentz-uniqueness}]
By $L^{3,\mathfrak r}\hookrightarrow L^{3,\infty}$ and
Lemma~\ref{app:lorentz-energy}, the reference solution $u$ is strongly
continuous into $L^2$ and satisfies
\eqref{app:lorentz-energy-identity} on $[0,T]$.
Set $M=\|u\|_{L^\infty(0,T;L^{3,\mathfrak r})}$.
For every $s\in[0,T]$, choose times $t_j\to s$ outside the exceptional
set in \eqref{app:lorentz-hypothesis}.
Strong $L^2$ convergence admits a spatially almost-everywhere
convergent subsequence, and Fatou's lemma in
\eqref{app:lorentz-norm} gives
\begin{equation}\label{app:lorentz-all-traces}
 u(s)\in L^2_\sigma\cap L^{3,\mathfrak r},\qquad
 \|u(s)\|_{L^{3,\mathfrak r}}\leq M,\qquad 0\leq s\leq T.
\end{equation}
In particular, the initial datum belongs to $L^{3,\mathfrak r}$.
Lemma~\ref{app:lorentz-local} therefore gives agreement of $u$ and
$v$ on an initial interval.

Let $\tau_\ast$ be the supremum of the times $s\in[0,T]$ for which
$u=v$ on $[0,s]$.  Weak $L^2$ continuity of both solutions gives
\begin{equation}\label{app:lorentz-restart-datum}
 u(\tau_\ast)=v(\tau_\ast)\in L^2_\sigma\cap L^{3,\mathfrak r},
\end{equation}
where the membership follows from \eqref{app:lorentz-all-traces}.
Suppose that $\tau_\ast<T$.
Subtract the energy equality for $u$ on $[0,\tau_\ast]$ from the
energy inequality \eqref{intro:energy-inequality} for $v$ starting
at zero.  Since the solutions and their weak gradients agree before
$\tau_\ast$, the result is
\begin{equation}\label{app:lorentz-restarted-energy}
 \frac12\|v(t)\|_2^2+
       \int_{\tau_\ast}^t\|\nabla v(s)\|_2^2\idiff s
   \leq\frac12\|u(\tau_\ast)\|_2^2,
 \qquad \tau_\ast\leq t\leq T.
\end{equation}
Weak continuity and \eqref{app:lorentz-restarted-energy} give
$v(t)\to u(\tau_\ast)$ strongly in $L^2$ as $t\downarrow\tau_\ast$.
The distributional equation with this initial value follows by
restricting its weak formulation, and the almost-everywhere
starting-time energy inequalities are inherited from $v$.
Thus both solutions, translated by $\tau_\ast$ in time, are
Leray--Hopf solutions with initial value $u(\tau_\ast)$.
Lemma~\ref{app:lorentz-local}, applied to
\eqref{app:lorentz-restart-datum}, extends their agreement beyond
$\tau_\ast$.  This contradicts its definition, so $\tau_\ast=T$
and \eqref{app:lorentz-uniqueness-conclusion} follows.
\end{proof}

\bibliographystyle{myalpha}
\bibliography{references}

\end{document}